\documentclass[12pt,a4paper,twoside,openright,titlepage]{book}
\usepackage[latin1]{inputenc}
\usepackage[italian,english]{babel}
\usepackage{graphicx,amsmath,amsfonts,verbatim,enumerate,subfigure}
\usepackage{amsthm}
\usepackage{yhmath}
\usepackage[norules]{frontespizio}

\theoremstyle{plain}
\usepackage{fancyhdr}
\newcommand{\fncyblank}{\fancyhf{}}
\newenvironment{abstract}%
{\cleardoublepage\fncyblank\null\vfill\begin{center}%
\bfseries\abstractname\end{center}}%
{\vfill\null}

\providecommand{\norm}[1]{\lVert#1\rVert}
\newcommand{\numberset}{\mathbb}
\newcommand{\R}{\numberset{R}}

\newcommand{\Z}{\numberset{Z}}
\newcommand{\N}{\numberset{N}}
\newcommand{\Poly}{\numberset{P}}

\usepackage{amsthm}
\theoremstyle{plain}
\newtheorem{theorem}{Theorem}[section]

\newtheorem{proposition}[theorem]{Proposition}
\newtheorem{corollary}[theorem]{Corollary}
\theoremstyle{definition}
\newtheorem{definition}{Definition}
\theoremstyle{remark}
\newtheorem{remark}{Remark}
\newtheorem{example}{Example}

\begin{document}

\frontmatter

\numberwithin{equation}{chapter}

\selectlanguage{english}
%\maketitle

%\begin{frontespizio}
%\Universita{Padova}
%\Facolta{Scienze Matematiche, Fisiche e Naturali}
%\Logo[4cm]{logo1.pdf}
%%\Filigrana[height=10cm,before=1,after=1]{logo3.pdf} %disegno sullo sfondo
%\Corso[Laurea]{Matematica}% e Applicazioni
%\Titoletto{Tesi di Laurea Magistrale}
%\Titolo{%Medical image reconstruction\\ using kernel based
%%methods
%MEDICAL IMAGE RECONSTRUCTION\\ USING KERNEL BASED METHODS}
%\Candidato{Amos Sironi}
%\Relatore{Prof. Stefano De Marchi}
%\Correlatore{Prof. Armin Iske \\ (University of Hamburg)}
%\Annoaccademico{2010-2011}
%%\Rientro{2cm}
%\Margini{2cm}{1.5cm}{2cm}{1cm}
%\Preambolo{\usepackage{type1cm}}
%\end{frontespizio}

%\subfile{./frontmatter/title_page}

%\subfile{./frontmatter/acknowledgments}
%
%abstract
\begin{abstract}
The image reconstruction problem consists in finding an approximation of a function $f$ starting from its Radon transform $Rf$. This problem arises in the ambit of medical imaging when one tries to reconstruct the internal structure of the body, starting from its X-ray tomography.
The classical approach to this problem is based on the Back-Projection Formula. This formula gives an analytical inversion of the Radon transform, provided that all the values of $Rf$ are known. In applications only a discrete set of values of $Rf$ is given, thus, one can only obtain an approximation of $f$.
Another class of methods, called ART, can be used to solve the reconstruction problem. Following the ideas contained in ART, we try to apply the Hermite-Birkhoff interpolation to the reconstruction problem. It turns out that, since the Radon transform of a kernel basis function can be infinity, a regularization technique is needed.
The method we present here is then based on positive definite kernel functions and it is very flexible thanks to the possibility to choose different kernels and parameters. We study the behavior of the methods and compare them with classical algorithms.
\end{abstract}

\tableofcontents

%content
\chapter{Introduction and content}
This thesis is the result of a three-months stage at the Univerit\"{a}t Hamburg during which I studied the problem of clinical image reconstruction, i.e. the problem of obtaining the image of the internal structure of a sample starting from its X-ray tomography. From a mathematical point of view this correspond to find a function $f$ knowing its Radon transform $Rf$.

In the first Chapter the problem of image reconstruction is defined, we formalize the concept of Computed Axial Tomography and the history behind it.

In Chapter 2 we discuss the mathematical aspect of the problem and its relation with the Radon transform. Then we follow the classical approach for solving the problem and deduce an inversion formula for the Radon transform: the \emph{Back-Projection Formula}. Finally, we adapt the Back-Projection Formula to be used in real applications and thus we obtain the classical Fourier-based discrete image reconstruction algorithms.

In Chapter 3 we introduce a different class of methods, called Algebraic Reconstruction Techniques (ART) and use them to solve our problem.

Following the ART approach, in Chapter 4 we describe kernel based methods and show how they can be used to solve the image reconstruction problem.

Chapters 5 and 6 are the original part of the work. In Chapter 5 we introduce a regularization technique that is necessary to implement kernel based image reconstruction and we realize such methods using specific positive definite kernel functions. In Chapter 6 we study from a numerical point of view the behavior of the methods in function of particular shape parameters and compare these methods with the Fourier based algorithms.

In order to use the algorithms in a simple way, we also realized a graphical user interface that allows the user to test the algorithms on a set of predefined mathematical phantoms, with the possibility to choose options.

%\subfile{./frontmatter/symbols}
\chapter{List of symbols}
\begin{tabular}{ll}
$l_{t,\theta}$ & line in the plane characterized by values $t$ and $\theta$\\
$Rf(t,\theta)$ & Radon transform of the function $f$ at a point $(t,\theta)$\\
$\mathcal{S}$ & Schwartz space of rapidly decreasing functions\\
$Bh(x,y)$ & back projection of the function $h$ at point $(x,y)$\\
$F_{n}f(\omega)$ & $n$-dimensional Fourier transform of the function $f$ at a point $\omega$\\
$Ff(\omega)$ & $1$-dimensional Fourier transform of the function $f$ at a point $\omega$\\
$R_{D}f$ & discrete Radon transform of the function $f$ \\
$B_{D}h$ & discrete back projection of the function $h$ \\
$F_{D}f$ & discrete Fourier transform of the function $f$ \\
$FWHM(\phi)$ & full width half maximum of the function $\phi$\\
$\lambda^{y}K(\cdot,y)$ & linear operator $\lambda$ applied to the function $K$ with respect to the variable $y$ \\
$R_{w}f$ & Radon transform of the function $f$ multiplied by the window function $w$\\
$\text{erf}(x)$ & error function evaluated at a point $x$\\
$\Poly^{d}_{k}$ & space of polynomial of degree lower or equal to $k$ on $\R^{d}$\\
$k(A)$ & 1-norm condition number of the matrix $A$
\end{tabular}

\mainmatter

%Chapter 1 - Definition of the problem
\chapter{Computed axial tomography}
\begin{figure}[htbp]
\centering%
\subfigure[A modern CT scanner \label{fig: toshiba_ct_scanner}]%
{\includegraphics[width=0.49\textwidth]{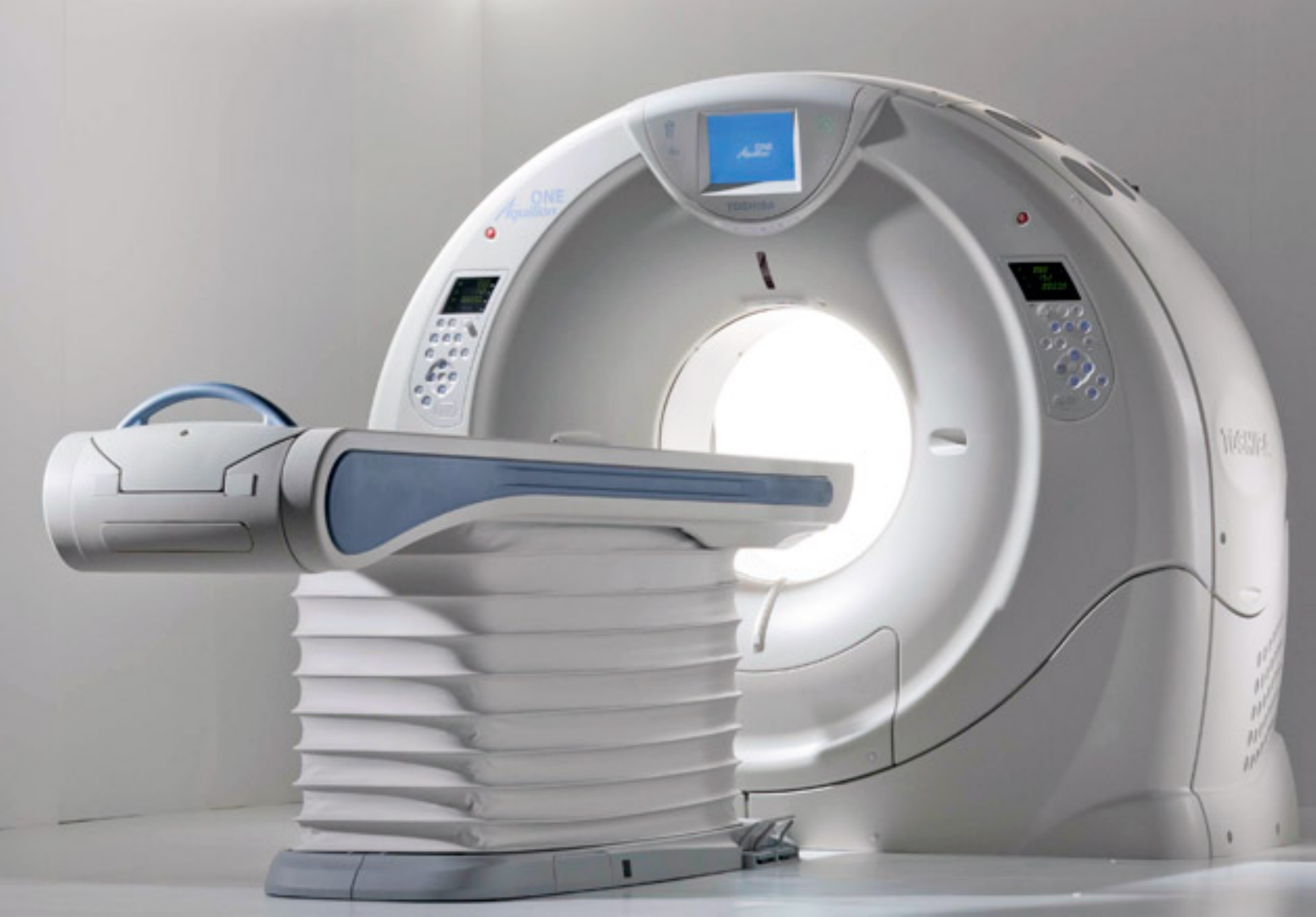}}
\subfigure[Brain CT scan \label{fig: ct_brain}]%
{\includegraphics[width=0.34\textwidth]{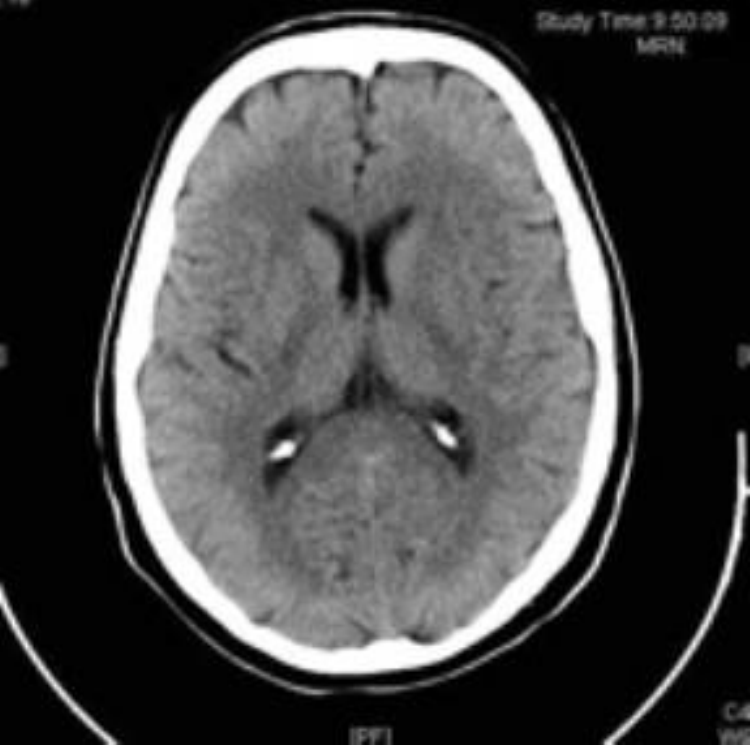}}
\caption{Computed axial tomography today}
\label{fig: ct_today}
\end{figure}

%\section{Computed axial tomography}
%\subsection{Introduction}
Computed axial tomography (CAT or CT) is a method that generates images of the interior of the body by digital computation applied to the measured transmission of X-rays tomography. In this process, an X-ray source and a set of aligned X-ray detectors are rotated around the patient (see Figure \ref{fig: toshiba_ct_scanner}). The word tomography is derived from the Greek \emph{tomos} (slice) and \emph{graphein} (to write).

\begin{figure}[htbp]
\centering%
\subfigure[Picture of Wilhelm R\"ontgen \label{fig: william_roentgen}]%
{\includegraphics[width=0.4\textwidth]{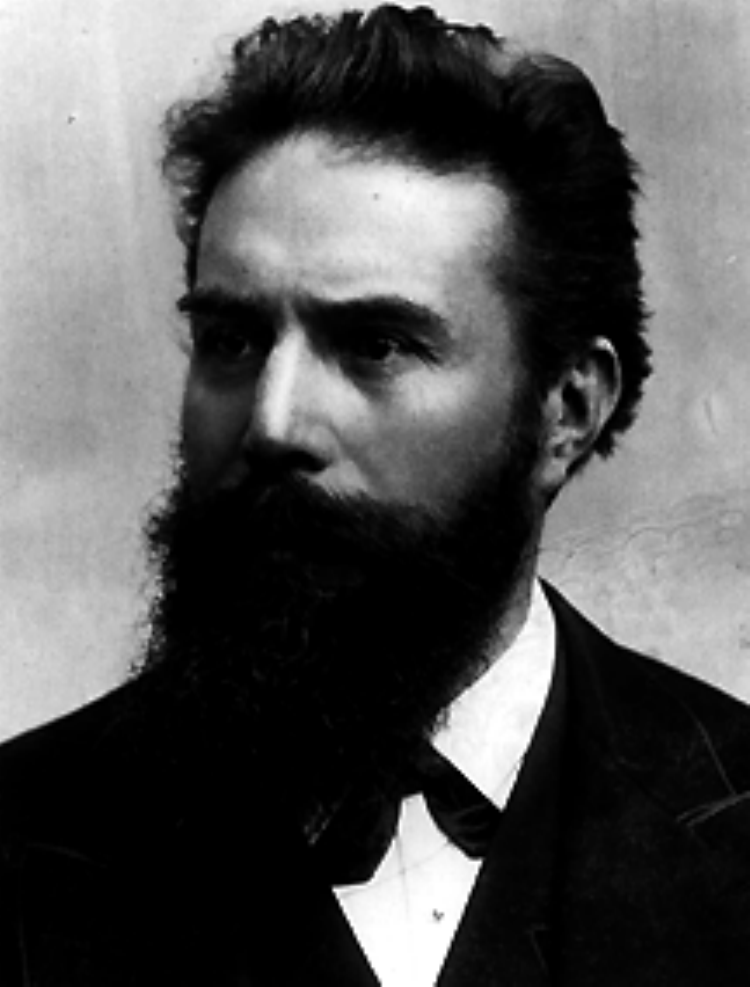}}
\subfigure[The first x-ray Frau R\"ontgen's left hand \label{fig: roentgen_first_xray}]%
{\includegraphics[width=0.36\textwidth]{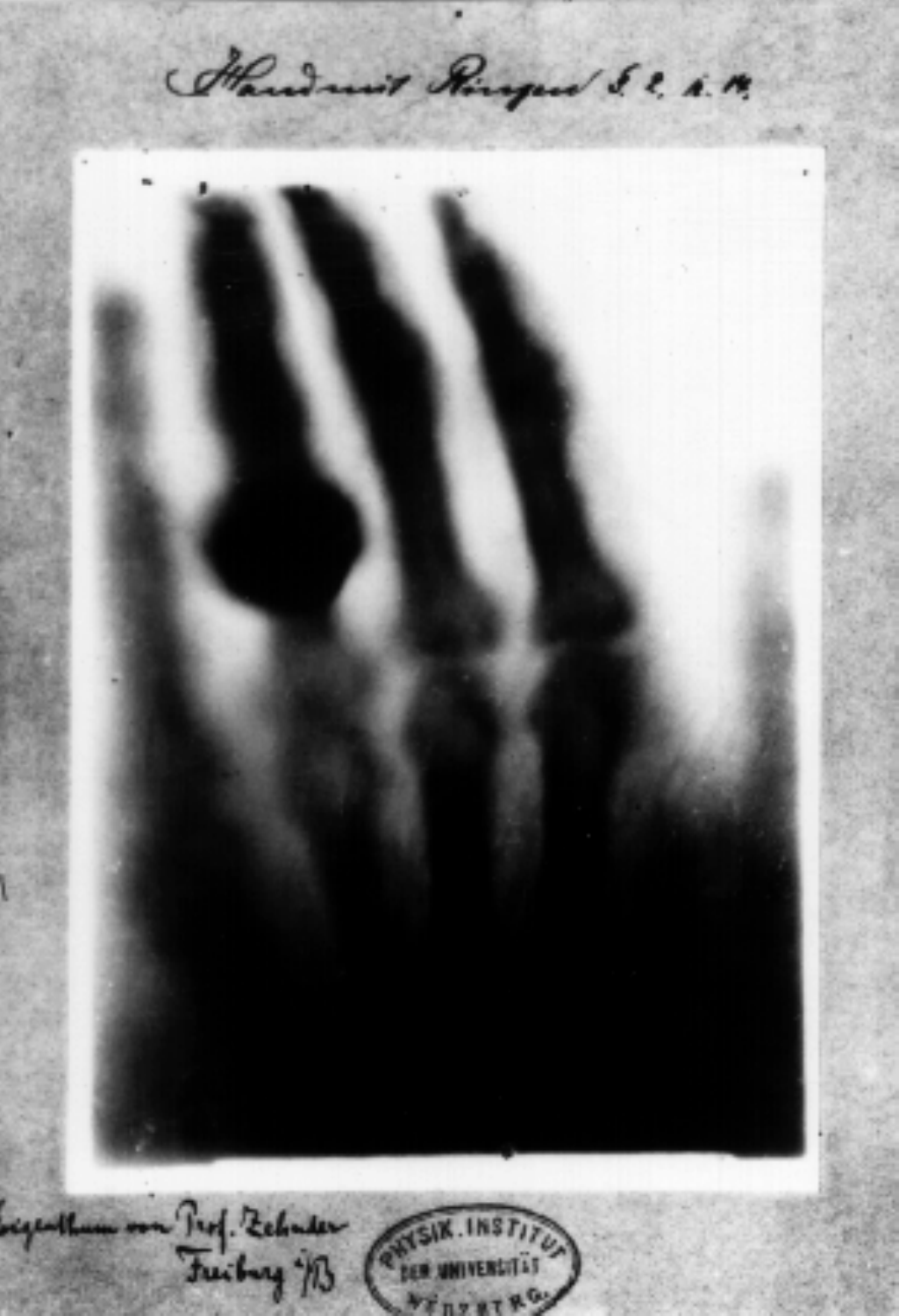}}
\caption{The discovery of X-ray}
\label{fig: dics_xray}
\end{figure}

The history of CT scan starts in Germany in 1895, when Wilhelm Conrad R\"ontgen (1859-1923; Figure \ref{fig: william_roentgen}) discovered a new type of radiation, which he called X-rays \cite{RON}. This type of electromagnetic radiation, which has shorter wavelength then visible light and the ability to penetrate matter, was immediately used to image the interior of the human body. Figure \ref{fig: roentgen_first_xray} shows one of the first X-ray images, this kind of images showed a two dimensional projection of the inner structures. In 1901 R\"ontgen received the first Nobel prize for physics. Basic to the CT technology are the theoretical principles of reconstruction of a three-dimensional object from multiple two-dimensional views relying on a mathematical model formulated by Johann Radon (1887-1956) in 1917 \cite{RAD17}.
%intro storica presa da phd th. Alfonso Isola

%da qui in poi testo "basics of madical imaging"
In 1979 the Nobel Prize for Medicine and Physiology was awarded jointly to Allan McLeod Cormack (1924-1998) and Godfrey Newbold Hounsfield (1919-2004), the two scientists primarily responsible for the development of computerized axial tomography in the 1960s and early 1970s. Cormack developed certain mathematical algorithms that could be used to create an image from X-ray data \cite{COR}. Working completely independently of Cormack and at about the same time, Hounsfield, a research scientist at EMI Central Research Laboratories in the United Kingdom, designed the first operational CT scanner, the first commercially available model and presented the first pictures of a patient's head \cite{HOU}. Compared to a plan X-ray image, the CT image showed remarkable contrast between tissues with small differences in X-ray attenuation coefficient (Figure \ref{fig: ct_brain} shows the CT scan of a section of the brain). Since 1980, the number of CT scans performed every year in the United States has risen from about 3 million to over 67 million (for further details about X-ray history one should refer to \cite{BASIC} or \cite{HISTORY}).

The problem behind CT scans is essentially mathematical: if we know the values of the integral of two- or three- dimensional function along all possible cross-sections, then how can we reconstruct the function itself? This is a particular case of what is called as an \emph{inverse problem} and it was studied by the Austrian mathematician Johann Radon  in the early part of the twentieth century. Radon's work incorporated a sophisticated use of theory of transform and integral operators.

The practical obstacles to implementing Radon's theories are several. First Radon's inversion methods assume knowledge of the behavior of the function along every cross-section, while in practice only a discrete set of cross-sections can be sampled. Thus it is possible to construct only an approximation of the solution. Second, the computation power needed to process a multitude of discrete measurements and obtain from them a good approximation of the solution has been available for just a few decades. In order to overcome these obstacles theoretical approaches and approximation methods have been developed.

\section{X-rays}
A CT scan is generated form a set of thousands of X-ray beams, consisting of 160 or more beams at each of 180 directions. When a single X-ray beam of known intensity passes through a medium, some of the energy present in the beam is absorbed by the medium and some passes through. The intensity of the beam as it emerges from the medium can be measured by a detector. The difference between the initial and final intensities tell us about the ability of the medium to absorb energy.

The idea behind the CT scan is that, by measuring the changes in the intensity of X-ray beams passing through the medium in different directions and by comparing the measurements, we can determine which location within the sample are more or less absorbent than others.

In our analysis of the X-rays behavior we will make some assumptions:
\begin{itemize}
\item X-ray beam is \emph{monochromatic}. That is each photon has the same energy level $E$ and the beam propagates at a constant frequency. If $N(x)$ denotes the number of photons per second passing through a point $x$, then the intensity of the beam at the point $x$ is 
\begin{equation*}
I(x)=E\cdot N(x);
\end{equation*}
\item X-ray beam has \emph{zero width};
\item X-ray beams are \emph{not subject to refraction or diffraction}.
\end{itemize}

Every substance has the property to absorbs a part of the photons that pass through it. To quantify this property we define the \emph{attenuation coefficient} of a material:

 \begin{definition} 
The \emph{attenuation coefficient} of a substance is the fractional number of photons removed from a beam of radiation per unit thickness of material through which it is passing due to all absorption and scattering processes.
\end{definition}

In radiology a variant of the attenuation coefficient is used: the \emph{Hounsfield unit}. Developed by Godfrey Hounfield, the Hounsfield unit represents a comparison of the attenuation coefficient of the medium with that of water. Specifically:

 \begin{definition} 
The \emph{Hounsfield unit} of a medium is 
\begin{equation*}
H_{\text{medium }}=\frac{A_{\text{medium}}-A_{\text{water}}}{A_{\text{water}}},
\end{equation*}
where $A$ denotes the attenuation coefficient.
\end{definition}

Suppose now an X-ray beam passes through some medium located between the position $x$ and the position $x+\Delta x$. Suppose $A(x)$ is the attenuation coefficient of the medium located there. Then the portion of all photons that will be absorbed in the interval $[x,x+\Delta x]$ is $p(x)=A(x)\Delta x$. The number of photons absorbed per second by the medium is then $p(x)N(x)=A(x)N(x)\Delta x$. Multiplying both sides by the energy level $E$ of each photon, we see that the loss of intensity of the X-ray over this interval is
\begin{equation*}
\Delta I\approx -A(x)I(x)\Delta x.
\end{equation*}
Let $\Delta x\rightarrow0$ to get the differential equation known as the \emph{Beer's law}:
\begin{equation}
\frac{dI}{dx}=-A(x)I(x).
\label{eq: beerLaw}
\end{equation}
In other words: \emph{The rate of change or intensity per millimeter of a nonrefractive, monochromatic, zero-width X-ray beam passing through a medium is jointly proportional to the intensity of the beam and to the attenuation coefficient of the medium.}

The differential equation \eqref{eq: beerLaw} is separable. If the beam starts at the point $x_{0}$ with initial intensity $I_{0}=I(x_{0})$ and is detected, after passing through the medium, at the point $x_{1}$ with final intensity $I_{1}=I(x_{1})$, we get
\begin{equation*}
\int_{x_{0}}^{x_{1}}{\frac{dI}{I}}=-\int_{x_{0}}^{x_{1}}{A(x)\,dx},
\end{equation*}
from which it follows that
\begin{equation}
\int_{x_{0}}^{x_{1}}{A(x)\,dx}=\ln\left( \frac{I_{0}}{I_{1}}\right).
\label{eq: averageA}
\end{equation}
Here we know the initial and final values of $I$ and we want to determine the coefficient function $A$. Thus, form the measured intensity of the X-ray we are able to compute not the values of $A$ itself, but the value of the integral of $A$ along the line of the X-ray.

From equation \eqref{eq: averageA} it is easy to see that we can not discriminate two functions that have the same value of the integral along the X-ray path $[x_{0},x_{1}]$. The fundamental question of image reconstruction asks if it is possible to do that knowing the value of the integral of $A$ along every line:

%\begin{remark}
\textbf{The fundamental question of image reconstruction}: \emph{Can we reconstruct the function $A(x,y,z)$ (within some finite region) if we know the average value of $A$ along every line that passes through the region?} (cfr. \cite{BASIC} pp. 7.)
%\end{remark}

In our study of CT scans, we will consider a two dimensional slice of the sample, obtained as the intersection of the sample and some plane, which we will generally assume coincides with the $xy$-plane. In this context, we interpret the attenuation coefficient function as a function $A(x,y)$ of two variables.

%Chapter 2 - Fourier based methods
\chapter{Fourier based methods}\label{chap: fourier_methods}
%\section{Introduction}
In this chapter we study the methods that are used nowadays in the CT scanner. The mathematical foundation of these methods is based on the work of J. Radon on an integral transform, called  in his honor \emph{Radon transform}, and its inverse. Roughly speaking we can think that sending a set of X-ray beams through a sample and measuring the intensity of the beams after their passage through it, correspond to compute the Radon transform of the sample's attenuation coefficient. Thus applying an inversion formula of the Radon transform gives us the value of the  attenuation coefficient within the sample.

In theory this is possible if we know the value of the Radon transform in every point of the sample. In practice only a discrete set of values can be recorded by a X-ray machine, that's why we can only obtain an approximation of the original attenuation coefficient function and we will have to consider problems that arise working with discrete functions, such as sampling, filtering and interpolation.

We start this chapter formalizing the concept of Radon transform. Since this operator involves the computation of the integral of a function along lines in the plane, we need first to define a suitable characterization of lines in $\R^{2}$.

\section{Characterization of lines in $\R^{2}$}
Consider again the equation \eqref{eq: averageA}:
\begin{equation}
\int_{x_{0}}^{x_{1}}{A(x)\,dx}=\ln\left( \frac{I_{0}}{I_{1}}\right).
\label{eq: averageA2}
\end{equation}
Suppose a sample of material occupies a finite region in space. At each point $(x,y,z)$ within the sample, the material there has an attenuation coefficient $A(x,y,z)$. An X-ray beam passing through the sample follows a line $l$ from an initial point $P$ (assumed to be outside the region) to a final point $Q$ (also assumed to be outside the region). The emission/detection machine measures the initial and final intensities of the beam at $P$ and $Q$, from which the value $\ln(I_{0}/I_{1})$ is calculated. According to \eqref{eq: averageA2} this is equal to the value of the integral $\int_{\overline{PQ}}{A(x,y,z)\,ds}$, where $ds$ represents arclength units along the segment $\overline{PQ}$ of the line $l$. Thus the measurement of each X-ray beam gives us information about the average value of $A$ along the path of the beam and it is fundamental to find a useful representation of lines that can help us in solving the image reconstruction problem.

For simplicity let assume that we are interested only in the cross-section of a sample that lies in the $xy$-plane. Each X-ray will follow a segment of a line in the plane and we look for a way of cataloging all such lines.

\begin{figure}[htbp]
\centering
\includegraphics[width=0.5\textwidth]{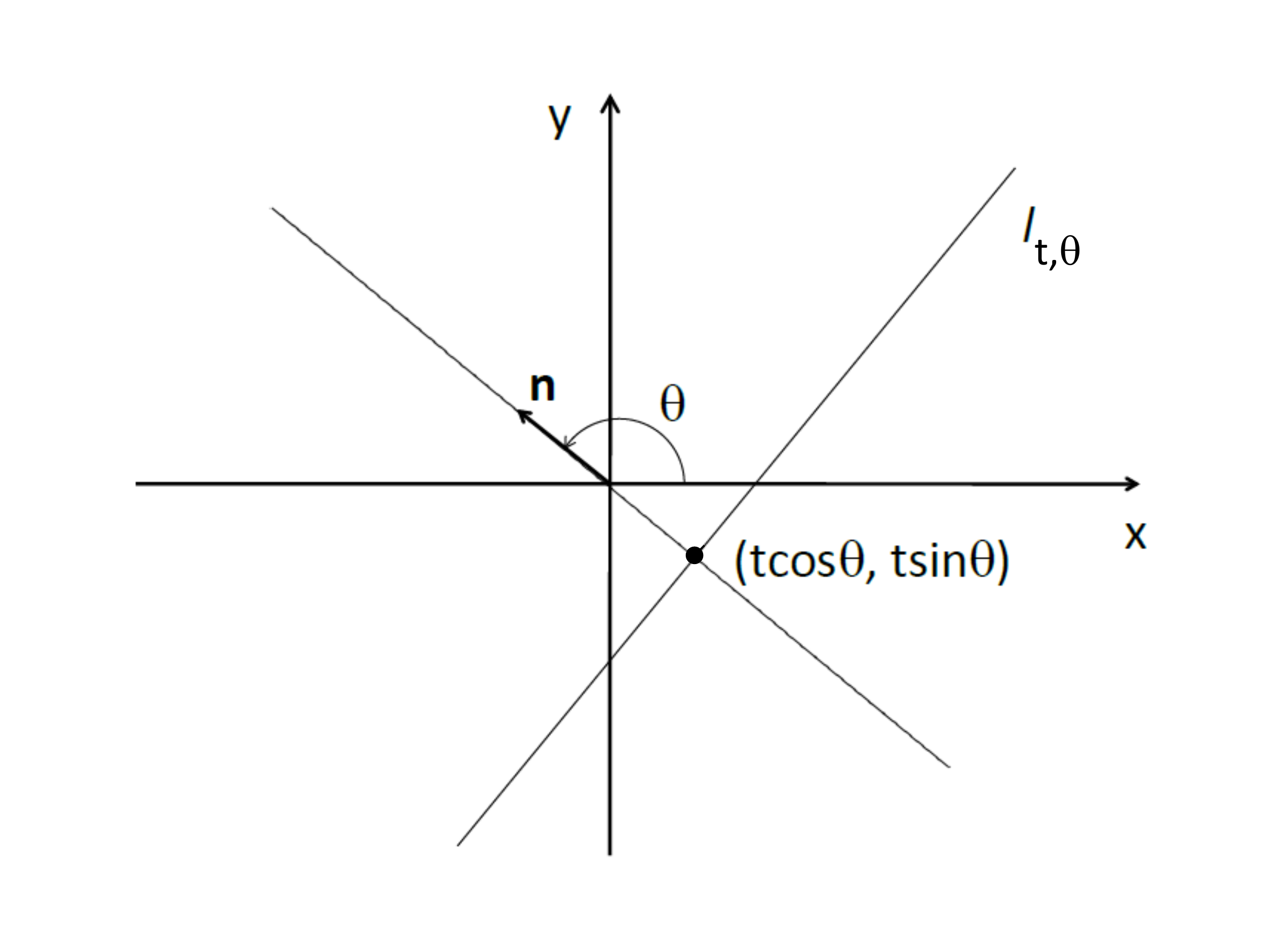} %[scale=1]width=0.7\textwidth width=1\textwidth, height=1\textwidth
\caption{A line in the plane can be characterized by two real numbers $t$, $\theta$}
\label{fig: charac_line}
\end{figure}

The approach we adopt is characterizing every line in the plane by a point that the line passes through and a normal vector to the line. Then, let $\vec{\textsf{n}}$ be a vector that is normal to a given line $l$, then there exists some angle $\theta$ such that $\vec{\textsf{n}}$ is parallel to the line radiating out from the origin at an angle $\theta$ measured counterclockwise from the positive $x$-axis (Figure \ref{fig: charac_line}). This line is also perpendicular to $l$ and thus intersects $l$ at some point whose coordinates in the plane have the form $(t\cos{\theta},t\sin{\theta})$ for some  real number $t$. The line $l$ is hence characterized by the values of $t$ and $\theta$ and so we denote $l=l_{t,\theta}$.
\begin{definition}
For any real numbers $t$ and $\theta$, the line $l_{t,\theta}$ is the line passing through the point  $(t\cos{\theta},t\sin{\theta})$ and perpendicular to the vector  $\vec{\textsf{n}}=(\cos{\theta},\sin{\theta})$.
\end{definition}

Because of the relationships $l_{t,\theta+2\pi}=l_{t,\theta}$ and $l_{t,\theta+\pi}=l_{-t,\theta}$ for all $t,\theta$, there is not a unique representation of the form $l_{t,\theta}$ for a line. For this reason we will consider only the set of lines
\begin{equation*}
\{l_{t,\theta}\ : \ t\in\R, \, 0\leq\theta<\pi\}.
\end{equation*}
If we consider the unit vector $(-\sin{\theta},\cos{\theta})$, perpendicular to $\vec{\textsf{n}}$, every point on $l_{t,\theta}$ can be written as 
\begin{equation*}
(t\cos{\theta},t\sin{\theta})+s(-\sin{\theta},\cos{\theta}),
\end{equation*}
for some number $s\in\R$. So we can parametrize a line $l_{t,\theta}$ as $(x(s),y(s))$, where $s\in\R$ and
\begin{equation*}
\begin{cases}
&x(s)=t\cos{\theta}-s\sin{\theta}\\
&y(s)=t\sin{\theta}+s\cos{\theta}
\end{cases}
\end{equation*}
Note that for every point $(x(s),y(s))\in l_{t,\theta}$, we have $x(s)^{2}+y(s)^{2}=t^{2}+s^{2}$.

With this parametrization the arclenght element along the line $l_{t,\theta}$ is given by
\begin{equation*}
\sqrt{\left( \frac{dx}{ds}\right)^{2}+\left( \frac{dy}{ds}\right)^{2}}ds=\sqrt{(-\sin{\theta})^{2}+(\cos{\theta})^{2}}ds=ds
\end{equation*}
Therefore for a given function $A(x,y)$ defined in the plane, we get
\begin{equation}
\int_{l_{t,\theta}}{A(x,y)}=\int_{\R}{A(t\cos{\theta}-s\sin{\theta},t\sin{\theta}+s\cos{\theta})\,ds}.
\label{eq: int_line}
\end{equation}
The value of this integral is exactly what an X-ray emission/detection machine measures when an X-ray is emitted along the line $l_{t,\theta}$.

Finally note that for an arbitrary point $(x_{0},y_{0})$ in the plane and for a given value $\theta$, there is a unique value of $t$ such that $(x_{0},y_{0})\in l_{t,\theta}$. The value of $t$ is given by the solution of the system
\begin{equation*}
\begin{cases}
&x_{0}=t\cos{\theta}-s\sin{\theta}\\
&y_{0}=t\sin{\theta}+s\cos{\theta},
\end{cases}
\end{equation*}
that is
\begin{equation*}
\begin{cases}
&t=x_{0}\cos{\theta}+y_{0}\sin{\theta}\\
&s=-x_{0}\sin{\theta}+y_{0}\cos{\theta}.
\end{cases}
\end{equation*}
This formula will be used in the next sections to operate some change of variables that will be used in finding an inversion formula of the Radon transform.

\section{The Radon transform}
\subsection{Definition and basic properties}
The fundamental question of image reconstruction is: is it possible to reconstruct a function $f$, representing the attenuation coefficient of a cross section of a sample, starting from the value of the integral of $f$ along \emph{every} line $l_{t,\theta}$ in the plane?

 We will consider the integral of $f$ for any values of $t$ and $\theta$, in other words, given a function $f$ we associate to every point $(t,\theta)$ a number representing the value of the integral $\int_{l_{t,\theta}}{f}$. This leads us to the definition of the Radon transform:

\begin{definition} 
For a given function $f:\R^{2}\rightarrow\R$, the Radon transform of $f$ is defined by
\begin{equation*}
Rf(t,\theta)=\int_{l_{t,\theta}}{f\,ds}=\int_{\R}{f(t\cos{\theta}-s\sin{\theta},t\sin{\theta}+s\cos{\theta})\,ds}, 
\end{equation*}
$\forall\,t\in\R,\,\theta\in[0,\pi).$
\end{definition}

So the Radon transform is an operator that, to a given function $f$ of the Cartesian coordinates $(x,y)$, associates a function $Rf$ of the polar coordinates $(t,\theta)$.

\begin{example}
\label{es: radonCirc}
%\subsubsection{Example}
%\begin{enumerate}
%\item
Consider a circle of radius $r>0$ and a function $r_{r}$ defined as follows:
\begin{equation*}
f_{r}(x,y)=\left\{
\begin{aligned}
&1 \qquad \text{if} \ x^{2}+y^{2}\leq r^{2}\\
&0 \qquad \text{otherwise},
\end{aligned}
\right.
\end{equation*}
since $x^{2}+y^{2}=t^{2}+s^{2}$, the Radon transform of $f$ is 
\begin{equation*}
Rf(t,\theta)=\int_{\R}{f(t\cos{\theta}-s\sin{\theta},t\sin{\theta}+s\cos{\theta})\,ds}=\int_{t^{2}+s^{2}\leq r^{2}}{1\,ds},
\end{equation*}
so we get %\int_{-\sqrt{r^{2}-t^{2}}}^{\sqrt{r^{2}-t^{2}}}{1\,ds}=
\begin{equation}
\label{eq: radonCirc}
Rf(t,\theta)=\left\{
\begin{aligned}
&2\sqrt{r^{2}-t^{2}}  \quad &\text{if} \ t\leq |r|\\
&0 \quad  &\text{if} \ t>|r|,
\end{aligned}
\right.
\end{equation}
\end{example}
%\item The Radon transform of a Gaussian $f(x,y)=\frac{1}{\sqrt{2\pi\sigma^{2}}}e^{-\frac{x^{2}+y^{2}}{2\sigma^{2}}}$ is
%\begin{equation*}
%Rf(t,\theta)=\frac{1}{\sqrt{1+t^{2}}}e^{-\frac{\theta^{2}}{2(i+t^{2})}\sigma^{2}}
%\end{equation*}
%\end{enumerate}

\begin{proposition}
The Radon transform is a linear operator: for two functions $f$ and $g$ and constants $\alpha$ and $\beta$,
\begin{equation*}
R(\alpha f+\beta b)=\alpha Rf+\beta Rg.
\end{equation*}
\end{proposition}
\begin{proof}
It follows from the linearity of the integral.
\end{proof}

\begin{example}
Consider the function
\begin{equation*}
f(x,y)=\left\{
\begin{aligned}
&\frac{1}{2} \qquad \text{if} \ x^{2}+y^{2}\leq r_{1}^{2}\\
&1  \qquad \text{if} \ r_{1}^{2}<x^{2}+y^{2}\leq r_{2}^{2}\\
&0 \qquad \text{otherwise},
\end{aligned}
\right.
\end{equation*}
with $0<r_{1}<r_{2}$. We observe that $f$ can be rewritten as $f=f_{r_{2}}-\frac{1}{2}f_{r_{1}}$, where $f_{r_{1}}$ and $f_{r_{2}}$ are defined as in Example \ref{es: radonCirc}, by equation \eqref{eq: radonCirc} and the linearity of the Radon transform, we get
\begin{equation*}
Rf(t,\theta)=Rf_{r_{2}}(t,\theta)-\frac{1}{2}Rf_{r_{1}}(t,\theta)
=\left\{
\begin{aligned}
&2\sqrt{r_{2}^{2}-t^{2}}-\sqrt{r_{1}^{2}-t^{2}} \ &\text{if} \ |t|\leq r_{1}\\
&2\sqrt{r_{2}^{2}-t^{2}} \ &\text{if} \ r_{1}<|t|\leq r_{2}\\
&0 \ &\text{if} \ |t|>r_{2}.
\end{aligned}
\right.
\end{equation*}
\end{example}

\subsubsection{Domain of Radon transform}
As we see from the definition, the Radon transform is a linear operator acting on functions and it involves improper integral on lines that can be infinity for some function. It is then natural to ask ourself for what kind of function is defined the Radon transform and in particular which is the domain of this operator, i.e. which space of functions is composed of all and only the functions that admit finite Radon transform. It can be proved (see \cite{SIGU}) that the space we are looking for is the Schwartz space 
\begin{equation*}
\mathcal{S}=\left\{f:\ \sup_{x}{\left|\ |x|^{m}P(\partial_{1},\partial_{2})f(x)\right|<\infty}, \ \forall\,m\in\N, P\ \text{polynomial} \right\}.
\end{equation*}
of rapidly decreasing functions, but for the moment we can not consider this problem since the functions involved in medical imaging correspond to attenuation coefficient of finite size samples and therefore are compact supported. In Chapter \ref{chap: kernelMethods} we will face the problem of how to compute Radon transforms of functions that do not belong to $\mathcal{S}$ and we will find there some expedient to overcome this obstacle.

\subsection{Back projection}
Our aim is to recover a function $f$, representing the attenuation-coefficient of a sample, from the values of its Radon transform $Rf$. 

We start by considering a point $(x_{0},y_{0})$ in the plane. For any values of $\theta$ there exists one and only one value of $t$ such that the line $l_{t,\theta}$ passes through the point $(x_{0},y_{0})$. In particular, the value of $t$ is $t=x_{0}\cos{\theta}+y_{0}\sin{\theta}$. 

In practice, any X-ray beam passing through a point $(x_{0},y_{0})$ follows the line $l_{(x_{0}\cos{\theta}+y_{0}\sin{\theta}),\theta}$ for some angle $\theta$. So the Radon transform $Rf ( x_{0}\cos{\theta}+y_{0}\sin{\theta},\theta )$ takes into account the value of the attenuation coefficient $f(x_{0},y_{0})$.

The first way one can try to recover $f(x_{0},y_{0})$ is to compute the average of the Radon transform along all lines passing through $(x_{0},y_{0})$, that is
\begin{equation*}
\frac{1}{\pi}\int_{0}^{\pi}{Rf{(x_{0}\cos{\theta}+y_{0}\sin{\theta},\theta)\,d\theta}}
\end{equation*}  

This leads us to the definition of the following transform, called \emph{back projection}:
\begin{definition}
Let $h=h(t,\theta)$ a function in polar coordinates. The \emph{back projection} of $h$ at the point $(x,y)$ is given by
\begin{equation*}
Bh(x,y)=\frac{1}{\pi}\int_{0}^{\pi}{h{(x\cos{\theta}+y\sin{\theta},\theta)\,d\theta}}
\end{equation*}  
\end{definition}
Back projection is a linear transform:

\begin{proposition}
The back projection is a linear transform, i.e. for all functions $h_{1}$ and $h_{2}$ and for all constants $c_{1}$ and $c_{2}$, we have
\begin{equation*}
B(c_{1}h_{1}+c_{2}h_{2})=c_{1}Bh_{1}+c_{2}Bh_{2}
\end{equation*}   
\end{proposition}

We observe that the back projection 
\begin{equation*}
BRf(x,y)=\frac{1}{\pi}\int_{0}^{\pi}{Rf{(x\cos{\theta}+y\sin{\theta},\theta)\,d\theta}}
\end{equation*}
of the Radon transform, does not give us the value of $f(x,y)$. Indeed, the value $Rf(x\cos{\theta}+y\sin{\theta},\theta)$ represents the total accumulation of the attenuation-coefficient $f$ along a particular line. The integral $BRf$ is computing the average values of those averages. Hence it gives us a smoothed version of $f$.

\begin{example} \label{ex: backPro}
Consider $f_{1}$ a function corresponding to a disc of radius $1/2$ centered at the origin with constant density 1, that is
\begin{equation*}
f_{1}(x,y)=\left\{
\begin{aligned}
&1 \quad \text{if} \ x^{2}+y^{2}<\frac{1}{4}\\
&0 \quad \text{otherwise}.
\end{aligned}
\right.
\end{equation*} 
Then, for each line passing through the origin, we have $Rf_{1}(0,\theta)=1$ and consequently $BRf_{1}(0,0)=1$.

Now suppose $f_{2}$ be defined by
\begin{equation*}
f_{2}(x,y)=\left\{
\begin{aligned}
&1 \quad \text{if} \ \frac{1}{4}<x^{2}+y^{2}<\frac{3}{4}\\
&0 \quad \text{otherwise}.
\end{aligned}
\right.
\end{equation*} 
 Again, for every line $l_{0,\theta}$ passing through the origin, we have $Rf_{2}(0,\theta)=1$ and $BRf_{2}(0,0)=1$.

Thus $BRf_{1}(0,0)=BRf_{2}(0,0)=1$, but $f_{1}(0,0)=1$ and $f_{2}(0,0)=0$. This shows the fact that the back projection of the Radon transform does not necessarily reproduce the original function.
\end{example}

\section{The Filtered Back-Projection Formula}
In this section we will discuss the relationships between the Radon transform, the back projection and the Fourier transform. Thanks to these formulas we will obtain the inversion of the Radon transform. In other words we will be able to get the values of a function $f$, representing for example an X-ray attenuation coefficient, starting form the values of its Radon transform.

In the next paragraphs we will consider successive transforms of a function, for example in the \emph{central slice theorem} in section \ref{subsec: centralSlice} we will consider the Fourier transform of the Radon transform. In all these cases we will assume that all the transforms are well defined, i.e. we will assume that a function $f$ belongs to the Schwartz space of rapidly decreasing functions $\mathcal{S}$. For example one can think to $f$ as a compact supported function.

\subsection{The Central Slice Theorem}\label{subsec: centralSlice}
The interaction between the Radon transform and the Fourier transform is given by the \emph{Central Slice Theorem}, also known as the \emph{Central Projection Theorem}.

We recall that the $n$-dimensional Fourier transform of a function $f:\R^{n}:\rightarrow\R$ is defined as 
\begin{equation*}
(F_{n}f(x))(\omega)=\int_{\R^{n}}{f(x)e^{-ix\cdot\omega}\,dx} \quad \forall\,\omega\in\R^{n},
\end{equation*}
where $i$ denotes the imaginary unit and $x\cdot\omega$ the standard inner product in $\R^{n}$. For a function in polar coordinates $f(t,\theta)$, we consider the 1-dimensional Fourier transform $F=F_{1}$, applied only to the variable $t$, i.e. 
\begin{equation*}
(Ff(t,\theta))(\omega)\int_{\R}{f(t)e^{-it\omega}\,dt}, \quad \omega\in\R.
\end{equation*}

We can now state the Central Slice Theorem in the case $n=2$:
\begin{theorem}[The Central Slice Theorem]
For a function $f$ defined in the plane and for all real numbers $r$, $\theta$,
\begin{equation*}
F_{2}f(r\cos{\theta},r\sin{\theta})=F(Rf)(r,\theta).
\end{equation*} 
\end{theorem}
\proof
The definition of the Fourier transform gives
\begin{equation}
F_{2}f(r\cos{\theta},r\sin{\theta})=\int_{-\infty}^{+\infty}{\int_{-\infty}^{+\infty}{f(x,y)e^{-ir(x\cos{\theta}+y\sin{\theta})}\,dx}\,dy}
\label{eq: centralSlice1}
\end{equation} 
Consider now the change of variables 
\begin{equation*}
\begin{cases}
&x=t\cos{\theta}-s\sin{\theta}\\
&y=t\sin{\theta}+s\cos{\theta}
\end{cases}
\qquad
\begin{cases}
&t=x\cos{\theta}+y\sin{\theta}\\
&s=-x\sin{\theta}+y\cos{\theta}.
\end{cases}
\end{equation*}
Note that the quantity $t=x\cos{\theta}+y\sin{\theta}$ is exactly the line $l_{t,\theta}$. Moreover $dxdy=dtds$, indeed
\begin{equation*}
\left|
\begin{array}{cc}
\frac{\partial x}{\partial t} & \frac{\partial x}{\partial s} \\
\frac{\partial y}{\partial t} & \frac{\partial y}{\partial s}
\end{array}
\right|=
\left|
\begin{array}{cc}
\cos{\theta} & -\sin{\theta} \\
\sin{\theta} & \cos{\theta}
\end{array}
\right|=
\cos^{2}{\theta}+\sin^{2}{\theta}=1.
\end{equation*}
The integral in \eqref{eq: centralSlice1} becomes then
\begin{align*}
&\int_{-\infty}^{+\infty}{\int_{-\infty}^{+\infty}{f(t\cos{\theta}-s\sin{\theta},t\sin{\theta}+s\cos{\theta})e^{-irt}\,ds}\,dt}=\\
&=\int_{-\infty}^{+\infty}{\left(\int_{-\infty}^{+\infty}{f(t\cos{\theta}-s\sin{\theta},t\sin{\theta}+s\cos{\theta})\,ds}\right)e^{-irt}\,dt},
\end{align*} 
where we have factored out the inner integral since the term $e^{-irt}$ does not depends on $s$. Now the inner integral in the last equation is exactly the definition of the Radon transform of the function $f$ evaluated at point $(t,\theta)$. Thus the last integral equals
\begin{equation*}
\int_{-\infty}^{+\infty}{Rf(t,\theta)e^{-irt}\,dt},
\end{equation*}
that is the definition of the 1-dimensional Fourier transform of $Rf$ at the point $(r,\theta)$.

In conclusion
\begin{equation*}
F_{2}f(r\cos{\theta},r\sin{\theta})=F(Rf)(r,\theta).
\end{equation*}
\endproof

\subsection{The Filtered Back-Projection}
Applying the back projection to the Radon transform gives a smoothed version of the original function. The following theorem, called \emph{Filtered Back-Projection Formula}, shows how to correct the smoothing effect and recover the original function.

\begin{theorem}[The Filtered Back-Projection Formula]\label{thm: filteredBP}
For all function $f$ and for all real number $x$,$y$,
\begin{equation}
f(x,y)=\frac{1}{2}B\{F^{-1}[|r|F(Rf(r,\theta))]\}(x,y).
\label{eq: filteredBP}
\end{equation}
\end{theorem}
\proof By the Fourier inversion theorem, for any function $f$ and any point in the plane $(x,y)$, we have
\begin{equation*}
f(x,y)=F_{2}^{-1}F_{2}f(x,y).
\end{equation*}
Applying the definition we have
\begin{equation*}
f(x,y)=\frac{1}{4\pi^{2}}\int_{-\infty}^{+\infty}{\int_{-\infty}^{+\infty}{F_{2}f(X,Y)e^{i(xX+yY)}\,dX}\,dY}.
\end{equation*}
We pass now from Cartesian coordinates $(X,Y)$ to polar coordinates $(r,\theta)$, where $X=r\cos{\theta}$ and $Y=r\sin{\theta}$, with $r\in\R$ and $\theta\in[0,\pi]$. Because of this change of coordinates in the integral, we have $dXdY=|r|drd\theta$ and so
\begin{equation*}
f(x,y)=\frac{1}{4\pi^{2}}\int_{0}^{\pi}{\int_{-\infty}^{+\infty}{F_{2}f(r\cos{\theta},r\sin{\theta})e^{ir(x\cos{\theta}+y\sin{\theta})}|r|\,dr}\,d\theta}.
\end{equation*}
Applying the central slice theorem to the factor $F_{2}f(r\cos{\theta},r\sin{\theta})=F(Rf(r,\theta))$, we get
\begin{equation*}
f(x,y)=\frac{1}{4\pi^{2}}\int_{0}^{\pi}{\int_{-\infty}^{+\infty}{F(Rf)(r,\theta)e^{ir(x\cos{\theta}+y\sin{\theta})}|r|\,dr}\,d\theta}.
\end{equation*}
In the last equation, the inner integral is by definition, $2\pi$ times the inverse Fourier transform of the function $|r|F(Rf)(r,\theta)$, evaluated at the point $(x\cos{\theta}+y\sin{\theta},\theta)$. So we can write
\begin{equation*}
f(x,y)=\frac{1}{2\pi}\int_{0}^{\pi}{F^{-1}[|r|F(Rf)(r,\theta)](x\cos{\theta}+y\sin{\theta},\theta)\,d\theta},
\end{equation*}
that is half of the back projection of the function $F^{-1}[|r|F(Rf)(r,\theta)]$. Hence we finally obtain the desired formula
\begin{equation*}
f(x,y)=\frac{1}{2}B\{F^{-1}[|r|F(Rf)(r,\theta)]\}(x,y).
\end{equation*}
\endproof

Observe that the factor $|r|$ in the formula \eqref{eq: filteredBP} is fundamental. Indeed without this factor, the Fourier transform and its inverse, would cancel out and the result would be simply the back projection of the Radon transform of $f$, that as shown in example \ref{ex: backPro} does not lead to recover $f$.

The Filtered Back-Projection formula is the basis for image reconstruction. However it assumes that the values of $Rf(t,\theta)$ are known for all possible values $(t,\theta)$. In practice only a finite number of X-ray samples are taken and we must approximate an image from the resulting data. 

\section{Filtering} \label{sec: filter}
Consider the Filtered Back-Projection formula in \eqref{eq: filteredBP}:
\begin{equation*}
f(x,y)=\frac{1}{2}B\{F^{-1}[|r|F(Rf(r,\theta))]\}(x,y)
\end{equation*}
and suppose there exists a function $\phi(t)$ such that $F\phi(r)=|r|$. In this case we could write
\begin{equation*}
|r|F(Rf)(r,\theta)=[F\phi\cdot F(Rf)](r,\theta).
\end{equation*} 
By the properties of the Fourier transform we would have
\begin{equation*}
|r|F(Rf)(r,\theta)=F(\phi\ast Rf)(r,\theta),
\end{equation*}
hence
\begin{align*}
F^{-1}[|r|F(Rf)(r,\theta)]&=F^{-1}[F(\phi\ast Rf)(r,\theta)]=\\
&=(\phi\ast Rf).
\end{align*}
We could then write equation \eqref{eq: filteredBP} as 
\begin{equation}
f(x,y)=\frac{1}{2}B(\phi\ast Rf)(x,y).
\label{eq: modFilteredBP}
\end{equation}
In this way the formula of the reconstruction of $f$ would be simpler. The problem is that such a function $\phi$ does not exist. However, the previous discussion will be useful if we consider data $Rf$ to be affected of noise, that is the case when we have to work with real data from the X-ray machine.

Consider the function $|r|F(Rf)(r,\theta)$. The variable $r$ represent a frequency that is present in a signal, so if the Radon transform has a component at high frequency, this component is magnified by the factor $|r|$. Since noise has high frequency, that means that the noise present in the image is amplified and this effect corrupt the reconstructed image.

In order to obtain a formula less sensitive to noise, instead of $|r|$ we use a function, actually a \emph{low-pass filter}, such that for $r$ close to 0, it is near to the absolute-value function $|r|$, but vanishes if the value of $|r|$ is large.
Moreover, in order to use the formula \eqref{eq: modFilteredBP}, in place of $|r|$ we consider a function of the form $A=F\phi$, where $A$ has compact support, or in other words, we consider $\phi$ \emph{band-limited function}.

In this way we obtain an approximation of $f$:
\begin{equation*}
f(x,y)\approx\frac{1}{2}B(F^{-1}A\ast Rf)(x,y).
\end{equation*}

Typically the function $A$ is of the form $A(\omega)=|\omega|F(\omega)\chi_{[-L,L]}(\omega)$, for some $L>0$, where $\chi_{I}$ represents the characteristic function of the set $I$. The function $F$ is even and $F(0)=1$ in order to have an approximation of the function $|\cdot|$ near the origin and $\phi$ real valued.

\begin{figure}[htbp]
\centering
\includegraphics[width=0.7\textwidth]{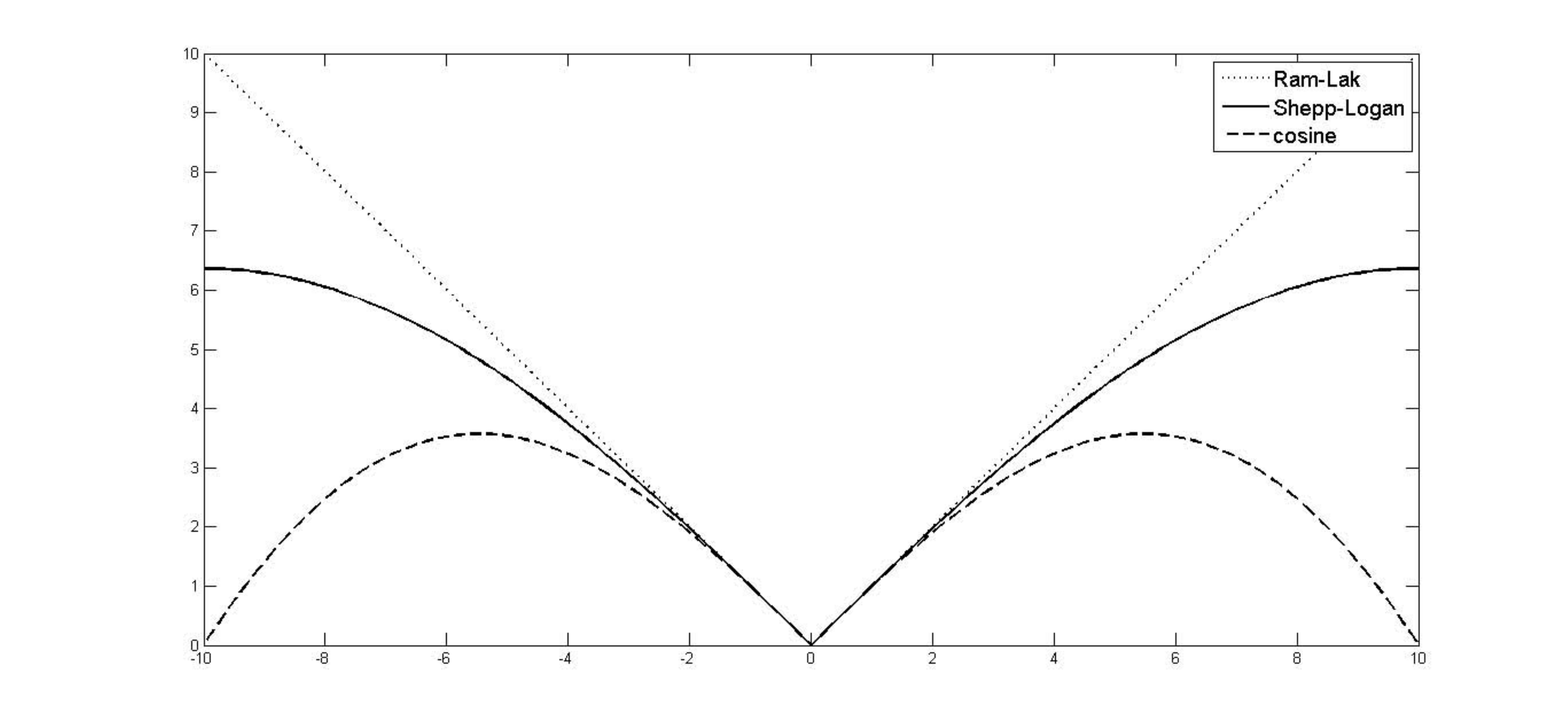} 
\caption{Main low pass filters used in medical imaging}
\label{fig: lowpass_filter}
\end{figure}

%\subsection{Low-pass filter}\label{subsec: lowPassFilter}
Typical low-pass filters used in medical imaging are:
\begin{itemize}
\item The \emph{Ram-Lak filter}:
\begin{equation*}
A_{1}(\omega)=|\omega|\chi_{[-L,L]}(\omega),
\end{equation*}
is simply the truncation of the absolute-value function to a finite interval. 
\item The \emph{Shepp-Logan filter}:
\begin{equation*}
\begin{aligned}
A_{3}(\omega)&=|\omega|\left(\frac{\sin(\pi\omega/(2L))}{\pi\omega/(2L)}\right)\chi_{[-L,L]}(\omega)=\\
&=\left\{
\begin{aligned}
&\frac{2L}{\pi}|\sin(\pi\omega/(2L))| & &\text{if} \ |\omega|\leq L\\
&0 & &\text{otherwise}.
\end{aligned}
\right.
\end{aligned}
\end{equation*}
\item The \emph{low-pass cosine filter}:
\begin{equation*}
A_{2}(\omega)=|\omega|\cos(\pi\omega/(2L))\chi_{[-L,L]}.
\end{equation*}
\end{itemize}
The plot of these filters in the case $L=10$ is shown in Figure \ref{fig: lowpass_filter}.

\subsection{Filter resolution}
Consider a function $\phi$, suppose $\phi\geq0$, with a single maximum value $M$ in $x=0$ and increasing for $x<0$, decreasing for $x>0$ (for example $\phi$ can be a Gaussian). For another function $f$, the filtered version of $f$ with $\phi$ is given by $f\ast\phi$.

Let now the numbers $x_{1}, x_{2}$ be such that $x_{1}<0<x_{2}$ and $\phi(x_{1})=\phi(x_{2})=M/2$, half of the maximum value of $\phi$. The distance $x_{2}-x_{1}$ is called \emph{full width half maximum} of the function $\phi$, in symbol $FWHM(\phi)$. The resolution of the filter defined by the convolution with $\phi$ is set to be equal to $FWHM(\phi)$.
\begin{figure}[htbp]
\centering
\includegraphics[width=0.7\textwidth]{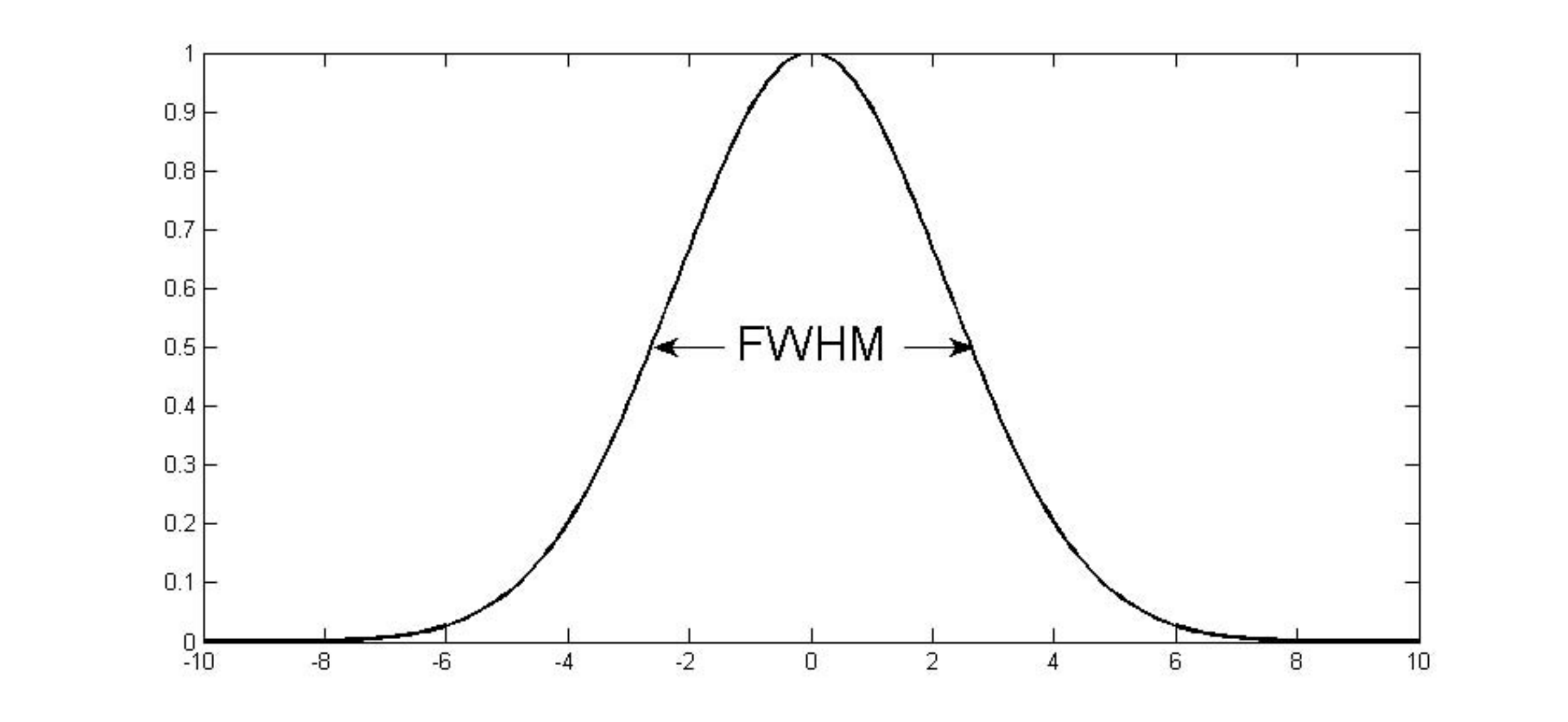} 
\caption{Full width half maximum of a Gaussian}
\label{fig: fwhm}
\end{figure}

To understand the reason of this definition, consider a function $f$ that consists of two unit impulses separated by a distance of $d$. It is easy to show that if $d>FWHM(\phi)$, then the graph of $f\ast\phi$ has two peaks, but if $d\leq FWHM(\phi)$, then the graph of $f\ast\phi$ has only one peak and so we have lost of details. So we conclude that the smallest distance between two different features of $f$ that can still be seen in the filtered signal $f\ast\phi$ is $FWHM(\phi)$.
One should choose the filter function $\phi$ in accordance with the resolution required.
Intuitively we can think that a function $\phi$ with a small $FWHM$ is spiker than a function having large $FWHM$ and has better resolution. The following examples can help to understand better

\begin{example}[$FWHM$ of a Gaussian]
Let $F(\omega)=e^{-B\omega^{2}}$, where $\omega\in\R$ and $B$ is a positive constant. The maximum value of $F$ is $F(0)=1$. Half maximum is hence achieved for $e^{-B\omega^{2}}=1/2$ i.e. $\omega=\pm\sqrt{\ln(2)/B}$. therefore
\begin{equation*}
FWHM=2\sqrt{\ln(2)/B}
\end{equation*}
\end{example}
\begin{example}[$FWHM$ of a the Lorentz signal]
The Lorentz signal is given by
\begin{equation*}
g(\omega)=\frac{T_{2}}{1+4\pi^{2}T_{2}^ {2}(\omega-\omega_{0})^{2}},
\end{equation*}
where $\omega\in\R$ and $T_{2},\omega_{0}$ are constants. The maximum of $g$ is given by $g(\omega_{0})=T_{2}$ and $g(\omega)=T_{2}/2$ if and only if $\omega=\pm1/(2\pi T_{2})$. Therefore
\begin{equation*}
FWHM=\frac{1}{\pi T_{2}}
\end{equation*}
 \end{example}
\begin{figure}[htbp]
\centering
\includegraphics[width=0.7\textwidth]{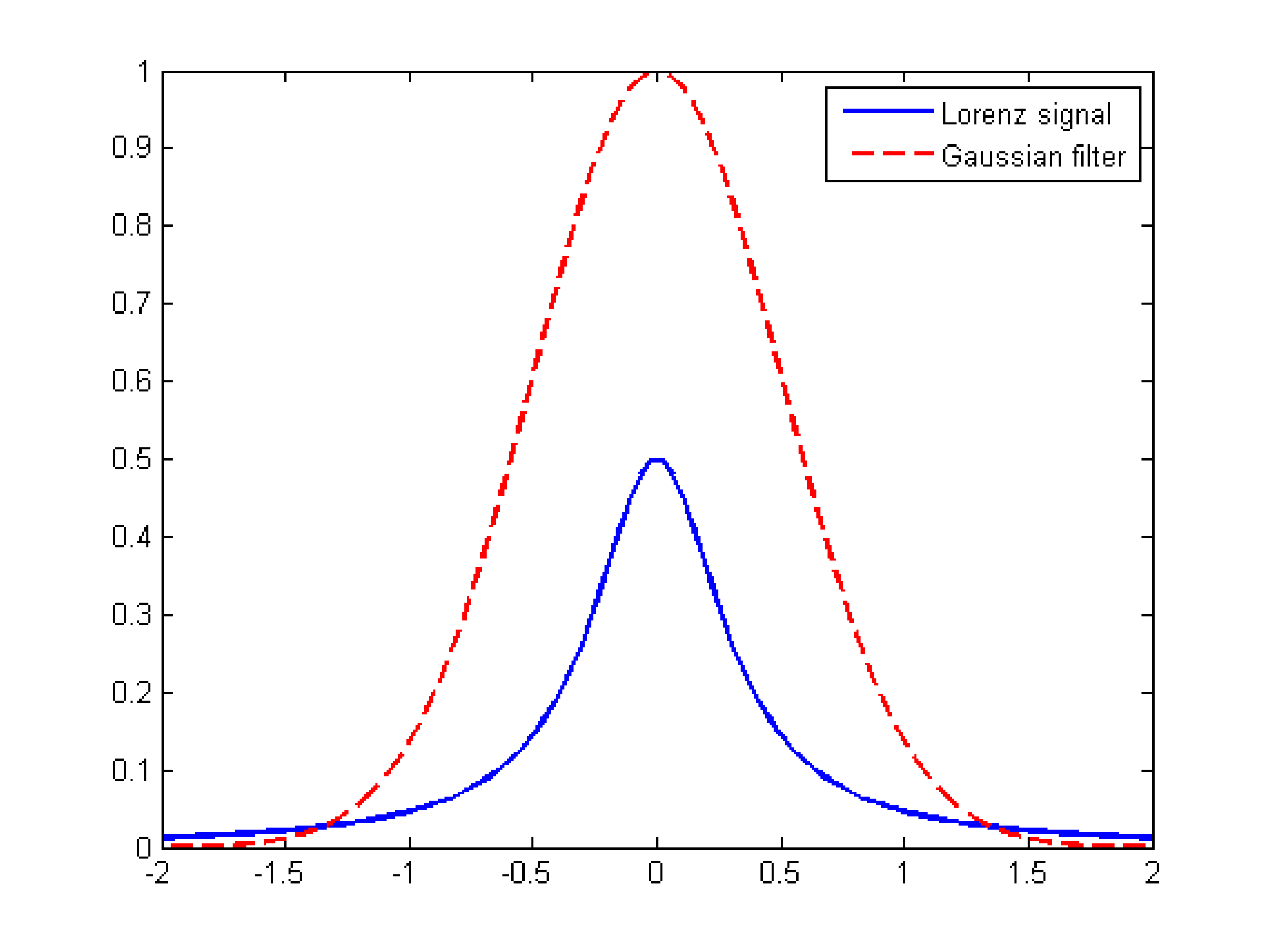} 
\caption{Gaussian filter ($B=2,\ FWHM=1.1774$) and Lorenz signal ($T_{2}=0.5,\ \omega_{0}=0,\ FWHM=0.6366$).}
\label{fig: example_fwhm}
\end{figure}

\section{Discrete problem}
By Theorem \ref{thm: filteredBP} we know that if completed continuous data are available, then we can exactly reconstruct a function $f$ starting from its Radon transform. In particular this is possible thanks to the back projection formula \eqref{eq: filteredBP}
\begin{equation*}
f(x,y)=\frac{1}{2}B\{F^{-1}[|r|F(Rf)(r,\theta)]\}(x,y).
\end{equation*} 

We have also seen, in Section \ref{sec: filter}, that in practice is convenient to replace the absolute-value function with a low-pass filter $A$. Thus, we may use the approximation
\begin{equation}
f(x,y)\approx\frac{1}{2}B[F^{-1}(A\ast Rf)](x,y).
\label{eq: approxBP2}
\end{equation}

In the practical implementation of this formula, we have to consider that only a finite number of values of $Rf(r,\theta)$ are measured by the X-ray machine. As a consequence of this fact we have to answer to some question about accuracy and computation.  
First of all we have to understand the sampling process, i.e. the process of computing only a discrete set of value of a continuous function; then we have to find the corresponding form of formula \eqref{eq: approxBP2} for discrete functions; and finally we will use the process of interpolation to obtain value of the function we can not directly measure. 

\subsection{Phantoms} \label{sec: phantoms}
Different choices of filters, interpolation methods, and other parameters, will give us different reconstruction of the same image, thus we need a technique for testing the accuracy of one particular image reconstruction algorithm.

In order to have a good accuracy test, we should know the original image we want to reconstruct. Moreover the method should be independent from the possible noise present in the data, but should depend only on the algorithm used in the reconstruction. To solve this problem, Shepp and Logan (\cite{SHEP}) introduced the concept of \emph{mathematical phantom}.

A mathematical phantom (or simply a phantom) is a simulated object whose structure is defined by mathematical formulas. Thus no errors occur in collecting the data from the object and when an algorithm is applied to produce a reconstruction of the phantom, all inaccuracies are due to the algorithm. In this way we can compare different algorithms meaningfully.

\begin{figure}[htbp]
\centering
\includegraphics[width=\textwidth]{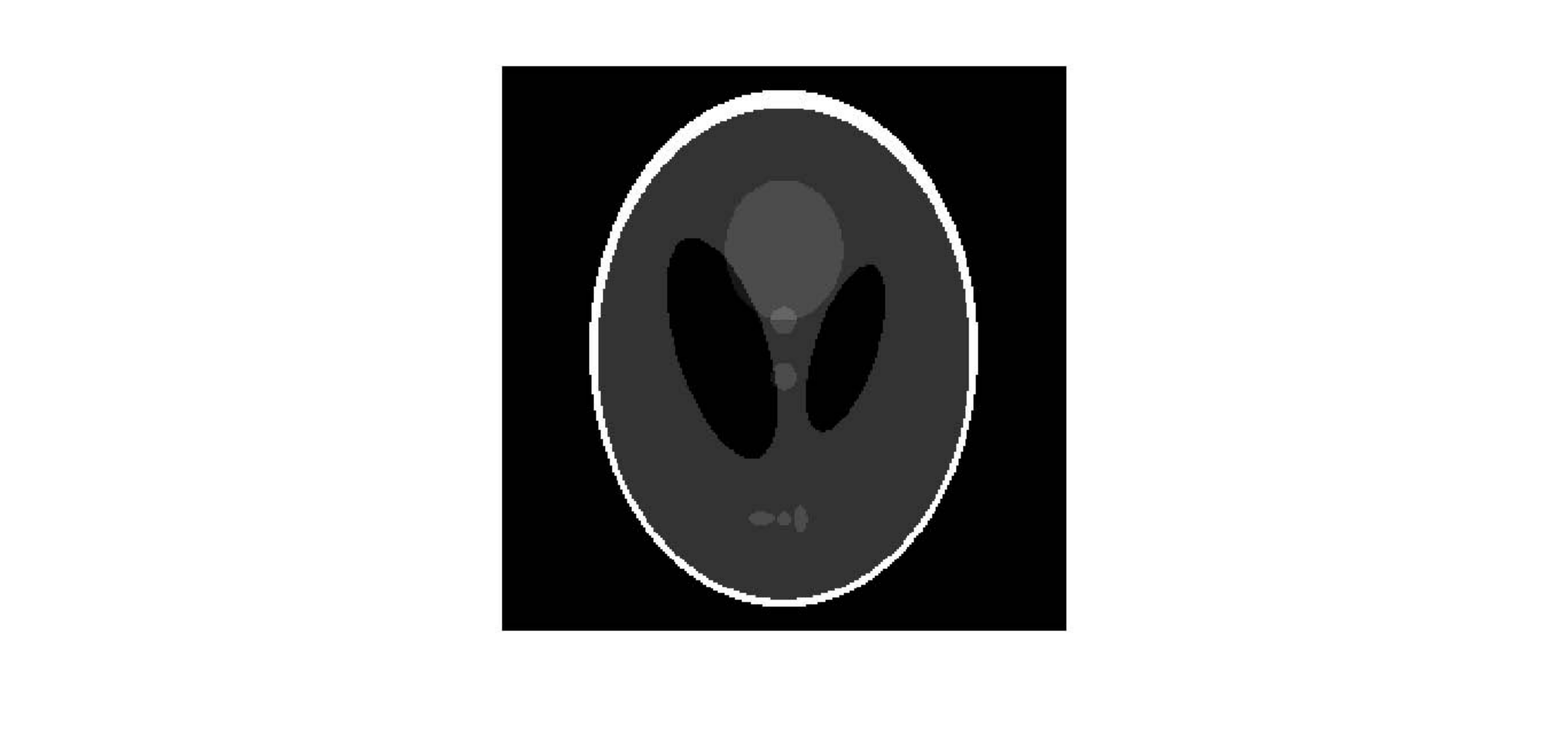} 
\caption{The Shepp-Logan phantom}
\label{fig: shepp_phantom}
\end{figure}

Figure \ref{fig: shepp_phantom} shows the well-known \emph{Shepp-Logan phantom}. This phantom is widely used to test the quality of an image reconstruction algorithm since it is a good imitation of the human brain. 

\subsection{Sampling}
Sampling is the process of computing the values of a function, or a signal defined in $\R$, only on a discrete set of points $\{x_{k}\}_{k\in\Z}$. For example, points $x_{k}$ can be taken with uniformly spacing, i.e. $x_{k}=k\cdot d$ for some positive number $d$, called the sampling spacing. The sampling spacing $d$ determines the smallest detail of $f$ that can be seen after sampling: if $d$ is small we have a better resolution, while bigger values of $d$ give us less resolution. On the other hand, small values of $d$ generate a bigger amount of data and make algorithm slower, so we want to find an optimal values of $d$ which is a compromise for this trade-off.

If we think to a signal as a sum of sinusoidal waves, the narrowest detail in the signal is given by the wave with the shortest wavelength (maximum frequency). If the signal is band limited, the Nyquist Theorem \ref{thm: nyquist} below tells us that the signal can be completely recovered starting from its sampled version, provided that the sampling spacing is small enough.

Suppose $f$ band limited, i.e. its Fourier transform is zero outside a finite interval: $Ff(\omega)=0$ for $|\omega|>L$. If we extend $Ff$ periodically out of $[-L,L]$, its Fourier coefficients are given by
\begin{equation*}
c_{n}=\frac{1}{2L}\int_{-L}^{L}{Ff(\omega)e^{-i\omega n\frac{\pi}{L}}\,d\omega}, \qquad n\in\Z,
\end{equation*} 
thus
\begin{align*}
2\pi f\left(n\frac{\pi}{L}\right)&=2\pi F^{-1}Ff(n\frac{\pi}{L})=\int_{\R}{Ff(\omega)e^{i\omega n\frac{\pi}{L}}\,d\omega}=\\
&=\int_{-L}^{L}{Ff(\omega)e^{i\omega n\frac{\pi}{L}}\,d\omega}=2Lc_{-n}.
\end{align*} 
Assuming $Ff$ continuous we have
\begin{equation*}
Ff(\omega)=\sum_{n=-\infty}^{\infty}{c_{-n}e^{-i\omega n\frac{\pi}{L}}}=\frac{\pi}{L}\sum_{n=-\infty}^{\infty}f\left(n\frac{\pi}{L}\right)e^{-i\omega n\frac{\pi}{L}}
\end{equation*}
and so
\begin{equation*}
f(x)=F^{-1}Ff(x)=\sum_{n=-\infty}^{\infty}f\left(\pi\frac{n}{L}\right)\frac{\sin{(Lx-n\pi)}}{Lx-n\pi}
\end{equation*}
that is $f$ can be exactly reconstructed from the values $f(n\pi/L)$, $n\in\Z$. The optimal sampling spacing is therefore $d=\frac{\pi}{L}$, since $L$ is the maximum value of $|\omega|$ in $Ff$, the smallest wavelength is $\frac{2\pi}{L}$, hence the optimal sampling distance is equal to half the size of the smallest detail present in the signal. This result is resumed in the following
\begin{theorem}[Nyquist Theorem]\label{thm: nyquist}
If $f$ is a square integrable and band limited function, i.e. $Ff(\omega)=0$ for all $|\omega|>L$, then for all $x\in\R$
\begin{equation}
f(x)=\sum_{n=-\infty}^{\infty}f\left(n\frac{\pi}{L}\right)\frac{\sin{(Lx-n\pi)}}{Lx-n\pi}.
\label{eq: nyquist}
\end{equation}
\end{theorem}

We observe that formula \eqref{eq: nyquist} involves an infinite series and that its general term $\sin{(Lx-n\pi)}/(Lx-n\pi)$ converges slowly. So we need a large number of samples for a good approximation. To address this, we can use a smaller sampling distance $\frac{\pi}{R},\ R>L$ to gain a series with better convergence. This process is called \emph{oversampling}.

\subsection{Discrete filters}
\begin{figure}[htbp]
\centering%
\subfigure[Ram-Lak filter \label{fig: ift_rlf}]%
{\includegraphics[width=0.30\textwidth]{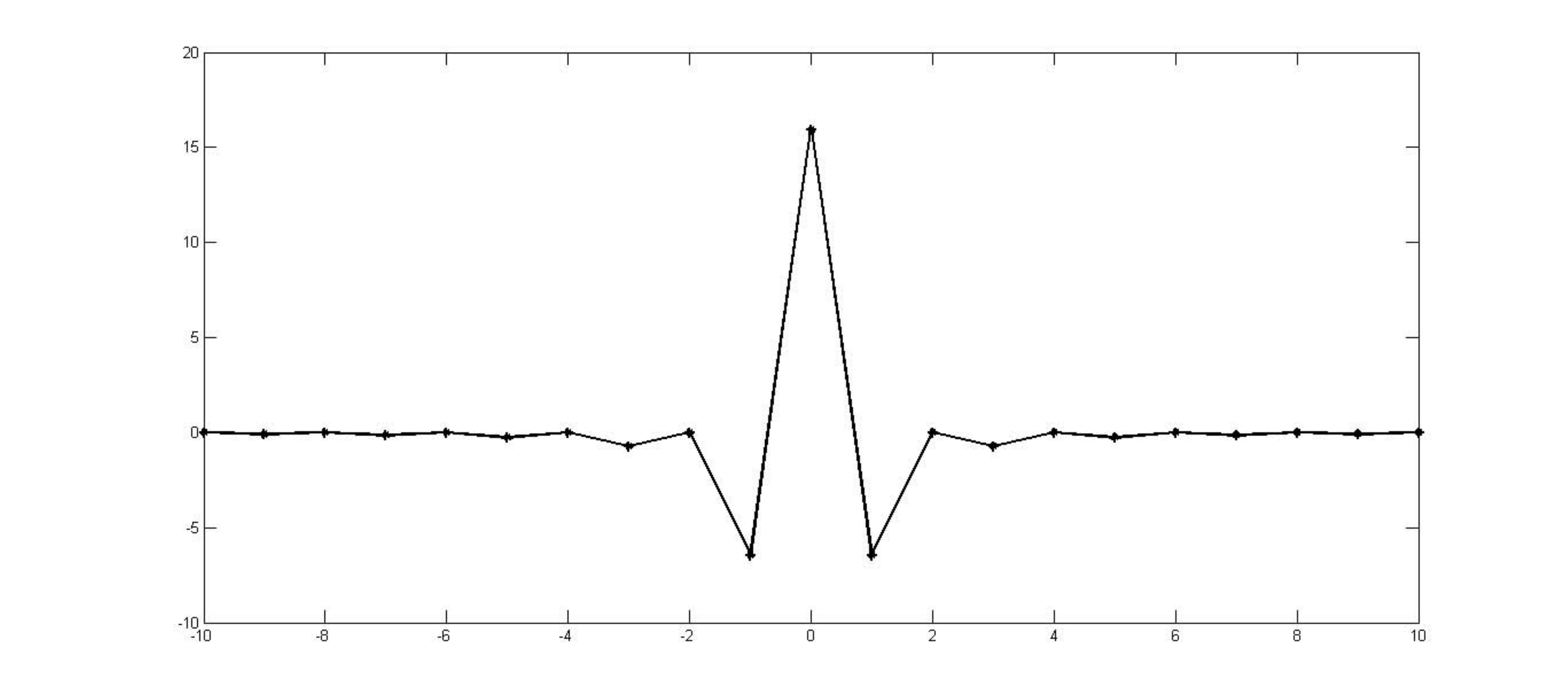}}
\subfigure[Shepp-Logan filter \label{fig: ift_slf}]%
{\includegraphics[width=0.30\textwidth]{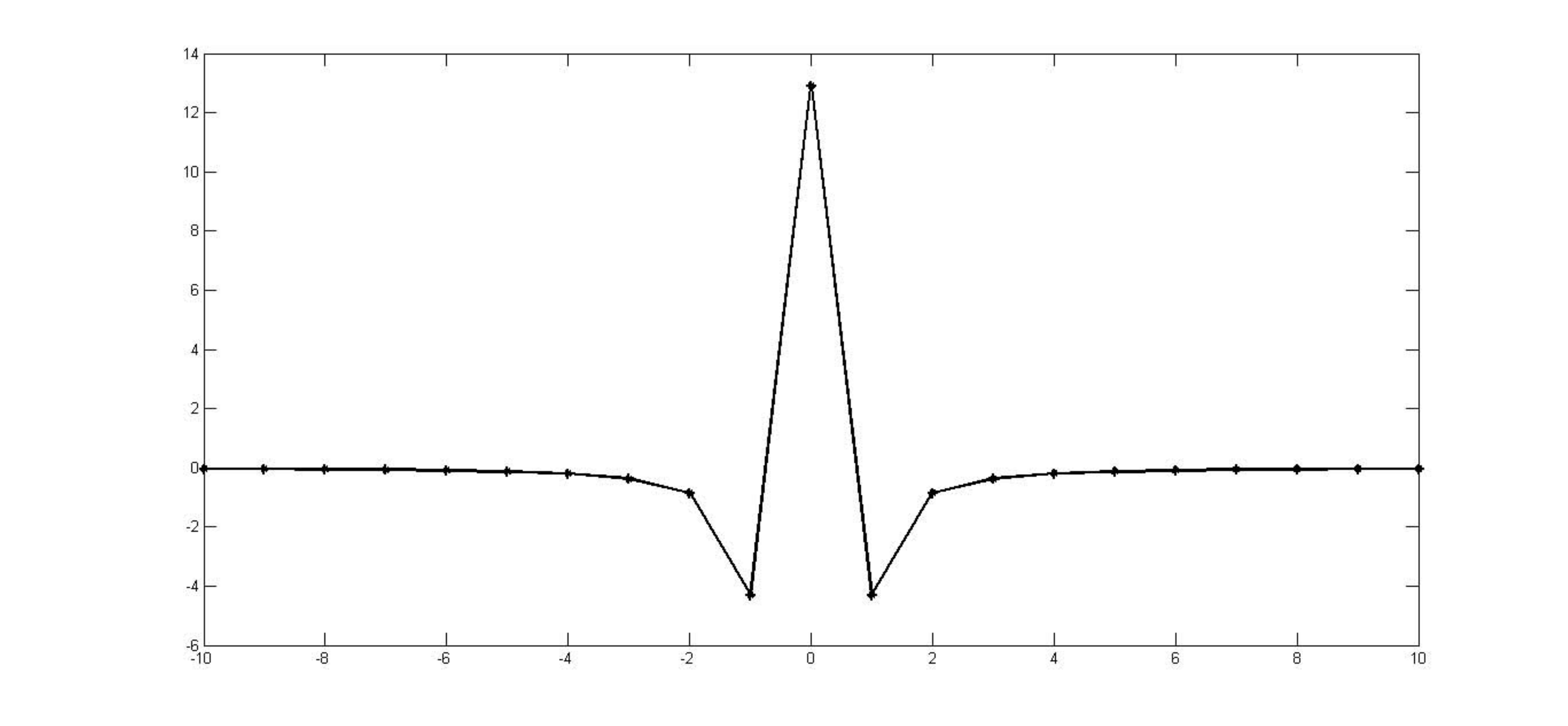}}
\subfigure[Cosine filter \label{fig: ift_cos}]%
{\includegraphics[width=0.30\textwidth]{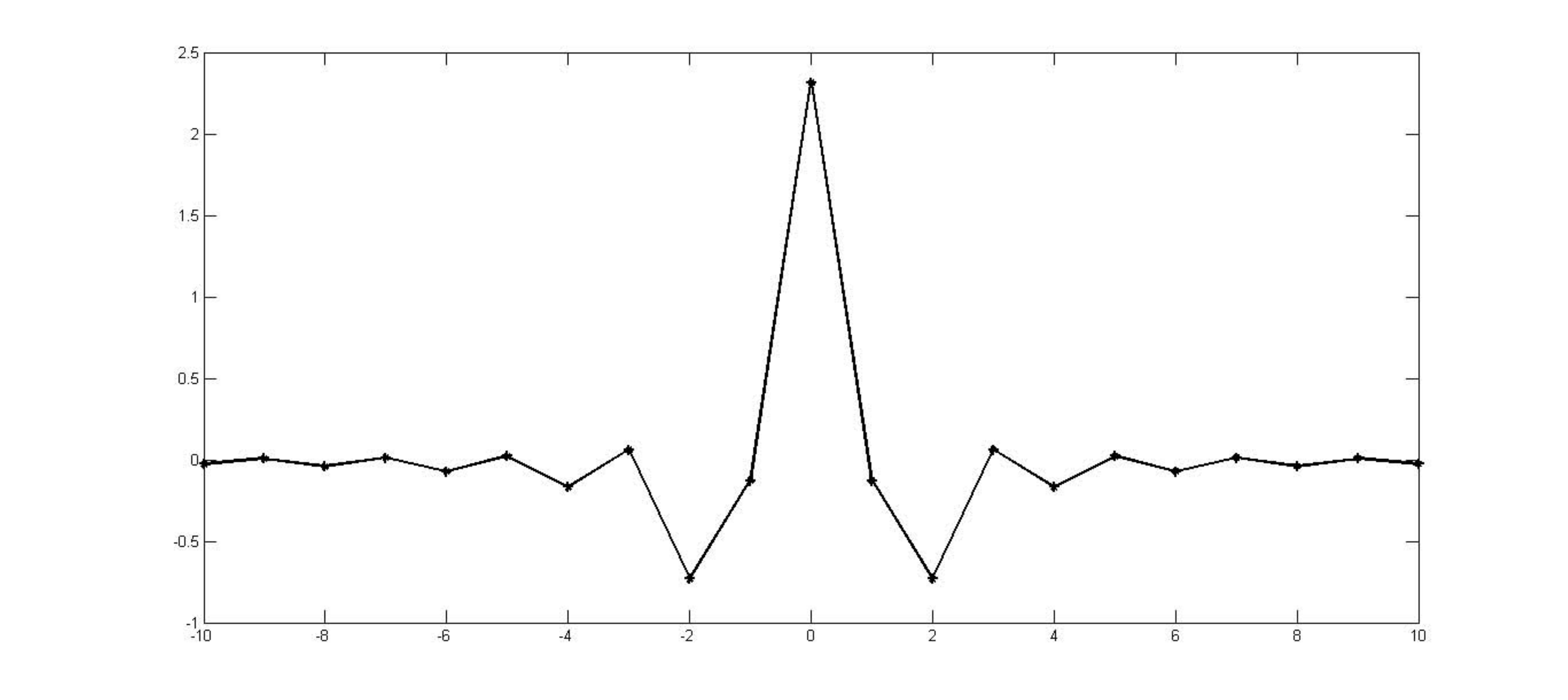}}
\caption{Sampled inverse Fourier transform}
\label{fig: ift}
\end{figure}

The image reconstruction formula \eqref{eq: approxBP2} involves the inverse Fourier transform of a low pass filter. In practice also this function will be sampled like the Radon transform. Since the filters we consider are band limited, we use the Nyquist theorem \ref{thm: nyquist} to know how many samples are needed to get an accurate representation of the filter. Here, we reconsider the filters introduced in section \ref{sec: filter}:
\begin{itemize}
\item The \emph{Shepp-Logan filter} is defined by
\begin{equation*}
\begin{aligned}
A_{3}(\omega)&=|\omega|\left(\frac{\sin(\pi\omega/(2L))}{\pi\omega/(2L)}\right)\chi_{[-L,L]}(\omega)=\\
&\left\{
\begin{aligned}
&\frac{2L}{\pi}|\sin(\pi\omega/(2L))| \qquad \text{if} \ |\omega|\leq L\\
&0 \qquad \text{otherwise}.
\end{aligned}
\right.
\end{aligned}
\end{equation*}
for some $L>0$. The inverse Fourier transform of $A_{3}$ is a band limited function and is given by 
\begin{align*}
F^{-1}A_{3}(x)&=\frac{1}{\pi}\int_{0}^{L}{\frac{2L}{\pi}\sin{(\pi\omega/(2L))}\cos{\omega}\,d\omega}=\\
&=\frac{L}{\pi^{2}}\left[\left(\frac{\cos{(Lx-\pi/2)}}{x-\pi/(2L)}-\frac{\cos{(Lx+\pi/2L)}}{x+\pi/(2L)}
 \right)- \right.\\
&\left.\left(\frac{1}{x-\pi/(2L)}-\frac{1}{x+\pi/(2L)}\right)
\right].
\end{align*}
According to Nyquist theorem, $F^{-1}A_{3}$ can be reconstructed exactly from its values taken at distance $\pi/L$. Setting $x=n\pi/L$, for $n\in\Z$, we get
\begin{equation*}
F^{-1}A_{3}(\pi n/L)=\frac{4L^{2}}{\pi^{3}(1-4m^{2})}.
\end{equation*}
\item The \emph{Ram-Lak filter} is given by
\begin{equation*}
A_{1}(\omega)=|\omega|\chi_{[-L,L]}(\omega).
\end{equation*}
Proceeding as in the previous case we find that the inverse Fourier transform of the Ram-Lak filter satisfies
\begin{equation*}
F^{-1}A_{1}(x)=\frac{1}{\pi}\left[\frac{Lx\sin{(Lx)}}{x^{2}}-\frac{2\sin^{2}{(Lx/2)}}{x^{2}} \right].
\end{equation*}
Setting again $x=\pi n/L$ we obtain
\begin{equation*}
F^{-1}A_{1}(\pi n/L)=\frac{L^{2}}{2\pi}\left[\frac{2\sin{(\pi n)}}{\pi n}-\left(\frac{\sin{(\pi n/2)}}{\pi n/2}\right)^{2} \right].
\end{equation*}
\item Finally we consider the  \emph{low-pass cosine filter}:
\begin{equation*}
A_{2}(\omega)=|\omega|\cos(\pi\omega/(2L))\chi_{[-L,L]}.
\end{equation*}
The inverse Fourier transform of $A_{2}$, evaluated at multiples of the Nyquist distance is
\begin{equation*}
F^{-1}A_{2}(\pi n/L)=\frac{2L^{2}}{\pi^{2}}\left[\frac{\pi\cos{(\pi n)}}{1-4n^{2}}-\frac{2(1+4n^{2})}{(1-4n^{2})^{2}} \right].
\end{equation*}
\end{itemize}
Figure \ref{fig: ift} shows the sampled version of the inverse Fourier transform of these filters.

\subsection{Discrete functions}
\subsubsection{Discrete convolution}
In order to implement formula \eqref{eq: approxBP2} we have to decide what is convolution of discrete functions.

A discrete function is a mapping from the integers into the set of real numbers. For a discrete function $g$, we write $g_{n}$ for $g(n)$, for all $n\in\Z$.
\begin{definition} 
The discrete convolution of two discrete functions $f$ and $g$ is defined by
\begin{equation*}
(f\ast g)_{m}=\sum_{j=-\infty}^{+\infty}{f_{j}g_{m-j}} \quad \forall m\in\Z.
\end{equation*}
\end{definition}
The discrete convolution satisfies all principal properties of the standard convolution (e.g. commutativity and linearity).

If only a finite set of values $\{f_{k}=f(dk): k=0,\ldots,N-1\}$ is known, like in real applications, there exist two different ways to extend the sequence to all integers:
\begin{enumerate}
\item Set $f_{k}=0$ for all $k\notin\{0,\ldots,N-1\}$;
\item Extend the sequence to be periodic with period $N$, $f_{m}=f_{m+nN}$, where for $m\in\Z$, $n$ is the only integer such that $m+nN\in\{0,\ldots,N-1\}$. We call such a function \emph{$N$-periodic discrete function}.
\end{enumerate} 

The convolution of two $N$-periodic discrete functions is also a $N$-periodic discrete function, defined by
\begin{equation*}
(f\ast g)_{m}=\sum_{j=0}^{N-1}{f_{j}g_{m-j}} \quad \forall m\in\Z.
\end{equation*}
 
Some problem can arise using discrete functions. For example if we are sampling a non periodic function, the periodic model is not the best to be used. But, even if the function is periodic, it may be not clear what the appropriate period is and so we might sample the function on a set of values that do not correspond to one period, then extending data to form a discrete periodic function, we have the wrong one. The solution to these problems is a technique called \emph{zero padding}. We take a finite set of values of a function $g$, then we pad the sequence with a lot of zeros and finally we form a periodic discrete function.

The following theorem tells us that the convolution between a zero padded function and another discrete function gives the same result of "true" discrete convolution at least at the points where the values has been sampled.

\begin{theorem} \label{thm: discrete_convolution}
Let $f,g$ be discrete functions and suppose that $\exists K\in\N$ such that $g_{k}=0$ for $k<0$ and $k\geq K$. Let $M\in\Z$, $M\geq K-1$ and let $\tilde{f},\tilde{g}$ $(2M+1)$-periodic discrete functions defined by $\tilde{f}_{m}=f_{m}$, $\tilde{g}_{m}=g_{m}$ for $-M\leq m\leq M$. Then for all $m$ such that $0\leq m \leq K-1$ we have $(f\ast g)_{m}=(\tilde{f}\ast\tilde{g})_{m}$.
\end{theorem}
%\proof
\begin{remark} 
The proof of this theorem is just an application of the definition of convolution for discrete functions. For details we refer the reader to \cite{BASIC}. 
%\endproof
\end{remark}

\subsubsection{Discrete Radon transform}
\begin{figure}[htbp]
\centering
\includegraphics[width=0.7\textwidth]{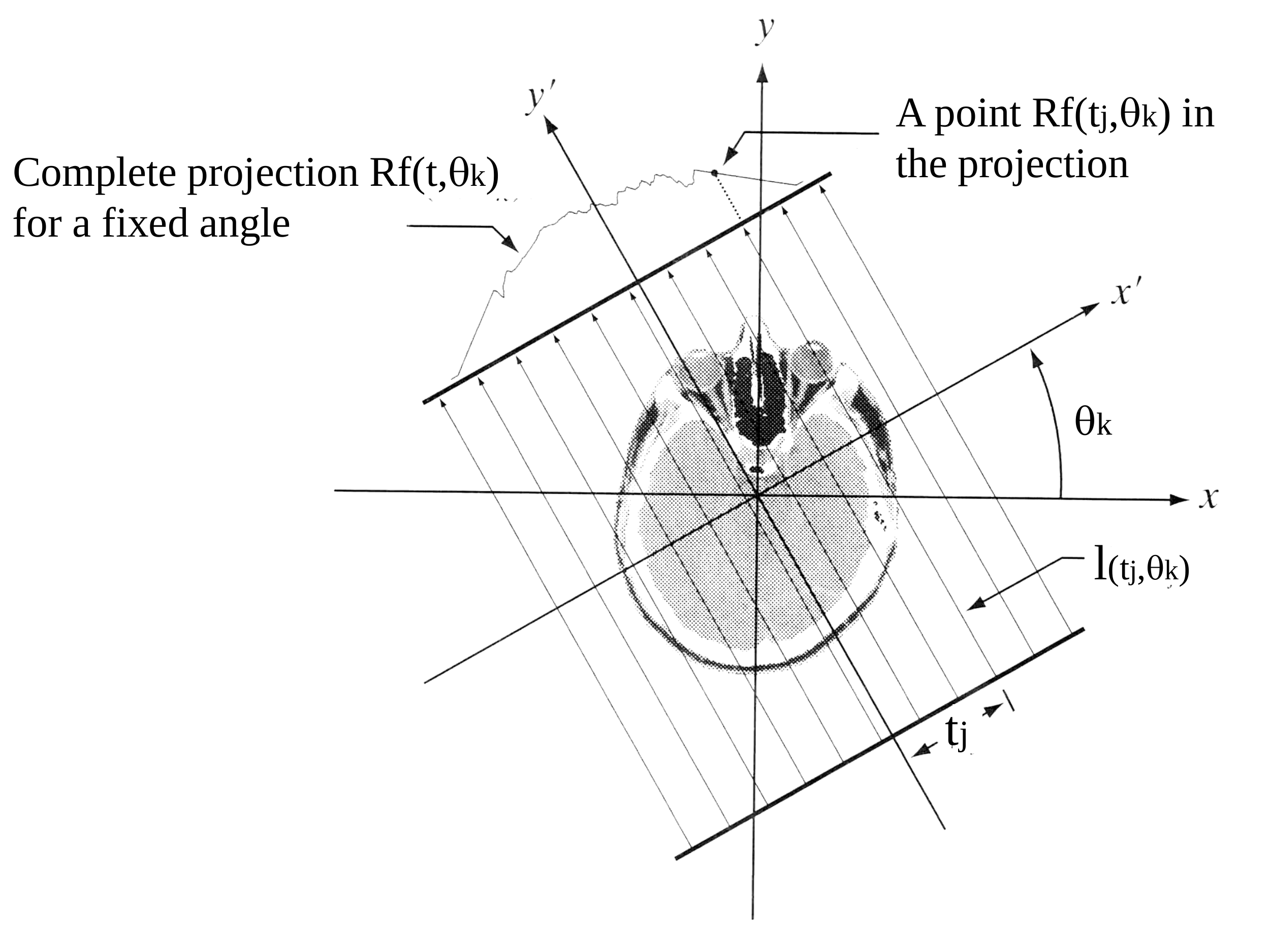} 
\caption{The parallel beam geometry}
\label{fig: parallel_beam}
\end{figure}

In the context of a CT scan, the X-ray machine does not access the attenuation coefficient along every line, but the Radon transform is sampled for finite number of angles $\theta\in[0,\pi)$ and, for each angle, for a finite number of values of $t$. Values of $\theta$ and $t$ are equally spaced and we consider the \emph{parallel beam geometry}: the X-ray machine rotates by a fixed angle and, at each angle, the beams form a set of parallel lines (Figure \ref{fig: parallel_beam}).

If $N$ is the number of angles at which the machine takes scans, then the values of $\theta$ that occur are $\{ k\frac{\pi}{N}, \ k=1,\ldots,N-1\}$. Assume that, at each angle, the set of parallel beams is composed of $2M+1$ equally spaced lines and let $d$ be the distance between two lines, with the object to be scanned centered at the origin. Then the corresponding values of $t$ are $\{jd: \ j=-M,\ldots,M\}$.

The continuous Radon transform $Rf$ is then replaced by its discrete counterpart $R_{D}f$, defined by
\begin{equation*}
R_{D}f_{j,k}=Rf(jd,k\pi/N)
\end{equation*}
for $j=-M,\ldots,M$ and $k=0,\ldots,N-1$.

Theorem \ref{thm: discrete_convolution} above applies to the discrete convolution of the sampled band-limited function $F^{-1}A$ and the sampled Radon transform $R_{D}f$. Since the scanned object has finite size, we can set $R_{D}f(j,\theta)=0$  for $|j|$ sufficiently large. Thus with enough zero padding the discrete Radon transform can be extended to be periodic in the radial variable $jd$.

For discrete function in polar coordinates, the discrete convolution is carried out in the radial variable only, so in the reconstruction formula \eqref{eq: approxBP2} we have:
\begin{equation*}
(F^{-1}A\ast R_{D}f)_{m,\theta}=\sum_{j=0}^{N-1}{F^{-1}A_{j}R_{D}}f_{m-j,\theta}.
\end{equation*}

\subsubsection{Discrete Fourier transform}
\begin{definition}
The \emph{discrete Fourier transform} of a $N$-periodic discrete function is another $N$-periodic discrete function $F_{D}f$ defined by
\begin{equation*} 
(F_{D}f)_{j}=\sum_{k=0}^{N-1}{f_{k}e^{-i2\pi k j/N}}, \quad \text{for}\ j=0,\ldots,N-1
\end{equation*}
and extended to be periodic for other values of $j$.

The \emph{discrete inverse Fourier transform} of $f$ is given by 
\begin{equation*} 
(F^{-1}_{D}f)_{j}=\frac{1}{N}\sum_{k=0}^{N-1}{f_{k}e^{i2\pi k j/N}}, \quad \text{for}\ j=0,\ldots,N-1
\end{equation*}
and extended to be periodic for other values of $j$.
\end{definition}
The following theorems show that the properties of the Fourier transform are still valid for its discrete version.
\begin{theorem}
For a discrete function $f$ with period $N$,
\begin{equation*}
F^{-1}_{D}(F_{D}f)_{n}=f_{n}, \quad \text{for all integers $n$}. 
\end{equation*}
\end{theorem}
\begin{theorem}
For two $N$ discrete functions $f=\{f_{k}:\ 0\leq k\leq N-1\}$ and $g=\{g_{k}:\ 0\leq k\leq N-1\}$, we have
\begin{itemize}
\item $F_{D}(f\ast g)=(F_{D}f)(F_{D}g)$;
\item $F_{D}(fg)=\frac{1}{N}(F_{D}f)\ast (F_{D}g)$;
\item $(F_{D}\bar{f})_{j}=(\overline{F_{D}f})_{-j}$;
\item Parceval equality:
\begin{equation*} 
\sum_{j=0}^{N-1}{|f_{j}|^{2}}=\frac{1}{N}\sum_{j=0}^{N-1}{|(F_{D}f)_{j}|^{2}}
\end{equation*}
\end{itemize}
\end{theorem}
For a proof of these facts we suggest to see \cite{BASIC}.

\subsubsection{Discrete back projection}
In the continuous setting the back projection has been defined by 
\begin{equation*}
Bh(x,y)=\frac{1}{\pi}\int_{0}^{\pi}{h(x\cos{\theta}+y\sin{\theta})\,d\theta}.
\end{equation*}
Now, in the discrete case, we replace the continuous variable $\theta$ with angles $k\pi/N$ for $k=0,\ldots,N-1$ and state the following 
\begin{definition}
The \emph{discrete back projection} of a function $h$ is defined by 
\begin{equation*}
B_{D}h(x,y)=\frac{1}{N}\sum_{k=0}^{N-1}{h(x\cos{k\frac{\pi}{N}}+y\sin{k\frac{\pi}{N}},k\pi/N)}.
\end{equation*}
\end{definition}

In our case, $B_{D}$ has to be applied to $h=(F^{-1}_{D}A)\ast (R_{D}f)$ and the reconstruction grid within which the final image is to be presented is a rectangular array of pixels located at $(x_{m},y_{n})$, each of which is to be assigned a color or a gray-scale value. Hence $B_{D}$ needs the values of $h$ at points $(x_{m}\cos{k\pi/N}+y_{n}\sin{k\pi/N},k\pi/N)$, while the Radon transform is sampled at points $(jd,k\pi/N)$ arranged in a polar grid. The solution to this problem is \emph{interpolation}.

\section{Interpolation}
The process to obtain a function $f(x)$, $x\in\R$ starting form a discrete set of values $f_{k}=f(x_{k})$, $k=1,\ldots,N+1$ is called interpolation. There exist several interpolation schemes. Here we give a short introduction to the most commonly used.
\begin{itemize}
\item \textbf{Nearest neighbor}: $f(x)=f_{k}$, where $x_{k}$ is the closest point to $x$. This is the simplest method but generates a discontinuous function;
\item \textbf{Linear}: $f$ is obtained connecting successive points $(x_{k},f_{k}),\ (x_{k+1},f_{k+1})$ with segment:
\begin{equation*}
f(x)=\frac{f_{k+1}-f_{k}}{x_{k+1}-x_{k}}(x-x_{k})+f_{k} \quad \text{for}\ x\in[x_{k},x_{k+1}];
\end{equation*} 
\item \textbf{Cubic polynomial spline}: successive points $(x_{k},f_{k}),\ (x_{k+1},f_{k+1})$ are connected by apiece of a cubic polynomial. The pieces are joint together asking for $C^{2}$ continuity of the resulting curve. Also values of $f'(x_{k})$ are prescribed;
\item \textbf{Lagrange interpolation}: $f$ is given by a polynomial of degree $N$:
\begin{equation*}
f(x)=\sum_{j=1}^{N+1}{f_{j}\frac{\prod_{k\neq j}{(x-x_{k})}}{\prod_{k\neq j}{(x_{j}-x_{k})}}}.
\end{equation*}
\end{itemize}

We notice that the nearest neighbor interpolation can be written as
\begin{equation*}
I_{f}(x)=\sum_{m}f_m\chi_{[-\frac{1}{2},\frac{1}{2})}\left(\frac{x}{d}-m\right),
\end{equation*}
where $\chi_{J}$ denotes the characteristic function of a set $J$ and $f_m$ is the value of the function $f$ at the sample point $md$. Similarly the linear interpolation can be written $I_{f}(x)=\sum_{m}f_m\Lambda\left(\frac{x}{d}-m\right)$, with
\begin{equation*}
\Lambda(x)=
\begin{cases}
1-|x| &\quad\text{if}\ |x|\leq1\\
0 &\quad \text{if}\ |x|>1.
\end{cases}
\end{equation*}
Generalizing this approach we define, for a weighting function $W$ satisfying certain conditions, the $W$-interpolation $I_{W}(f)$ of a discrete function $f$ is
\begin{equation*}
I_{W}(f)=\sum_{m}{f_mW\left(\frac{x}{d}-m\right)} \quad \ x\in\R.
\end{equation*}
We want $I_{W}f(kd)=f_k$. Then we choose $W$ such that $W(0)=1$ and $W(m)=0$ for all $m\in\Z$, $m\neq0$. Moreover, if we want to preserve also the integral, we ask $W$ to be such that
\begin{equation*}
\int_{\R}{I_{W}(f)(x)\,dx}=d\sum_{m}{f_{m}}.
\end{equation*}
Then, $W$ should satisfy
\begin{equation*}
\int_{\R}{{W}(u)\,du}=1.
\end{equation*}

\begin{remark}[Interpolation and convolution]\label{rk: intconv}
Suppose that a discrete function $g$ is given by the discrete convolution $g=\phi\ast f$ and let $W$ be a weighting function. Then the $W$-interpolation
\begin{equation*}
I_{W}g(x)=I_{W}(\phi\ast f)(x)\sum_{k}{\sum_{m}{\phi(m-k)W\left(\frac{x-kd}{d}-(m-k)\right)f(k)}}
\end{equation*}
\end{remark}
can be approximated as 
\begin{equation*}
I_{W}(\phi\ast f)(x)\approx \sum_{k}{I_{W}(\phi)(x-kd)f(k)},
\end{equation*}
that is, we can approximate the interpolation of $\phi\ast f$ with a weighted sum of values $f(k)$ and the interpolation $I_{W}(\phi)$ of $\phi$ at points $x-kd$ (cfr. \cite{BASIC}, pages 82-86). 

\section{Discrete image reconstruction: Algorithms}
Having examined the discrete version of all elements in the formula \eqref{eq: approxBP2}, we have now all the necessary tools for approximating $f$ starting from a discrete set of samples of its Radon transform.
\begin{enumerate}
\item \textbf{Image reconstruction algorithm I}. Let $I$ be the interpolation of $(F^{-1}_{D}A)\ast(R_{D}f)$, so that $I(t,k\pi/N)$ is interpolated from the computed values $(F^{-1}_{D}A)\ast(R_{D}f)(jd,k\pi/N)$. Then for all points $(x_{m},y_{n})$ in the grid, we approximate
\begin{align*}
f(x_{m},y_{n})&\approx\frac{1}{2}B_{D}I(x_{m},y_{n})=\\
&=\frac{1}{2N}\sum_{k=0}^{N-1}{I\left(x_{m}\cos{\left(k\frac{\pi}{N}\right)}+y_{n}\sin{\left(k\frac{\pi}{N}\right)},k\frac{\pi}{N}\right)}.
\end{align*}
\item \textbf{Image reconstruction algorithm II}. Instead of interpolating the filtered Radon transform, we interpolate the filter and then, as shown in remark \ref{rk: intconv}, we form a weighted sum of the sampled Radon transform:
\begin{align*}
&W(k)=\sum_{j}{I_{F^{-1}A}\left(x_{m}\cos{\left(k\frac{\pi}{N}\right)}+y_{n}\sin{\left(k\frac{\pi}{N}\right)}-jd,k\frac{\pi}{N}\right)}R_{D}f(jd,k\frac{\pi}{N})
\end{align*}
\begin{align*}
&f(x_{m},y_{n})\approx\frac{1}{2N}\sum_{k=0}^{N-1}{W(k)}.
\end{align*}
\end{enumerate}

We conclude this chapter applying the reconstruction formula in a particular case.

\subsubsection{Example: crescent-shaped phantom} \label{subsec: crescent_shape}
\begin{figure}[htbp]
\centering%
\subfigure[Phantom $f$ \label{fig: phantom}]%
{\includegraphics[width=0.45\textwidth]{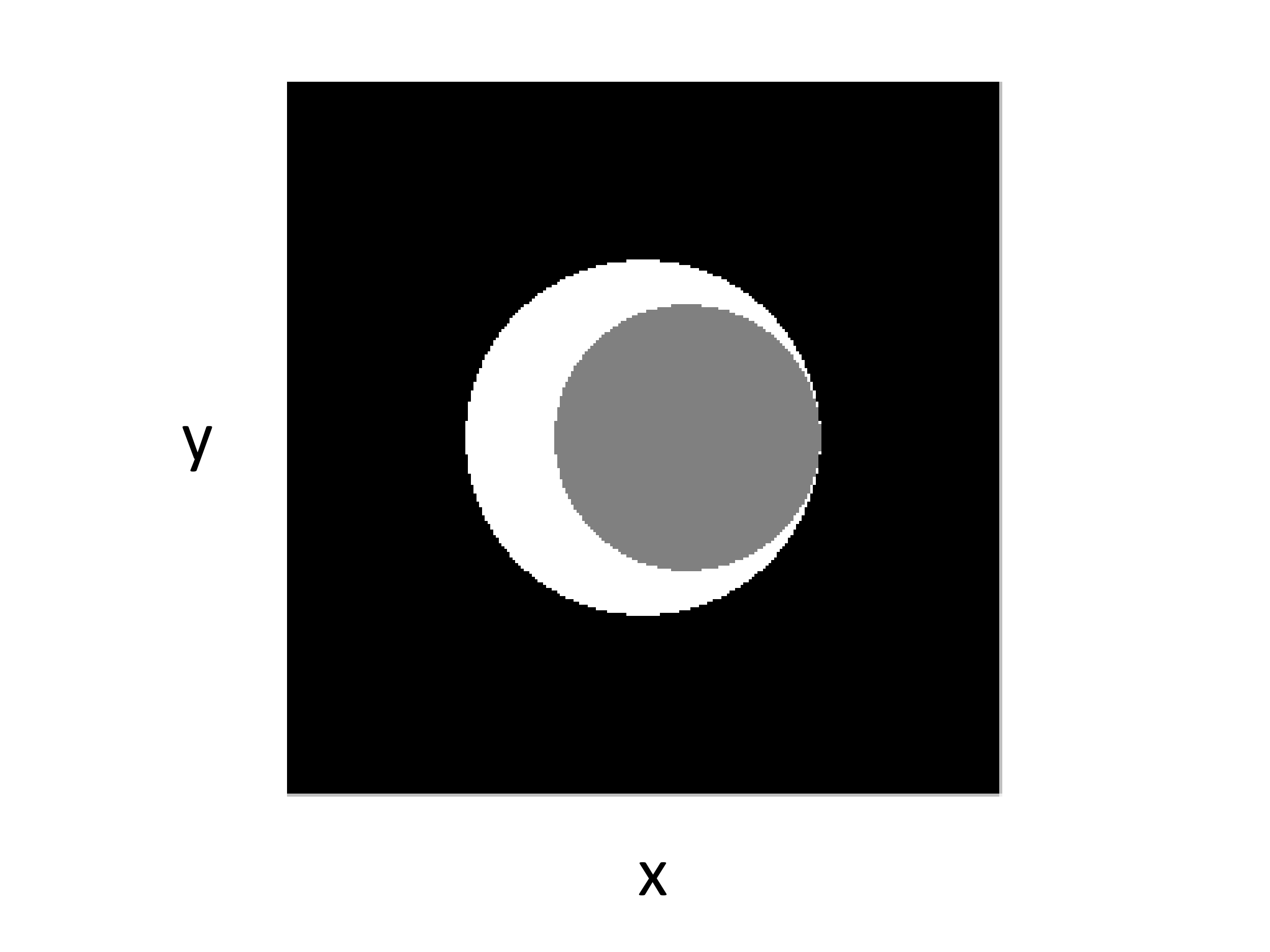}}
\subfigure[Radon transform $Rf$ \label{fig: radon}]%
{\includegraphics[width=0.45\textwidth]{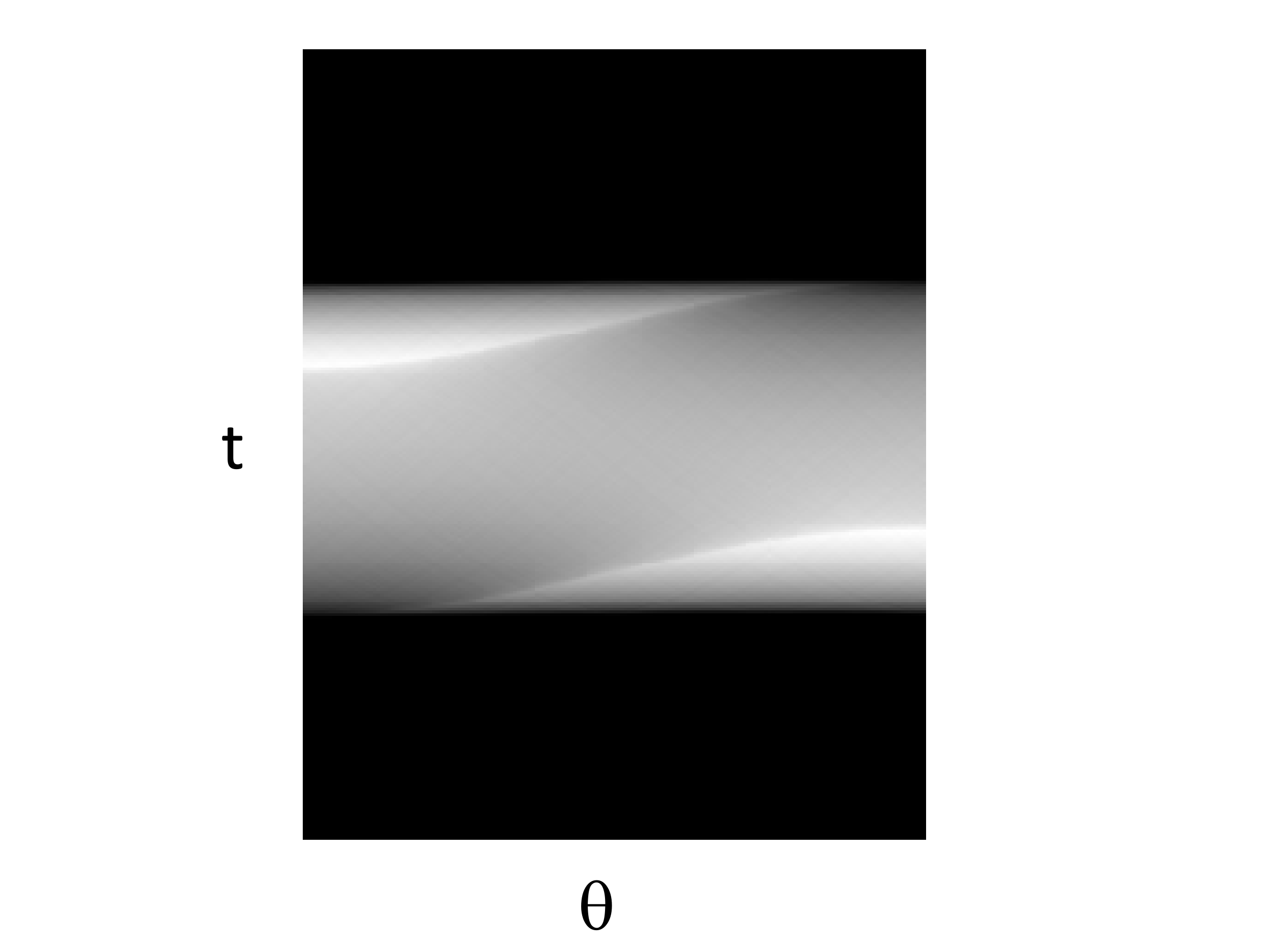}}
\caption{Crescent-shaped phantom and its Radon transform }
\label{fig: phantomRadon}
\end{figure}

We want to apply the reconstruction algorithm introduced in the previous section to a particular phantom called \emph{crescent-shaped phantom} (Figure \ref{fig: phantom}) whose analytic expression is
\begin{equation*}
f(x,y)=\left\{
\begin{aligned}
&1 & &\text{if}\ x^{2}+y^{2}\leq\frac{1}{4}\,\wedge\,(x-\frac{1}{8})^{2}+y^{2}>\frac{9}{64}\\
&\frac{1}{2} & &\text{if}\ (x-\frac{1}{8})^{2}+y^{2}\leq\frac{9}{64}\\
&0 & &\text{if}\ x^{2}+y^{2}>\frac{1}{4}.
\end{aligned}
\right.
\end{equation*}
In order to compute samples of the Radon transform, we calculate $Rf$ analytically.

We observe that $f$ can be written as a sum of two functions: $f=f_{1}-\frac{1}{2}f_{2}$, where $f_{1}$ and $f_{2}$ are given by
\begin{equation*}
f_{1}(x,y)=\left\{
\begin{aligned}
&1 & &\text{if}\ x^{2}+y^{2}\leq\frac{1}{4}\\
&0 & &\text{otherwise}
\end{aligned}
\right.\qquad
f_{2}(x,y)=\left\{
\begin{aligned}
&1 & &\text{if}\ (x-\frac{1}{8})^{2}+y^{2}\leq\frac{9}{64}\\
&0 & &\text{otherwise}
\end{aligned}
\right.
\end{equation*}
for all $(x,y)\in\R^{2}$. By the linearity of the Radon transform, we have
\begin{equation}
\boxed{
Rf=Rf_{1}-\frac{1}{2}Rf_{2}.
}
\label{eq: linearity}
\end{equation}

We know (see example \ref{es: radonCirc}) that for all fixed value $r>0$, the Radon transform of the function
\begin{equation*}
f_{r}(x,y)=\left\{
\begin{aligned}
&1 & &\text{if}\ x^{2}+y^{2}\leq r^{2}\\
&0 & &\text{otherwise},
\end{aligned}
\right.\qquad
\end{equation*}
is given by
\begin{equation}
Rf_{r}(t,\theta)=\left\{
\begin{aligned}
&2\sqrt{r^{2}-t^{2}}& &\text{if}\ |t|\leq r\\
&0 & &\text{if}\ |t|>r
\end{aligned}
\right.
\label{eq: radon_fr}
\end{equation}
and the function $f_{1}$ equals $f_{r_{1}}$ for $r_{1}=\frac{1}{2}$. So we can use equation \eqref{eq: radon_fr} to compute its Radon transform.

Function $f_{2}$ is not of the form $f_{r}$ for some $r$, but can be obtained shifting such a function. More precisely $f_{2}(x,y)=f_{r_{2}}(x-c,y)$, where $r_{2}=\frac{3}{8}$ and $c=\frac{1}{8}$. So we can use the shift property of the Radon transform for computing $Rf_{2}$:
\begin{theorem}[Shift property of the Radon transform] \label{thm: shiftProp}
Let $g:\R^{2}\rightarrow\R$ a function and let $G(t,\theta)=Rg(t,\theta)$ it's Radon transform. If 
\begin{equation*}
h(x,y)=g(x-c_{x},y-c_{y}),
\end{equation*}
then the Radon transform $H(t,\theta)=Rh(t,\theta)$ of $h$ is given by 
\begin{equation}
H(t,\theta)=G(t-c_{x}\cos{\theta}-c_{y}\sin{\theta},\theta).
\label{eq: shift}
\end{equation}  
\end{theorem}
See \cite{PEY} for more details about Radon transform shifting properties.

Thus, from \eqref{eq: radon_fr} and \eqref{eq: shift}, we gain $Rf_{2}(t,\theta)=Rf_{r_{2}}(t-c\cos{\theta},\theta)$, that is
\begin{equation*}
Rf_{2}(t,\theta)=\left\{
\begin{aligned}
&2\sqrt{r_{2}^{2}-(t-c\cos{\theta})^{2}}& &\text{if}\ |t-c\cos{\theta}|\leq r_{2}\\
&0 & &\text{if}\ |t-c\cos{\theta}|>r_{2}.
\end{aligned}
\right.
\end{equation*}
By equation \eqref{eq: linearity} we conclude that
\begin{equation*}
Rf(t,\theta)=\left\{
\begin{aligned}
&2\sqrt{r_{1}^{2}-t^{2}} & &\text{if} \ |t|\leq r_{1} \, \wedge \, |t-c\cos{\theta}|>r_{2}\\
&2\sqrt{r_{1}^{2}-t^{2}}-\sqrt{r_{2}^{2}-(t-c\cos{\theta})^{2}}& &\text{if}\ |t-c\cos{\theta}|\leq r_{2}\\
&0 & &\text{if}\ |t|>r_{1}.
\end{aligned}
\right.
\end{equation*}
Figure \ref{fig: radon} shows the spectra of this function.

Suppose now $M=20$ and $N=18$. We sample the domain $[-1,1]\times[0,\pi)$ with values $t_{k}=kd$, $k=-M,\ldots,M$ and $\theta_{j}=j\frac{\pi}{N}$, for $j=0,\ldots,N-1$, where $d=0.05$, obtaining the discrete Radon transform $Rf(kd,j\pi/N)$.

We consider \emph{Shepp-Logan filter} as low pass filter and we consider $\frac{1}{2L}$ as sampling spacing. Indeed, when we used the continuous Fourier transform, we considered $\frac{\pi}{L}$ as sampling spacing, in accordance with Nyquist theorem. Now, to compensate the additional factor $2\pi$ in the definition of the discrete inverse Fourier transform, we use $\frac{\pi}{2\pi L}=\frac{1}{2L}$. To match this spacing with that of the Radon transform, we want $\frac{1}{2L}=d=0.05$ and so $L=10$, then
\begin{equation*}
A(\omega)=
\begin{cases}
\frac{\pi}{20}|\sin{(0.05\pi\omega)}|  &\quad \ |\omega|\leq10\\
0 &\quad \ |\omega|>10
\end{cases}
\end{equation*}
and
\begin{equation*}
(F_{D}^{-1}A)_{n}=\frac{400}{\pi^{3}(1-4n^{2})}.
\end{equation*}

Next we compute the discrete convolution $\gamma=F_{D}^{-1}A\ast Rf$: for $-20\leq m\leq20$, $0\leq j\leq17$
\begin{equation*}
\gamma(m,j\pi/N)=\sum_{k=-20}^{20}{(F_{D}^{-1}A)_{m-k}Rf(0.05k,j\pi/N)}.
\end{equation*}
Applying \emph{linear interpolation} to the variable $t$ of $\gamma$, we obtain
\begin{equation*}
h(t,j\pi/N)=\sum_{m=-20}^{20}{\gamma(m,j\pi/N)\Lambda(20t-m)}, \quad t\in[-1,1].
\end{equation*}

\begin{figure}[htbp]
\centering
\includegraphics[width=0.7\textwidth]{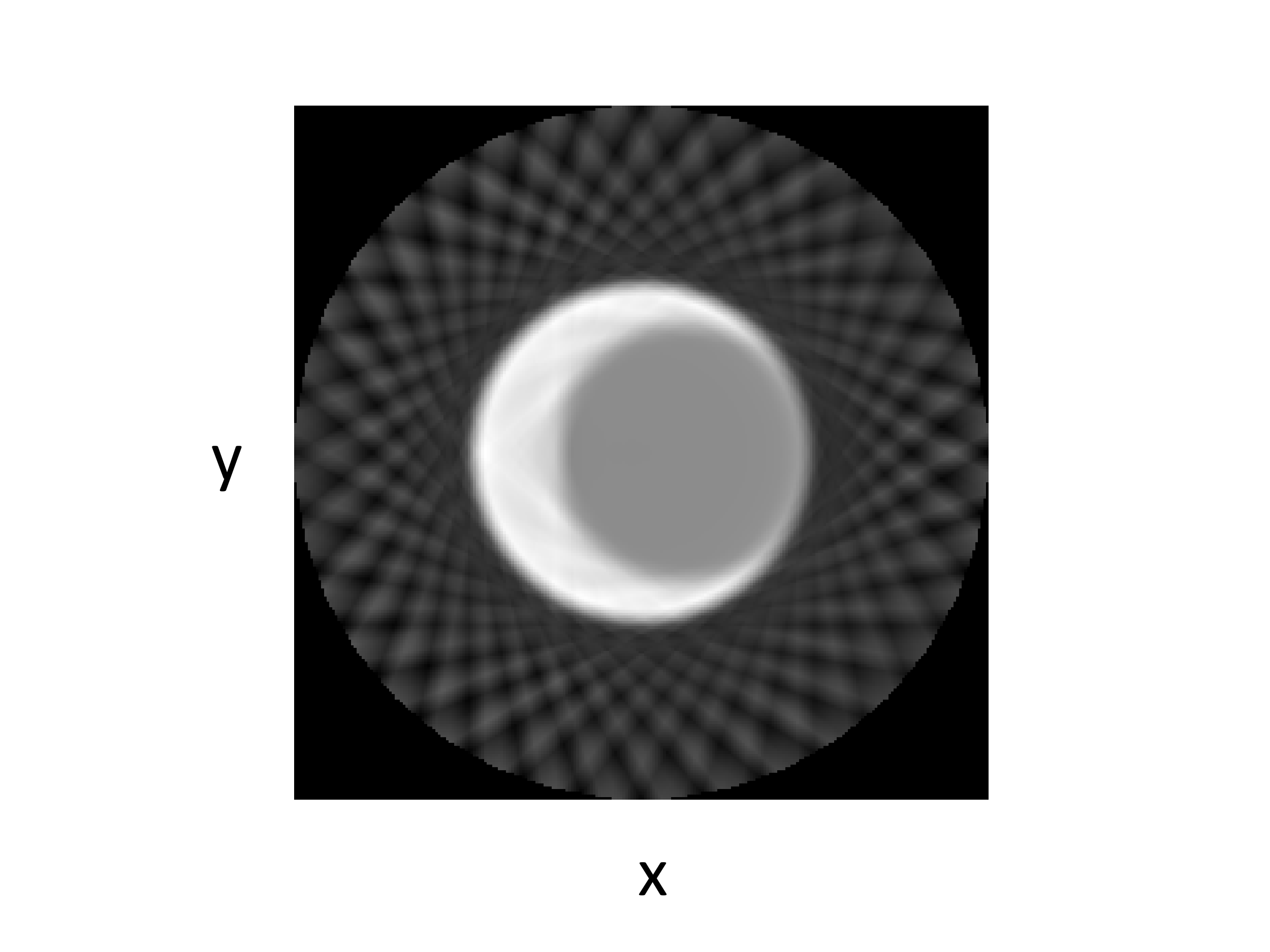} 
\caption{Reconstruction of $f$}
\label{fig: reconstCshape}
\end{figure}

Finally we use the \emph{reconstruction algorithm I} and we have the approximation
\begin{equation*}
f(x,y)\approx\frac{1}{36}\sum_{j=0}^{17}{h(x\cos{(j\pi/N)}+y\sin{(j\pi/N)},j\pi/N)}.
\end{equation*}
Figure \ref{fig: reconstCshape} shows the reconstructed function.

%Chapter 3 - Algebraic Reconstruction Techniques

\chapter{Algebraic Reconstruction Techniques}\label{chap: art}
%\section{Introduction}
The Fourier based methods we have seen so far are the algorithms used in modern CT scan. Another approach to image reconstruction is based on linear algebra. Algorithms that use this approach are known as \emph{algebraic reconstruction techniques}, or ART. For example, the first CT scanner designed in the late 1960s by Godfrey Hounsfield used these methods.

While the Fourier transform approach solves the continuous problem and then passes to the discrete one, ART considers the discrete problem from the beginning.

Let us start reminding that an image is given by a grid of pixels (picture elements) and at each pixel is assigned a color (or a gray scale value) that represents the value of the attenuation coefficient in the region of the given pixel.

Suppose that our image is formed by $K\times K$ pixels, each of them representing a small square in the plane. Define the pixel basis functions $b_{1},b_{2},\ldots,b_{K^{2}}$ as 
\begin{equation*}
b_{i}(x,y)=\left\{
\begin{aligned}
&1 & &\text{if $(x,y)$ lies in pixel number $i$}\\
&0 & &\text{otherwise,}
\end{aligned}
\right. 
\end{equation*} 
and let $x_{i}$ the color value of the $i$-th pixel. Then the resulting image can be written as
\begin{equation*}
I(x,y)=\sum_{i=1}^{K^{2}}{b_{i}(x,y)x_{i}}.
\end{equation*}
Applying the Radon transform to both sides, we get 
\begin{equation*}
RI(t,\theta)=\sum_{i=1}^{K^{2}}{Rb_{i}(x,y)x_{i}}.
\end{equation*}

The X-ray machine gives us the value of the attenuation coefficient function $f$ for some finite set of lines $l_{t_{j},\theta_{j}}$, $j=1,\ldots,J$. Let us denote by $p_{j}=Rf(t_{j},\theta_{j})$ these values. We want to approximate the attenuation coefficient $f$ with image $I$, so we set $p_{j}=RI(t_{j},\theta_{j})$ and $r_{j.i}=Rb_{i}(t_{j},\theta_{j})$, for $j=1,\ldots,J$ and $i=1,\ldots,K^{2}$, and we ask that
\begin{equation}
p_{j}=\sum_{i=1}^{K^{2}}{x_{i}r_{j,i}}, \quad j=1,\ldots,J.
\label{eq: system}
\end{equation} 
Thus we obtain a system of $J$ linear equations and $K^{2}$ unknowns. This system is very large but spare and typically overdetermined or underdetermined. We need then specific techniques for the solution of such a system. Before looking to these methods, let us see in detail how to generate the linear system \eqref{eq: system}.

\section{Generation of the linear system}
In this section we consider the problem of generate the linear system $Ax=p$, i.e. we want to compute $A$ and $p$ starting from the values of the Radon transform $Rf$ of an attenuation coefficient function $f$ obtained from a X-ray machine working with parallel beam geometry.

We know the values $Rf(t_{k},\theta_{l})$ with $t_{k}=kd$, $k=-M,\ldots,M$ and $\theta_{l}=l\frac{\pi}{N}$, $j=0,\ldots,N$. We want to compute $A=(r_{j,i})\ i=1,\ldots,K^{2}\ j=1,\ldots,J$ where $K^2$ is the dimension of the reconstructed gray-scale image $I=\{PX_{i}\}_{i=1,\ldots,K^{2}}$, with the components $x_{i}$ of the solution of the system representing the color of pixel $PX_{i}$; $J=(2M+1)N$ is the number of samples $(t_{j},\theta_{j})$ on which $Rf$ is measured and $r_{j,i}=Rb_{i}(t_{j},\theta_{j})$ is the Radon transform of the $i$-th pixel-basis function $b_{i}$, computed at point $(t_{j},\theta_{j})$, with $b_{i}$ defined by
\begin{equation*}
b_{i}(x,y)=\left\{
\begin{aligned}
&1 &\text{if} \ (x,y)\in PX_{i}\\
&0 &\text{if} \ (x,y)\notin PX_{i}
\end{aligned}
\right. .
\end{equation*} 

In order to solve this problem we assume that:
\begin{enumerate}
\item The support of the function $f$ and the samples points $(t_{j},\theta_{j})$ are contained in the unit square $[-1,1]\times[-1,1]$, this implies that $d=\frac{1}{M}$;
\item The reconstructed image $I$ also lies in $[-1,1]\times[-1,1]$ and its center is at the origin $(0,0)$. If we consider $I$ as a matrix $I(r,s)$, $r=1,\ldots,K$, $s=1,\ldots,K$, whose components are the values $x_{i}$ of pixels $PX_{i}$, the center is the pixel of indexes $r=\lfloor\frac{K+1}{2}\rfloor$, $s=\lfloor\frac{K+1}{2}\rfloor$. We denote $c=\lfloor\frac{K+1}{2}\rfloor$;
\item The $K^{2}$ pixels in $I$ are ordered as follows:
\begin{equation*}
I=\left(
\begin{array}{cccc}
x_{1} & x_{2} & \cdots &x_{K}\\
x_{K+1} & x_{K+2} & \cdots & x_{2K}\\
\vdots & & & \vdots\\
x_{K(K-1)+1} & x_{K(K-1)+2} & \cdots &x_{K^{2}}\\
\end{array}
\right);
\end{equation*}
\item Considering $I$ as a function $I:\R^2\rightarrow\R$, i.e. $I(x,y)=\sum_{i}{b_{i}(x,y)x_{i}}$, the Cartesian coordinates of pixels are $PX_{i}=[x_{i},x_{i+1})\times(y_{i+1},y_{i}]$. Thus, we are considering a top-down, left-right enumeration of vertexes, in accord with the matrix indexing. Note that $PX_{i}$ includes the top horizontal side and the left vertical side, but not the right and the bottom sides (see Figure \ref{fig: pixel}), exception are the pixels in the last row and in the last column of $I$ that include all sides. We identify a pixel with the coordinates of its top-left vertex $(x_{i},y_{i})$;%, and with an abuse of notation we will also write $PX_{i}=(x_{i},y_{i})$;
\begin{figure}[htbp]
\centering
\includegraphics[width=0.3\textwidth]{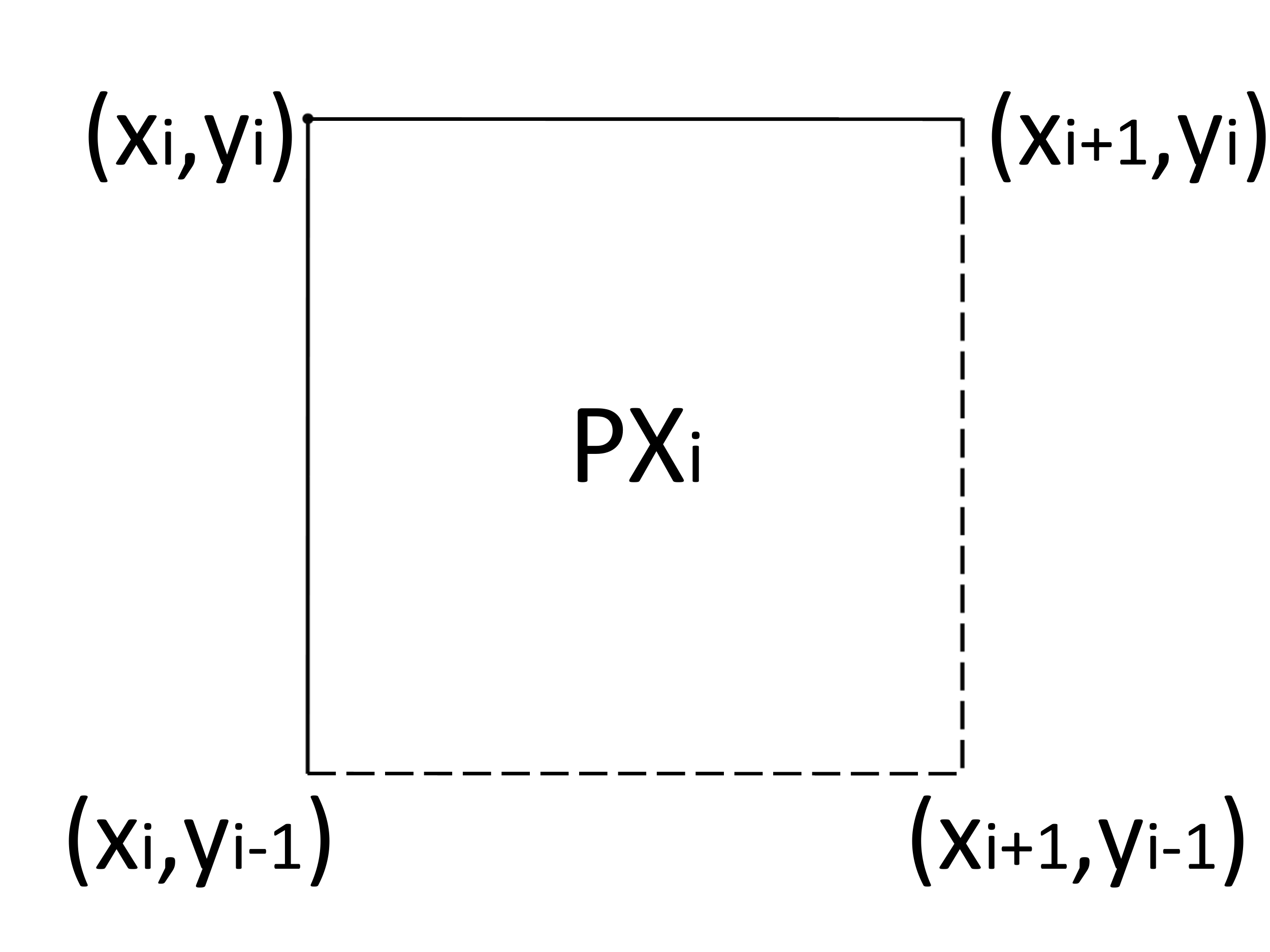} 
\caption{Coordinates of a pixel}
\label{fig: pixel}
\end{figure}
\item X-ray beams has zero width.
\end{enumerate}

Using these assumptions, we find that pixel $PX_{i}$, determinate by coordinates $(x_{i},y_{i})$, given by
\begin{equation*}
\boxed{
x_{i}=\frac{\bar{x}_{i}}{c}-1; \qquad y_{i}=\frac{-\bar{y}_{i}}{c}+1.
}
\end{equation*}
In particular 
\begin{itemize}
\item $PX_{i}=[x_{i},x_{i+1})\times(y_{i+1},y_{i}]$, $x_{i+1}=x_{i}+c^{-1}$, $y_{i+1}=y_{i}-c^{-1}$;
\item If $i\in\{K,2K,\ldots,K^{2}\}\ \Rightarrow\ PX_{i}=[x_{i},x_{i+1}]\times(y_{i+1},y_{i}]$;
\item if $i\in\{K(K-1)+1,\ldots,K^{2}\}\ \Rightarrow\ PX_{i}=[x_{i},x_{i+1})\times[y_{i+1},y_{i}]$.
\end{itemize}

\section{Construction of $A$ and $p$}
Assume that we know $(t_{j},\theta_{j})$ and pixel coordinates $PX_{i}=[x_{i},x_{i+1})\times(y_{i+1},y_{i}]$. What we want to do now is to compute $r_{j,i}=Rb_{i}(t_{j},\theta_{j})$. Let 
\begin{equation*}
A=(r_{j,i})=\left(
\begin{array}{cccc}
r_{11} & r_{12} & \cdots & r_{1K^{2}}\\
r_{21} & r_{12} & \cdots & r_{2K^{2}}\\
\vdots & & & \vdots\\
r_{J1} & r_{J2} & \cdots & r_{JK^{2}}
\end{array}
\right)=(A^{1},A^{2},\ldots,A^{K^{2}}).
\end{equation*}
and $p=Rf(t_{j},\theta_{j})$.

For a fixed $i\in\{1,\ldots,K^{2}\}$ we define $\tilde{A}^{i}\in\R^{2M+1}\times\R^{ N}$ the matrix such that $\tilde{A}^{i}(r,s)=Rb_{i}(rd,s\frac{\pi}{N})$, with $r=-M,\ldots,M,\ s=0,\ldots,N-1$, i.e. the columns of $\tilde{A}^{i}$ represent values of $Rb_{i}$ for a fixed value of $\theta$:
\begin{equation*}
\tilde{A}^{i}=\left(
\begin{array}{ccc}
Rb_{i}(-Md,0)  & \cdots & Rb_{i}(-Md,(N-1)\frac{\pi}{N})\\
\vdots & & \vdots\\

Rb_{i}(Md,0)  & \cdots & Rb_{i}(Md,(N-1)\frac{\pi}{N})
\end{array}
\right).
\end{equation*}
If $R\in\R^{2M+1}\times\R^{ N}$ is the matrix containing data $Rf(t_{k},\theta_{l})$ and if we set
\begin{equation*}
p=R(:)
\end{equation*}
then we have that the $i$-th column of $A$ is  
\begin{equation*}
A^{i}:=\tilde{A}^{i}(:),
\end{equation*}
where the operator $(:)$ indicates the analogous Matlab operator (see \cite{MATLAB}).
%RIF MATLAB
% that applied to a $m\times n$ matrix $B$, gives a column vector of length $mn$ obtained putting in sequence columns of $B$.

The problem can therefore be reduced to computation of the columns of $\tilde{A}^{i}$ for a fixed $i$.% In other words, we want to compute $Rb_{i}(t,\theta)$ for fixed $i\in\{1,\ldots,K^{2}\}$ and fixed $\theta=l\frac{\pi}{N}$, for all values $t=t_{k}=kd$, $k\in\{-M,\ldots,M\}$.

%\subsection{Computation of $r_{j,i}$}
Let $i\in\{1,\ldots,K^{2}\}$ and $\theta\in[0,\pi)$ fixed. For $t\in\R$ let $r(t)=Rb_{i}(t,\theta)$. We observe that since $b_{i}\equiv1$ inside pixel $PX_{i}$ and $b_{i}\equiv0$ outside (and since we assume X-ray beams to have zero width), the value of $r$ is the length of the intersection between line $l_{t,\theta}$ and $PX_{i}$. Indeed:
\begin{align*}
Rb_{i}(t,\theta)&=\int_{\R}{b_{i}(t\cos{\theta}-s\sin{\theta},t\sin{\theta}+s\cos{\theta})\,ds}=\\
&=\int_{\{s\in\R: (t\cos{\theta}-s\sin{\theta},t\sin{\theta}+s\cos{\theta})\in PX_{i}\}}{ds}=m(C),
\end{align*}
where $m$ denotes the Lebesgue measure on $\R$ and $C=\{s\in\R:\ t\cos{\theta}-s\sin{\theta}\in[x_{i},x_{i+1}),\ t\sin{\theta}+s\cos{\theta}\in(y_{i+1},y_{i}] \}=l_{t,\theta}\cap PX_{i}$.

\begin{figure}[htbp]
\centering%
\subfigure[$\theta\in[0,\frac{\pi}{4})$ \label{fig: calcolo_t1}]%
{\includegraphics[width=0.4\textwidth]{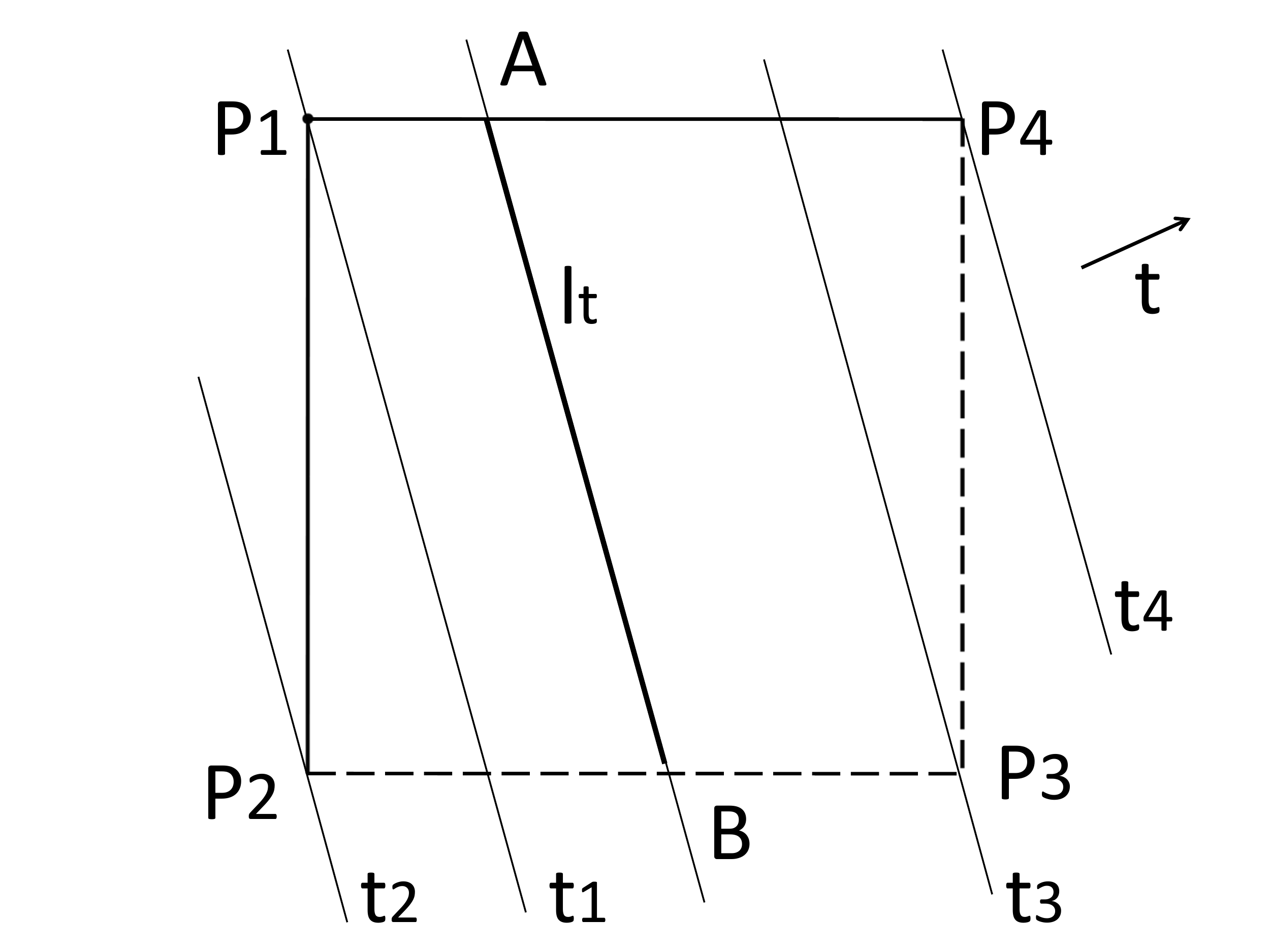}}
\subfigure[$\theta\in[\frac{\pi}{4},\frac{\pi}{2})$ \label{fig: calcolo_t2}]%
{\includegraphics[width=0.4\textwidth]{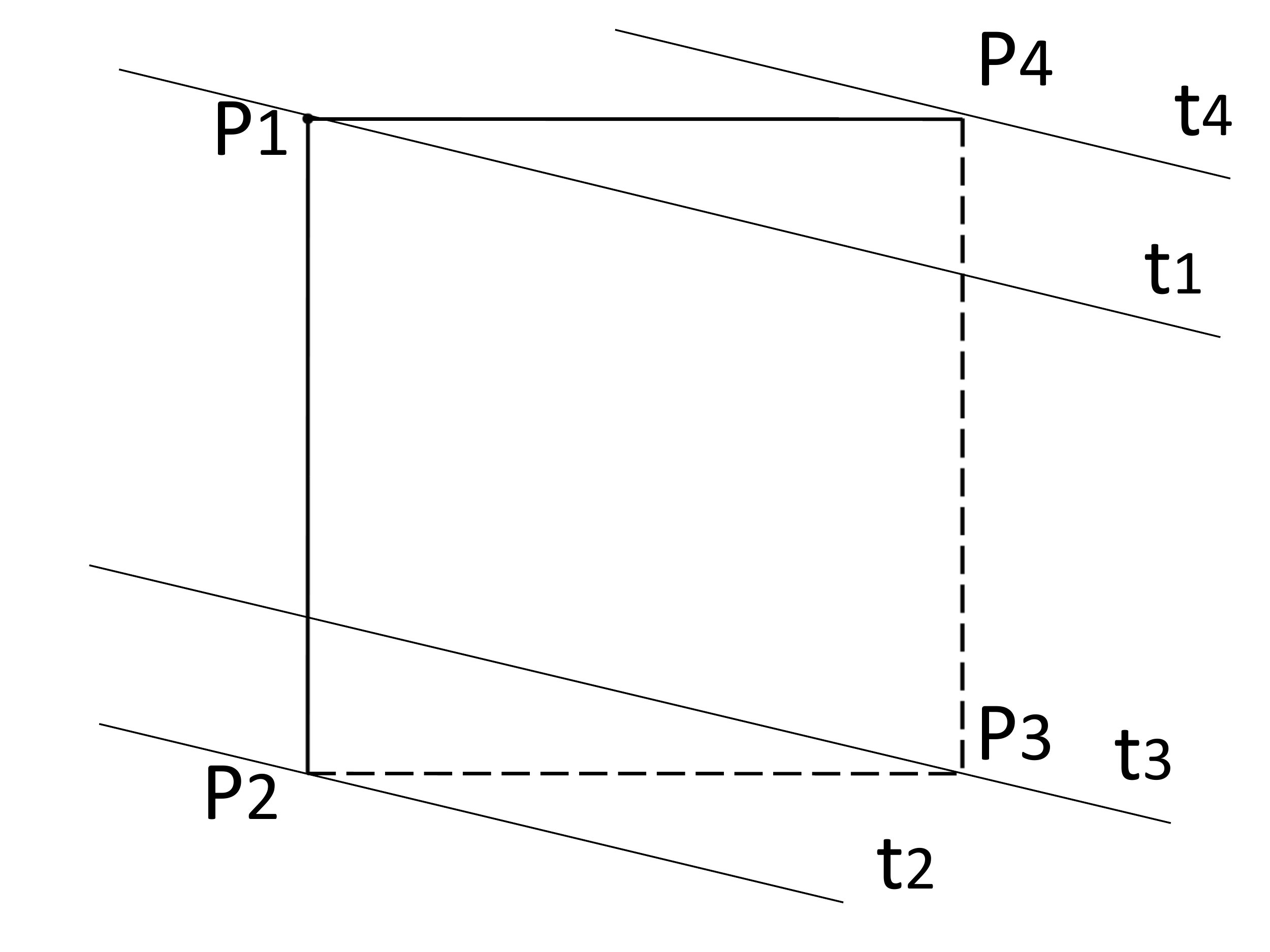}}
\caption{Lines $t_{h}$ for $\theta\in[0,\frac{\pi}{2})$}
\label{fig: calcolo_t0}
\end{figure}

To determine for which values of $t$ the line $l_{t,\theta}$ lies in $PX_{i}$, we consider lines that pass through the vertexes of $PX_{i}$. Let 
\begin{align*}
&P_{1}=(x_{i},y_{i}) & &P_{2}=(x_{i},y_{i+1}) & &P_{3}=(x_{i+1},y_{i+1}) & &P_{4}=(x_{i+1},y_{i})
\end{align*}
and let $t_{h}$ be such that the line $l_{t_{h},\theta}$ passes through point the $P_{h},\ h=1,2,3,4$ (Figure \ref{fig: calcolo_t0}). Since $(x_{0},y_{0})\in l_{x_{0}\cos{\theta}+y_{0}\sin{\theta},\theta}$, we have 
\begin{align*}
&t_{1}=x_{i}\cos{\theta}+y_{i}\sin{\theta} & &t_{2}=x_{i}\cos{\theta}+y_{i+1}\sin{\theta} \\ &t_{3}=x_{i+1}\cos{\theta}+y_{i+1}\sin{\theta} & &t_{4}=x_{i+1}\cos{\theta}+y_{i}\sin{\theta}
\end{align*}

Moreover, to determine the length of the intersection $l_{t,\theta}\cap PX_{i}$, we need to know the intersections between $l_{t,\theta}$ and the sides of $PX_{i}$. Let 
\begin{align*}
&E_{12}=l_{t,\theta}\cap P_{1}P_{2} & &E_{23}=l_{t,\theta}\cap P_{2}P_{3} & &E_{34}=l_{t,\theta}\cap P_{3}P_{4} & &E_{14}=l_{t,\theta}\cap P_{1}P_{4}.
\end{align*}
Let us compute for example $E_{12}$:
\begin{equation*}
l_{t,\theta}= (t\cos{\theta}-s\sin{\theta},t\sin{\theta}+s\cos{\theta})=(x(s),y(s))
\end{equation*}
the line through $P_{1}P_{2}$ is $x=x_{i}$, so we want $x(s)=x_{i}$, $\Rightarrow$
\begin{align*}
&s=\frac{t\cos{\theta}-x_{i}}{\sin{\theta}},
&y(s)=t\sin{\theta}+\frac{t\cos{\theta}-x_{i}}{\sin{\theta}}\cos{\theta}=\frac{t-x_{i}\cos{\theta}}{\sin{\theta}}.
\end{align*}
In a similar way we find $E_{23},E_{34},E_{41}$:
\begin{align*}
&E_{12}=\left(x_{i},\frac{t-x_{i}\cos{\theta}}{\sin{\theta}}\right) & &E_{23}=\left(\frac{t-y_{i+1}\sin{\theta}}{\cos{\theta}},y_{i+1}\right)\\
&E_{34}=\left(x_{i+1},\frac{t-x_{i+1}\cos{\theta}}{\sin{\theta}}\right) &  &E_{14}=\left(\frac{t-y_{i}\sin{\theta}}{\cos{\theta}},y_{i}\right)\\
\end{align*}
Of course this values are different in the limit case $\cos{\theta}=0$ or $\sin{\theta}=0$, i.e. for $\theta=0,\pi/2$. In these cases  intersections between $l_{t,\theta}$ and $PX_{i}$ coincide with the vertexes $P_{h}$.

Depending on the values of $\theta$, the behavior of $l_{t,\theta}$ and $l_{t_{h},\theta}$ changes (see Figure \ref{fig: calcolo_t0}).%\ref{fig: calcolo_tpi}
Therefore one should distinguish the cases $\theta\in[0,\frac{\pi}{4}),\ \theta\in[\frac{\pi}{4},\frac{\pi}{2}),\ \theta\in[\frac{\pi}{2},\frac{3}{4}\pi)\ \text{and}\ \theta\in[\frac{3}{4}\pi,\pi)$. Consider for example
%\begin{enumerate}
%\item 
$\theta\in[0,\frac{\pi}{4})$. In this case we have $t_{2}\leq t_{1}<t_{3}\leq t_{4}$, therefore
\begin{align*}
&\text{if} \ t<t_{2}\, \vee\, t\geq t_{4} \ \Rightarrow \ l_{t,\theta}\cap PX_{i}=\emptyset\\ 
&\text{if} \ t_{2}\leq t<t_{4} \ \Rightarrow \ l_{t,\theta}\cap PX_{i}=AB
\end{align*}
where $AB$ is given by
\begin{enumerate}
\item $AB=E_{23}E_{12} \quad$ if $t_{2}\leq t<t_{1}$;
\item\label{case: 1b} $AB=E_{23}E_{14}\quad$ if $t_{1}\leq t\leq t_{3}$;
\item $AB=E_{34}E_{14}\quad$ if $t_{3}<t<t_{4}$;
\item\label{case: 1d} In the limit case $i\in\{K,2K,\ldots,K^{2}\}$, i.e. if we are considering a pixel in the last column of $I$, we have to account that side $P_{3}P_{4}\in PX_{i}$, so
\begin{equation*}
\text{if} \ x_{i}=x_{K^{2}}\, \wedge \, \theta=0 \, \wedge \, t=t_{4} \ \Rightarrow r=\overline{AB}=\overline{P_{3}P_{4}}=c^{-1}.
\end{equation*}
\end{enumerate}
We notice that in case \ref{case: 1b}. $AB(t)=AB(t_{1})=AB(t_{3})=const$. Moreover for $\theta=0$, only cases \ref{case: 1b}. and \ref{case: 1d}. are possible, then we do not need to compute $1/\sin{\theta}$.
The determination of $AB$ in the others 3 cases is similar. What remains to do now is to compute the length of $AB$.

\begin{figure}[htbp]
\centering%
\subfigure[$\theta\in[0,\frac{\pi}{2})$ \label{fig: calcolo_AB1}]%
{\includegraphics[width=0.4\textwidth]{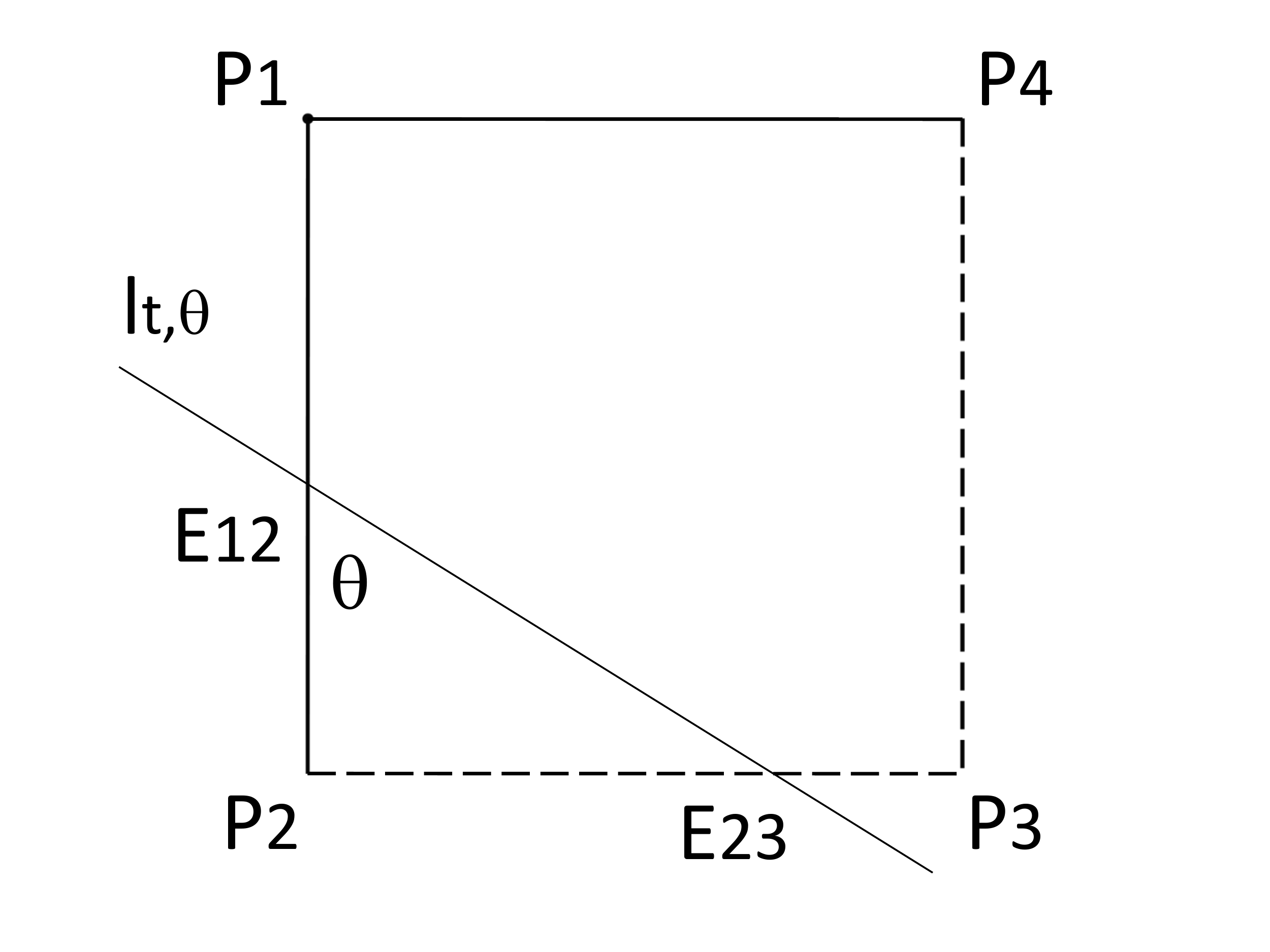}}
\subfigure[$\theta\in[\frac{\pi}{2},\pi)$ \label{fig: calcolo_AB2}]%
{\includegraphics[width=0.4\textwidth]{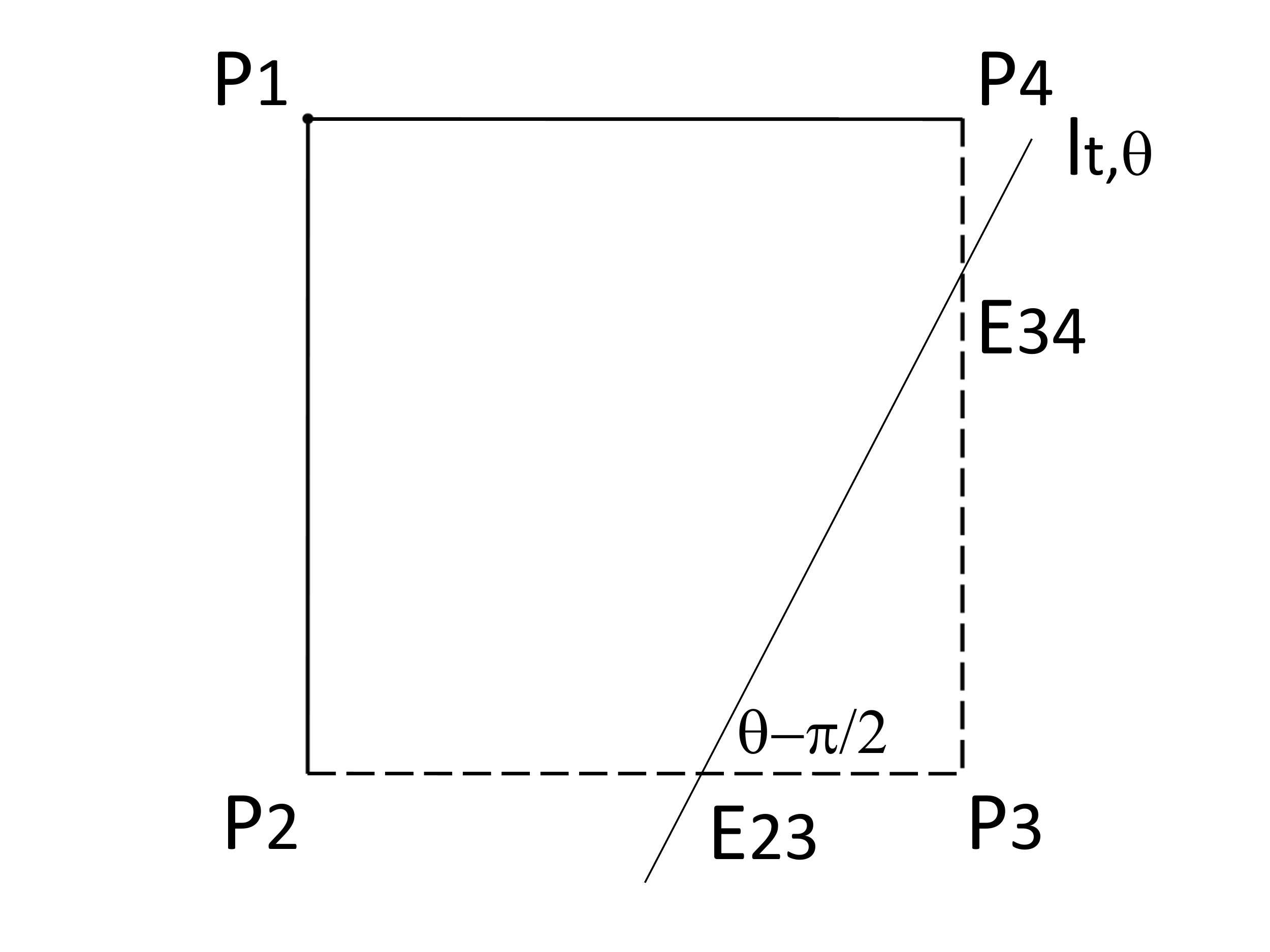}}
\caption{Computation of $AB$ }
\label{fig: calcoloAB}
\end{figure}

%\subsection{Computation of $\overline{AB}$}
Let us start considering $AB=E_{12}E_{23}$  (see Figure \ref{fig: calcoloAB}). The coordinates of the points are 
\begin{align*}
&E_{12}=\left(x_{i},\frac{t-x_{i}\cos{\theta}}{\sin{\theta}}\right) & &E_{23}=\left(\frac{t-y_{i+1}\sin{\theta}}{\cos{\theta}},y_{i+1}\right)
\end{align*}
We notice that the triangle $E_{12}\widetriangle{P}_{2}E_{23}$ is rectangle in $P_{2}$ and that, by definition, $\theta=P_{2}\hat{E}_{12}E_{23}$. 
By Pythagoras theorem 
\begin{equation*}
\overline{E_{12}E_{23}}=\frac{\overline{P_{2}E_{23}}}{\sin{\theta}},
\end{equation*} 
thus
\begin{align*}
P_{2}E_{23}=&x_{E_{23}}-x_{P_{2}}=\frac{t-y_{i+1}\sin{\theta}}{\cos{\theta}}-x_{i}=\\
&=\frac{t-y_{i+1}\sin{\theta}-x_{i}\sin{\theta}}{\cos{\theta}}=\\
&=\frac{t-t_{2}}{\cos{\theta}}.
\end{align*}
We conclude that 
\begin{equation*}
\overline{E_{12}E_{23}}=\frac{t-t_{2}}{\sin{\theta}\cos{\theta}}, \qquad \text{for}\ \theta\in(0,\frac{\pi}{2}).
\end{equation*}

We now consider $AB=E_{23}E_{14}$: as stated before, for all $t\in[t_{1},t_{3}]$, $\theta\in[0,\frac{\pi}{4})$, $E_{23}E_{14}(t)=E_{23}E_{14}(t_{1})$, hence 
\begin{equation*}
\overline{E_{23}E_{14}}=\frac{\overline{P_{1}P_{2}}}{\cos{\theta}}=\frac{c^{-1}}{\cos{\theta}}.
\end{equation*}

%In a similar way one computes all the other values.
%\begin{align*}
%&\overline{E_{12}E_{23}}=\frac{t-t_{2}}{\sin{\theta}\cos{\theta}}, & &\text{for} \ \theta\in(0,\frac{\pi}{4})\ t_{2}\leq t<t_{1} \, \vee \, \theta\in[\frac{\pi}{4},\frac{\pi}{2})\ t_{2}< t<t_{3}\\
%%
%&\overline{E_{23}E_{14}}=\frac{c^{-1}}{\cos{\theta}}, & &\text{for} \ \theta\in[0,\frac{\pi}{4})\ t_{1}\leq t\leq t_{3}\\
%%
%&\overline{E_{34}E_{14}}=\frac{t_{4}-t}{\sin{\theta}\cos{\theta}}, & &\text{for} \ \theta\in(0,\frac{\pi}{4})\ t_{3}< t<t_{4} \, \vee \, \theta\in[\frac{\pi}{4},\frac{\pi}{2})\ t_{1}< t<t_{4}\\
%%
%&\overline{E_{34}E_{12}}=\frac{c^{-1}}{\sin{\theta}}, & &\text{for} \ \theta\in[\frac{\pi}{4},\frac{\pi}{2})\ t_{3}\leq t\leq t_{1}\\
%\end{align*}
%\begin{align*}
%&\overline{E_{23}E_{34}}=\frac{t_{3}-t}{\sin{\theta}\cos{\theta}}, & &\text{for} \ \theta\in(\frac{\pi}{2},\frac{3}{4}\pi)\ t_{3}< t<t_{2} \, \vee \, \theta\in[\frac{3}{4}\pi,\pi)\ t_{3}< t<t_{4}\\
%%
%&\overline{E_{12}E_{34}}=\frac{c^{-1}}{\sin{\theta}}, & &\text{for} \ \theta\in[\frac{\pi}{2},\frac{3}{4}\pi)\ t_{2}\leq t\leq t_{4}\\ 
%%
%&\overline{E_{12}E_{14}}=\frac{t-t_{1}}{\sin{\theta}\cos{\theta}}, & &\text{for} \ \theta\in(\frac{\pi}{2},\frac{3}{4}\pi)\ t_{4}< t\leq t_{1} \, \vee \, \theta\in[\frac{3}{4}\pi,\pi)\ t_{2}< t<t_{1}\\
%%
%&\overline{E_{23}E_{14}}=-\frac{c^{-1}}{\cos{\theta}}, & &\text{for} \ \theta\in[\frac{3}{4}\pi,\pi)\ t_{4}\leq t\leq t_{2}\\
%\end{align*}

We can reduce the number of cases if we order $t_{h}$, $h=1,2,3,4$ in increasing order: $t_{\min}\leq t_{\min2}\leq t_{\max2}\leq t_{\max}$. Thus, we conclude that
\begin{itemize}
\item for $\theta\neq0,\frac{\pi}{2}$ we have
\begin{equation}
Rb_{i}(t,\theta)=\left\{
\begin{aligned}
&r_{1} & &\text{if}\ t_{\min}<t<t_{\min2}\\
&r_{2} & &\text{if}\ t_{\min2}\leq t\leq t_{\max2}\\
&r_{3} & &\text{if}\ t_{\max2}<t<t_{\max},
\end{aligned}
\right.
\label{eq: Rgen}
\end{equation}
where
\begin{align*}
r_{1}&=\left| \frac{t-t_{\min}}{\sin{\theta}\cos{\theta}}\right| & r_{2}&=c^{-1}\min{\left(\frac{1}{|\cos{\theta}|},\frac{1}{|\sin{\theta}|}\right)} & r_{3}&=\left| \frac{t-t_{\max}}{\sin{\theta}\cos{\theta}}\right|.
\end{align*}
\item for $\theta=0$
\begin{equation}
Rb_{i}(t,0)=\left\{
\begin{aligned}
&c^{-1} & &\text{if}\ t_{\min}\leq t<t_{\max} \, \vee \, (x_{i}=x_{K^{2}}\,\wedge\,t=t_{\max})\\
&0 & &\text{otherwise}
\end{aligned}
\right.
\label{eq: R0}
\end{equation}
\item for $\theta=\frac{\pi}{2}$
\begin{equation}
Rb_{i}(t,\frac{\pi}{2})=\left\{
\begin{aligned}
&c^{-1} & &\text{if}\ t_{\min}< t\leq t_{\max} \, \vee \, (y_{i}=y_{K^{2}}\,\wedge\,t=t_{\min})\\
&0 & &\text{otherwise}
\end{aligned}
\right. .
\label{eq: Rpi2}
\end{equation}
\end{itemize}

\subsection{Algorithm}
We now summarize in a pseudo-algorithm the principal steps involved in the computation of $A$ and $p$. Operations are indicated in Matlab language.
\begin{enumerate}
\item Input: 
\begin{verbatim}
R: (2M+1)N matrix representing the Radon data
K: dimension of the output image
\end{verbatim} 
\item Initialization:
\begin{verbatim}
A=zeros((2M+1)*N,K^{2});
c=floor((K+1)/2);
\end{verbatim}
\item Compute pixels coordinates:
\begin{verbatim}
Y=repmat([0:K,K+1,1]);
X=Y';
x=X/c-1; y=-Y/c+1;
\end{verbatim}
\item Compute values $t_{h}=t_{h}(x_{i},y_{i},\theta)$, $h=1,2,3,4$:
\begin{verbatim}
T=zeros(K+1,K+1,N);
for j=1:N,
     T(:,:,j)=x*cos(theta(j))+y*sin(theta(j));
end
\end{verbatim}
\item Compute $A$:
\begin{itemize}
\item For all \texttt{i=1:K*K} extract sub-matrix \texttt{Ti} of \texttt{T} containing $t_{h}$ values corresponding to pixel $PX_{i}$;
\item Sort \texttt{Ti} in increasing order;
\item for all \texttt{j=1:N}, for all \texttt{k=1:2M+1}, compute $\tilde{A}^{i}(t_{k},\theta_{j})$ using equations \eqref{eq: Rgen},\eqref{eq: R0},\eqref{eq: Rpi2};
\item Fill $i$-th column of \texttt{A}:
\begin{verbatim}
A(:,i)=Atilde_i(:);
\end{verbatim}
\end{itemize}
\item Compute $p$:
\begin{verbatim}
p=R(:);
\end{verbatim}
\end{enumerate} 

\section{Solving the system}
In the previous section we saw how to compute matrix $A$ and the r.h.s. $p$ of the linear system generated from the algebraic approach to the image reconstruction problem. In this section we discuss other methods useful for the solution of this system.

We start observing that the matrix $A$ can be very large. In fact, every sampling of the Radon transform produces an equation, while at every pixel in the output image is associated an unknown. Moreover the system can be typically both underdetermined (more unknowns then equations) or overdetermined (more equations then unknowns). Another important property of the system $Ax=p$ is that the matrix $A$ is sparse. Indeed every particular line $l_{t_{j},\theta_{j}}$ passes through relatively few pixels in the grid. Thus most of the values $r_{jk}$ are equal to zero.

In order to solve the system we will use two different methods depending on whether the system is overdetermined or underdetermined. In the first case we will use the least square approximation, that means that the solution $x$ will be given by $x=\text{argmin}_{y}\norm{Ay-b}$. In the second case we will use an iterative method called \emph{Kaczmarz's method} that will be discussed in the next paragraph. 
%We will not discuss here in detail the lsqr method, for more information we remand you at CITE.

\subsection{Kaczmarz's method}
The Kaczmarz's method \cite{KACZ} is an iterative procedure for approximating a solution of a linear system $Ax=p$, $A\in\R^{m\times n}$, $p\in\R^{m}$. If we denote $r_{i}$ the $i$-th row of $A$ and $p_{i}$ the $i$-th component of $p$, we can say that a vector $x$ is solution of $Ax=p$ if and only if
\begin{equation*} 
r_{i}\cdot x=p_{i} \quad \forall i=1,\ldots,m.
\end{equation*}
We also notice that the set $L_{i}=\{ x\in\R^{n}:\, r_{i}\cdot x=p_{i} \}$ is an affine subspace of $\R^{n}$. The idea of Kaczmarz's method is to project an initial approximated solution $x_{0}$ on all these affine spaces, generating in this way a sequence of vectors, each of them satisfies one of the equations $r_{i}\cdot x=p_{i}$.% The following definition describe what an affine projection is:
\begin{definition}
Let $L_{p,r}=\{ x\in\R^{n}:\, r\cdot x=p \}$, for $r\in\R^{n}$ and $p\in\R$, an affine space, let $u\in\R^{n}$. The \emph{affine projection} of $u$ on $L_{p,r}$ is the vector $\bar{u}\in L_{p,r}$ such that 
\begin{equation*}
\norm{\bar{u}-u}_{2}=\min_{x\in L_{p.r}}{\norm{x-u}_{2}}.
\end{equation*}
\end{definition}
\begin{proposition}
The affine projection $\bar{u}$ of a vector $u$ in the affine space $L_{p,r}$ is given by
\begin{equation*}
\bar{u}=u-\frac{r\cdot u-p}{\norm{r}_{2}^{2}}r.
\end{equation*}
\end{proposition}

The Kaczmarz's method proceeds as following. From an initial guess it computes its affine projection on the first affine space. This projection is then projected on the next affine space in our list and so on until the last affine space. These operations consist of one iteration and the result of this iteration become the starting point of the next one. In detail the algorithm proceed as follow:
\begin{enumerate}
\item Select $x_{0}$;
\item for $k=1,\ldots,K_{\max}$, $x_{k-1}^{0}=x_{k-1}$ (where $K_{\max}$ is the maximum number of iteration allowed);
\item for $i=1,\ldots,m$, 
\begin{equation}
x_{k}^{i}=x_{k}^{i-1}-\frac{r_{i}\cdot x_{k}^{i-1}-p_{i}}{\norm{r_{i}}^{2}}r_{i}
\label{eq: kaczmarz}
\end{equation}
\item $x_{k}=x_{k-1}^{m}$.
\end{enumerate}
The sequence $x_{0},x_{1},x_{2}\ldots$ generated by the method converges to a vector $x$ that satisfies $Ax=p$ (see Theorem 9.14 in \cite{BASIC} and references there).
However the convergence can be slow and a lot of steps are needed to get a good approximation. Moreover if the system has no solution, like in many image reconstruction applications, then the behavior of the sequence is not clear and can be chaotic.

In the field of medical imaging the size of the system can be a serious problem, but, as we know, the matrix $A$ is also sparse. This means that when we compute $x_{k}^{i}$ from $x_{k}^{i-1}$, we only change the components of $x_{k}^{i-1}$ that correspond to non zero entries of $r_{i}$. So we can increase efficiency storing the location of these entries.

Another fact that is connected to the nature of the reconstruction problem is that adjacent X-ray beams transmitted along line $l_{t,\theta}$, for similar values of $t$ and $\theta$, will intersect many of the same pixels, thus the corresponding affine spaces will be almost parallel. As a consequence, the convergence is slow and a lot of iteration is needed to reach a good approximate image.

We conclude this section introducing a variation of the Kaczmarz's method that involves the introduction of a relaxation parameters in the formula \eqref{eq: kaczmarz}. Let $\lambda_{i,k}$ be such that $0<\lambda_{i,k}<2$, then we replace formula \eqref{eq: kaczmarz} with 
\begin{equation*}
x_{k}^{i}=x_{k}^{i-1}-\lambda_{i,k}\frac{r_{i}\cdot x_{k}^{i-1}-p_{i}}{\norm{r_{i}}^{2}}r_{i}.
\end{equation*}
The parameter $\lambda_{i,k}$ can accelerate the convergence of an indeterminate system. Note that if $\lambda_{i,k}=2$, then the vector $x_{k}^{i}$ is just the reflection of $x_{k}^{i-1}$ across $L_{i}$ and there is no improvement in the proximity to a solution. That's why we consider $0<\lambda_{i,k}<2$.

%Chapter 4 - Kernel based methods\documentclass[../thesis01.tex]{subfiles} 
\chapter{Kernel based methods} \label{chap: kernelMethods}
%\section{Introduction}
In this chapter we present another approach for solving the image reconstruction problem based on kernel functions. Reproducing kernels have already been used in image reconstruction \cite{REI}. Here we use a different approach.

As usual our data are the discrete Radon transform of a function $f:\R^2\rightarrow\R$ $\{Rf(t_{j},\theta_{j})\}_{j=1}^{n}$, from which we want to find an approximation of the function $f$. 

The basic idea is to seek for the approximation $s$ of $f$ in a functions space $S$ with finite dimension $n$, that is $S=span\{s_{1},s_{2},\ldots, s_{n}\}$. Thus a function $s\in S$ can be written as $s=\sum_{j=1}^{n}{c_{j}s_{j}}$ for some $c_{j}\in\R$.

Then we ask that $Rs$ coincides with $Rf$ on the points $(t_{j},\theta_{j})$ for all $j=1,\ldots,n$, i.e.
\begin{equation}
(Rs)(t_{j},\theta_{j})=(Rf)(t_{j},\theta_{j}) \quad j=1,\ldots,n.
\label{eq: kerProblem}
\end{equation} 
By linearity of $R$ the coefficients $c_{j},\ j=1,\ldots,n$ are given by the solution of the linear system $Ac=b$, $c\in\R^{n}$, where 
\begin{align*}
A&=(Rs_{k}(t_{j},\theta_{j}))_{j,k=1,\ldots,n} & b&=(Rf(t_{j},\theta_{j}))_{j=1,\ldots,n}.
\end{align*}

\section{Hermite-Birkhoff interpolation}
We generalize the image reconstruction problem \eqref{eq: kerProblem} considering the problem of finding a function $s\in S$ such that $f|_{\Lambda}=s|_{\Lambda}$ for some function $f$, where $\Lambda=\{\lambda_{1},\ldots,\lambda_{n}\}$ is a set of linearly independent linear functionals (see \cite{ISKE3,ISKE4}). In our specific case we will consider $\lambda_{j}f=Rf(t_{j},\theta_{j})$. We also assume $S=span\{s_{1},s_{2},\ldots, s_{n}\}$ with $|\Lambda|=n$. By linearity, the problem is equivalent to the linear system \begin{equation*}
Ac=f_{\Lambda}, 
\label{eq: HermBirk}
\end{equation*}
 where
\begin{align*}
A&=(\lambda_{j}s_{k})_{j,k=1,\ldots,n} & f_{\Lambda}&=(\lambda_{j}f)_{j=1,\ldots,n} & s=\sum_{j=1}^{n}{c_{k}s_{k}}.
\end{align*}

\begin{theorem}[Mairhuber-Curtis \cite{MAIR}, \cite{CURTIS}]
Let $\Omega\subseteq\R^{d}$, $d\geq2$, suppose $\Omega$ contains a interior point, then there is no Chebjichev system $s_{1},\ldots,s_{n}$, $n\geq2$ on $\Omega$, i.e. for all $s_{1},\ldots,s_{n}$ real valued functions on $\Omega$, exists $\Lambda=\{\lambda_{1},\ldots,\lambda_{n}\}$, $|\Lambda|=n$ such that the matrix $(\lambda_{k}(s_{j}))_{j,k=1,\ldots,n}$ is singular, where $\lambda_{k}=\delta_{\xi_{k}}$ for pairwise distinct points $\xi_{k}\in\Omega$ are the functionals "evaluation at $\xi_{k}$".  
\end{theorem}
This theorem tells us that if we want to find a basis $s_{1},\ldots,s_{n}$ of $S$ such that the system \eqref{eq: HermBirk} has a unique solution for all data $\Lambda$, the basis should depends on the location of the data, i.e. on $\Lambda$ itself. For this reason we will choose 
\begin{equation*}
s_{j}=\lambda_{j}^{y}K(\cdot,y),
\end{equation*}
where $K:\R^{d}\times\R^{d}\rightarrow\R$ and $\lambda_{j}^{y}$ indicates that operator $\lambda_{j}$ is applied to variable $y$.

\subsection{Positive definite kernels}
The problem is then solving $f|_{\Lambda}=s|_{\Lambda}$ with $s=\sum_{j}{c_{j}\lambda_{j}^{y}K(\cdot,y)}$. It can be also written as a linear system $A_{K,\Lambda}c=f_{\Lambda}$, where $c\in\R^{n}$ and $A_{K,\Lambda}=(\lambda_{k}^{x}\lambda_{j}^{y}K(x,y))_{j,k=1}^{n}$.

In order to have a unique solution for every choice of $\Lambda$, $A_{K,\Lambda}$ must be non-singular. This is certain true if we assume $K$ symmetric, i.e. $K(x,y)=K(y,x)$ for all $x,y\in\R^{d}$, and positive definite, that means that $A_{K,\Lambda}$ is positive definite for all $\Lambda$. The following theorem gives us a characterization of positive definite functions.
\begin{theorem}[Bochner]
Assume $\Phi:\R^{d}\rightarrow\R$ even and continuous and that its Fourier transform $\hat{\Phi}$ is such that the Fourier inversion theorem holds:
\begin{equation*}
\Phi(x)=\frac{1}{(2\pi)^{d}}\int_{\R^{d}}{\hat{\Phi}(\omega)e^{ix\omega}\,d\omega},
\end{equation*}
then, if $\hat{\Phi}(\omega)\geq0$, $\forall\omega\in\R^{d}$, $K(x,y)=\Phi(\norm{x-y})$ is positive definite.
\end{theorem}

\begin{example}[Gaussian]
$\Phi(x)=e^{-\norm{x}^{2}}$ is positive definite since $\hat\Phi(\omega)=e^{-\frac{\norm{\omega}^{2}}{4}}>0$;
\end{example}

\begin{example}[Inverse multiquadric]
$\Phi(x)=\frac{1}{\sqrt{1+\norm{x}^{2}}}$ is positive definite since $\hat\Phi(\omega)=K_{\frac{d-1}{2}}(\norm{\omega})\norm{\omega}^{-\frac{d-1}{2}}>0$, where $K_{\nu}$ denotes the Bessel function of 2nd kind of order $\nu$.
\end{example}

\section{Conditionally positive definite kernels}
\begin{definition}
A set of functional $\Lambda$ is said to be $k$-unisolvent (or unisolvent w.r.t $\Poly^{k}_{d}$) if for $p\in\Poly_{k}^{d}$ we have
\begin{equation*}
p|_{\Lambda}=0 \quad \Rightarrow \quad p\equiv0,
\end{equation*}
where $\Poly_{k}^{d}$ denotes the set of all polynomial of degreed less or equal of $k$ on $\R^{d}$.
\end{definition}
\begin{definition}
The radial kernel $\Phi$  is conditionally positive definite of order $k$ and we write $\Phi\in cdp(k)$, if for $K(x,y)=\Phi(\norm{x-y})$ the quadratic form
\begin{equation*}
c^{T}A_{K,\Lambda}c=\sum_{i,j=1}^{n}{c_{i}c_{j}\lambda_{i}^{x}\lambda_{j}^{y}\Phi(\norm{x-y})}
\end{equation*}
is positive for all possible $\Lambda$, $|\Lambda|=n$ and vectors $c\in\R^{n}\setminus\{0\}$ satisfying 
\begin{equation*}
\sum_{j=1}^{n}{c_{j}\lambda_{j}(p)}=0, \qquad \forall\,p\in\Poly_{k}^{d}.
\end{equation*}
\end{definition}
\subsection{Reconstruction by conditionally positive kernel functions}
We introduce a polynomial part in the interpolant $s$, so what we have to do now is to solve $f|_{\Lambda}=s|_{\Lambda}$ with $s$ of the form
\begin{equation*}
s=\sum_{j=1}^{n}{c_{j}\lambda_{j}^{y}\Phi(\norm{\cdot-y})}+p,
\end{equation*}
where $p\in\Poly^{d}_{k-1}$ ($k=k(\Phi)$) and vector $c\in\R^{n}$ satisfying the \emph{vanishing moment condition}
\begin{equation*}
\sum_{j=1}^{n}{c_{j}\lambda_{j}(p)}=0 \quad \forall\, p\in\Poly_{k-1}^{d}.
\end{equation*}
\begin{theorem}[Michelli,Wu\label{thm: unicity}]
The reconstruction problem $f|_{\Lambda}=s|_{\Lambda}$ has under vanishing moment condition a unique solution $s$, provided that $\Phi\in cdp(k)$ and the functionals $\Lambda$ are $(k-1)$-unisolvent.
\end{theorem}
\proof
Let $p_{1},\ldots,p_{m}\in\Poly_{k-1}^{d}$ a basis of $\Poly_{k-1}^{d}$, where $m=\text{dim}(\Poly_{k-1}^{d})=\left(\begin{array}{c}
(k-1)+d\\
d
\end{array}\right)$. Then $s$ can be written as
\begin{equation*}
s=\sum_{j=1}^{n}{c_{j}\lambda_{j}^{y}\Phi(\norm{\cdot-y})}+\sum_{l=1}^{m}{d_{l}p_{l}},
\end{equation*}
for some $c\in\R^{n},\ d\in\R^{m}$. Condition $f_{\Lambda}=s|_{\Lambda}$ and the vanishing moment condition, are equivalent to the linear system
\begin{equation*}
\left\{
\begin{aligned}
&\lambda_{i}^{x}s=\lambda_{i}^{x}f & &\forall\, i=1,\ldots,n\\
&\sum_{j=1}^{m}{c_{j}\lambda_{j}^{x}(p_{l})}=0 & &\forall\,l=1,\ldots,m
\end{aligned}
\right.
\end{equation*}
that in matricial form is
\begin{equation}
\left(
\begin{array}{cc}
\Phi & P\\
P^{T} & O
\end{array}
\right)
\left(
\begin{array}{c}
c\\
d
\end{array}
\right)=\left(
\begin{array}{c}
f|_{\Lambda}\\
\textbf{0}
\end{array}
\right)
\label{eq: recSyst}
\end{equation}
where $\Phi_{i,j}=\lambda_{i}^{x}\lambda_{j}^{y}\Phi(\norm{x-y})$, $i,j=1,\ldots,n$, $P_{j,l}=\lambda_{j}(p_{l})$, $j=1,\ldots,n$, $l=1,\ldots,m$ and $O\in\R^{m}\times\R^{m}$ and $\textbf{0}\in\R^{m}$ are a matrix and a vector with all components equal to zero.

We consider the homogeneous system
\begin{equation*}
\left\{
\begin{aligned}
&\Phi c+Pd=\textbf{0}\\
&P^{T}c=\textbf{0}
\end{aligned}
\right.\Leftrightarrow
\left\{
\begin{aligned}
&c^{T}\Phi c+c^{T}Pd=0\\
&c^{T}P=\textbf{0}
\end{aligned}
\right.
\end{equation*}
substituting the second equation in the first one, we have $c^{T}\Phi c=0$. Since $\Phi\in cdp(k)$, $c=\textbf{0}$. The first equation becomes then $Pd=\textbf{0}$, but $\Lambda$ is $(k-1)$-unisolvent, that implies $d=\textbf{0}$. So the unique solution to the homogeneous system is the null solution.
\endproof

\begin{example} Conditionally positive functions:
\begin{enumerate}
\item Polyharmonic splines:
\begin{equation*}
\varphi(r)=\left\{
\begin{aligned}
&r^{2k-d}\log{r} & &\text{if $d$ is even}\\
&r^{2k-d} & &\text{if $d$ is odd}
\end{aligned}
\right.
\end{equation*}
$2k>d$;
\item Gaussian: $\varphi(r)=e^{-r^{2}}$, $k=0$;
\item Multiquadrics: $\varphi(r)=(1+r^{2})^{\nu}$, $\nu>0,\ \nu\notin\N$, $k=\lceil\nu\rceil$;
\item Inverse multiquadrics: $\varphi(r)=(1+r^{2})^{\nu}$, $\nu<0$, $k=0$;
\item Power function: $\varphi(r)=r^{\beta}$, $0<\beta\notin2\N$, $k=\lceil\frac{\beta}{2}\rceil$.
\end{enumerate}
\end{example}

\section{Native function spaces}
In this section we will show that the solution of the Hermite-Birkhoff interpolation is optimal in the sense that it is the function of minimum norm among all functions that interpolate data $f|_{\Lambda}$, where the norm is taken in a suitable Hilbert space.

For a fixed positive definite function $K:\R^{d}\times\R^{d}\rightarrow\R$, we define the function spaces
\begin{align*}
S_{\Lambda}&=span\{\lambda{y}K(.,y):\ \lambda\in\Lambda\} & S&=\{s\in S_{\Lambda}: \ |\Lambda|<\infty\}
\end{align*} 
and the dual space
\begin{equation*}
L=\{\lambda\equiv\lambda_{c,\Lambda}=\sum_{j=1}^{n}{c_{j}\lambda_{j}: \ c\in\R^{n}, \ |\Lambda|<\infty}\}.
\end{equation*}
We observe that, for all $s\in S$ there exists $\lambda\in L$ such that $s\equiv s_{\lambda}=\lambda^{y}K(\cdot,y)$. Indeed
\begin{equation*}
s=\sum_{j=1}^{n}{c_{j}\lambda_{j}^{y}K(\cdot,y)=\sum_{j=1}^{n}{c_{j}\lambda_{j}^{y}K(\cdot,y)}  =\left(\sum_{j=1}^{n}{c_{j}\lambda_{j}}\right)^{y}K(\cdot,y)}=\lambda^{y}K(\cdot,y),
\end{equation*}
with $\lambda=\sum{c_{j}\lambda_{j}}\in L$.

We define an inner product on $L$:
\begin{equation*}
(\lambda,\mu)_{K}=\lambda^{x}\mu^{y}K(x,y)=\sum_{j=1}^{n}{\sum_{k=1}^{n}{c_{j}d_{k}\lambda_{j}^{x}\mu_{k}^{y}K(x,y)}},
\end{equation*}
where  $\lambda=\sum{c_{j}\lambda_{j}}$ and $\mu=\sum{d_{k}\mu_{k}}$, and the norm
\begin{equation*}
\norm{\lambda}_{K}=(\lambda,\lambda)_{K}^{1/2}.
\end{equation*}
Thanks to the duality relation between $L$ and $S$, we introduce a topology also on $S$ so that $(s_{\lambda},s_{\mu})_{K}=(\lambda,\mu)_{K}$, $\norm{\cdot}_{K}=(\cdot,\cdot)_{K}^{1/2}$

\begin{remark}
$L\cong S$: $L$ and $S$ are isometric with respect to the norm $\norm{\cdot}_{K}$. Moreover, for all $\mu\in L$, $\mu$ is continuous on $S$. In fact
\begin{equation*}
|\mu(s_{\lambda})|=|\mu^{x}\lambda^{y}K(x,y)|=|(\mu,\lambda)_{K}|\leq\norm{\mu}_{K}\norm{\lambda}_{K}=\norm{\mu}_{K}\norm{s_{\lambda}}_{K}.
\end{equation*}
\end{remark}

We now set $D=\bar{L}$ and $F=\bar{S}$ the topological closures with respect to $\norm{\cdot}_{K}$. The following theorem holds:
\begin{theorem}[Madych-Nelson, 1983\label{thm: reprDF}]
For all $s_{\mu}\in S$, for all $\lambda\in D$ we have 
\begin{equation}
(\lambda^{y}K(\cdot,y),s_{\mu})_{K}=(s_{\lambda},s_{\mu})_{K}=(\lambda,\mu)_{K}=\lambda^{x}\mu^{y}K(x,y)=\lambda(s_{\mu}).
\label{eq: reprDF}
\end{equation}
\end{theorem}
\proof
The statements holds for $\lambda,\mu\in L$ and the representation \eqref{eq: reprDF} follows by continuity. 
\endproof

\begin{definition}\label{def: rep_kernel}
Let $H=\{f:\ \Omega\subseteq\R^{d}\rightarrow\R\}$ a Hilbert space, $K:\Omega\times\Omega\rightarrow\R$ is a \emph{reproducing kernel} for $H$ if:
\begin{enumerate}
\item $K(\cdot,x)\in H$ for all $x\in\Omega$;
\item $f(x)=(f,K(\cdot,x))_{H}$, for all $f\in H,\ x\in\Omega$.
\end{enumerate}
\end{definition}
\begin{corollary}
$K:\R^{d}\times\R^{d}\rightarrow\R$ is the reproducing kernel of the Hilbert space $F$.
\end{corollary}
\proof We prove properties 1. and 2. of Definition \ref{def: rep_kernel}
\begin{enumerate}
\item For $\delta_{z}\in L, \ z\in\R^{d}$, $\delta_{z}^{y}K(\cdot,y)=K(\cdot,z)\in F$, for all $z\in\R^{d}$;
\item For $\lambda=\delta_{z}\in L$, by theorem \ref{thm: reprDF} $(K(\cdot,z),f)_{K}=f(z)$, for all $f\in F,\ z\in\R^{d}$
\end{enumerate}
\endproof
\begin{corollary}
The point evaluation $\delta_{z}:F\rightarrow\R$ are continuous on $F$.
\end{corollary}
\proof
$|\delta_{z}(f)|=|f(z)|\leq\norm{\delta_{z}}_{K}\norm{f}_{K}$ for all $z,f$.
\endproof
\begin{corollary}
Let $s\equiv s_{f,\Lambda}\in S$ denote the unique interpolation to $f\in F$ on $\Lambda$: $s|_{\Lambda}=f|_{\Lambda}$, then the Pythagoras theorem
\begin{equation*}
\norm{f}^{2}_{K}=\norm{s}^{2}_{K}+\norm{f-s}^{2}_{K}
\end{equation*}  
holds.
\end{corollary}
\proof
For $s=\lambda^{y}K(\cdot,y)\in S$, $\lambda\in L$, we find $(s,g)_{K}=0$ for all $g$ such that $\lambda(g)=0$, i.e. $s\equiv s_{f,\Lambda}$ is orthogonal to the kernel of the functional $\lambda\in\Lambda$. Hence  
\begin{equation*}
(s_{f,\Lambda},f-s_{f,\Lambda})_{K}=0,
\end{equation*}
i.e.the interpolant $s_{f,\Lambda}\in S$ is the orthogonal projection of $f\in F$ onto $S$.   
\endproof
\subsection{Optimality of the interpolation method}
The following results are consequences of Theorem \ref{thm: reprDF} and its corollaries.
\begin{theorem}
The interpolant $s\equiv s_{f\Lambda}\in S$ is the unique minimizer of the energy functional $\norm{\cdot}_{K}$ among all interpolants to data $f|_{\Lambda}$, i.e. 
\begin{equation*}
\norm{s}_{K}\leq \norm{g}_{K} \ \forall g\in S \ \text{s.t.}\ g|_{\Lambda}=f|_{\Lambda}.
\end{equation*}
In this sense the interpolation scheme is optimal.
\end{theorem}
\begin{corollary}
The interpolant $s\equiv s_{f\Lambda}\in S$ is the unique best approximation to $f\in F$ from $S$ with respect to $\norm{\cdot}_{K}$. 
\end{corollary}

%Chapter 5 - Kernel based image reconstruction

\chapter{Kernel based image reconstruction}\label{chap: kernelRec}
In this chapter we apply the kernel based methods saw in Chapter \ref{chap: kernelMethods} to the problem of image reconstruction. We will see that the Hermite-Birkhoff interpolation can not be applied to the original reconstruction problem because the Radon transform of a kernel basis function can be infinity. We will then overcome this obstacle introducing a regularization of the integrals involved in the computation of the Radon transform. Thanks to this regularization it is possible to generate a liner system, solving the linear system one can find an approximation of the image to reconstruct.

This technique can be used with both parallel beam geometry and scattered data. This second case is useful when one wants to reduce the dosage of X-rays passing through the sample. Scattered data can then be interpolated using suitable methods. For example a radial functions method was recently introduced by Beatson and zu Castell \cite{ZUC} to obtain new values of the Radon transform. This technique can be combined with the methods introduced below to obtain a reconstruction using less initial data.
\newline

Let $f:\R^{2}\rightarrow\R$ be a function. Consider again the problem $s|_{\Lambda}=f|_{\Lambda}$, where 
\begin{align*}
f|_{\Lambda}&=\{\lambda_{j}f\}_{j=1}^{n}, & \lambda_{j}f&=R[f(x_{1},x_{2})](t_{j},\theta_{j}), \ j=1,\ldots,n,
\end{align*}
and $s$ is an approximation of $f$ belonging to the space
\begin{equation*}
S=\left\{\sum_{j=1}^{n}{c_{j}\lambda_{j}^{y}K(\cdot,y)}:\ c_{j}\in\R,\ K(x,y)=\Phi(\norm{x-y})\, \text{positive definite}\right\}.
\end{equation*}
Let us denote $b_{j}(x)=\lambda_{j}^{y}(K(x,y))=R[K(x,y)](t_{j},\theta_{j})$, $x\in\R^{2}$, $1\leq j\leq n$, the basis of $S$, so that $s(x)=\sum_{j=1}^{n}{c_{j}b_{j}(x)}$ for some $c\in\R^{n}$. The interpolation conditions $s|_{\Lambda}=f|_{\Lambda}$ are equivalent to $\lambda_{k}s=\lambda_{k}f\ \forall k=1,\ldots,n$. By linearity of the Radon transform we obtain 
\begin{equation*}
\sum_{j=1}^{n}{c_{j}\lambda^{x}_{k}\lambda_{j}^{y}K(x,y)}=\lambda_{k}f, \ k=1,\ldots,n,
\end{equation*} 
or in matrix form
\begin{align}
Ac&=f_{\Lambda}
\label{eq: sist_kernel}
\end{align}
with $ f_{\Lambda}=(\lambda_{1}f,\ldots,\lambda_{n}f)^{t}$ and $A=(a_{k,j})_{1\leq k,j \leq n}$ given by
\begin{equation*}
a_{k,j}=\lambda_{k}^{x}\lambda_{j}^{y}K(x,y)=R^{x}\left\{R^{y}[K(x,y)](t_{j},\theta_{j})\right\}(t_{k},\theta_{k}).
\end{equation*}
Thus, to determine $c$ and then solution $s$, we have to solve the linear system \eqref{eq: sist_kernel}.

The first step in solving system \eqref{eq: sist_kernel} is of course computing the matrix $A$. We start by considering a generic basis function $b_{j}(x)=R^{y}[K(x,y)](t_{j},\theta_{j})$ (for simplicity of notation we omit index $j$ and so we denote $(t_{j},\theta_{j})=(t,\theta)$).

We notice that since the kernel $K$ is of the form $K(x,y)=\Phi(\norm{x-y})$, we can use the shift property of the Radon transform to simplify the computation of $b_{j}$. Indeed, if we set $k(y)=K(0,y)=\Phi(\norm{y})$, then $K(x,y)=\Phi(\norm{x-y})=k(x-y)=k(y-x)$. Hence, by theorem \ref{thm: shiftProp} (shift property), if $g(t,\theta)=R[k(y)](y,\theta)$, we have
\begin{equation*}
R^{y}[K(x,y)](t,\theta)=R^{y}[k(y-x)](t,\theta)=g(t-x\cdot v,\theta), 
\end{equation*}
where $v=(\cos{\theta},\sin{\theta})$. So, in order to obtain $b_{j}$ we have only to compute $R[k(y)](t,\theta)=R[K(0,y)](t,\theta)$. Notice that this property is independent of the particular kind of kernel (Gaussian, multiquadrics, etc.) used and so is applicable with any kernel function of the form $K(x,y)=\Phi(\norm{x-y})$.

\section{Gaussian kernel reconstruction} \label{sec: gaussRec}
We start considering the Gaussian kernel 
\begin{equation*}
K(x,y)=e^{-\norm{x-y}^2}.
\end{equation*}
\begin{align*}
R[k(y)](t,\theta)&=\int_{\R}{k(t\cos{\theta}-s\sin{\theta},t\sin{\theta}+s\cos{\theta})\,ds}=\\
&=\int_{\R}{\exp{(-({\theta}-s\sin{\theta})^{2}-(t\sin{\theta}+s\cos{\theta})^{2})}\,ds}=\\
&=\int_{\R}{e^{-(t^{2}+s^{2})}\,ds}=e^{-t^2}\int_{\R}{e^{-s^{2}}\,ds}=\\
&=\sqrt{\pi}e^{-t^{2}}=g(t,\theta).
\end{align*}
Thus
\begin{align*}
R^{y}[K(x,y)](t,\theta)&=g(t-x\cdot v,\theta)=\sqrt{\pi}e^{-(t-x\cdot v)^{2}}
\end{align*}
and we conclude that
\begin{align*}
&\boxed{
b_{j}(x)=\sqrt{\pi}e^{-(t_{j}-x\cdot v_{j})^{2}} } & &\text{where}\ v_{j}=(\cos{\theta_{j},\sin{\theta_{j}}}).
\end{align*}

We now want to compute $a_{k,j}=R[b_{j}](t_{k},\theta_{k})$. Again for simplicity of notation, we write $(t_{j},\theta_{j})=(t,\theta)$ and $(t_{k},\theta_{k})=(r,\varphi)$, then
\begin{align*}
R[\sqrt{\pi}e^{-(t-x\cdot v)^{2}}](r,\varphi)&=\sqrt{\pi}\int_{\R}\exp\left(-\left[t-(r\cos{\varphi}-s\sin{\varphi})\cos{\theta}\right.\right.&\\
&\qquad\left.\left.-(r\sin{\varphi}+s\cos{\varphi})\sin{\theta}\right]^{2}\right)\,ds=\\
&=\sqrt{\pi}\int_{\R}\exp\left(-\left[t-r(\cos{\varphi}\cos{\theta}+\sin{\varphi}\sin{\theta})+\right.\right.\\
&\qquad\left.\left.+s(\sin{\varphi}\cos{\theta}-\cos{\varphi}\sin{\theta})\right]^{2}\right)\,ds=\\
&=\sqrt{\pi}\int_{\R}{\exp{(-[t-r\cos{(\varphi-\theta)}+s\sin{(\varphi-\theta)}]^{2})}\,ds}.
\end{align*}
If we set $a=\sin{(\varphi-\theta)}$ and $b=t-r\cos{(\varphi-\theta)}$, we can write
\begin{equation*}
a_{k,j}=R[\sqrt{\pi}e^{-(t-x\cdot v)^{2}}](r,\varphi)=\sqrt{\pi}\int_{\R}{e^{-(as+b)^{2}}\,ds}.
\end{equation*}
Hence, if $a\neq0$ we have
\begin{equation*}
a_{k,j}=\frac{\sqrt{\pi}}{a}\int_{\R}{e^{-(as+b)^{2}}\,d(as+b)}=\frac{\sqrt{\pi}}{a}\int_{\R}{e^{-u^{2}}\,du}=\frac{\pi}{a},
\end{equation*}
while, in the case in which $a=0$,
\begin{equation*}
a_{k,j}=\sqrt{\pi}\int_{\R}{e^{-b^{2}}\,ds}=\infty,
\end{equation*}
since both $\varphi$ and $\theta$ are in $[0,\pi)$, $a=0$ if and only if $\varphi=\theta$, so we conclude
\begin{equation*}
\boxed{
a_{k,j}=\left\{
\begin{aligned}
&\frac{\pi}{\sin{(\theta_{k}-\theta_{j})}} & &\text{if}\ \theta_{k}\neq\theta_{j}\\
&+\infty & &\text{if}\ \theta_{k}=\theta_{j}.
\end{aligned}
\right.
}
\end{equation*}

\subsection{Regularization}\label{subsec: regularization}
We saw that matrix $A$ can have infinity entries. More precisely $R[b_{j}](t_{k},\theta_{k})=+\infty$ for some values of $j$ and $k$ (those values s.t. $\theta_{k}=\theta_{j}$), that means that for these values $b_{j}(x)$ is not integrable on line $l_{t_{k},\theta_{k}}$. %From the graphical point of view we see that  ... FIGURE K(x,y) b_{j}(x) with lines

To overcome this obstacle we must find some regularization technique so that the value of the Radon transform of the basis elements $b_{j}$ is finite for all $k,j$. The simplest choice is to consider a truncation of the integral, i.e. computing 
\begin{equation*}
\int_{-\bar{L}}^{\bar{L}}{b_{j}(x(s))\,ds} \qquad \bar{L}\gg0,
\end{equation*}
instead of the integral on the whole real line. This approach is equivalent to compute $R[b_{j}(x)\chi_{[-L,L]}(\norm{x})](t_{k},\theta_{k})$, where
\begin{equation*}
\chi_{[-L,L]}(r)=\left\{
\begin{aligned}
&1 & &\text{if}\ -L\leq r\leq L\\
&0 & &\text{otherwise}
\end{aligned}
\right.
\end{equation*}
is the characteristic function of the set $[-L,L]$ for some $L>0$. 

Moreover, in general, we can multiply $b_{j}$ for a window function $w$, where $w$ is such that 
\begin{equation*}
\int_{l_{t_{k},\theta_{k}}}{b_{j}w}<\infty \quad \forall\,j,k.
\end{equation*}
Possible choices of $w$ are:
\begin{itemize}
\item the characteristic function of a compact set $w(x)=\chi_{[-L,L]}(\norm{x})$;
\item the Gaussian function $w(x)=e^{-\varepsilon^{2}\norm{x}^{2}}$;
\item the cosine window $w(x)=\cos{\frac{\pi\norm{x}}{2L}}\chi_{[-L,L]}{\norm{x}}$;
%\item hamming window
\end{itemize}
This approach can be interpreted also as substituting the operator $R$ with another operator, say $R_{w}$, defined by 
\begin{equation*}
R_{w}[f]=R[fw], \qquad \text{for all} \ f:\R^{2}\rightarrow\R.
\end{equation*}
Note that, since $R$ is a linear operator, also $R_{w}$ is so. Indeed, for all $f,g$ functions, for all $\alpha,\beta$ constants
\begin{align*}
R_{w}[\alpha f+\beta g]&=R[(\alpha f+\beta g)w]=R[\alpha fw+\beta gw]=\\
&=\alpha R[fw]+\beta R[gw]=\\
&=\alpha R_{w}[f]+\beta R_{w}[g].
\end{align*}
Then if we approximate $f_{k}=R[f](t_{k},\theta_{k})$ with $f_{w,k}=R_{w}[f](t_{k},\theta_{k})$, we can consider the interpolation problem 
\begin{equation*}
f_{k}\approx f_{w,k}=R_{w}[s](t_{k},\theta_{k})\qquad \forall\,k=1,\ldots,n.
\end{equation*} 
By linearity of $R_{w}$
\begin{equation*}
f_{k}\approx \sum_{j=1}^{n}{c_{j}R_{w}[b_{j}](t_{k},\theta_{k})}\qquad \forall\,k=1,\ldots,n,
\end{equation*} 
that leads us to the linear system $A_{w}c=f$, where $(A_{w})_{k,j}=R_{w}[b_{j}](t_{k},\theta_{k})=R[b_{j}w](t_{k},\theta_{k})$.

We notice that for all $k$, the difference between $f_{k}$ and $f_{w,k}$ is bounded by
\begin{align*}
|f_{w,k}-f_{k}|&=|R_{w}[f](t_{k},\theta_{k})-R[f](t_{k},\theta_{k})|=\\
&=\left| \int_{\R}{f(x(s))w(x(s))\,ds}-\int_{\R}{f(x(s))\,ds}\right|\leq\int_{\R}{|f||w-1|\,ds}\leq\\
&\leq\norm{w-1}_{\infty}\norm{f}_{L^{1}(\R^{2})}
\end{align*} 
Thus, for $w\rightarrow1$, $|f_{w,k}-f_{k}|\rightarrow0$ but also $R[b_{j}w](t_{k},\theta_{k})\rightarrow R[b_{j}](t_{k},\theta_{k})$ and this quantity can be infinity. For $w\rightarrow0$, $R[b_{j}w](t_{k},\theta_{k})\rightarrow0$ but the difference \begin{equation*}
|f_{w,k}-f_{k}|\rightarrow\left|\int{f(x(s))\,ds}\right|\leq\norm{f}_{L^{1}}.
\end{equation*}

Before starting on computing $A_{w}$ consider the following example. Let $K$ be the inverse multiquadric kernel given by
\begin{equation*}
K(x,y)=\frac{1}{\sqrt{1+\norm{x-y}^{2}}}.
\end{equation*}
As before we can just consider $R[k(y)]=R[K(0,y)]$ because of the relation $k(y-x)=K(x,y)$ and the shift property of the Radon transform. What we obtain is
\begin{equation*}
R[k(y)](t,\theta)=\int_{\R}{\frac{1}{\sqrt{1+t^{2}+s^{2}}}\,ds}=+\infty,
\end{equation*}
that means that in this case, not only $b_{j}(x)$ is not integrable on some line (as in the Gaussian kernel case), but even $K(x,y)$ is not integrable on any line $l_{t,\theta}$. In this case we have to consider a further regularization of the integral.

The remedy we adopt is to multiply the function $K$ itself by a window function $w$ such that $R^{y}[K(x,y)w](t,\theta)$ exists finite for all $(t,\theta)$. In choosing the function $w$, we consider that we would like to still use the shift property of the Radon transform, therefore we take $w$ of the form $w=w(x,y)=\tilde{w}(\norm{x-y})$ so that, if we set now $k(y)=K(0,y)w(0,y)$, it is still true that $k(y-x)=K(x,y)w(x,y)$. Moreover we will choose $w$ to be positive definite, in this way, since the product of positive definite function is positive definite, also $K\cdot w$ is so. 

Notice that this kind of regularization does not correspond, as in the first case, to replace the operator $R$ with another operator. In fact now the function $w$ depends on $x$. What we are doing now is simply substituting the positive definite kernel $K(x,y)$ with another positive definite kernel given by $K(x,y)w(x,y)$ that is integrable on every line in the plane $(y_{1},y_{2})$.

\subsection{Regularization by truncation}
We first consider the regularization of the Gaussian reconstruction problem using $w_{L}(x)=\chi_{[-L,L]}(\norm{x})$ as window function. With this choice $R[b_{j}](t_{k},\theta_{k})$ is replaced by $R_{L}[b_{j}](t_{k},\theta_{k})=R[b_{j}w_{L}](t_{k},\theta_{k})$, where $b_{j}(x)=R^{y}[K(x,y)](t_{j},\theta_{j})$ and K the Gaussian kernel.

In applications is useful to use kernels depending on a shape parameter $\varepsilon$ so that choosing suitable values of $\varepsilon>0$ one can obtain system matrix with a better condition number. We will consider 
\begin{equation*}
K(x,y)=e^{-\varepsilon^2\norm{x-y}^{2}}.
\end{equation*}
Basis $b_{j}$ then becomes
\begin{equation}
\boxed{
b_{j}(x)=\frac{\sqrt{\pi}}{\varepsilon}e^{-\varepsilon^{2}(t_{j}-x\cdot v_{j})^{2}}.
}
\label{eq: basis_gauss}
\end{equation}
Indeed 
\begin{align*}
R[k(y)]=\int_{\R}{e^{-\varepsilon^2(t^{2}+s^{2})}\,ds}=e^{-\varepsilon^{2}t^{2}}\int_{\R}{e^{-\varepsilon^{2}s^{2}}\,ds}=\frac{\sqrt{\pi}}{\varepsilon}e^{-\varepsilon^{2}t^{2}}.
\end{align*}
Components of matrix $A_{L}$ are given by $a_{k,j}=R_{L}[b_{j}](t_{k},\theta_{k})$, so we have
\begin{align*}
R_{L}[b_{j}](r,\varphi)=\int_{\R}{\frac{\sqrt{\pi}}{\varepsilon}e^{-\varepsilon^{2}(t-x(s)\cdot v)^{2}}\chi_{[-L,L]}(\norm{x(s)})\,ds},
\end{align*}
where $x(s)=(r\cos{\varphi}-s\sin{\varphi},r\sin{\varphi}+s\cos{\varphi})$ and $v=(\cos{\varphi},\sin{\varphi})$. If we set again
\begin{equation}
a=\sin{(\varphi-\theta)} \qquad b=t-r\cos{(\varphi-\theta)}
\label{eq: def_ab}
\end{equation}
we have
\begin{align*}
R_{L}[b_{j}](r,\varphi)&=\frac{\sqrt{\pi}}{\varepsilon}\int_{\R}{e^{-\varepsilon^{2}(as+b)^{2}}\chi_{[-L,L](\sqrt{r^{2}+s^{2}})}\,ds}=\\
&=\frac{\sqrt{\pi}}{\varepsilon}\int_{-\sqrt{L^{2}-r^{2}}}^{\sqrt{L^{2}-r^{2}}}{e^{-\varepsilon^{2}(as+b)^{2}}\,ds},
\end{align*}
where we are assuming $L\gg0$ so that $|r|<L$. Now we distinguish two cases: 
\begin{enumerate}
\item if $a=0$ then
\begin{equation*}
R_{L}[b_{j}](r,\varphi)=\frac{\sqrt{\pi}}{\varepsilon}\int_{-\sqrt{L^{2}-r^{2}}}^{\sqrt{L^{2}-r^{2}}}{e^{-\varepsilon^{2}b^{2}}\,ds}=\frac{\sqrt{\pi}}{\varepsilon}e^{-\varepsilon^{2}b^{2}}2\sqrt{L^{2}-r^{2}};
\end{equation*}
\item if $a\neq0$ we set $u=\varepsilon(as+b)$ so that 
\begin{equation*}
R_{L}[b_{j}](r,\varphi)=\frac{\sqrt{\pi}}{\varepsilon^{2}a}\int_{c_{1}}^{c_{2}}{e^{-u^{2}}\,du}
\end{equation*}
where $c_{1}=\varepsilon(-\sqrt{L^{2}-r^{2}}+b)$ and $c_{2}=\varepsilon(\sqrt{L^{2}-r^{2}}+b)$ and $\int_{c_{1}}^{c_{2}}{e^{-u^{2}}\,du}=\frac{\sqrt{\pi}}{2}\text{erf}{(c_{2})}-\text{erf}{(c_{1})}$, where \emph{erf} is the usual error function
\begin{equation*}
\text{erf}(x)=\frac{2}{\sqrt{\pi}}\int_{0}^{x}{e^{-u^{2}}\,du}.
\end{equation*}
\end{enumerate}

\begin{figure}[htbp]
\centering%
\subfigure[ \label{fig: orig_ph}]%
{\includegraphics[width=0.49\textwidth]{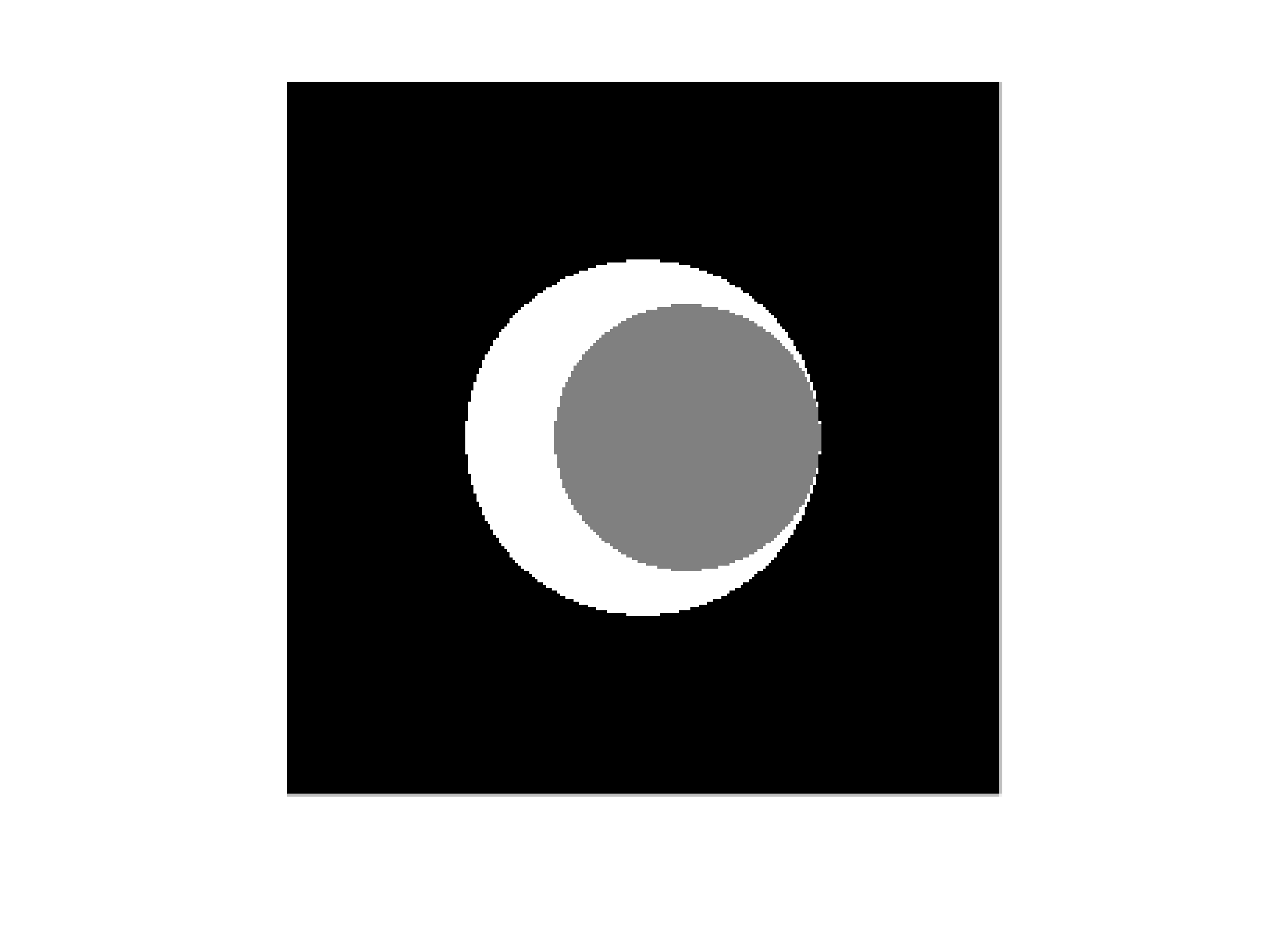}}%,height=0.3\textwidth
\subfigure[ \label{fig: gauss_trunc1}]%
{\includegraphics[width=0.49\textwidth]{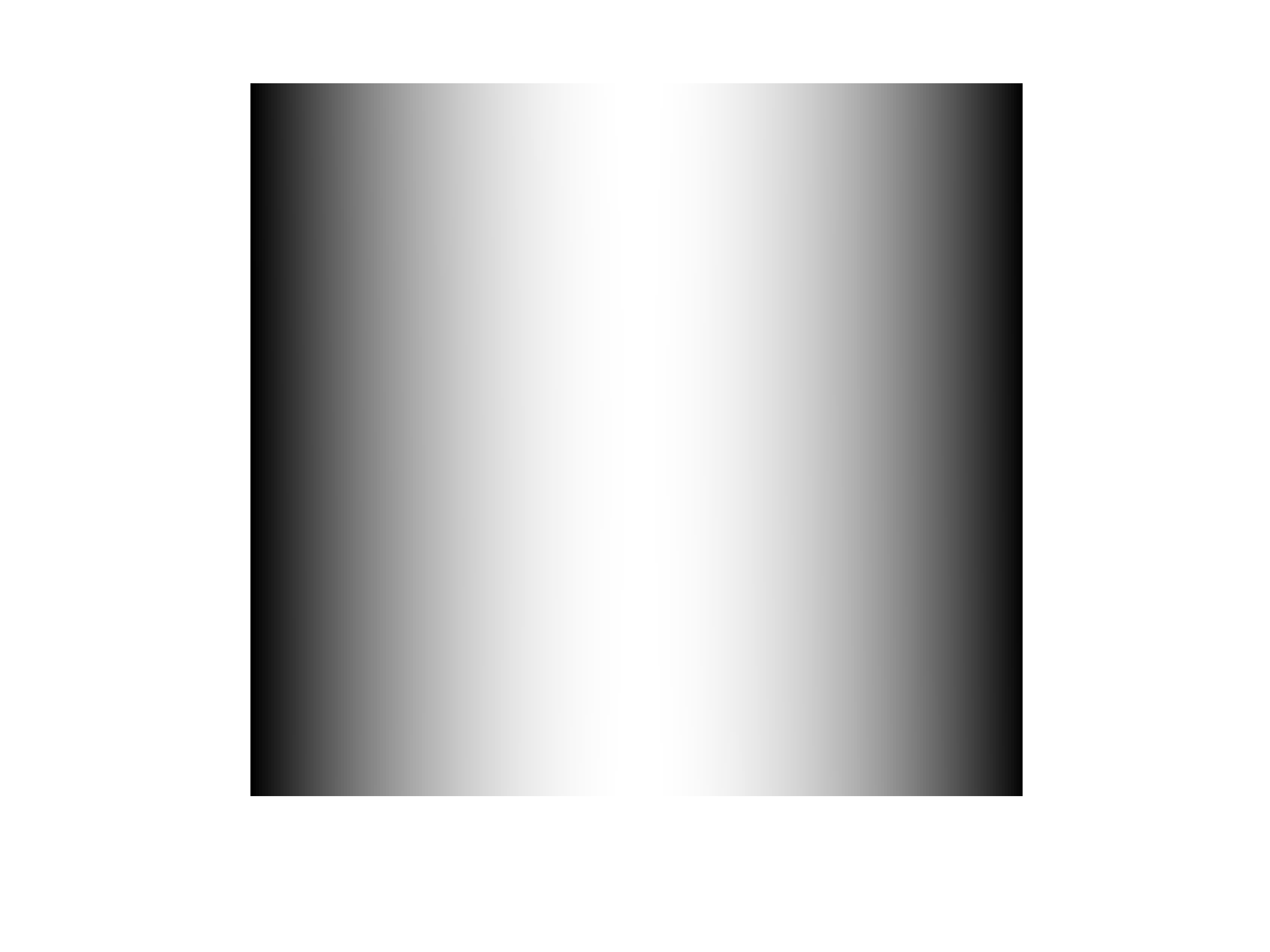}}
\subfigure[ \label{fig: gauss_trunc2}]%
{\includegraphics[width=0.49\textwidth]{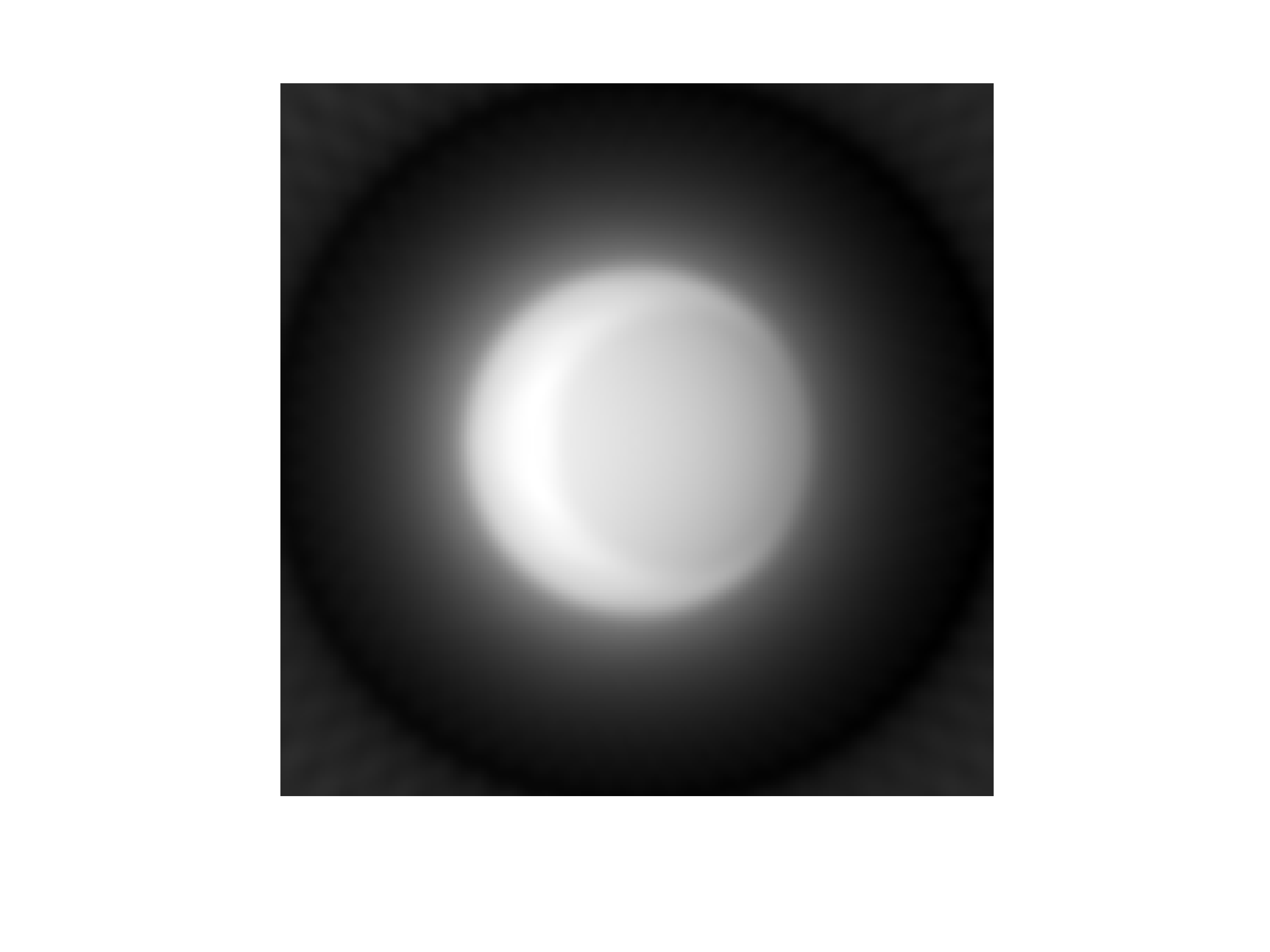}}%,height=0.3\textwidth
\subfigure[ \label{fig: mq_time3}]%
{\includegraphics[width=0.49\textwidth]{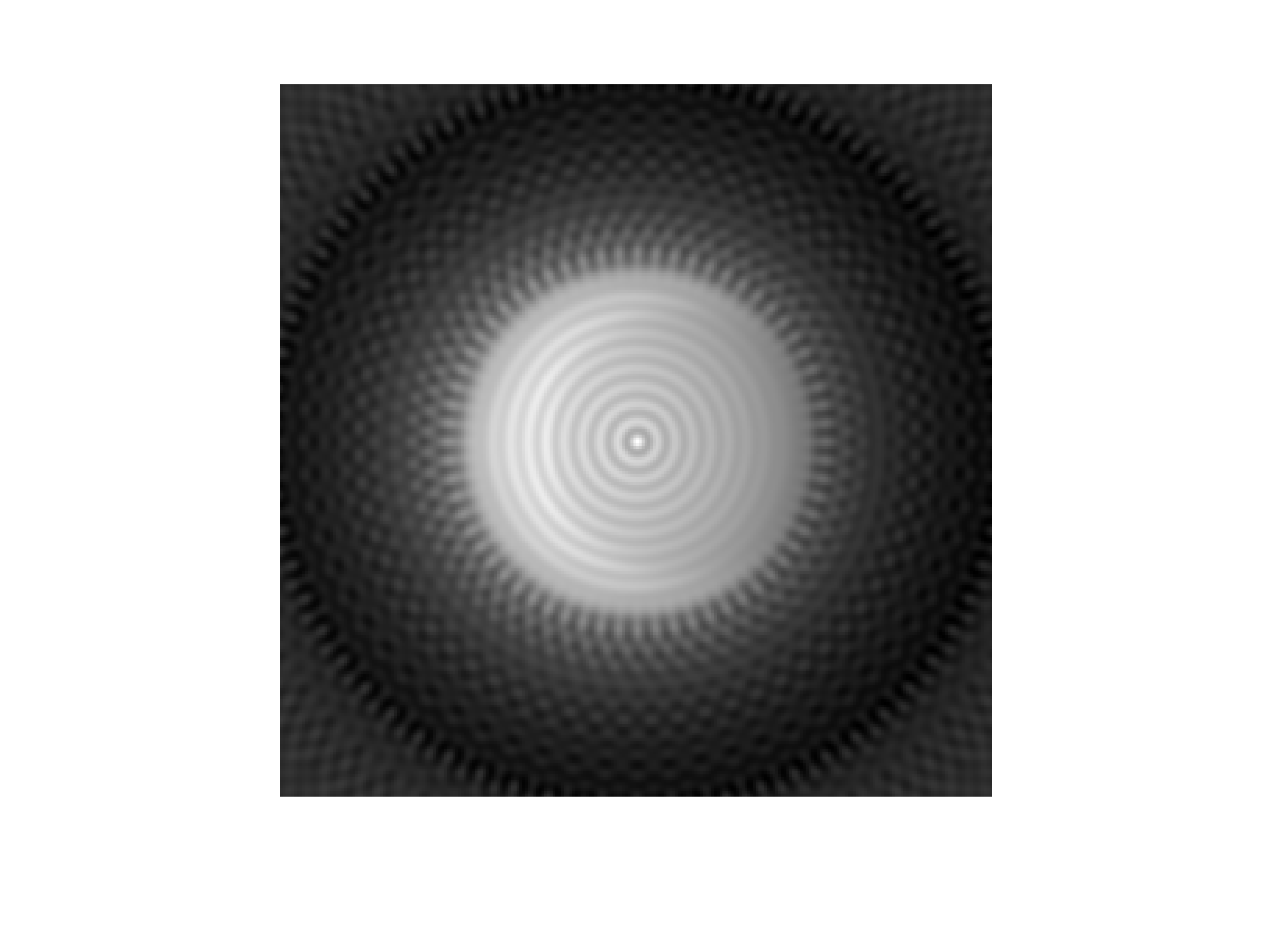}}
\caption{Reconstruction of the crescent-shaped phantom with Gaussian kernel and truncation regularization : (a) Original phantom; (b) $\varepsilon=1$, $L=10$, $k_{1}^{-1}(A)=4.04307\cdot10^{-28}$ ; (c) $\varepsilon=25$, $L=10$, $k_{1}^{-1}(A)=2.3734\cdot10^{-5}$; (d) $\varepsilon=50$, $L=10$, $k_{1}^{-1}(A)=4.4333\cdot10^{-5}$.}
\label{fig: gauss_trunc}
\end{figure}

In conclusion we have to solve the linear system $A_{L}c=f$ with components of $A_{L}$ given by
\begin{equation*}
\boxed{
a_{k,j}=\left\{
\begin{aligned}
&\frac{\sqrt{\pi}}{\varepsilon^{2}a}\int_{c1}^{c2}{e^{-u^{2}}\,du} & &\text{if}\ \theta_{k}\neq\theta_{j}\\
&\frac{2\sqrt{\pi}}{\varepsilon}e^{-\varepsilon^{2}b^{2}}\sqrt{L^{2}-t_{k}^{2}} & &\text{if}\ \theta_{k}=\theta_{j}
\end{aligned}
\right.
}
\end{equation*}
where
\begin{align*}
&a=\sin{(\theta_{k}-\theta_{j})} & &b=t_{j}-t_{k}\cos{(\theta_{k}-\theta_{j})}\\
&c_{1}=\varepsilon(-\sqrt{L^{2}-t_{k}^{2}}+b) & &c_{2}=\varepsilon(\sqrt{L^{2}-t_{k}^{2}}+b).
\end{align*}
Solving this system we obtain $c$ and we can then evaluate the solution $s$ by
\begin{equation*}
s(x)=\sum_{j=1}^{n}{c_{j}b_{j}(x)}
\end{equation*}
with $b_{j}$ given by \eqref{eq: basis_gauss}. Figure \ref{fig: gauss_trunc} shows the results of applying this method to the crescent-shaped phantom for suitable values of $\varepsilon$ and $L$ in the case of parallel beam geometry, where samples of the Radon transform of $f$ are taken at angles $\theta_{p}=p\pi/N$, $i=1,\ldots,N-1$ and $t_{q}=q/M$, $q=-M,\ldots,M$. We also observe that without using the shape parameter $\varepsilon$, i.e. using $\varepsilon=1$, the matrix $A_{L}$ would have been highly ill-conditioned and the result very different.
%\begin{figure}[htbp]
%\centering
%\includegraphics[width=0.7\textwidth]{gauss_trunc.pdf} 
%\caption{Reconstruction  with Gaussian kernel and truncation regularization: $\varepsilon=25$, $L=10$.}
%\label{fig: gauss_trunc}
%\end{figure} 

\subsection{Regularization by Gaussian filtering}\label{subsec: reg_gauss}
\begin{figure}[htbp]
\centering%
\subfigure[Option 1 \label{fig: gauss_gauss_1pm}]%
{\includegraphics[width=0.49\textwidth]{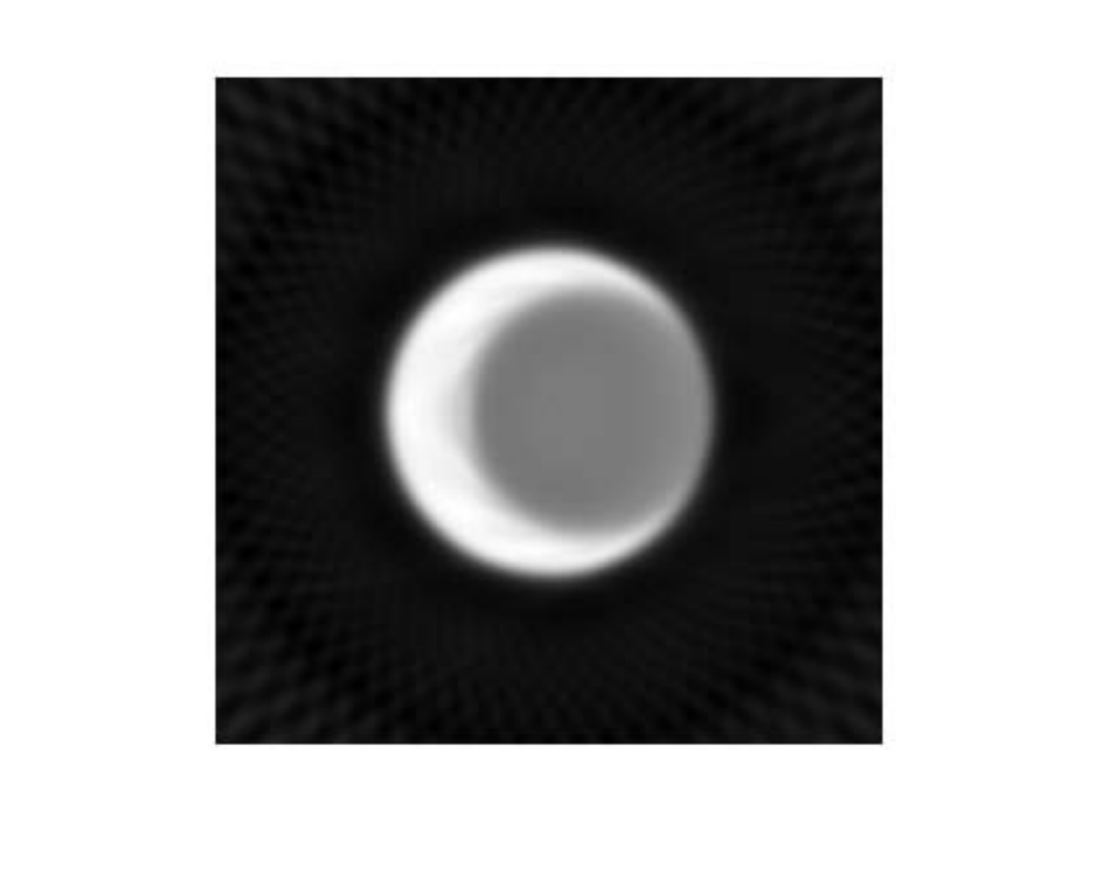}}
\subfigure[Option 2 \label{fig:  gauss_gauss_2pm}]%
{\includegraphics[width=0.49\textwidth]{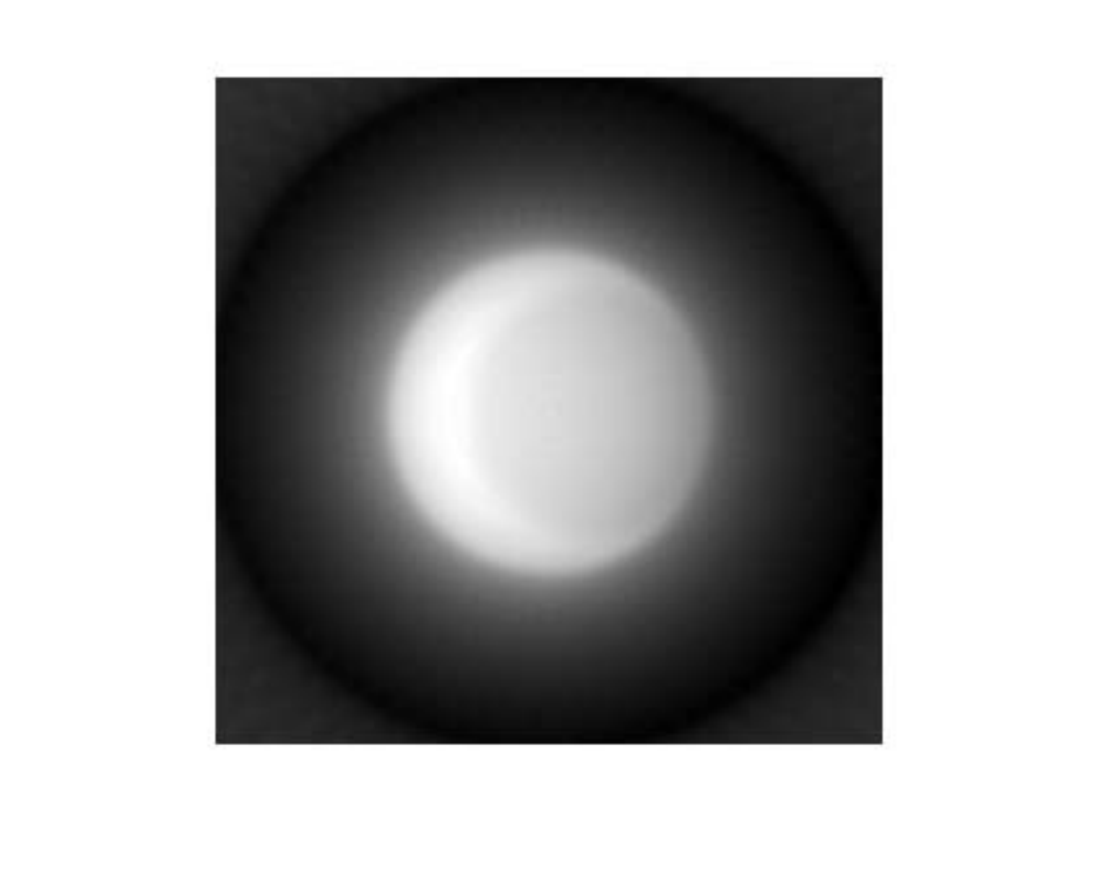}}
\caption{Gaussian reconstruction, parallel beam geometry data}
\label{fig: gauss_gauss_pm}
\end{figure}

\begin{figure}[htbp]
\centering%
\subfigure[Option 1 \label{fig: gauss_gauss_1s}]%
{\includegraphics[width=0.49\textwidth]{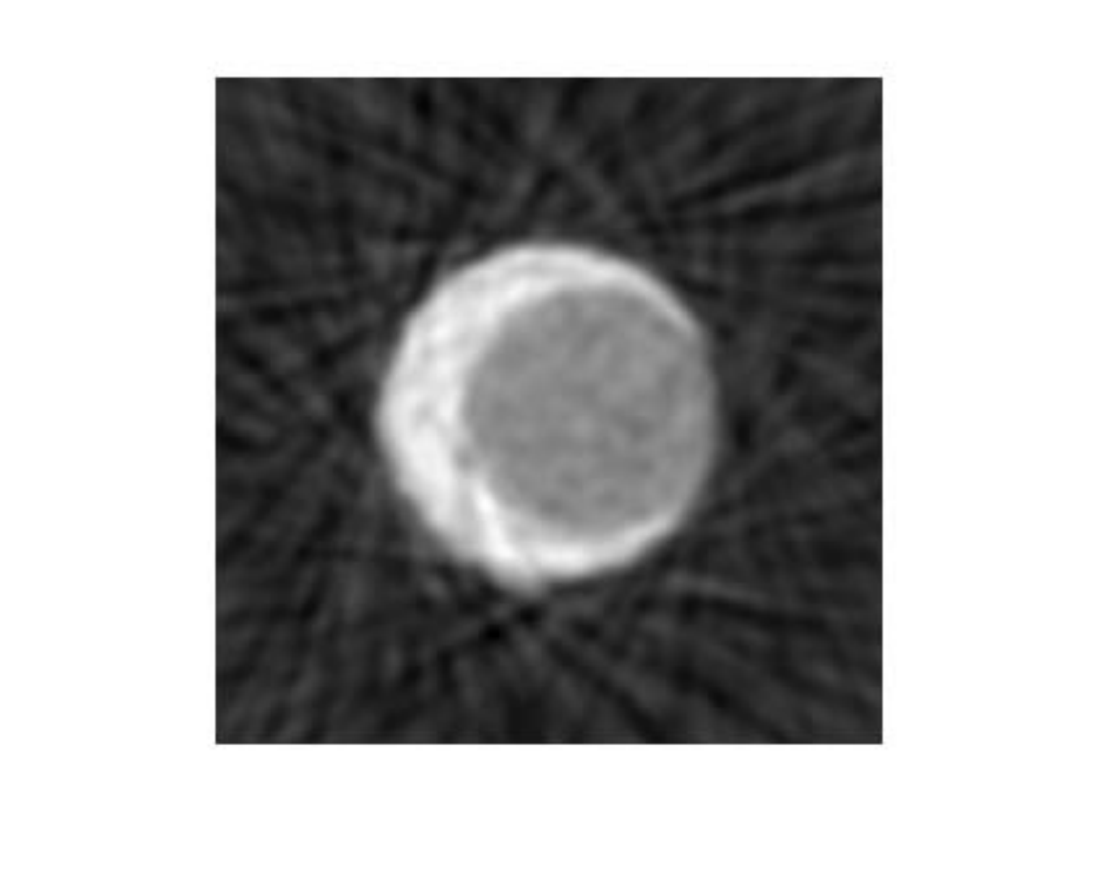}}
\subfigure[Option 2 \label{fig:  gauss_gauss_2pm}]%
{\includegraphics[width=0.49\textwidth]{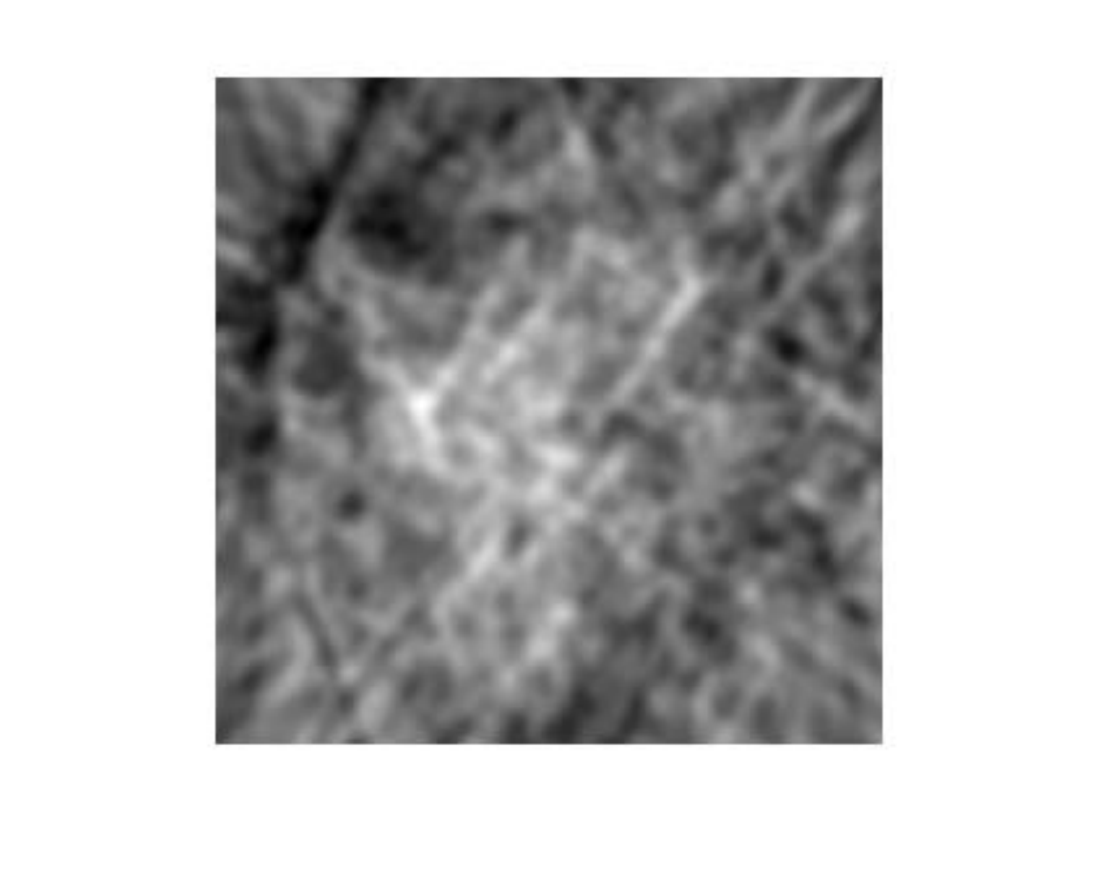}}
\caption{Gaussian reconstruction, scattered data}
\label{fig: gauss_gauss_s}
\end{figure}

Consider again the Gaussian kernel $K(x,y)=e^{-\varepsilon^{2}\norm{x-y}^{2}}$ and the associated basis of the space $S$, defined in \eqref{eq: basis_gauss}. In this section we will use another window function in order to regularize the integral $R[b_{j}]$, i.e. we will multiply the function $b_{j}$ by another Gaussian function
\begin{equation*}
w_{\nu}(x)=e^{-\nu^{2}\norm{x}^{2}}.
\end{equation*}
In other words we will consider the operator $R_{\nu}$ given by $R_{\nu}[f]=R[fw_{\nu}]$ instead of the classical Radon transform. We have
\begin{align*}
R_{\nu}[b_{j}](r,\varphi)&=\int_{\R}{b_{j}(x(s))e^{-\nu^{2}\norm{x(s)}^{2}}\,ds}=\frac{\sqrt{\pi}}{\varepsilon}\int_{\R}{e^{-\varepsilon^{2}(as+b)^{2}}e^{-\nu^{2}(r^{2}+s^{2})}\,ds}%=\\
%&=\frac{\sqrt{\pi}}{\varepsilon}\int_{\R}{e^{-\varepsilon^{2}(as+b)^{2}}
\end{align*}
where $a,b$ are defined by \eqref{eq: def_ab}. Then,
\begin{align*}
R_{\nu}[b_{j}](r,\varphi)&=\frac{\sqrt{\pi}}{\varepsilon}e^{-\varepsilon^{2}b^{2}-\nu^{2}r^{2}}\int_{\R}{\exp{(-[(\varepsilon^{2}a^{2}+\nu^{2})s^{2}+2abs\varepsilon^{2}])}\,ds}=\\
&=\frac{\sqrt{\pi}}{\varepsilon}\exp{\left(-\varepsilon^{2}b^{2}-\nu^{2}r^{2}+\frac{a^{2}b^{2}\varepsilon^{4}}{a^{2}\varepsilon^{2}+\nu^{2}}\right)}\cdot\\
&\qquad\int_{\R}{\exp{\left(-\left[\sqrt{\varepsilon^{2}a^{2}+\nu^{2}}s+\frac{ab\varepsilon^{2}}{\sqrt{\varepsilon^{2}a^{2}+\nu^{2}}}\right]^{2}\right)}\,ds}=\\
&=\frac{\pi}{\varepsilon\sqrt{a^{2}\varepsilon^{2}+\nu^{2}}}\exp{\left[-\nu^{2}\left(r^{2}+\frac{\varepsilon^{2}b^{2}}{a^{2}\varepsilon^{2}+\nu^{2}}\right)\right]}.
\end{align*}
We have now two options:
\begin{enumerate}
\item The regularization $R_{\nu}[b_{j}](t_{k},\theta_{k})$ for all values of $k,j$, that leads to a linear system with matrix $A_{1}^{\nu}$ whose components are
\begin{equation*}
\boxed{
a_{k,j}=\frac{\pi\exp{\left[-\nu^{2}\left(r^{2}+\frac{\varepsilon^{2}b^{2}}{a^{2}\varepsilon^{2}+\nu^{2}}\right)\right]}}{\varepsilon\sqrt{a^{2}\varepsilon^{2}+\nu^{2}}}.
}
\end{equation*}
\item The regularization $R_{\nu}[b_{j}](t_{k},\theta_{k})$ only for those values of $k,j$ for which $b_{j}$ has not finite Radon transform, i.e. only when $a=0$, while for $a\neq0$ we consider the usual Radon transform. This corresponds to the matrix $A_{2}^{\nu}$ whose elements are
\begin{equation*}
a_{k,j}=\left\{
\begin{aligned}
&\frac{\pi}{\varepsilon^{2}a}	& &\text{if}\ a\neq0\\
&\frac{\pi\exp{[-(\nu^{2}r^{2}+\varepsilon^{2}b^{2})]}}{\varepsilon\nu} & &\text{if}\ a=0.
\end{aligned}
\right.
\end{equation*}
\end{enumerate}
Numerical experiments (Figures \ref{fig: gauss_gauss_pm} and \ref{fig: gauss_gauss_s}) show that the first option gives better results, provided that the value of $\varepsilon$ is relatively big ($\approx30$) and value of $\nu$ is quite small $(\approx0.5)$ so that the condition number of the matrix $A^{\nu}_{1}$ is small.

\section{Inverse multiquadrics reconstruction}\label{sec: imqRec}
We now consider the same reconstruction problem by using the inverse multiquadrics as kernel function:
\begin{equation*}
K(x,y)=\frac{1}{\sqrt{1+\varepsilon^{2}\norm{x-y}^2}}.
\end{equation*}
As we saw in section \ref{subsec: regularization} $K$ does not admit finite Radon transform and so we have to multiply $K$ by another window function $w$ of the form $w=w(\norm{x-y})$, so that $R^{y}[K(x,y)w](t,\theta)<\infty$ for all $t,\theta$. 

The window function we consider is the characteristic function of a compact set:
\begin{equation*}
w=\chi_{[-L,L]}(\norm{x-y}), \qquad L\gg0.
\end{equation*}
As in the Gaussian case we use the shift property of the Radon transform and we first compute
\begin{align*}
R_{L}[k(y)](t,\theta)&=R[K(0,y)\chi_{[-L,L]}(\norm{y})](t,\theta)=\\
&=\int_{\R}{\frac{1}{\sqrt{1+\varepsilon^{2}(t^{2}+s^{2})}}\chi_{[-L,L]}(\sqrt{t^{2}+s^{2}})\,ds}=\\
&=\int_{-\sqrt{L^{2}-t^{2}}}^{\sqrt{L^{2}-t^{2}}}{\frac{1}{\sqrt{c^{2}+(\varepsilon s)^{2}}}\,ds}
\end{align*}
where $c=\sqrt{1+\varepsilon^{2}t^{2}}$ and we assume $|t|<L$. We then apply the substitution $u=\varepsilon s$ in the integral
\begin{align}\label{eq: asinh_sqrt}
\frac{1}{\varepsilon}\int_{-\varepsilon\sqrt{L^{2}-t^{2}}}^{\varepsilon\sqrt{L^{2}-t^{2}}}{\frac{1}{\sqrt{c^{2}+s^{2}}}\,du}&=\left[\frac{2}{\varepsilon}\text{asinh}\left(\frac{u}{c}\right)\right]_{0}^{\varepsilon\sqrt{L^{2}-t^{2}}}=\frac{2}{\varepsilon}\text{asinh}\left(\varepsilon\sqrt{\frac{L^{2}-t^{2}}{1+\varepsilon^{2}t^{2}}}\right),
\end{align}
where %asinh is the inverse function of the hyperbolic sinus and is given by
\begin{equation*}
\text{asinh}(x)=\log{(x+\sqrt{1+x^{2}})},
\end{equation*}
so we can also write integral \eqref{eq: asinh_sqrt} as 
\begin{align*}
\frac{2}{\varepsilon}\left[\text{asinh}\left(\frac{u}{c}\right)\right]_{0}^{\varepsilon\sqrt{L^{2}-t^{2}}}&=\frac{2}{\varepsilon}\left[\log(u+\sqrt{c^{2}+u^{2}})\right]_{0}^{\varepsilon\sqrt{L^{2}-t^{2}}}=\\
&=\frac{2}{\varepsilon}\left(\log{(\varepsilon\sqrt{L^{2}-t^{2}}+\sqrt{1+\varepsilon^{2}L^{2}})}-\frac{1}{2}\log{(1+\varepsilon^{2}t^{2})} \right).
\end{align*}
We conclude that the basis associated to the kernel $Kw$ is 
\begin{equation*}
\boxed{
b_{j}(x)=\frac{2}{\varepsilon}\text{asinh}\left(\varepsilon\sqrt{\frac{L^{2}-(t_{j}-x\cdot v_{j})^{2}}{1+\varepsilon^{2}(t_{j}-x\cdot v_{j})^{2}}}\right)\chi_{[-L,L]}(t_{j}-x\cdot v_{j}),
}
\end{equation*}
where, as usual, $v_{j}=(\cos{\theta_{j}},\sin{\theta_{j}})$.

\begin{figure}[htbp]
\centering%
\subfigure[Parallel beam geometry \label{fig: imq_trunc_pb}]%
{\includegraphics[width=0.49\textwidth]{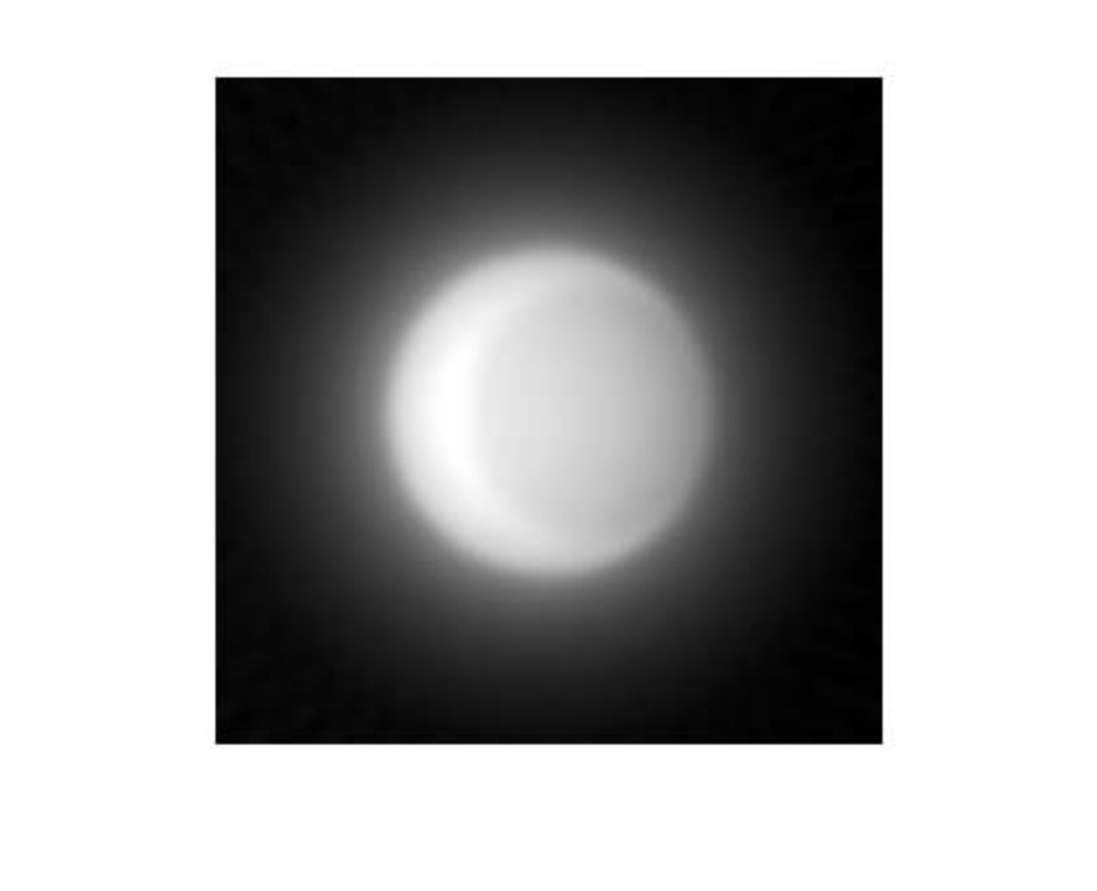}}
\subfigure[Scattered data \label{fig: imq_trunc_s}]%
{\includegraphics[width=0.49\textwidth]{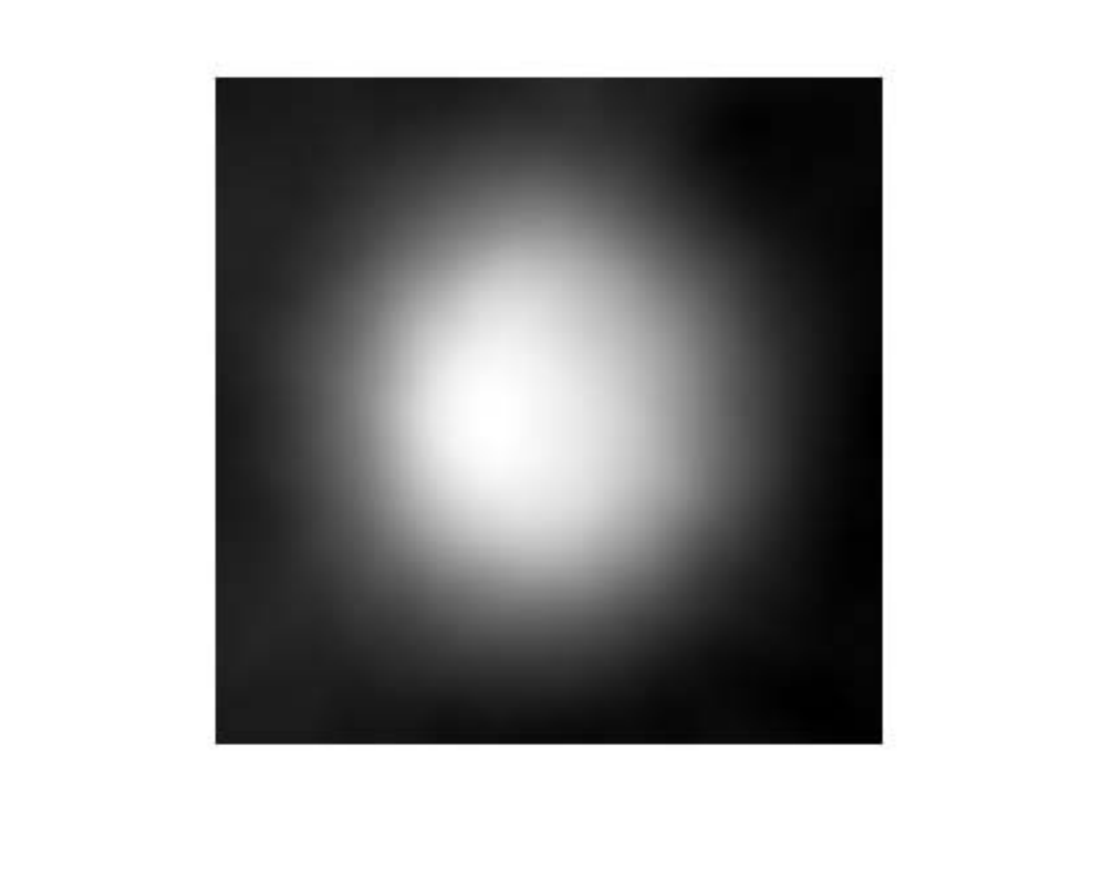}}
\caption{Inverse multiquadric reconstruction of the crescent-shaped phantom using $\varepsilon=30$, $H=L=20$.}
\label{fig: imq_trunc}
\end{figure}

We can now compute the matrix $A=(a_{k,j})_{k,j=1}^{n}$:
\begin{equation*}
a_{k,j}=R[b_{j}(x)](r,\varphi)=\int_{\R}{\frac{2}{\varepsilon}\text{asinh}\left(\varepsilon\sqrt{\frac{L^{2}-(as+b)^{2}}{1+\varepsilon^{2}(as+b)^{2}}}\right)\chi_{[-L,L]}(as+b)\,ds},
\end{equation*}
where $a,b$ are defined as usual and $r=t_{k}$, $\varphi=\theta_{k}$. We observe that if $a=0$, then 
\begin{equation*}
a_{k,j}=\frac{2}{\varepsilon}\text{asinh}\left(\varepsilon\sqrt{\frac{L^{2}-b^{2}}{1+\varepsilon^{2}b^{2}}}\right)\int_{\R}{\chi_{[-L,L]}(b)\,ds}=+\infty
\end{equation*}
when $|b|<L$, that is our case. So we have to consider a further regularization of $R$. We choose to truncate, i.e. we compute
\begin{equation*}
a_{k,j}=R_{H}[b_{j}(x)](r,\varphi)=R[b_{j}(x)\chi_{[-H,H]}(\norm{x})],\qquad H>0.
\end{equation*}
As for the Gaussian case, we have two options: consider the regularization $R_{H}$ for all values of $k$ and $j$ or use it only when $\theta_{k}=\theta_{j}$ i.e. when $a=0$. As before we consider the first option since the resulting matrix has a better condition number. Thus, we obtain
\begin{equation*}
a_{k,j}=\frac{2}{\varepsilon}\int_{\R}{\text{asinh}\left(\varepsilon\sqrt{\frac{L^{2}-(as+b)^{2}}{1+\varepsilon^{2}(as+b)^{2}}}\right)\chi_{[-L,L]}(as+b)\chi_{[-H,H]}(\sqrt{s^{2}+r^{2}})\,ds}
\end{equation*}
that in the case $a=0$ becomes
\begin{equation*}
a_{k,j}=\frac{4}{\varepsilon}\text{asinh}\left(\varepsilon\sqrt{\frac{L^{2}-b^{2}}{1+\varepsilon^{2}b^{2}}}\right)\sqrt{H^{2}-r^{2}},
\end{equation*}
provided that $|b|<L$ and $|r|<H$. In the case $a\neq0$ we consider the substitution $u=\varepsilon(as+b)$ that leads us to the integral
\begin{align}
\boxed{
a_{k,j}=\frac{2}{\varepsilon^{2}a}\int_{c_{1}}^{c_{2}}{\text{asinh}\left(\sqrt{\frac{\varepsilon^{2}L^{2}-u^{2}}{1+u^{2}}}\right)\,du}
}
\label{eq: int_asinh}
\end{align}
\begin{align*}
&c_{1}=\varepsilon\max{(-L,-|a|\sqrt{H^{2}-r^{2}}+b)} & &c_{2}=\varepsilon\min{(L,|a|\sqrt{H^{2}-r^{2}}+b)}.
\end{align*}
All what we have to do now is to compute integrals \eqref{eq: int_asinh}. The computation of these integrals can be found in appendix A. 
Applying this method again to the crescent-shaped phantom, choosing $\varepsilon=30$, $H=L=20\max|t_{j}|$ we obtain the reconstruction shown in Figure \ref{fig: imq_trunc}. We observe that we have acceptable reconstruction only with parallel beam geometry data, to obtain good results also with scattered data we should consider another window function instead of the characteristic function.

\section{Multiquadrics reconstruction}
We now consider the multiquadric kernel
\begin{equation*}
K(x,y)=\sqrt{1+\rho^{2}\norm{x-y}^2}, \qquad \rho>0.
\end{equation*}
As in the case of inverse multiquadrics, $K$ is not integrable on any line in the $(y_{1},y_{2})$ plane. The approach we follow to regularize the problem is to consider a Gaussian weighting function for computing  both basis functions $b_{j}$ and the matrix $A$.

We start with the Gaussian-filtered basis
\begin{equation*}
b_{j}(x)=R^{y}[K(x,y)e^{-\varepsilon^{2}\norm{x-y}^{2}}](t_{j},\theta_{j}), \quad \varepsilon>0.
\end{equation*}
We recall that this operation corresponds to use kernel $\tilde{K}(x,y)=K(x,y)e^{-\varepsilon^{2}\norm{x-y}^{2}}$ in place of $K$.
Proceeding as usual, we first compute $b_{j}(x)$ then the coefficients $a_{j,k}$.

Now, 
\begin{align*}
b_{j}(0)&=\int_{\R}{\sqrt{1+\rho^{2}(t^{2}+s^{2})}e^{-\varepsilon^{2}(t^{2}+s^{2})}\,ds}=\\
&=2\rho e^{\frac{\varepsilon^{2}}{\rho^{2}}}\int_{0}^{+\infty}{\sqrt{\frac{1}{\rho^{2}}+t^{2}+s^{2}}\ e^{-\varepsilon^{2}(\frac{1}{\rho^{2}}+t^{2}+s^{2})}\,ds}.
\end{align*}
Setting $c=\sqrt{\frac{1}{\rho^{2}}+t^{2}}$ and then integrating by parts, we obtain
\begin{align*}
b_{j}(0)&=2\rho e^{\frac{\varepsilon^{2}}{\rho^{2}}}\left\{c^{2}\left[\frac{1}{4}\sinh{(2s)}+\frac{s}{2}\right]_{s=0}^{+\infty}+\right.\\
&\left. -\int_{0}^{+\infty}{c^{2}\left[\left(\frac{1}{4}\sinh{(2s)}+\frac{s}{2}\right)e^{-\varepsilon^{2}(c^{2}+s^{2})}\right](-2\varepsilon^{2}s)e^{-\varepsilon^{2}(c^{2}+s^{2})}\,ds}\right\}=\\
&=c^{2}\rho e^{\frac{\varepsilon^{2}}{\rho^{2}}}\left[\left(\frac{1}{2}\sinh{(2s)}+s\right) e^{-\varepsilon^{2}(c^{2}+s^{2})}\right]_{s=0}^{+\infty}+\\
&+2c^{2}\rho e^{\frac{\varepsilon^{2}}{\rho^{2}}}\varepsilon^{2}\int_{0}^{+\infty}{\left(\frac{1}{2}\sinh{(2s)}+s\right) se^{-\varepsilon^{2}(c^{2}+s^{2})}\,ds}=\\
&=c^{2}\rho e^{\frac{\varepsilon^{2}}{\rho^{2}}}I_{1}+2c^{2}\rho e^{\frac{\varepsilon^{2}}{\rho^{2}}}\varepsilon^{2}I_{2},
\end{align*}
where 
\begin{align*}
I_{1}=\lim_{r\rightarrow+\infty}{\left[\left(\frac{1}{2}\sinh{(2r)}+r\right) e^{-\varepsilon^{2}(c^{2}+r^{2})}\right]}-e^{-\varepsilon^{2}c^{2}}(\frac{1}{2}\sinh{(0)}+0)=0
\end{align*}
and
\begin{align*}
I_{2}&=\frac{1}{2}\int_{0}^{+\infty}{s\sinh{(2s)}e^{-\varepsilon^{2}(c^{2}+s^{2})}\,ds}+\int_{0}^{+\infty}{s^{2}e^{-\varepsilon^{2}(c^{2}+s^{2})}\,ds}=\\
&=\frac{1}{2}e^{-\varepsilon^{2}c^{2}}\frac{\sqrt{\pi}}{2}\frac{e^{\varepsilon^{-2}}}{\varepsilon^{3}}+e^{-\varepsilon^{2}c^{2}}\frac{\sqrt{\pi}}{4}\frac{1}{\varepsilon^{3}}=\frac{\sqrt{\pi}}{4}\frac{e^{-\varepsilon^{2}c^{2}}}{\varepsilon^{3}}e^{\varepsilon^{-2}}.
\end{align*}
We conclude that
\begin{equation*}
b_{j}(0)=\frac{\sqrt{\pi}\rho}{2\varepsilon}\left(\frac{1}{\rho^{2}}+t^{2}\right)e^{-\varepsilon^{2}t^{2}+\varepsilon^{-2}}
\end{equation*}
and so 
\begin{equation*}
\boxed{
b_{j}(x)=\frac{\sqrt{\pi}\rho}{2\varepsilon}\left(\frac{1}{\rho^{2}}+(t_{j}-x\cdot v_{j})^{2}\right)e^{-\varepsilon^{2}(t_{j}-x\cdot v_{j})^{2}+\varepsilon^{-2}}, \quad v_{j}=(\cos{\theta_{j}},\sin{\theta_{j}}).
}
\end{equation*}

\begin{figure}[htbp]
\centering%
\subfigure[Parallel beam geometry \label{fig: mq_gauss_pb}]%
{\includegraphics[width=0.49\textwidth]{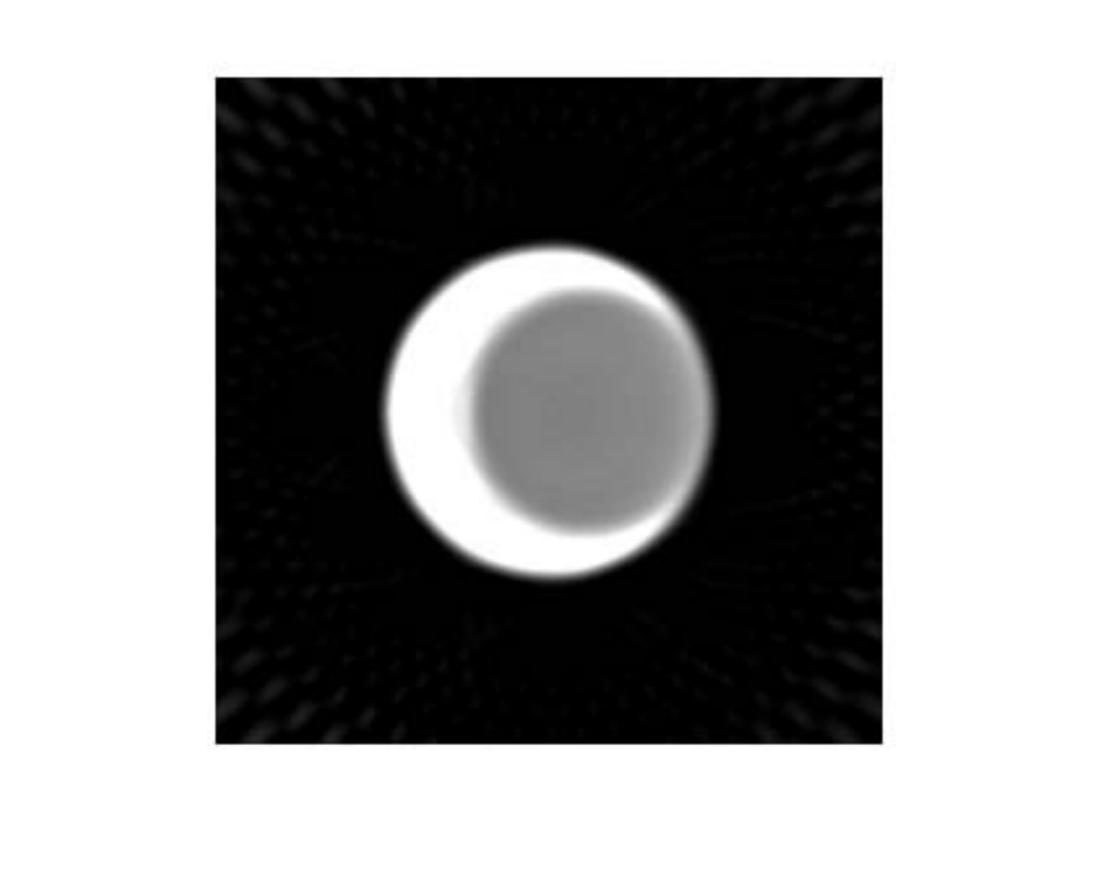}}
\subfigure[Scattered data \label{fig: mq_gauss_s}]%
{\includegraphics[width=0.49\textwidth]{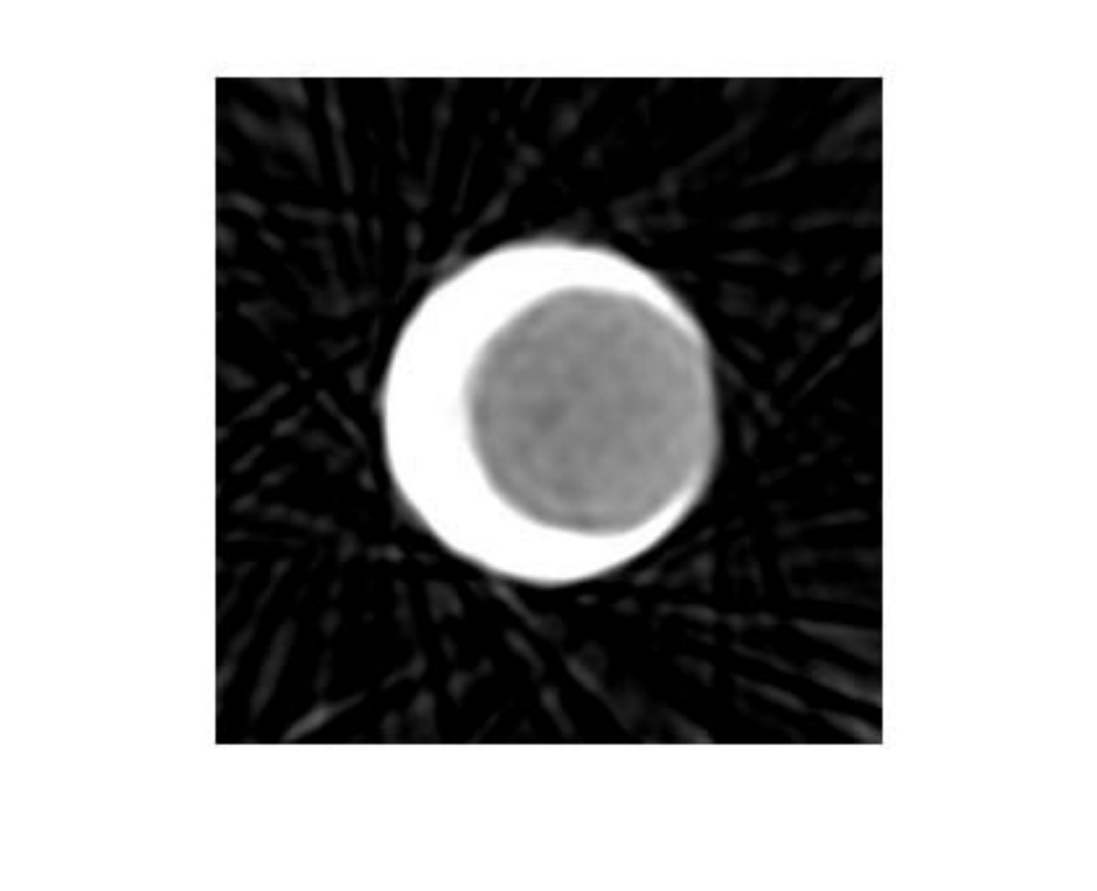}}
\caption{Multiquadric reconstruction of the crescent-shaped phantom using $\phi=1,\ \varepsilon=30,\ \nu=0.8.$}
\label{fig: mq_gauss}
\end{figure}

In order to compute the matrix $A$, we consider $R[b_{j}(x)](r,\varphi)$ . Setting $x_{s}=(r\cos{\theta}-s\sin{\theta},r\sin{\theta}+s\cos{\theta})$, we have
\begin{align*}
a_{k,j}&=R[b_{j}(x)](r,\varphi)=\int_{\R}{\frac{\sqrt{\pi}\rho}{2\varepsilon}\left(\frac{1}{\rho^{2}}+(t-x_{s}\cdot v)^{2}\right)e^{-\varepsilon^{2}(t-x_{s}\cdot v)^{2}+\varepsilon^{-2}}\,ds}=\\
&=\frac{\sqrt{\pi}\rho}{2\varepsilon}e^{\varepsilon^{-2}}\left[ \int_{\R}{\rho^{-2}e^{-\varepsilon^{2}(as+b)^{2}}\,ds}+\int_{\R}{(as+b)^{2}e^{-\varepsilon^{2}(as+b)^{2}}\,ds}
\right],
\end{align*}
where $a=\sin{(\varphi-\theta)}$ and $b=t-r\cos{(\varphi-\theta)}$. If $a=0$ this integral is infinite. To avoid this case we consider the regularization $R_{\nu}$ of $R$, that is
\begin{align*}
R_{\nu}[b_{j}(x)](r,\varphi)&=\frac{\sqrt{\pi}\rho}{2\varepsilon}e^{\varepsilon^{-2}}\int_{\R}{\left(\frac{1}{\rho^{2}}+(t-x_{s}\cdot v)^{2}\right)e^{-\varepsilon^{2}(t-x_{s}\cdot v)^{2}}e^{-\nu^{2}\norm{x_{s}}^{2}}\,ds}=\\
&=\frac{\sqrt{\pi}\rho}{2\varepsilon}e^{\varepsilon^{-2}}\int_{\R}{\left(\frac{1}{\rho^{2}}+(as+b)^{2}\right)e^{-\varepsilon^{2}(as+b)^{2}-\nu^{2}(r^{2}+s^{2})}\,ds}=\\
&=C_{\varepsilon,\nu}\left[ \int_{\R}{\rho^{-2}e^{-(cs+d)^{2}}\,ds}+\int_{\R}{(as+b)^{2}e^{-(cs+d)^{2}}\,ds} \right],
\end{align*}
where 
\begin{align*}
c&=\sqrt{\nu^{2}+\varepsilon^{2}a^{2}}, & d&=\frac{ab\varepsilon^{2}}{c}, & C_{\varepsilon,\nu}&=\frac{\sqrt{\pi}\rho}{2\varepsilon}\exp{(\varepsilon^{-2}+d^{2}-\varepsilon^{2}b^{2}-\nu^{2}r^{2})}.
\end{align*}
Since
\begin{equation*}
\int_{\R}{(as+b)^{2}e^{-(cs+d)^{2}}\,ds}=\frac{\sqrt{\pi}}{2|c|^{3}}(a^{2}(2d^{2}+1)-4abcd+2b^{2}c^{2}),
\end{equation*}
we conclude that
\begin{equation*}
\boxed{
a_{k,j}=\frac{\pi\exp{\left(\frac{1}{\varepsilon^{2}}-\frac{\nu^{2}\varepsilon^{2}b^{2}}{\nu^{2}+\varepsilon^{2}a^{2}}-\nu^{2}r^{2}\right)}}{2\varepsilon\sqrt{\nu^{2}+\varepsilon^{2}a^{2}}}\left[\frac{1}{\rho} +\frac{\rho}{2}\frac{a^{2}(\nu^{2}+\varepsilon^{2}a^{2})+2b^{2}\nu^{4}}{(\nu^{2}+\varepsilon^{2}a^{2})^{2}}\right]
}
\end{equation*}
Figure \ref{fig: mq_gauss} shows the result of using the multiquadric kernel with Gaussian filtering for the reconstruction of the crescent-shaped phantom.

\section{Compactly supported radial basis functions}
Another important class of positive definite functions we consider are the compactly supported radial basis functions. If a function as compact support, then it is automatically strictly positive definite and it is strictly positive definite on $\R^{d}$ only for a fixed maximum value of the dimension $d$. Moreover one can show that there not exist compactly supported radial functions that are strictly conditionally positive definite of order $m>0$ (see \cite{WU} for more informations about compactly supported radial functions).

\subsubsection{Wendland's compactly supported functions}
A popular family of compactly supported functions was introduced by Wendland \cite{WEND}. Wendland starts with the truncated power function $\varphi_{l}(r)=(1-r)^{l}_{+}$, which is striclty positive definite and radial on $\R^{d}$ for $d\leq2l-1$, and then applies repeatedly the integral operator $\mathcal{I}$ defined as follow:

\begin{definition}
Let $\varphi:[0,\infty)\rightarrow\R$ such that $t\varphi(t)\in L^{1}[0,\infty)$, then we define 
\begin{equation*}
\mathcal{I}\varphi(r)=\int_{r}^{\infty}{t\varphi(t)\,dt}, \quad r\geq0.
\end{equation*}
\end{definition}
We can now define the Wendland's compactly supported functions:
\begin{definition}
With $\varphi_{l}(r)=(1-r)^{l}_{+}$, we define 
\begin{equation*}
\varphi_{d,k}=\mathcal{I}^{k}\varphi_{\lfloor d/2\rfloor+k+1}
\end{equation*}
\end{definition}
\begin{example}\label{ex: wend}
The explicit representation of $\varphi_{d,k}$ for $k=0,1,2,3$ are:
\begin{align*}
\varphi_{d,0}(r)&=(1-r)^{l}_{+},\\
\varphi_{d,1}(r)&=(1-r)^{l}_{+}[(l+1)r+1],\\
\varphi_{d,2}(r)&=(1-r)^{l}_{+}[(l^2+4l+3)r^2+(3l+6)r+3],\\
\varphi_{d,3}(r)&=(1-r)^{l}_{+}[(l^3+9l^2+23l+15)r^3+(6l^2+36l+45)r^2+(15l+45)r+15],
\end{align*}
where $l=\lfloor d/2\rfloor +k+1$ and equalities are up to a multiplicative constant.
\end{example}
We observe that all the functions in Example \ref{ex: wend} are compactly supported and have a polynomial representation on their support. This is true in general as stated in the following
\begin{theorem}
The functions $\varphi_{d,k}$ are strictly positive definite and radial on $\R^{d}$ and are of the form 
\begin{equation*}
\varphi_{d,k}(r)=\left\{
\begin{aligned}
&p_{d,k}(r) & &\text{if}\ r\in[0,1]\\
&0 & &\text{if}\ r>1,
\end{aligned}
\right.
\end{equation*}
where $p_{d,k}$ is a polynomial of degree $\lfloor d/2\rfloor+3k+1$. Moreover $\varphi_{d,k}\in C^{2k}(\R^{d})$ are unique up to a constant factor and the polynomial degree is minimal for given space dimension $d$ and smoothness $2k$.
\end{theorem}
The proof of this theorem can be found in \cite{WEND}.

\begin{figure}[htbp]
\centering
\includegraphics[width=0.7\textwidth]{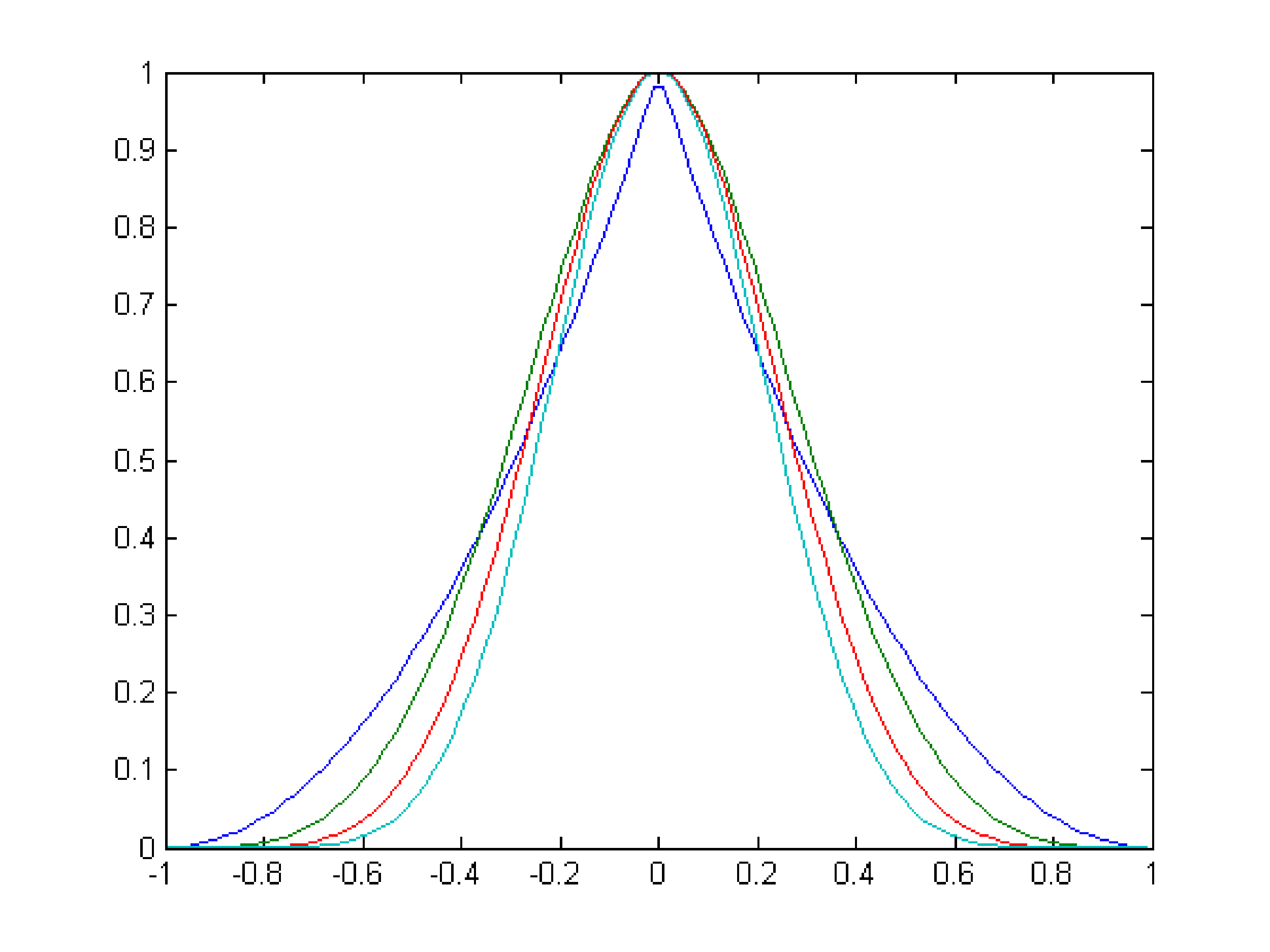} 
\caption{Plot of Wendland's compactly supported functions.}
\label{fig: compactlySupported}
\end{figure} 

\subsection{Compactly supported kernel reconstruction}\label{subsec: compSupp}
Compactly supported radial basis functions can be used as well in image reconstruction: if $\varphi(r)$ is a compactly supported function, one sets $K(x,y)=\varphi(\varepsilon\norm{x-y}), \ \varepsilon>0$, and uses $K$ as kernel function. The property of compact support can be useful in the computation of the Radon transform, in particular if $\varphi$ is compactly supported and of class at least $C^{0}$ on its support, then the Radon transform $Rf$ is well defined, indeed 
\begin{align*}
|Rf|=&\left|\int_{\R}{f(x_{s})\,ds}\right|=\left|\int_{Supp(f)}{f(x_{s})\,ds}\right|\leq\int_{Supp(f)}{|f(x_{s})|\,ds}\leq\\
&\leq\sup_{x\in Supp(f)}{|f(x)|}\cdot m(Supp(f))<\infty.
\end{align*}
For example Wendland's and Wu's compactly supported functions, having a polynomial representation on their domain, admit finite Radon transform. This means that for this class of functions the basis $b_{j}(x)=R^{y}[K(x,y)](t_{j},\theta_{j})$ is well defined.

However the fact that the kernel $K$ is compactly supported does not imply that $b_{j}(x)$ has finite Radon transform and then matrix $A=(a_{j,k})$ can have non finite entries as shown in the following example.

\begin{example}\label{ex: radon_wend}
Consider the Wendland's function $\varphi_{2,1}(r)=(1-r)_{+}^{4}(4r+1)$ and set 
\begin{equation*}
K(x,y)=\varphi_{2,1}(\varepsilon\norm{x-y})=(1-\varepsilon\norm{x-y})_{+}^{4}(4\varepsilon\norm{x-y}+1)
\end{equation*}
The support $\{(x,y):\ \varepsilon\norm{x-y}\leq1\}$ of $K$ is compact.

We compute $b_{j}(x)$ using the shift property of the Radon transform:
\begin{align*}
R^{y}[K(0,y)](t,\theta)&=\int_{\R}{(1-\varepsilon\sqrt{t^{2}+s^{2}})^{4}_{+}(4\varepsilon\sqrt{t^{2}+s^{2}}+1)\,ds}=\\
&=\int_{\sqrt{t^{2}+s^{2}}\leq\frac{1}{\varepsilon}}{4\varepsilon(1-\varepsilon\sqrt{t^{2}+s^{2}})^{4}\sqrt{t^{2}+s^{2}}\,ds}+\\
&+\int_{\sqrt{t^{2}+s^{2}}\leq\frac{1}{\varepsilon}}{(1-\varepsilon\sqrt{t^{2}+s^{2}})^{4}\,ds},
\end{align*}
thus
\begin{equation*}
R^{y}[K(0,y)](t,\theta)=\left\{
\begin{aligned}
&I_{1}+I_{2} & &\text{if}\  |t|\leq\frac{1}{\varepsilon}\\
&0 & &\text{if} \ |t|>\frac{1}{\varepsilon},
\end{aligned}
\right.
\end{equation*}
with
\begin{align*}
I_{1}&=\int_{|s|\leq\sqrt{\frac{1}{\varepsilon^{2}}-t^{2}}}{4\varepsilon\sqrt{t^{2}+s^{2}}(1-\varepsilon\sqrt{t^{2}+s^{2}})^{4}\,ds}=\\
&=\frac{1}{\varepsilon}\left[p_{1}(\varepsilon t)\text{acosh}\left(\frac{1}{\varepsilon|t|}\right)-p_{2}(\varepsilon t)\sqrt{1-\varepsilon^{2}t^{2}}\right],
\end{align*}
where
\begin{align*}
&p_{1}(r)=\frac{r^2}{2}(5r^4+36r^2+8), & &p_{2}(r)=\frac{1}{30}(437r^4+360r^2-8)\\
 &\text{acosh}(r)=\log(r+\sqrt{r^2-1}).
\end{align*}
While the second integral is
\begin{align*}
I_{2}&=\int_{|s|\leq\sqrt{\frac{1}{\varepsilon^{2}}-t^{2}}}{(1-\varepsilon\sqrt{t^{2}+s^{2}})^{4}\,ds}=\\
&=\frac{1}{\varepsilon}\left[-p_{3}(\varepsilon t)\text{acosh}\left(\frac{1}{\varepsilon|t|}\right)+p_{4}(\varepsilon t)\sqrt{1-\varepsilon^{2}t^{2}}\right],
\end{align*}
with
\begin{align*}
&p_{3}(r)=r^2(3r^2+4), & &p_{4}(r)=\frac{1}{15}(16r^4+83r^2+6).
\end{align*}
Then we obtain
\begin{align*}
I_{1}+I_{2}&=\frac{1}{\varepsilon}\left[(p_1-p_{3})(\varepsilon t)\text{acosh}\left(\frac{1}{\varepsilon|t|}\right)+(p_{4}-p_{2})(\varepsilon t)\sqrt{1-\varepsilon^{2}t^{2}}\right]=\\
&=\frac{1}{\varepsilon}\left[\frac{5}{2}\varepsilon^{4}t^{4}(\varepsilon^2t^2+6)\text{acosh}\left(\frac{1}{\varepsilon|t|}\right)-\frac{1}{6}(81\varepsilon^4t^4+28\varepsilon^2t^2-4)\sqrt{1-\varepsilon^{2}t^{2}}\right].
\end{align*}
We observe that $I_{1}$ and $I_{2}$ are not well defined for $t=0$, in this particular case it is easy to prove that the value of the integral is $\frac{2}{3\varepsilon}$, moreover $\lim_{t\rightarrow0}{(I_{1}+I_{2})}=\frac{2}{3\varepsilon}$ so we can define the continuous function
\begin{equation*}
g_{\varepsilon}(t)=\left\{
\begin{aligned}
&I_{1}+I_{2} & &\text{if} \ \varepsilon|t|\leq1, \ t\neq0\\
&\frac{2}{3\varepsilon} & &\text{if} \ t=0,
\end{aligned}
\right.
\end{equation*}
then $b_{j}(x)$ is given by
\begin{equation*}
b_{j}(x)=\left\{
\begin{aligned}
&g_{\varepsilon}(t_j-x\cdot v_j) & &\text{if} \ \varepsilon|t_j-x\cdot v_j|\leq1\\
&0 & &\text{if} \ \varepsilon|t_j-x\cdot v_j|>1.
\end{aligned}
\right.
\end{equation*}
For fixed $x$, $R^{y}[K(x,y)](t,\theta)$ as function of $(t,\theta)$ has compact support, but if we consider fixed values $(t_j,\theta_j)$, then $R^{y}[K(x,y)](t,\theta)=b_{j}(x)$ as function of $x$ has $\{x\in\R^{2}: \ \varepsilon|t_j-x\cdot v_j|\leq1\}$ as support, that is the strip between lines $\varepsilon(t_j-x\cdot v_j)=\pm1$ that is not limited and thus not compact.

If we now compute $a_{k,j}=R[b_{j}(x)](t_{k},\theta_{k})$ we see that this quantity can be infinity:
\begin{equation*}
a_{k,j}=R[b(x)](r,\varphi)=\int_{\R}{b(x_{s})\,ds},
\end{equation*}
where $x_{s}=(r\cos{\varphi}-s\sin{\varphi},r\sin{\varphi}+s\cos{\varphi})$. We set as usual $t-x_{s}\cdot v=as+b$, then
\begin{equation*}
a_{k,j}=\int_{\varepsilon|as+b|\leq1}{g_{\varepsilon}(as+b)\,ds}.
\end{equation*}
If $a\neq0$ we can set $u=\varepsilon(as+b)$ obtaining
\begin{align*}
a_{k,j}&=\frac{1}{\varepsilon^{2}a}\int_{|u|\leq1}{\left[\frac{5}{2}u^{4}(u^2+6)\text{acosh}\left(\frac{1}{|u|}\right)-\frac{1}{6}(81u^4+u^2-4)\sqrt{1-u^{2}}\right]\,du}=\\
&=-\frac{9}{112}\frac{\pi}{\varepsilon^2a}
\end{align*}

But if $a=0$ we have
\begin{align*}
a_{k,j}&=\frac{1}{\varepsilon^{2}}\left[\frac{5}{2}\varepsilon^{4}b^{4}(\varepsilon^2b^2+6)\text{acosh}\left(\frac{1}{\varepsilon|b|}\right)+\right.\\
&\left.-\frac{1}{6}(81\varepsilon^4b^4+28\varepsilon^2b^2-4)\sqrt{1-\varepsilon^{2}b^{2}}\right]\int_{\varepsilon|b|\leq1}{\,ds}
\end{align*}
that is 0 if $\varepsilon|b|>1$ but is $\infty$ if $\varepsilon|b|\leq1$.
\end{example}

In order to have a matrix $A$ with all finite entries, we again consider the regularization $R_{w}$ of $R$ for some weighting function $w$. Working with compactly supported radial basis functions it's natural to use another compactly supported function as weighting function, in this way we are sure that $R_{w}[b_{j}]=R[b_{j}w]$ is finite (indeed $bw$ is compactly supported) and it is possible to compute analytically  the value of $a_{k,j}$.

We consider the following case: 
\begin{align*}
K(x,y)&=\varphi_{2,0}=(1-\varepsilon\norm{x-y})_{+}^{2}, \quad \varepsilon>0,\\
w(x)&=(1-\nu^{2}\norm{x}^{2})_{+},\quad \nu>0.
\end{align*}

We start by computing 
\begin{align*}
R[K(0,y)](t,\theta)&=\int_{\R}{(1-\varepsilon\sqrt{t^2+s^2})_{+}^{2}\,ds}=\int_{\varepsilon\sqrt{t^2+s^2}\leq1}{(1-\varepsilon\sqrt{t^2+s^2})^2\,ds}=\\
&=\left\{
\begin{aligned}
&g(t) & &\text{if}\ |t|\leq\frac{1}{\varepsilon}\\
&0 & &\text{if}\ |t|>\frac{1}{\varepsilon}
\end{aligned}
\right.
\end{align*}
where
\begin{align*}
g(t)&=\int_{|s|\leq\sqrt{\frac{1}{\varepsilon^2}-t^2}}{(1-\varepsilon\sqrt{t^2+s^2})^2\,ds}=\\
&=\left\{
\begin{aligned}
&\frac{2}{\varepsilon}\left[\frac{\sqrt{1-\varepsilon^2t^2}}{3}(2\varepsilon^2t^2+1)-\varepsilon^2t^2\text{acosh}\left(\frac{1}{\varepsilon|t|}\right)\right] & &\text{if}\ t\neq0\\
&\frac{2}{3\varepsilon} & &\text{if}\ t=0
\end{aligned}
\right.
\end{align*}
We conclude that
\begin{align*}
b_{j}(x)=R^{y}[K(x,y)](t_j,\theta_j)=\left\{
\begin{aligned}
&g(t_j-x\cdot v_j) & &\text{if}\ |t_j-x\cdot v_j|\leq\frac{1}{\varepsilon}\\
&0 & &\text{if}\ |t_j-x\cdot v_j|>\frac{1}{\varepsilon}
\end{aligned}
\right.
\end{align*}

As in Example \ref{ex: radon_wend}, it is possible to show that $a_{k,j}=\infty$ if $a=0$ and $\varepsilon|b|\leq1$, where $a=\sin{(\theta_{k}-\theta_{j})}$ and $b=t_{j}-t_k\cos{(\theta_{k}-\theta_{j})}$, so we introduce the regularization of $R$ by multiplication for the weighting function $w=(1-\nu^{2}\norm{x}^{2})_{+}$ and we have
\begin{equation*}
a_{k,j}=R_{w}[b_{j}(x)](t_{k},\theta_{k}).
\end{equation*}

Using simpler notation,
\begin{equation*}
R[wb(x)](r,\varphi)=\int_{\R}{b(x_{s})(1-\nu^{2}\norm{x_{s}}^2)_{+}\,ds}.
\end{equation*}
The computation of this integral can be found in appendix B.

\section{Scaled problem}\label{sec: scaled_problem}
The Hermite-Birkhoff reconstruction problem can be expressed as finding a function $s\in S$ satisfying 
\begin{equation}
\lambda_{j}(s)=\lambda_{j}(f) \quad \forall\ j=1,\ldots,n,
\label{eq: hb_prob}
\end{equation} 
where $\lambda_{1},\ldots,\lambda_{n}$ are linearly independent linear operators and
\begin{equation*}
S=\left\{\sum_{j=1}^{n}{c_{j}\lambda_{j}^{y}K(\cdot,y)}:\ c_{j}\in\R,\ K \text{positive definite kernel}\right\}.
\end{equation*}
By linearity equation \ref{eq: hb_prob} can be written as a linear system $Ac=f|_{\Lambda}$, where the elements of the matrix $A$ are given by $a_{k,j}=\lambda^{x}_{k}\lambda^{y}_{j}[K(x,y)]$ and $(f|_{\Lambda})_{j}=\lambda_{j}(f)$.

Since the matrix $A$ can be highly ill-conditioned, it can be convenient to consider the scaled problem (\cite{ISKE1},\cite{ISKE2}). For $h>0$ the scaled reconstruction problem is  
\begin{equation*}
\lambda_{j}(s^{h}(h\cdot))=\lambda_{j}(f(h\cdot)) \quad \forall\ j=1,\ldots,n,
\end{equation*} 
where
\begin{equation*}
s^{h}\in S^{h}=\left\{\sum_{j=1}^{n}{c_{j}\lambda_{j}^{y}K(\cdot,hy)}:| c_{j}\in\R\right\}.
\end{equation*} 
In this way one obtains the linear system $A^{h}c=f^{h}|_{\Lambda}$ given by
\begin{align*}
a_{k,j}^{h}&=\lambda_{k}^{x}\lambda_{j}^{y}[K(hx,hy)], & &(f^{h}|_{\Lambda})_{j}=\lambda_{j}(f(h\cdot)).
\end{align*}

In the case of image reconstruction the operator $\lambda_{j}$ represents the Radon transform evaluated at point $(t_{j},\theta_{j})$. Thus, in order to compute $a_{k,j}^{h}$ and $f_{j}^{h}$ one has to understand which is the relationship between the Radon transform of a function $f$ and the Radon transform of the scaled function $f(h\cdot)$. This relationship is given by the following

\begin{theorem}[Dilatation-property of the Radon transform]\label{thm: dilatation_radon}
Let $f:\R^{2}\rightarrow\R$ be such that $R[f(x)](t,\theta)=g(t,\theta)$, then for all $h>0$
\begin{equation*}
R[f(hx)](t,\theta)=\frac{1}{h}g(ht,\theta).
\end{equation*}
\proof
\begin{align*}
R[f(hx)](t,\theta)&=\int_{\R}{f(ht\cos{\theta}-hs\sin{\theta},ht\sin{\theta}+hs\cos{\theta})\,ds}=\\
&=\int_{\R}{f(ht\cos{\theta}-r\sin{\theta},ht\sin{\theta}+r\cos{\theta})\,\frac{dr}{h}}=\\
&=\frac{1}{h}R[f(x)](th,\theta).
\end{align*}
\endproof
\end{theorem}
Thanks to this property it's possible to compute $a_{k,j}^{h},\ f_{k}^{h}$ and $b_{j}^{h}(x)=\lambda_{j}^{y}K(x,hy)$:
\begin{align*}
&f_{k}^{h}=\lambda_{k}^{x}(f(hx))=\frac{1}{h}R[f(x)](ht_{k},\theta_{k})\\
&b_{j}^{h}(x)=\frac{1}{h}R^{y}[K(\cdot,y)](ht_{j},\theta_{j})\\
&a_{k,j}^{h}=\lambda_{k}^{x}\lambda_{j}^{y}[K(hx,hy)]=\frac{1}{h}\lambda_{k}^{x}[b_{j}^{h}(hx)]=\frac{1}{h^2}R[b_{j}(x)](ht_{k},\theta_{k})
\end{align*}

We note that in order to compute $f_{k}^{h}$ we need to know the Radon transform of the unknown function $f$ at points $(ht_{k},\theta_{k})$ that is not possible if the X-ray machine gives us data only at points $(t_{k},\theta_{k})$. However, for our aim, we assume to know the analytical expression of $R[f]$.

The solution of the scaled reconstruction problem is 
\begin{equation*}
s^{h}(hx)=\sum_{j=1}^{n}{c_{j}b_{j}^{h}(hx)},
\end{equation*}
with $c=(c_{1},\ldots,c_n)^{T}$ solution of the linear system $A^{h}c=f^{h}$. It is important to notice that thanks to the dilatation-property \ref{thm: dilatation_radon} we don't have to compute the Radon transform for every different value of $h$, but we only need to compute it in the case $h=1$, multiply for $\frac{1}{h}$ and then scaling evaluation points $(t_{k},\theta_{k})$ to $(ht_{k},\theta_{k})$.

%Chapter 6 - Discussion

\chapter{Numerical results}
In the previous chapters we saw some theoretical tools that can be used to obtain the value of a function $f:\R^{2}\rightarrow\R$ starting from a sampling of its Radon transform.

In chapter \ref{chap: fourier_methods} we studied the continuous problem and we found an analytical inversion formula for the Radon transform: the back projection formula. Then, in order to use this formula in real applications, we introduced the process of linear filtering and interpolation.

In chapter \ref{chap: art} and \ref{chap: kernelMethods} we followed a different approach: starting form the discrete problem for finding an approximation of a function, belonging to a particular finite dimensional space of functions, such that its Radon transform coincides with the measured Radon transform of the unknown function $f$. The Kaczmarz's method consider pixel basis functions to determine an approximation of $f$, while kernel-based methods use positive definite functions to generate a functions space where to find a solution. We also saw that in this second case the problem need some kind of regularization so that the Radon transform of the kernel functions is well defined.

What we want to do now is to compare all these methods from a numerical point of view, studying the behavior of the solution and the approximation error in function of the data and the parameters involved in the algorithms.

\section{Optimal parameters}
In chapter \ref{chap: kernelRec} we introduced a regularization technique for solving  the Hermite-Birkhoff interpolation problem of image reconstruction using kernel based methods. In particular, the original problem was to find $s=\sum_{j=1}^{n}{c_{j}\lambda_{j}K(\cdot,y)}$ such that $\lambda_{k}f=\lambda_{k}s$ for all $k=1,\ldots,n$, where $\lambda_{k}g=R[g(x)](t_{k},\theta_{k})$ and $K$ is a positive definite kernel. By linearity this problem is equivalent to solve the linear system 
\begin{equation*}
\lambda_{k}f=\sum_{j=1}^{n}{c_{j}\lambda_{k}^{x}\lambda_{j}^{y}K(x,y)}, \quad k=1,\ldots,n.
\end{equation*}
 
The main problem found in applying this method is that the Radon transform $\lambda_{k}^{x}\lambda_{j}^{y}K(x,y)$ or $\lambda_{j}^{y}K(x,y)$ can be infinity.

The solution we adopted was to consider kernel functions $K$ such that $b_{j}(x)=\lambda_{j}^{y}K(x,y)$ is well defined (e.g. multiplying any kernel $K(x,y)$ for a suitable function $\phi(\norm{x-y})$) and to substitute operator $R$ with another linear operator so that,  when computing matrix $A=(a_{k,j})=(\lambda_{k}b_{j}(x))$, we have $a_{k,j}<\infty$ for all $k,j=1,\ldots,n$.

In our discussion we chose operator $R_{w}$ defined by $R_{w}[g]=R[gw]$ where $w$ is an appropriate window function.

Both kernel $K$ and window function $w$ depend on one or more parameters; in this section we will discuss, thanks to numerical experiment, how these parameters influence the quality of the reconstructed image. In order to do that we will apply our methods on predefined phantoms and by varying a parameter we will see the behavior of the solution. The error measure we will use to determine the quality of the result is the \emph{root mean square error}
\begin{equation*}
RMSE=\sqrt{\frac{\sum_{j=1}^{m}{(x_{i}-\hat{x}_{i})^{2}}}{m}},
\end{equation*} 
where $m$ is the dimension of the image, $x_{i}$ and $\hat{x}_{i}$ the gray scale value of pixel $i$ of the original and reconstructed image respectively. More the $RMSE$ is close to zero, more the solution will be considered accurate.

\subsection{Window function parameters}
We begin our analysis considering the parameter that influence operator $R_{w}$ and the window functions introduced in sections \ref{sec: imqRec} and \ref{sec: gaussRec}, i.e. $w_{\nu}(x)=\exp{(-\nu^{2}\norm{x}^{2})}$ and $w_{L}=\chi_{[-L,L]}(\norm{x})$.

Let us start with the truncated inverse multiquadric kernel 
\begin{equation*}K(x,y)=(1+\varepsilon^{2}\norm{x-y}^{2})^{-1/2}\chi_{[-L_{1},L_{1}]}(\norm{x-y})
\end{equation*}
and the characteristic window function $w_{L_{2}}=\chi_{[-L_{2},L_{2}]}(\norm{x})$, with $L_{1},L_{2}>2\max{|t_{j}|}$ (as we saw in section \ref{sec: imqRec}). 

\begin{figure}[htbp]
\centering%
\subfigure[Crescent-shaped phantom \label{fig: mse_L2_imq5}]%
{\includegraphics[width=0.49\textwidth]{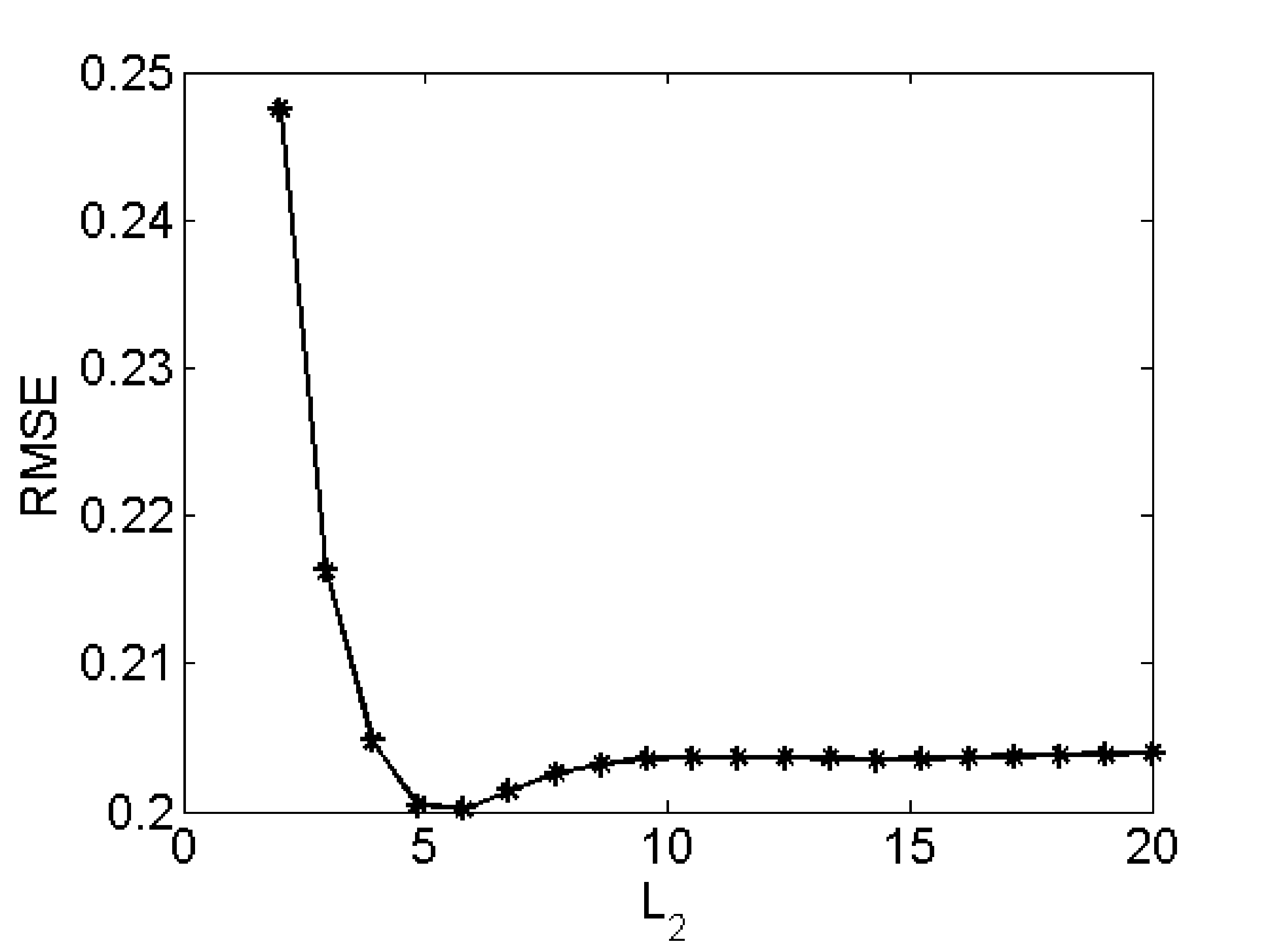}}
\subfigure[Bull's eye phantom \label{fig: mse_L2_imq6}]%
{\includegraphics[width=0.49\textwidth]{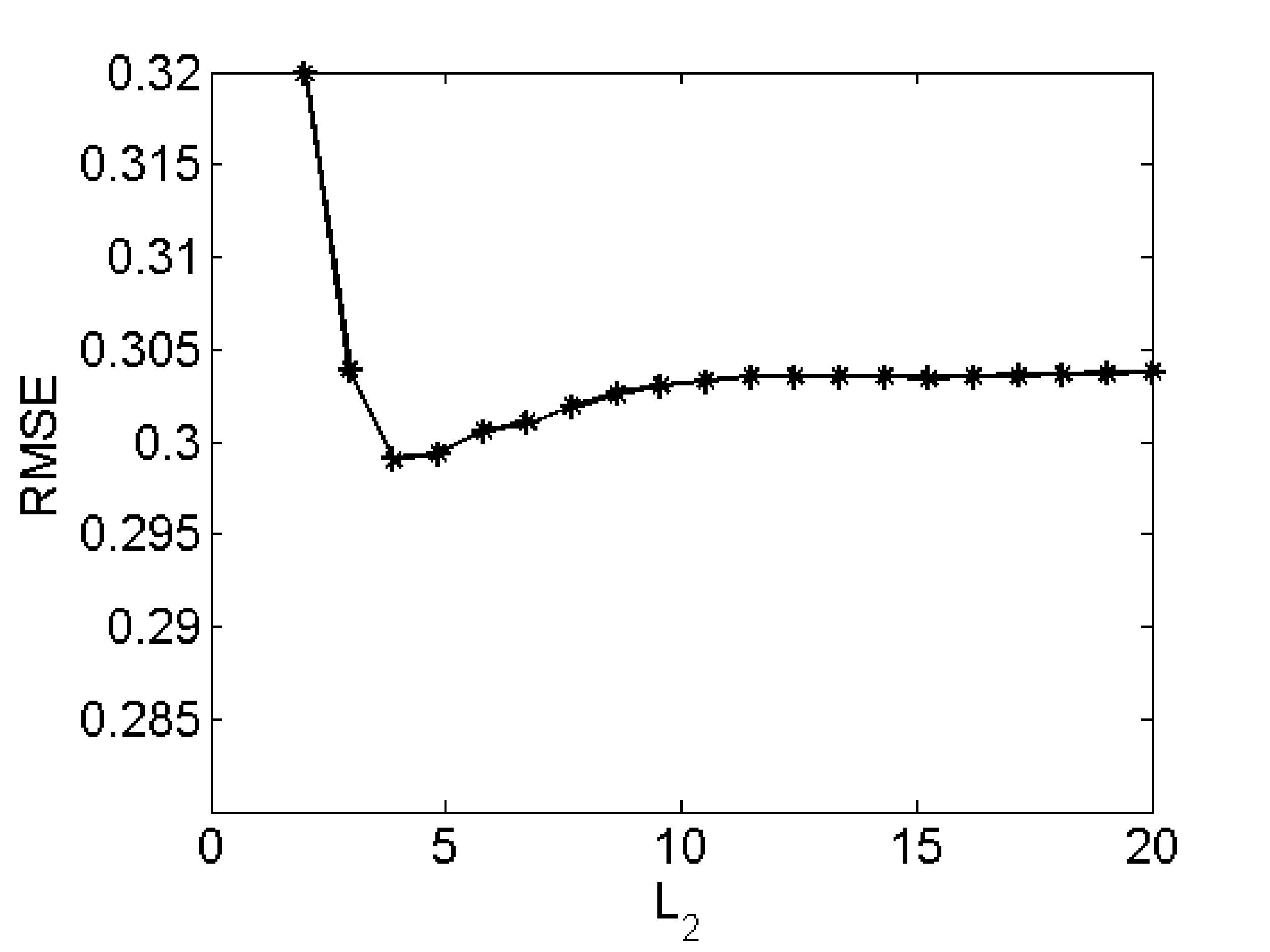}}
\caption{$RMSE$ of inverse multiquadric reconstruction in function of parameter $L_{2}$, with fixed number of samples $N=30,\ M=20,\ K=256$ and parameters $ \varepsilon=30,\ L_{1}=50.$}
\label{fig: mse_L2_imq56}
\end{figure}

\begin{figure}[htbp]
\centering%
\subfigure[$N=30,\ M=20$ \label{fig: mse_L2_imq1}]%
{\includegraphics[width=0.49\textwidth]{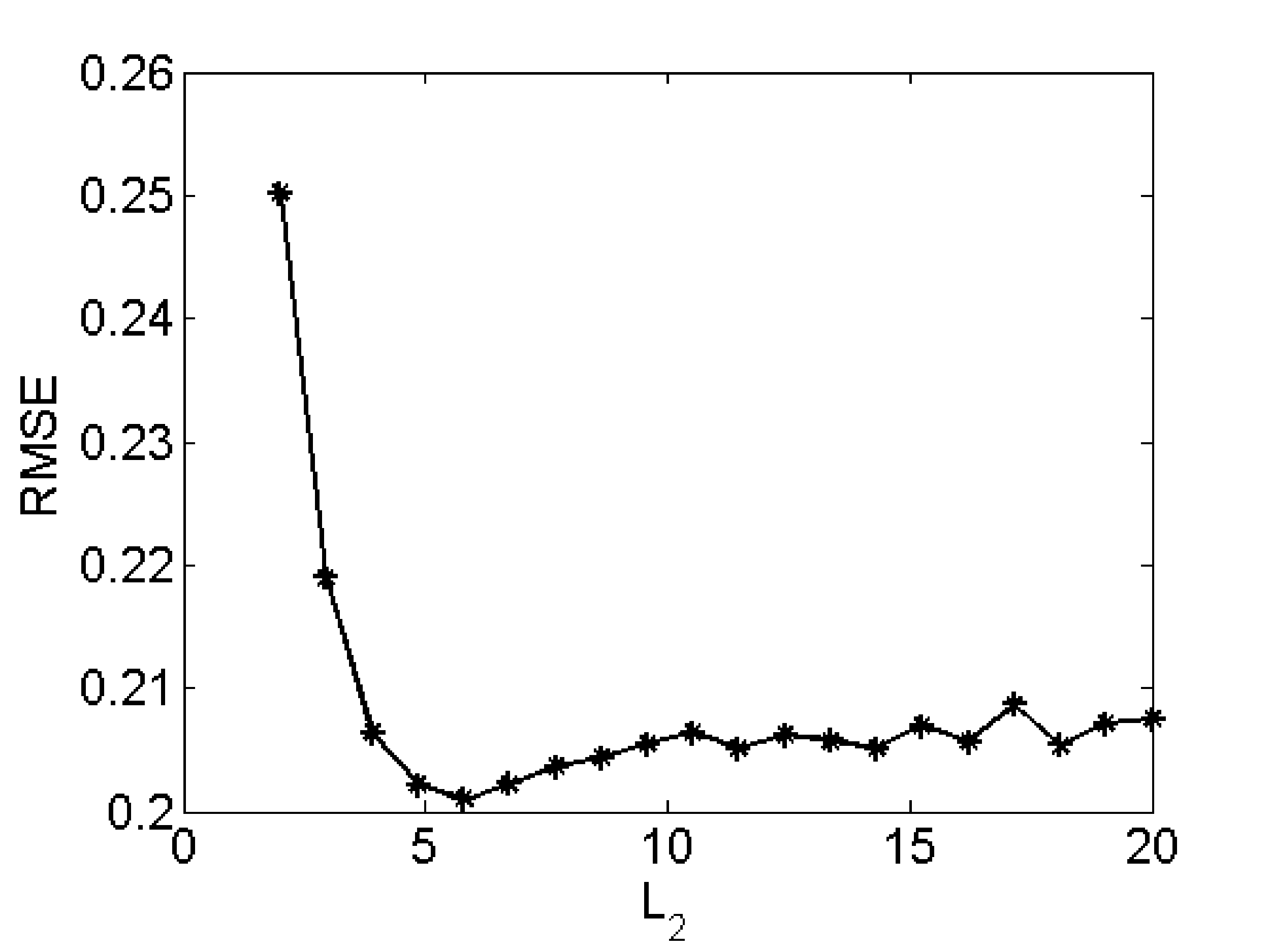}}%,height=0.3\textwidth
\subfigure[$N=50,\ M=40$ \label{fig: mse_L2_imq3}]%
{\includegraphics[width=0.49\textwidth]{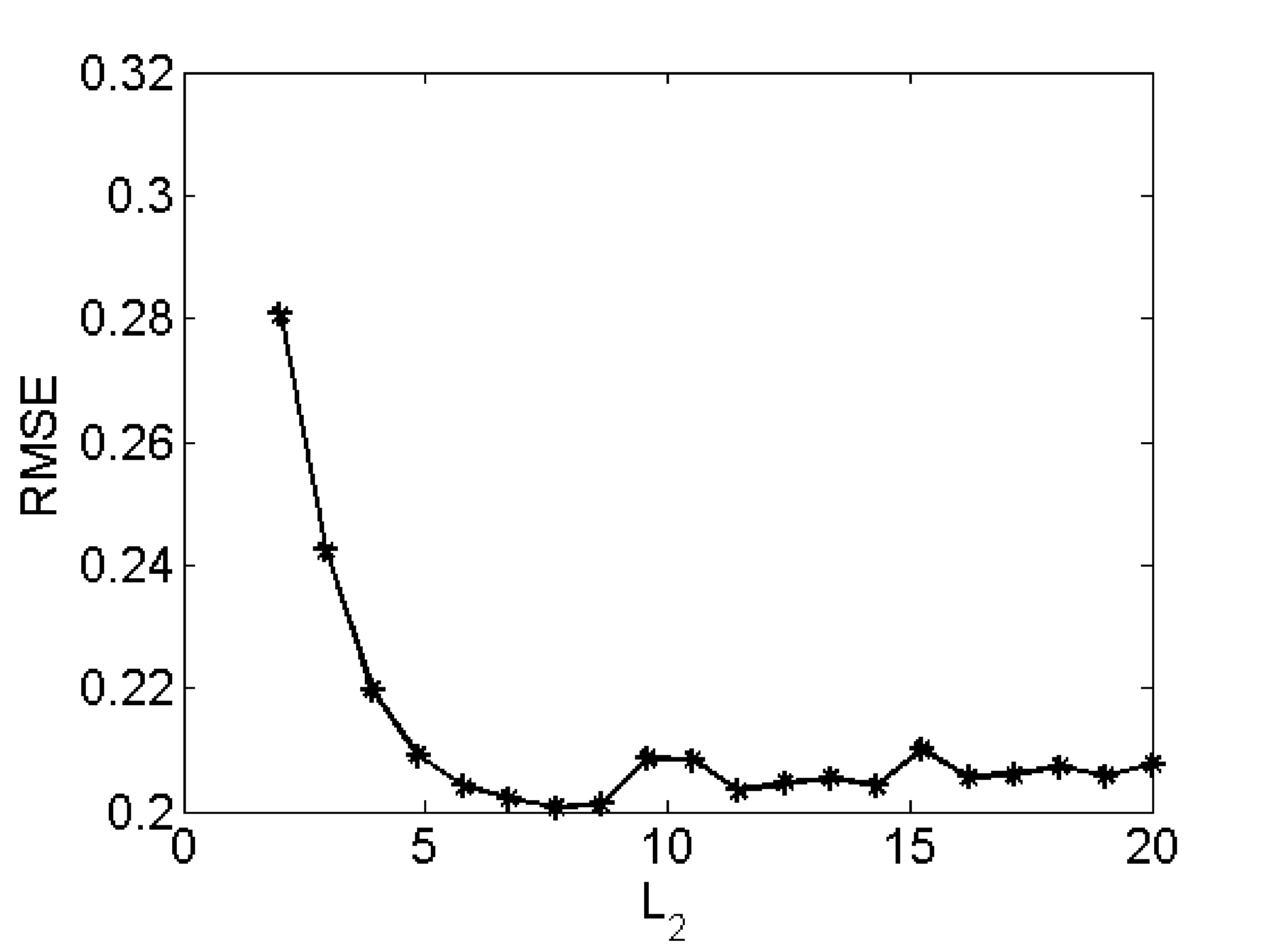}}
\caption{$RMSE$ of inverse multiquadric reconstruction as a function of the parameter $L_{2}$, with fixed output dimension $K=64$ and parameters $ \varepsilon=30,\ L_{1}=10$ for the crescent-shaped phantom.}
\label{fig: mse_L2_imq13}
\end{figure}

Varying the parameter $L_{2}$ for fixed values of $\varepsilon$, $L_{1}$ and data\footnote{in this chapter we will always consider the parallel beam geometry as acquisition method of data.} $\{(t_{j},\theta_{j})\}_{j=1}^{n}$ and applying this reconstruction technique to three different phantoms we can see that the $RMSE$ presents a minimum for a particular value $L_{2opt}$ (Figure \ref{fig: mse_L2_imq56}). This optimal value is influenced by $\varepsilon$, $L_{1}$, $n$ (in particular if the number of data increases, also $L_{2opt}$ increases - cfr. Figures \ref{fig: mse_L2_imq1} and \ref{fig: mse_L2_imq3}). We also notice that the $RMSE$ decrease very rapidly for $L_{2}<L_{2opt}$ but for $L_{2}\geq L_{2opt}$ the $RMSE$ is increasing with a very small rate, so one should choose $L_{2}$ in a way to be sure that $L_{2}\geq L_{2opt}$.

\begin{figure}[htbp]
\centering%
\subfigure[Parameters and data as in Figure \ref{fig: mse_L2_imq1} \label{fig: rcond_L2_imq1}]%
{\includegraphics[width=0.49\textwidth]{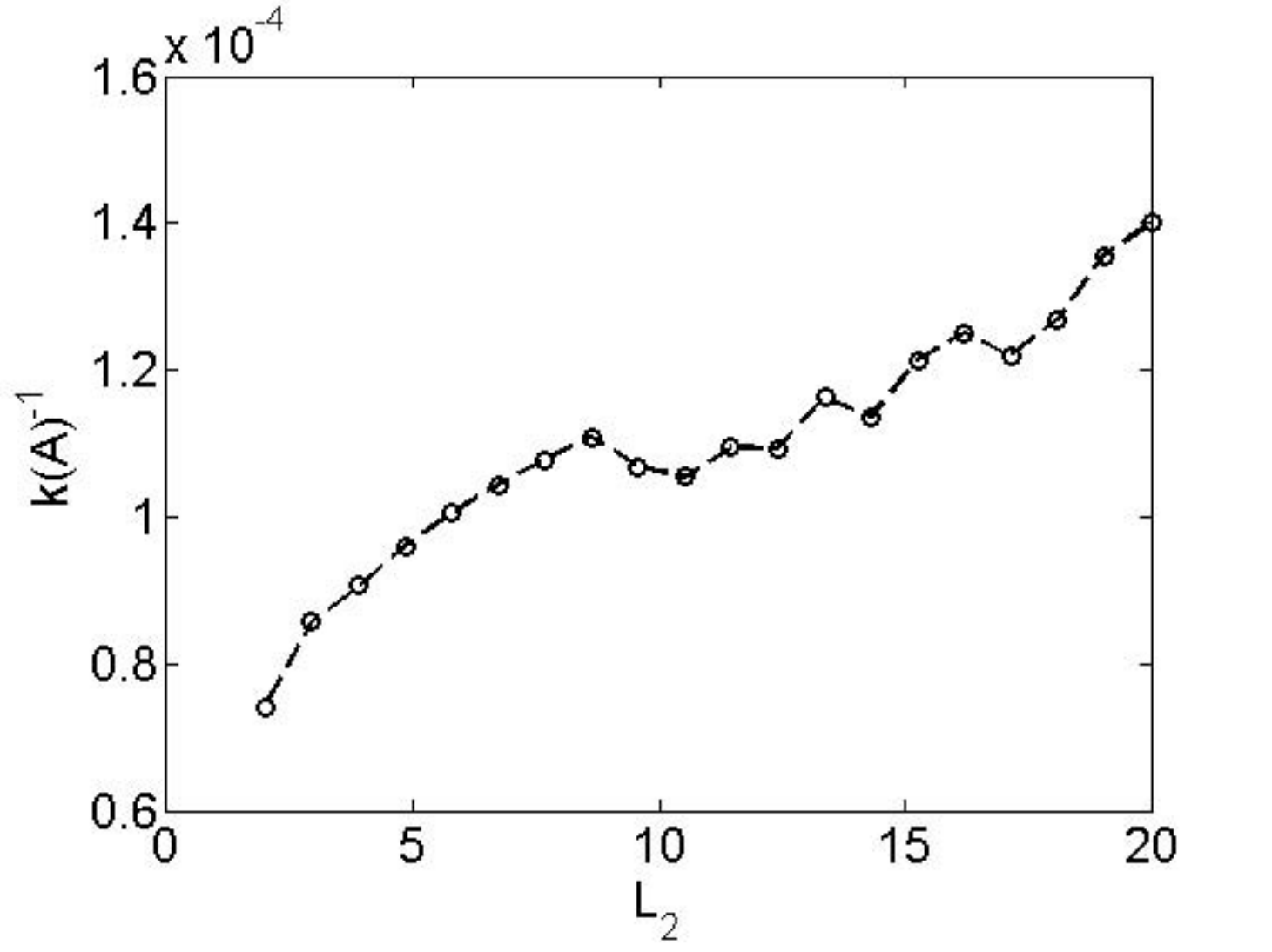}}%,height=0.3\textwidth
\subfigure[Parameters and data as in Figure \ref{fig: mse_L2_imq6} \label{fig: mse_L2_imq6}]%
{\includegraphics[width=0.49\textwidth]{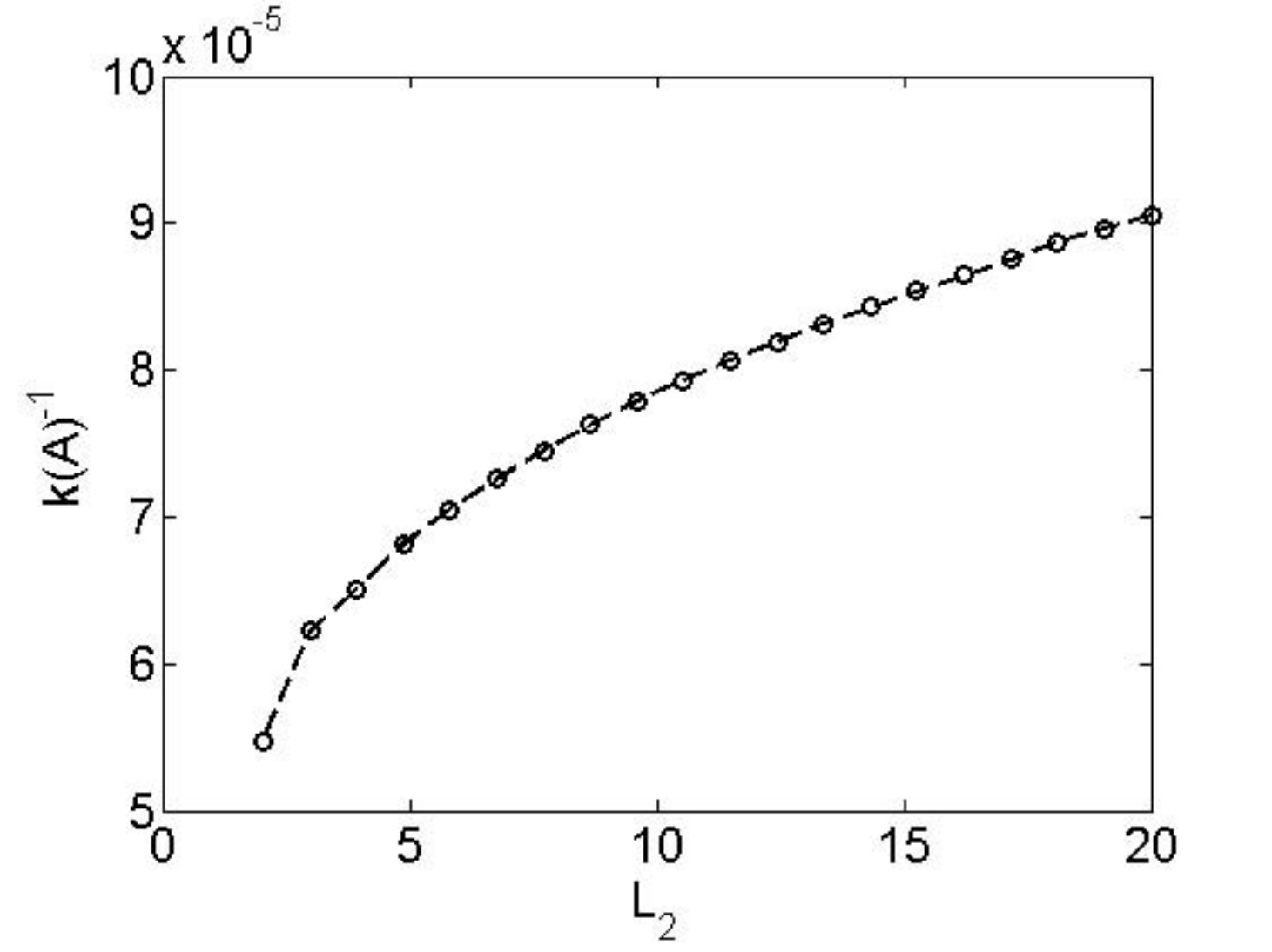}}
\caption{Reciprocal of the condition number of $A$ for inverse multiquadric reconstruction as a function of the parameter $L_{2}$}
\label{fig: rcond_L2_imq16}
\end{figure}

Another advantage in choosing $L_{2}$ large is that the condition number of the matrix $A$ is smaller (see Figure \ref{fig: rcond_L2_imq16}).

However the value of $L_{2}$ does not determine so drastically the behavior of the solution, whose quality remains acceptable for large values of $L_{2}$.

\begin{figure}[htbp]
\centering%
\subfigure[ \label{fig: mse_nu_gauss1}]%
{\includegraphics[width=0.49\textwidth]{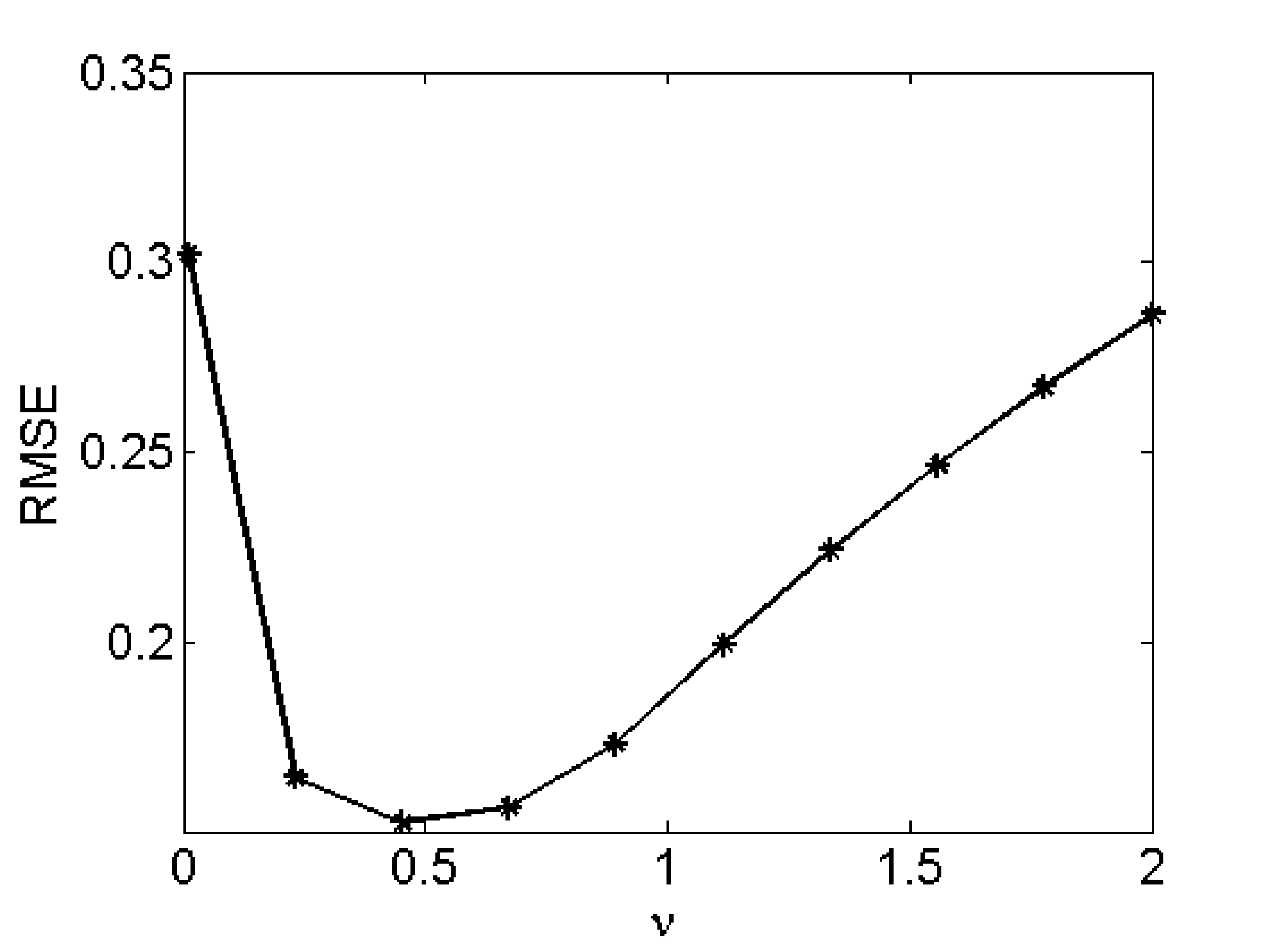}}%,height=0.3\textwidth
\subfigure[ \label{fig: mse_nu_gauss2}]%
{\includegraphics[width=0.49\textwidth]{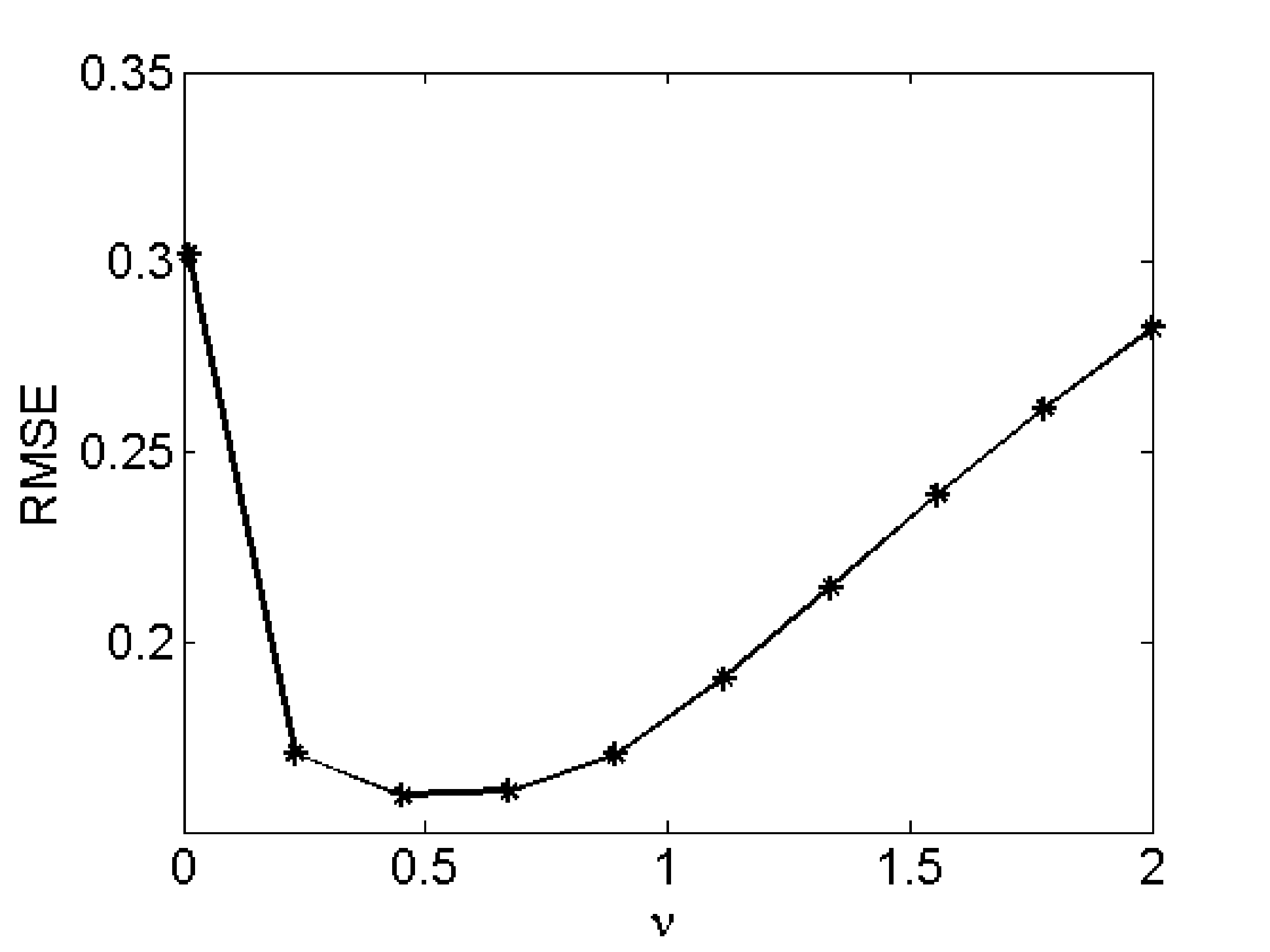}}
\subfigure[\label{fig: mse_nu_gauss3}]%
{\includegraphics[width=0.49\textwidth]{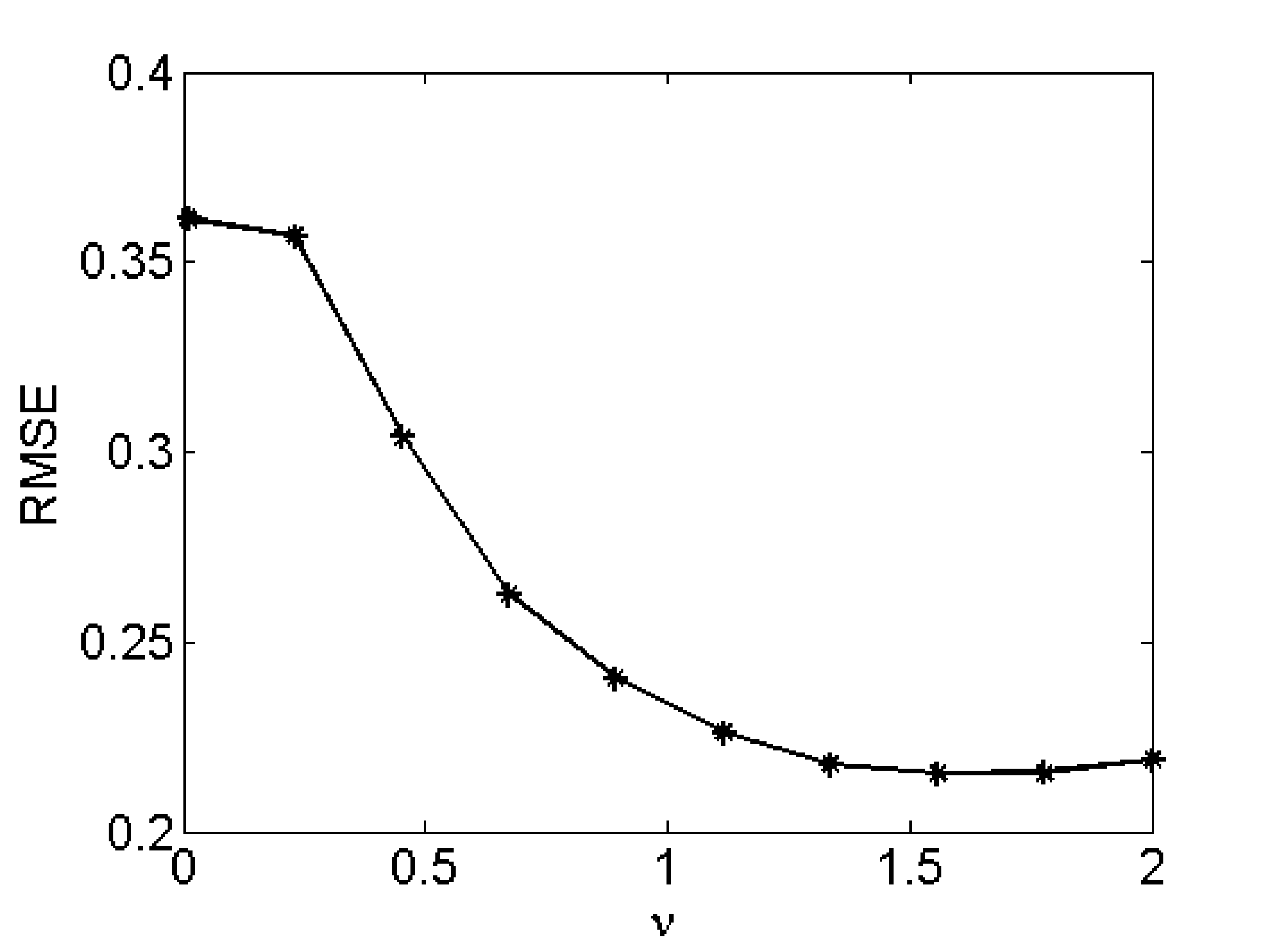}}%,height=0.3\textwidth
\subfigure[  \label{fig: mse_nu_mq2}]%
{\includegraphics[width=0.49\textwidth]{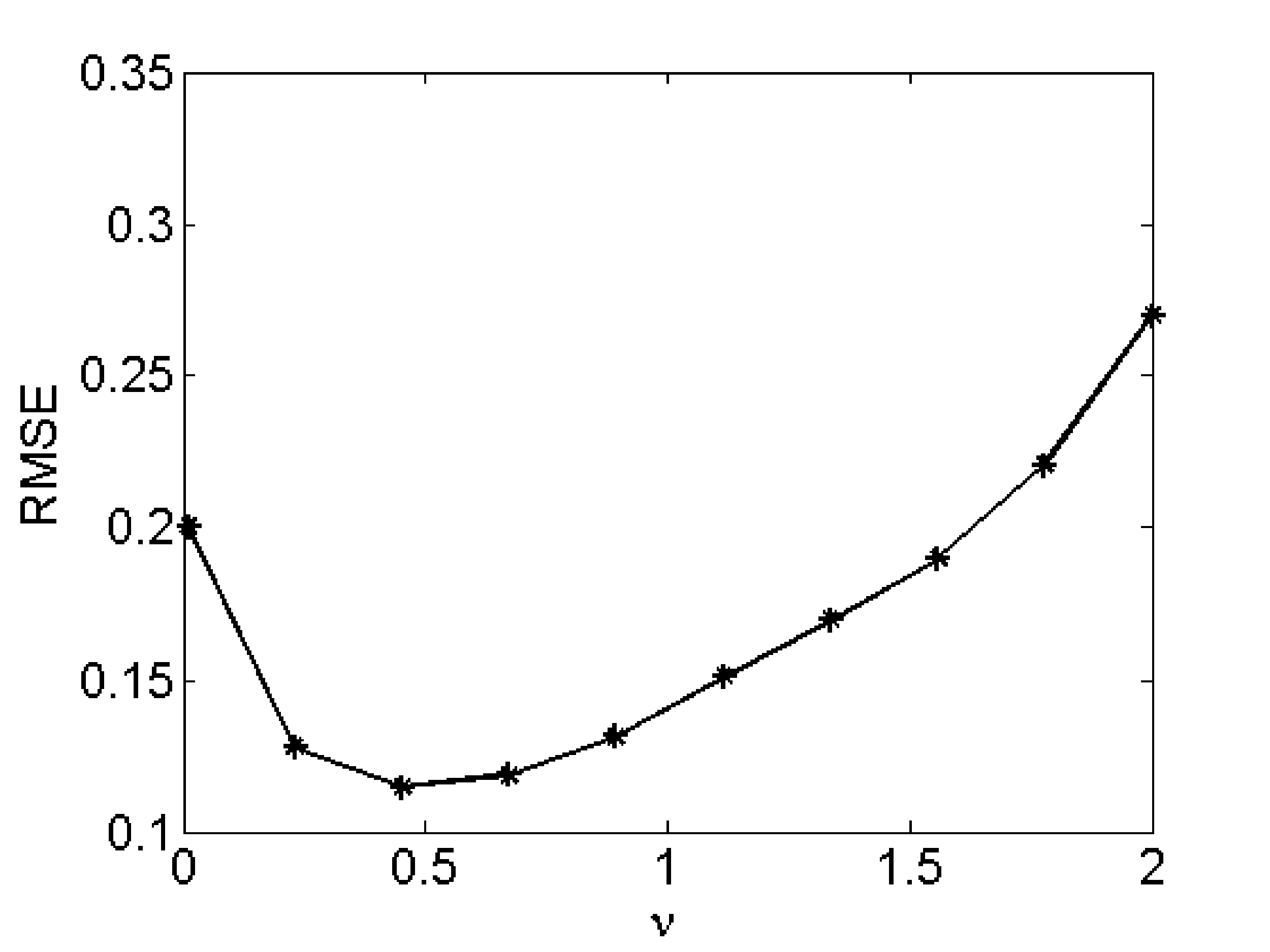}}
\caption{$RMSE$ of Gaussian and multiquadric  reconstruction as a function of the parameter $\nu$. (a) Gaussian kernel, bull's eye phantom, $N=50,\ M=40,\ K=64,\ \varepsilon=30$; (b) Gaussian kernel, bull's eye phantom, $N=50,\ M=40,\ K=64,\ \varepsilon=60$; (c) Gaussian kernel, Shepp-Logan phantom, $N=30,\ M=20,\ K=256,\ \varepsilon=50$  (d) Multiquadric kernel, crescent-shaped phantom, $N=30,\ M=20,\ K=64,\ \rho=1,\ \varepsilon=30$.}
\label{fig: mse_nu}
\end{figure}

More interesting is the case of the Gaussian window function $w(x)=e^{-\nu^{2}\norm{x}^{2}}$. Also in this case there exists an optimal value $\nu_{opt}$ such that $RMSE$ is minimum. But now for $\nu>\nu_{opt}$ the $RMSE$ increases with a fast rate and so the quality of the reconstruction becomes worse (see for example Figures \ref{fig: mse_nu_gauss3} and \ref{fig: mse_nu_mq2}).

The value $\nu_{opt}$ depends on the phantom used, i.e. on data. This is not surprising because we know that the approximation error depends on $|f_{w,k}-f_{k}|$ (see section \ref{subsec: regularization}). On the other hand $\nu_{opt}$ has only small variation w.r.t. the changing of other shape parameters (e.g. is independent on $\varepsilon$ in the case of Gaussian kernel - see Figures \ref{fig: mse_nu_gauss1} and \ref{fig: mse_nu_gauss2}).

\begin{figure}[htbp]
\centering%
\subfigure[ \label{fig: rcond_nu_gauss1}]%
{\includegraphics[width=0.49\textwidth]{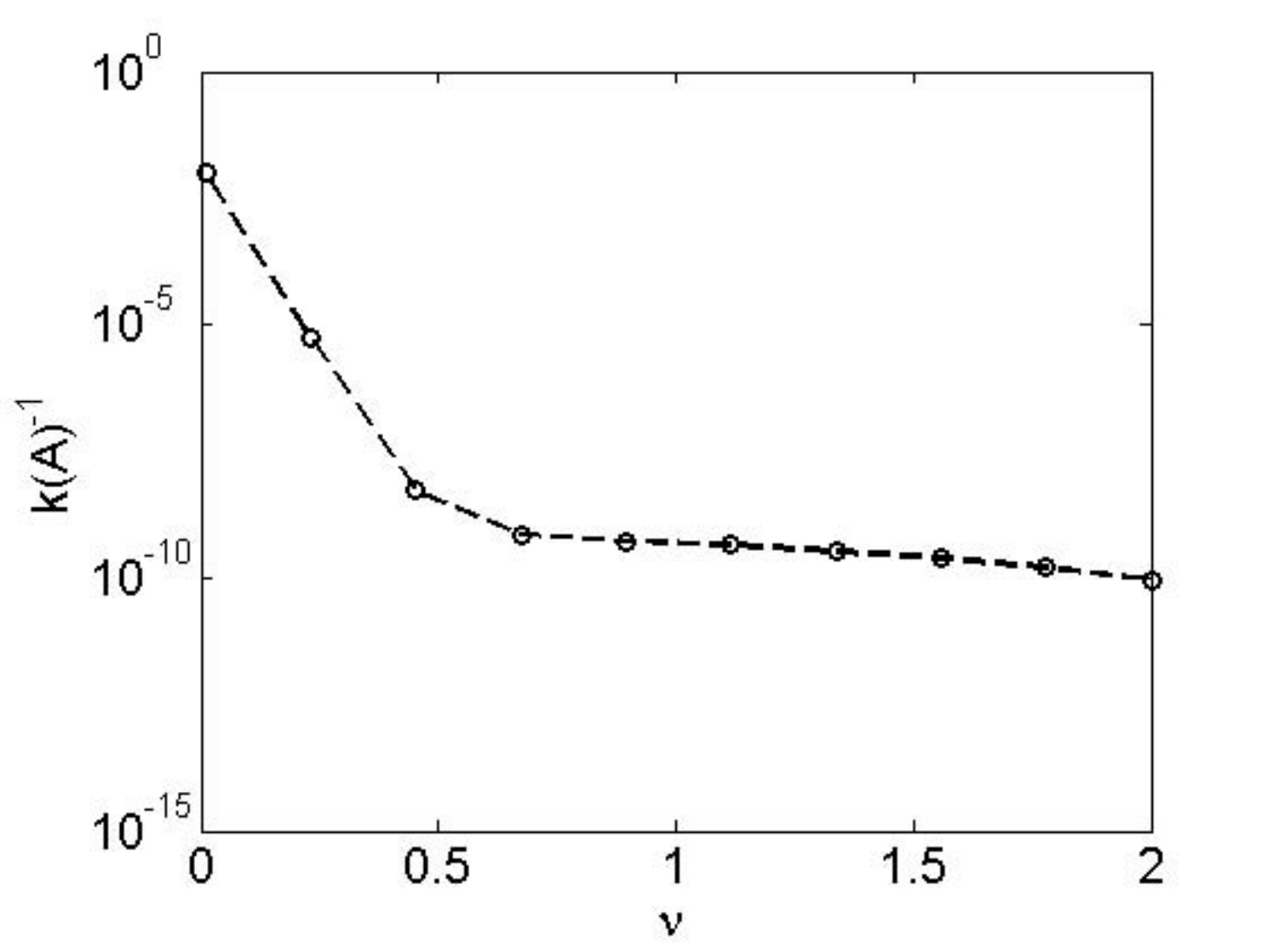}}%,height=0.3\textwidth
\subfigure[\label{fig: rcond_nu_gauss3}]%
{\includegraphics[width=0.49\textwidth]{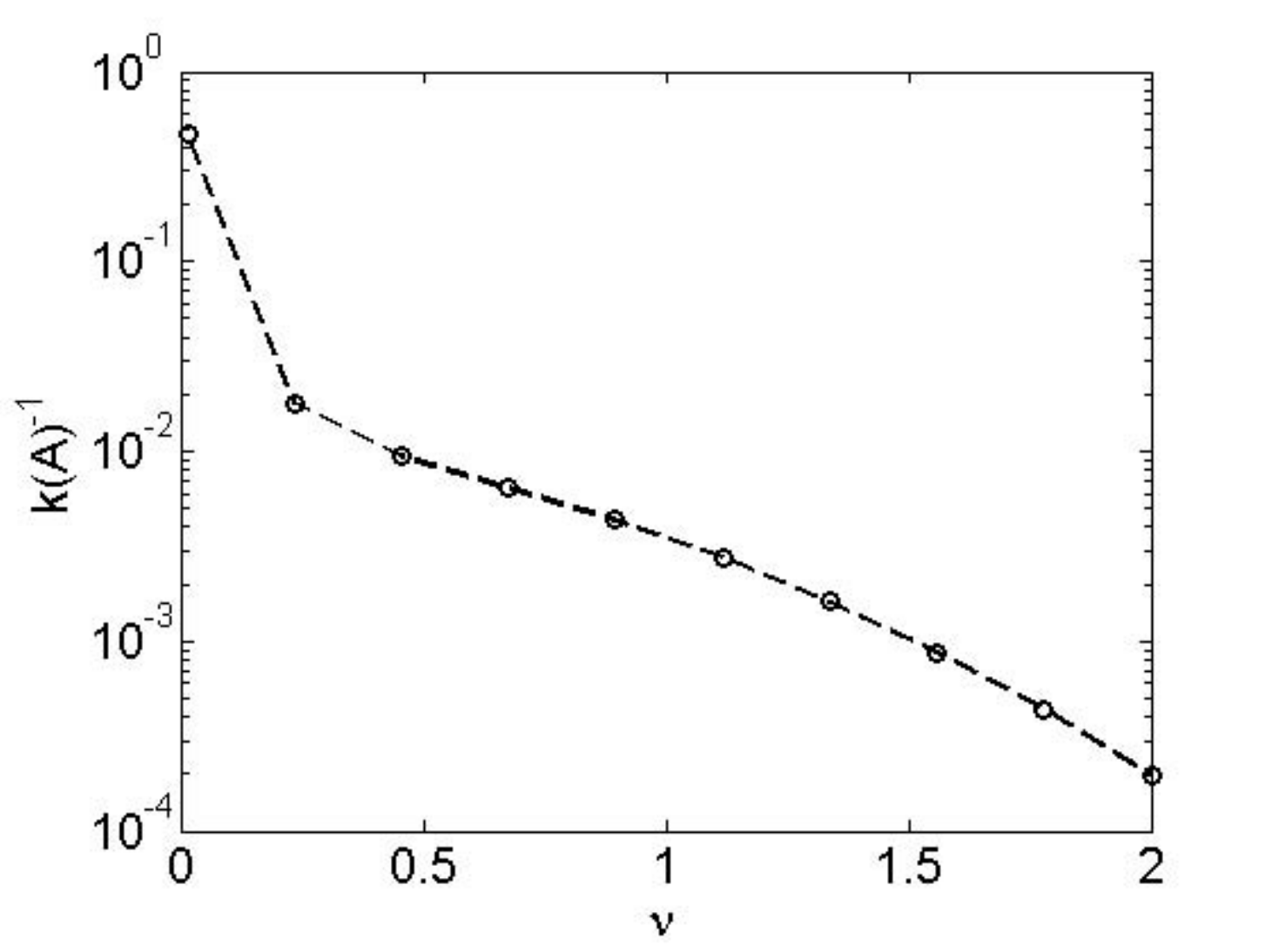}}
\caption{Reciprocal of the condition number of $A$ for Gaussian kernel reconstruction as a function of the parameter $\nu.$ (a) Parameters and data as in Figure \ref{fig: mse_nu_gauss1}; (b) Parameters and data as in Figure \ref{fig: mse_nu_gauss3}.}
\label{fig: rcond_nu}
\end{figure}

The fact that the $RMSE$ is increasing for $\nu>\nu_{opt}$ can be explained considering the condition number $k(A)$ of the system matrix $A$. In fact, $k(A)$ increases with $\nu$ (see Figure \ref{fig: rcond_nu} where the reciprocal of $k(A)$ is plotted in function of $\nu$). The quantity $k(A)=k_{1}(A)$ is the 1-norm condition number of the matrix $A$. Its inverse is estimated using the MATLAB function \texttt{rcond} (see \cite{MATLAB}).

%\begin{figure}[htbp]
%\centering
%\includegraphics[width=0.7\textwidth]{mse_L2_imq1.pdf} 
%\caption{Prova 3}
%\label{fig: mse_L2_imq}
%\end{figure} 

%\begin{center}
%\begin{tabular}[htbp]{c|*{5}{r}}
%$L_{2}$ & 2.0000  &   2.9474   & 3.8947 &   4.8421 &   5.7895\
%%\hline
%$MSE$ & 0.0626 &   0.0479  &  0.0426 &   0.0409 &   0.0404 \\
%$k(A)^{-1}$ & 0.0742  &  0.0857  &  0.0906  &  0.0961 &   0.1006\\
%\hline\\
%\hline
%%\end{tabular}
%
%%\begin{tabular}[htbp]{c|*{5}{r}}
%$L_{2}$ & 6.7368   & 7.6842  &  8.6316 &   9.5789 &  10.5263\\
%\hline
%$MSE$ & 0.0409  &  0.0415   & 0.0418  &  0.0422 &   0.0426 \\
%$k(A)^{-1}$ & 0.1044  &  0.1077  &  0.1108  &  0.1068 &   0.1054\\
%\hline\\
%\hline
%%\end{tabular}
%
%%\begin{tabular}[htbp]{c|*{5}{r}}
%$L_{2}$ & 11.4737 &   12.4211  & 13.3684 &  14.3158 &  15.2632\\
%\hline
%$MSE$ & 0.0421  &  0.0425 &   0.0423 &   0.0421 &   0.0428 \\
%$k(A)^{-1}$ & 0.1094 &   0.1093   & 0.1163  &  0.1136  &  0.1214\\
%\hline\\
%\hline
%%\end{tabular}
%
%%\begin{tabular}[htbp]{c|*{5}{r}}
%$L_{2}$ & 16.2105  & 17.1579 &  18.1053  & 19.0526 &  20.0000\\
%\hline
%$MSE$ & 0.0423  &  0.0435 &   0.0422 &   0.0429 &   0.0431 \\
%$k(A)^{-1}$ & 0.1250  &  0.1219 &   0.1270 &   0.1356 &   0.1400\\
%\hline
%\end{tabular}
%%
%\end{center}

The case of multiquadrics $K(x,y)=\sqrt{1+\rho^{2}\norm{x-y}}e^{-\varepsilon^{2}\norm{x-y}}$ with Gaussian window function is similar to the Gaussian kernel case, provided that $\rho$ is small enough, as explained in the next paragraph.

\begin{figure}[htbp]
\centering%
\subfigure[Crescent-shaped phantom \label{fig: mse_nu_compact_cre}]%
{\includegraphics[width=0.49\textwidth]{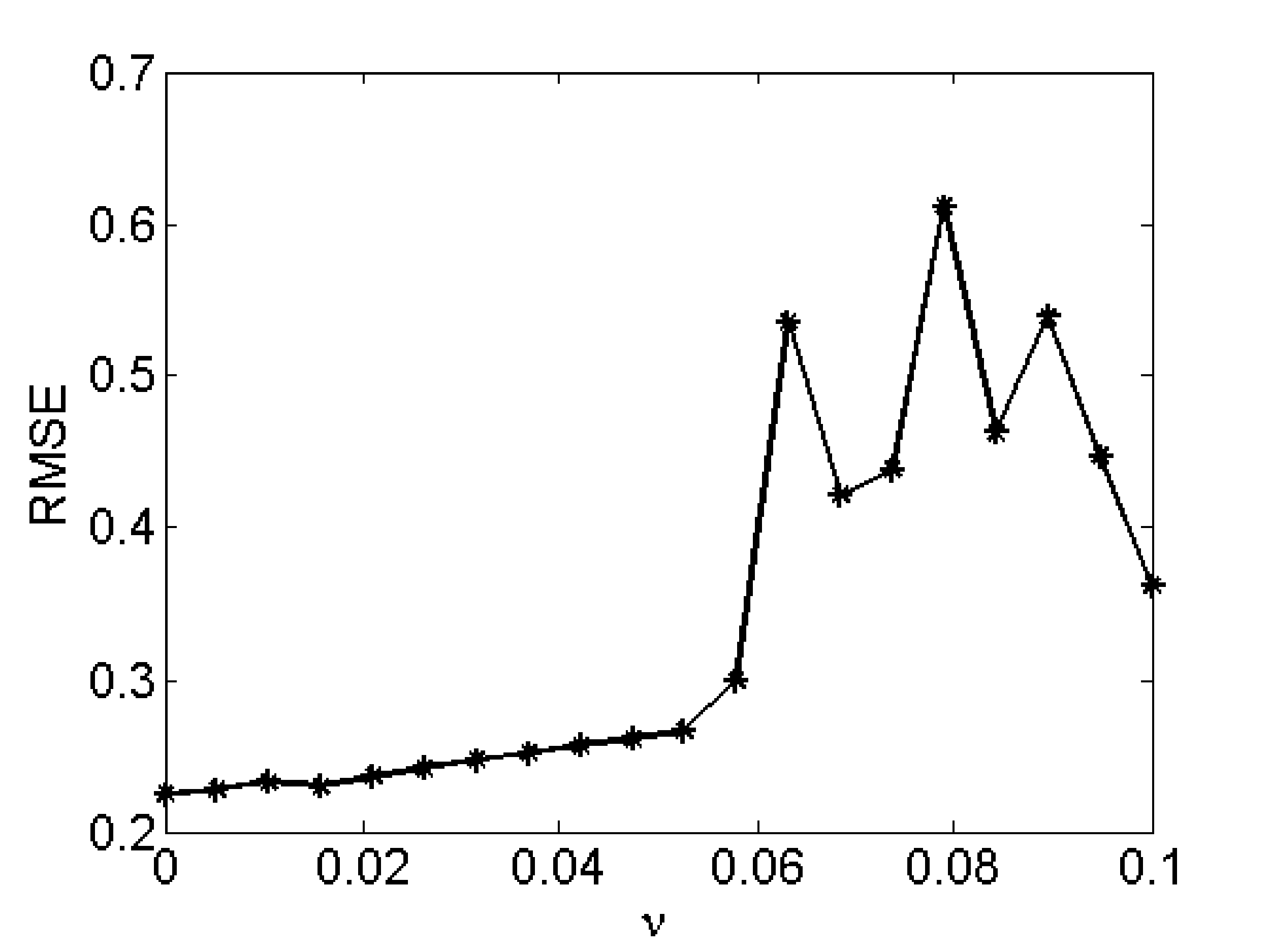}}%,height=0.3\textwidth
\subfigure[Bull's eye phantom \label{fig: mse_nu_compact_bull}]%
{\includegraphics[width=0.49\textwidth]{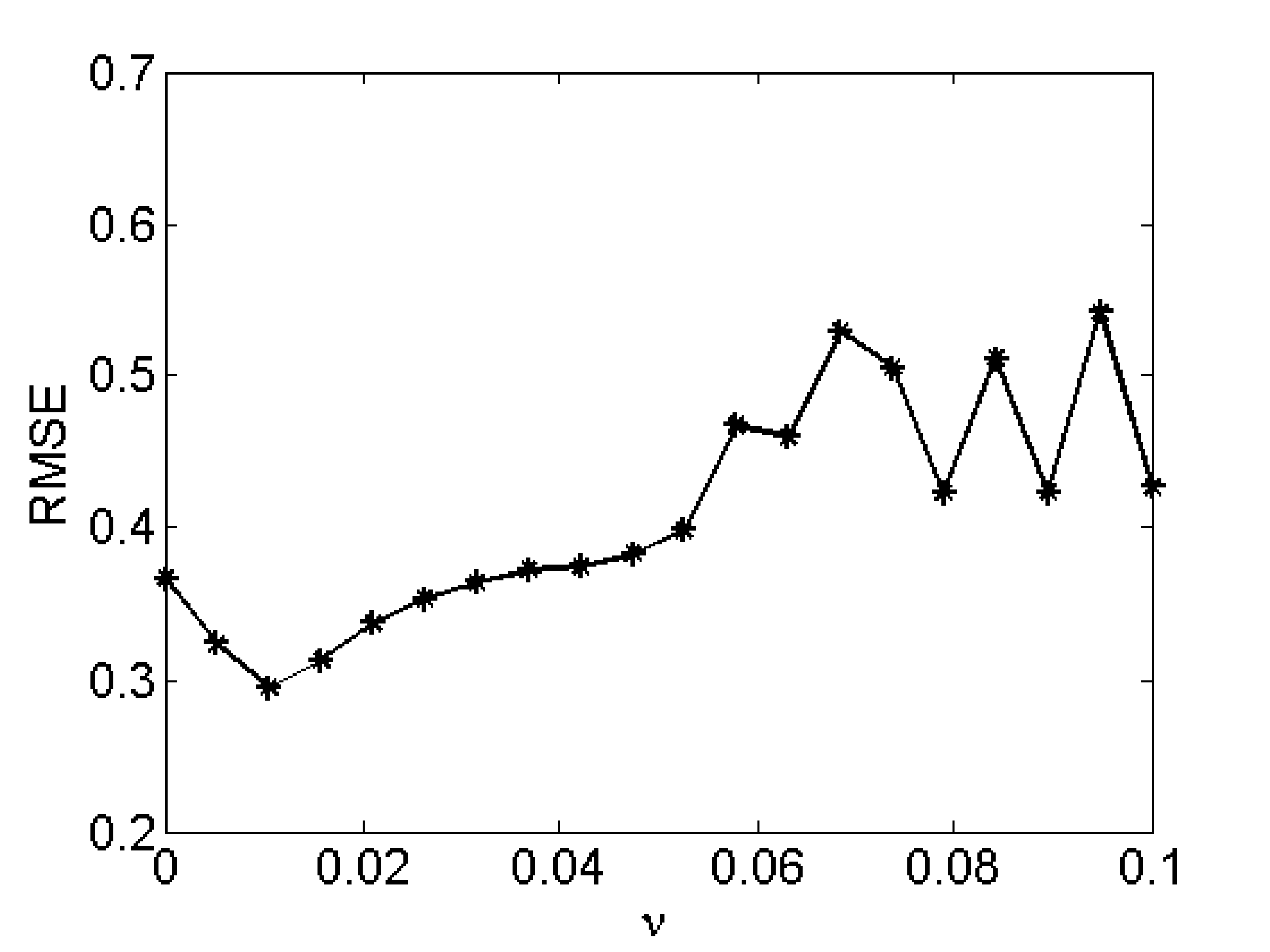}}
\caption{$RMSE$ of compactly supported reconstruction as a function of the parameter $\nu$.  $N=30$, $M=20,\ K=64,\ \varepsilon=1.1.$}
\label{fig: mse_nu_compact}
\end{figure}

\begin{figure}[htbp]
\centering%
\subfigure[Parameters and data as in Figure \ref{fig: mse_nu_compact} \label{fig: rcond_nu_cre}]%
{\includegraphics[width=0.49\textwidth]{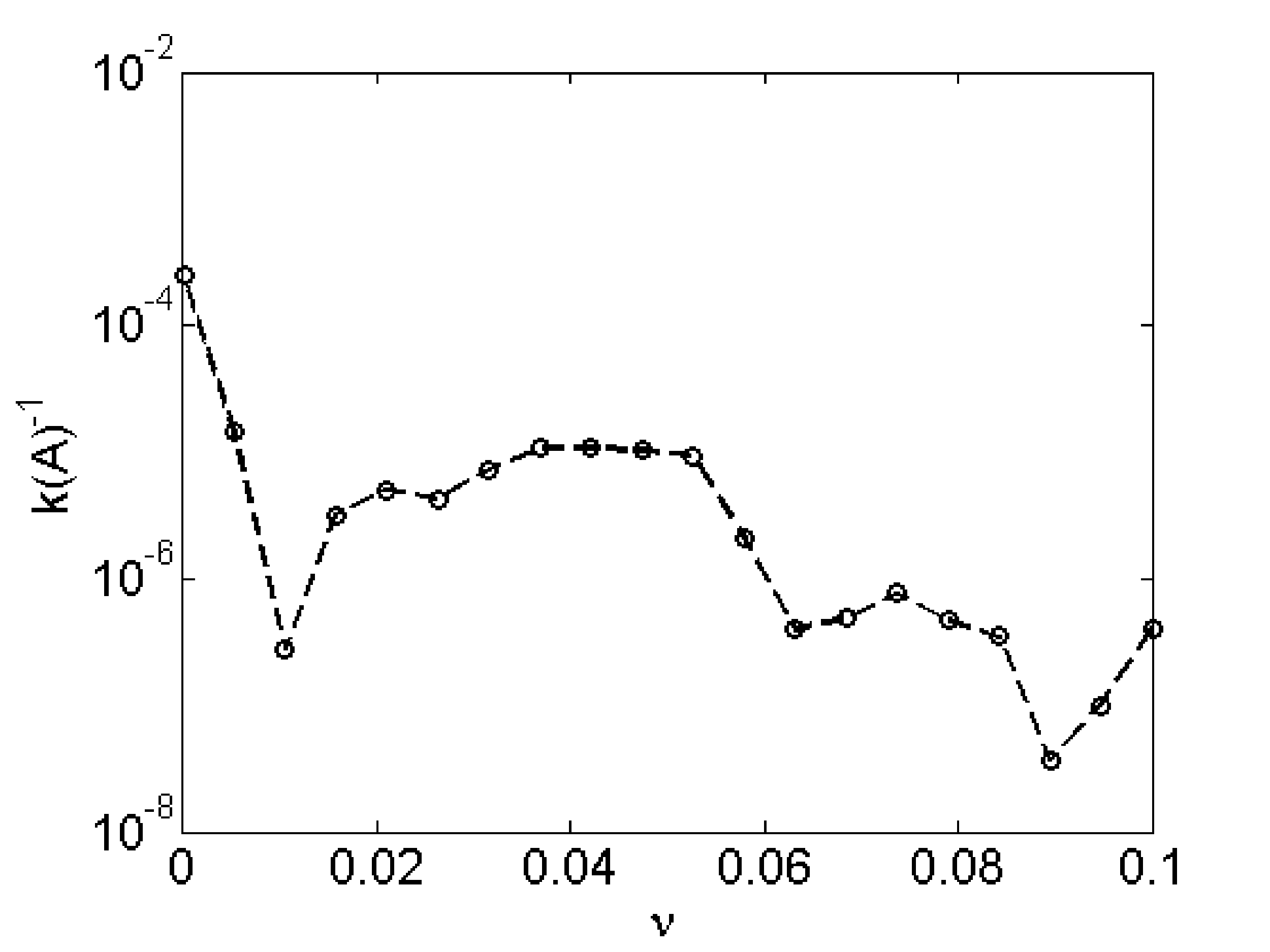}}%,height=0.3\textwidth
\subfigure[Parameters and data as in Figure \ref{fig: mse_nu_compact} \label{fig: rcond_nu_bull}]%
{\includegraphics[width=0.49\textwidth]{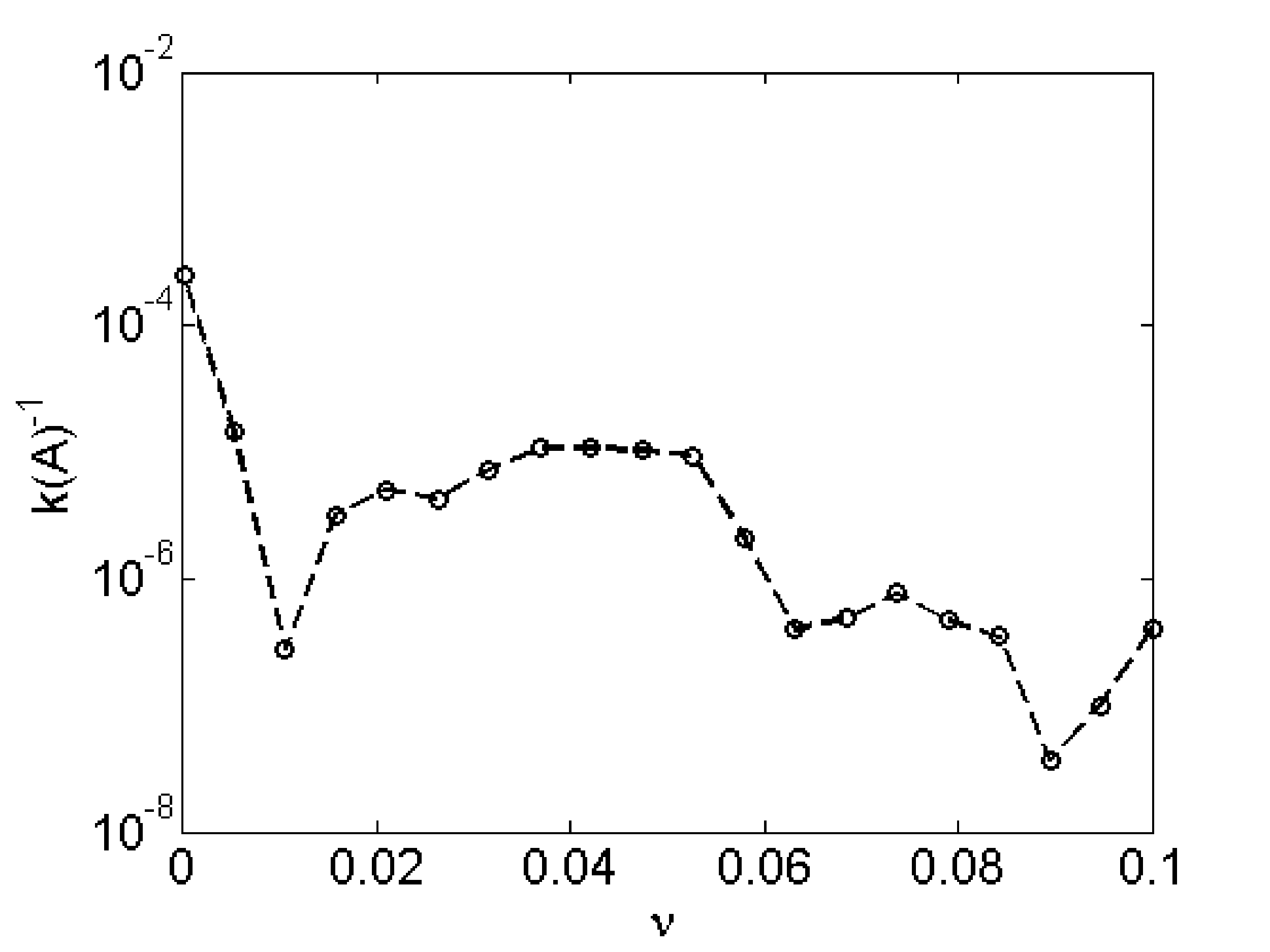}}
\caption{Reciprocal of the condition number of $A$ for compactly supported kernel reconstruction as a function of the parameter $\nu.$}
\label{fig: rcond_nu_compact}
\end{figure}

%INVIARE QUESTA ANALISI A DE MARCHI
At last we consider the case of compactly supported kernel. Let $K=(1-\varepsilon\norm{x-y})^{2}_{+}$ and $w(x)=(1-\nu^2\norm{x}^{2})_{+}$. Assume $\varepsilon\approx 1$, then it turns out that the $RMSE$ is minimal for $\nu\approx0$ (Figure \ref{fig: mse_nu_compact}). For small values of $\nu$ also the condition number of $A$ is larger (Figure \ref{fig: rcond_nu_compact}), so it is convenient to use $\nu\approx0$. 

We observe that using small values of $\nu$ corresponds to use a compactly supported window function $w$ with wide support, in this way one loses less information when filters basis $b_{j}$ with $w$. 

\begin{figure}[htbp]
\centering%
\subfigure[Crescent-shaped phantom \label{fig: mse_nu_compact_cre2}]%
{\includegraphics[width=0.49\textwidth]{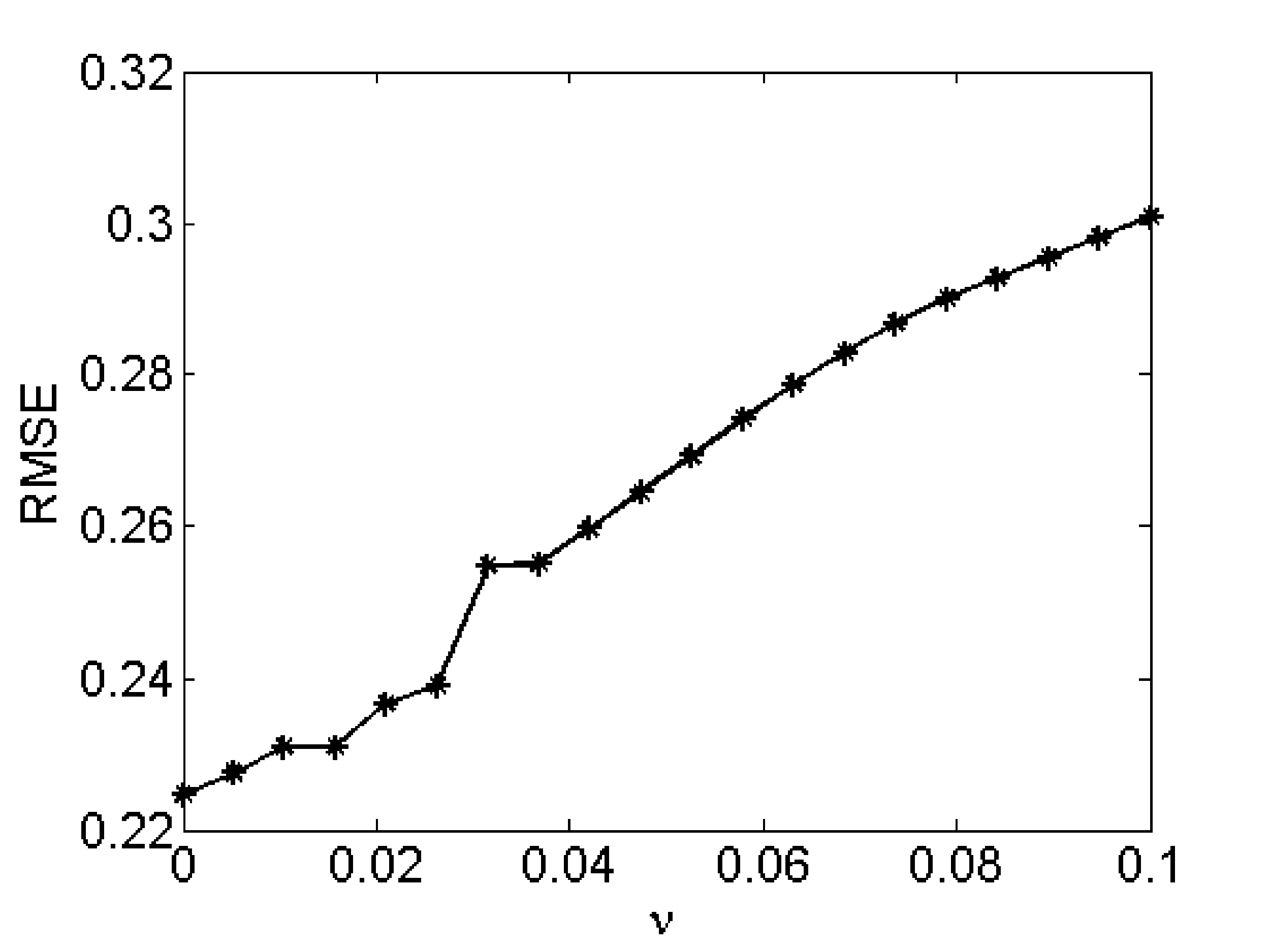}}%,height=0.3\textwidth
\subfigure[Bull's eye phantom \label{fig: mse_nu_compact_bull2}]%
{\includegraphics[width=0.49\textwidth]{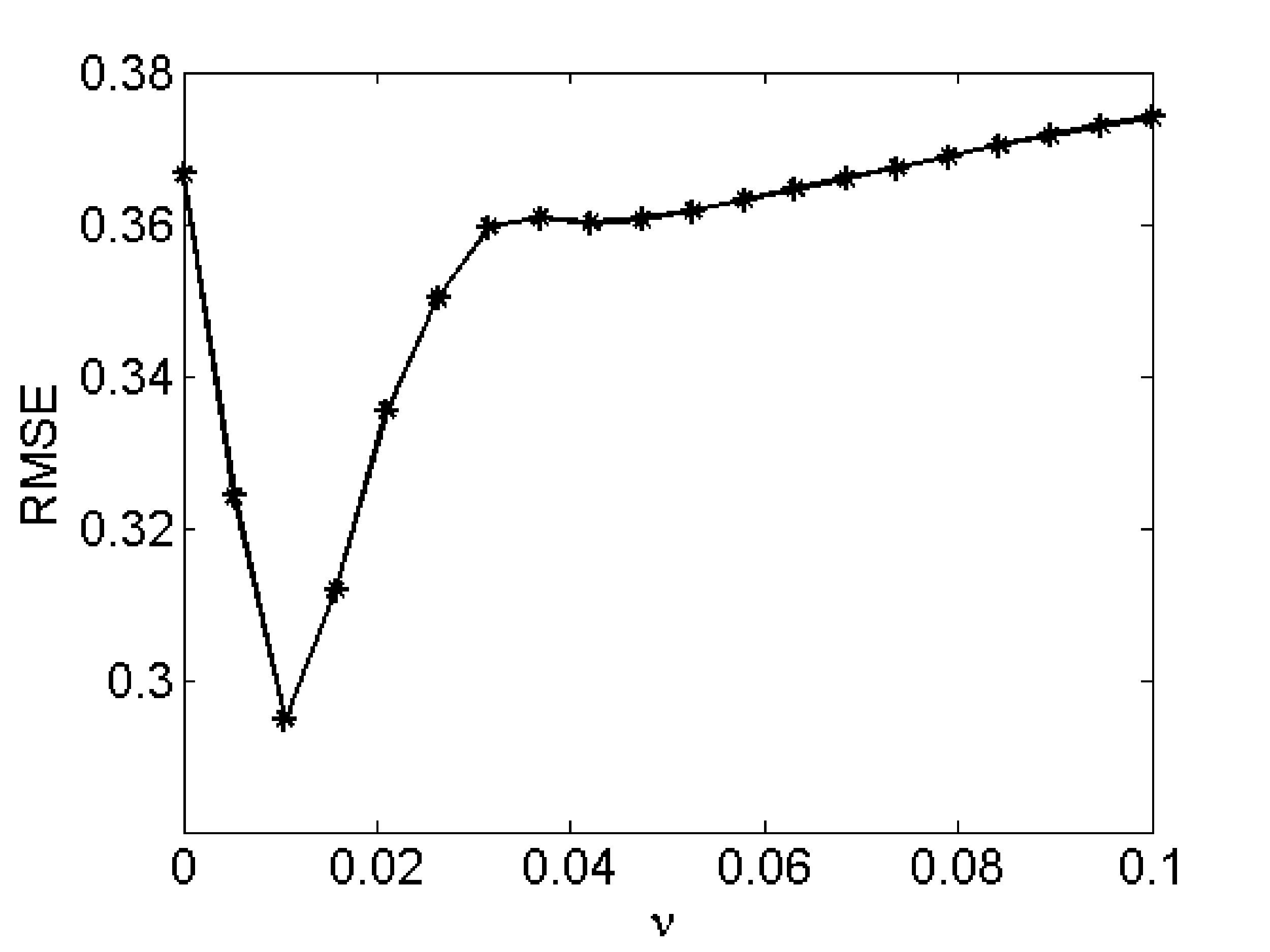}}
\caption{$RMSE$ of compactly supported reconstruction as a function of the parameter $\nu$. Regularization only for $a=0$. Parameters and data as in Figure \ref{fig: mse_nu_compact}.}
\label{fig: mse_nu_compact2}
\end{figure}

We can use the information that the approximation in optimal for $\nu\approx0$ to simplify the expression of the matrix $A$. Indeed, when $a\neq0$, if $\nu\rightarrow0$, then $a_{k,j}\rightarrow\frac{\pi}{6\varepsilon^2 a}$ and we can use this simpler expression of $a_{k,j}$ instead of \eqref{eq: matrix_cs} (see appendix B). This option is equivalent to consider the regularization $R_{w}$ only when $a=0$, while using the original operator $R$ when $a\neq0$ (that is what we called option 2 in the section \ref{subsec: reg_gauss}). 
Using this second option the behavior of $RMSE$ becomes more regular (see Figure \ref{fig: mse_nu_compact2}).

\begin{figure}[htbp]
\centering
\includegraphics[width=0.5\textwidth]{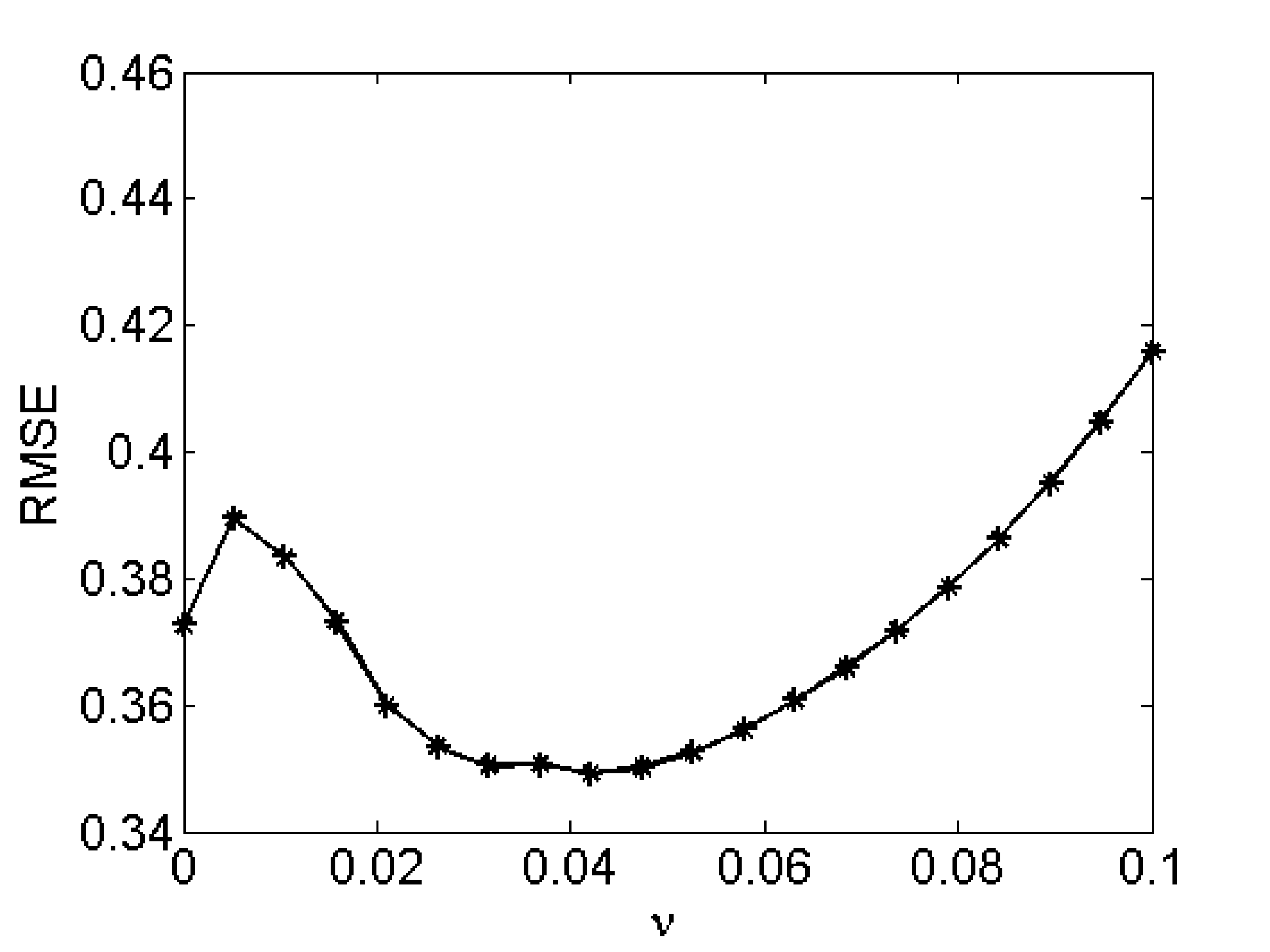}
\caption{$RMSE$ of compactly supported reconstruction as a function of the parameter $\nu$ for the Bull's eye phantom with $N=50$, $M=40,\ K=64,\ \varepsilon=1.1.$}
\label{fig: mse_nu_compact_data}
\end{figure}

Finally we observe that using more data, one should use a bigger value of $\nu$, as shown in Figure \ref{fig: mse_nu_compact_data}.

\subsection{Kernel shape parameters}
\begin{figure}[htbp]
\centering%
\subfigure[\label{fig: mse_ep_gauss1}]%
{\includegraphics[width=0.49\textwidth]{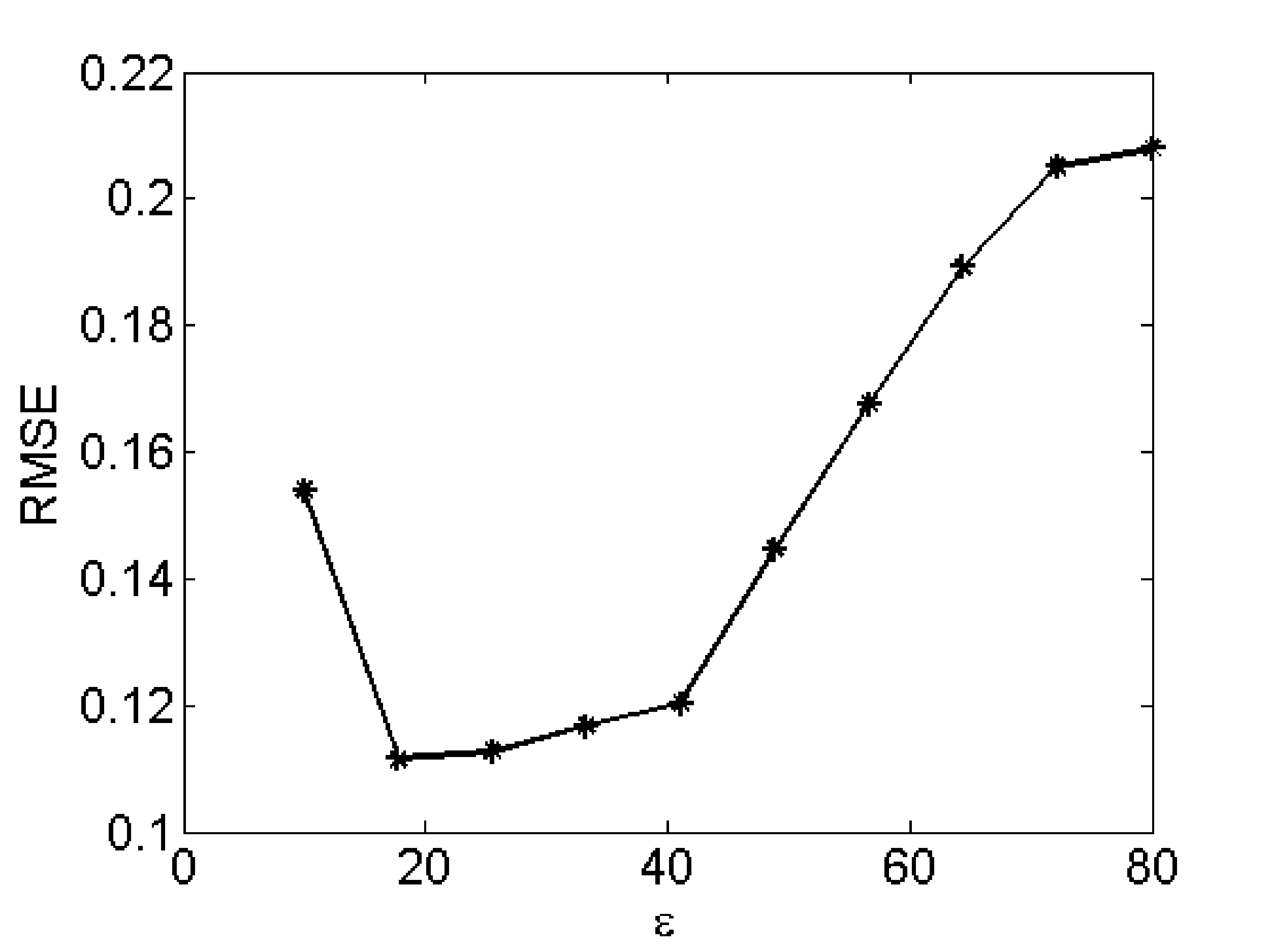}}%,height=0.3\textwidth
\subfigure[ \label{fig: mse_ep_gauss2}]%
{\includegraphics[width=0.49\textwidth]{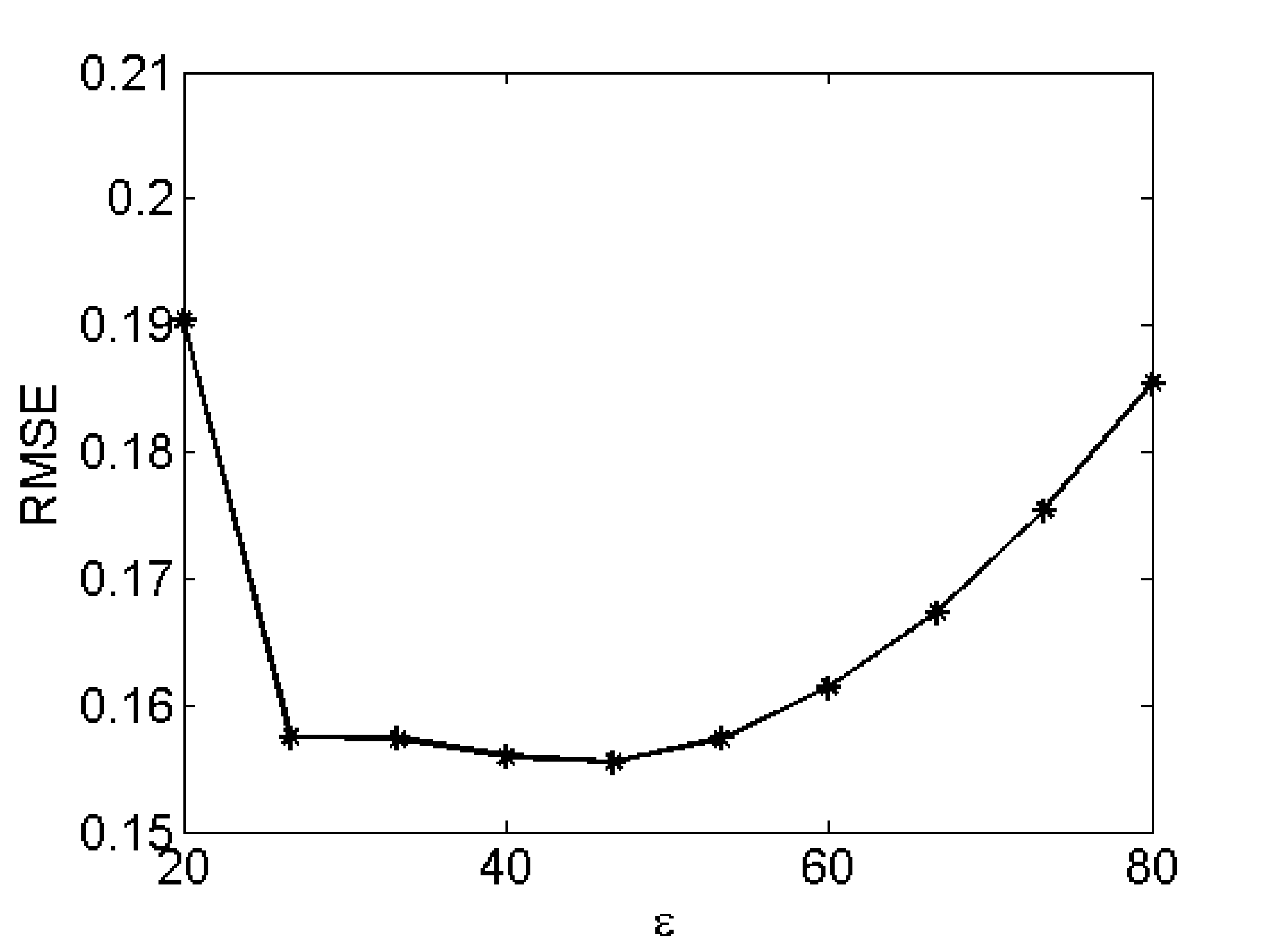}}
\subfigure[ \label{fig: mse_ep_imq1}]%
{\includegraphics[width=0.49\textwidth]{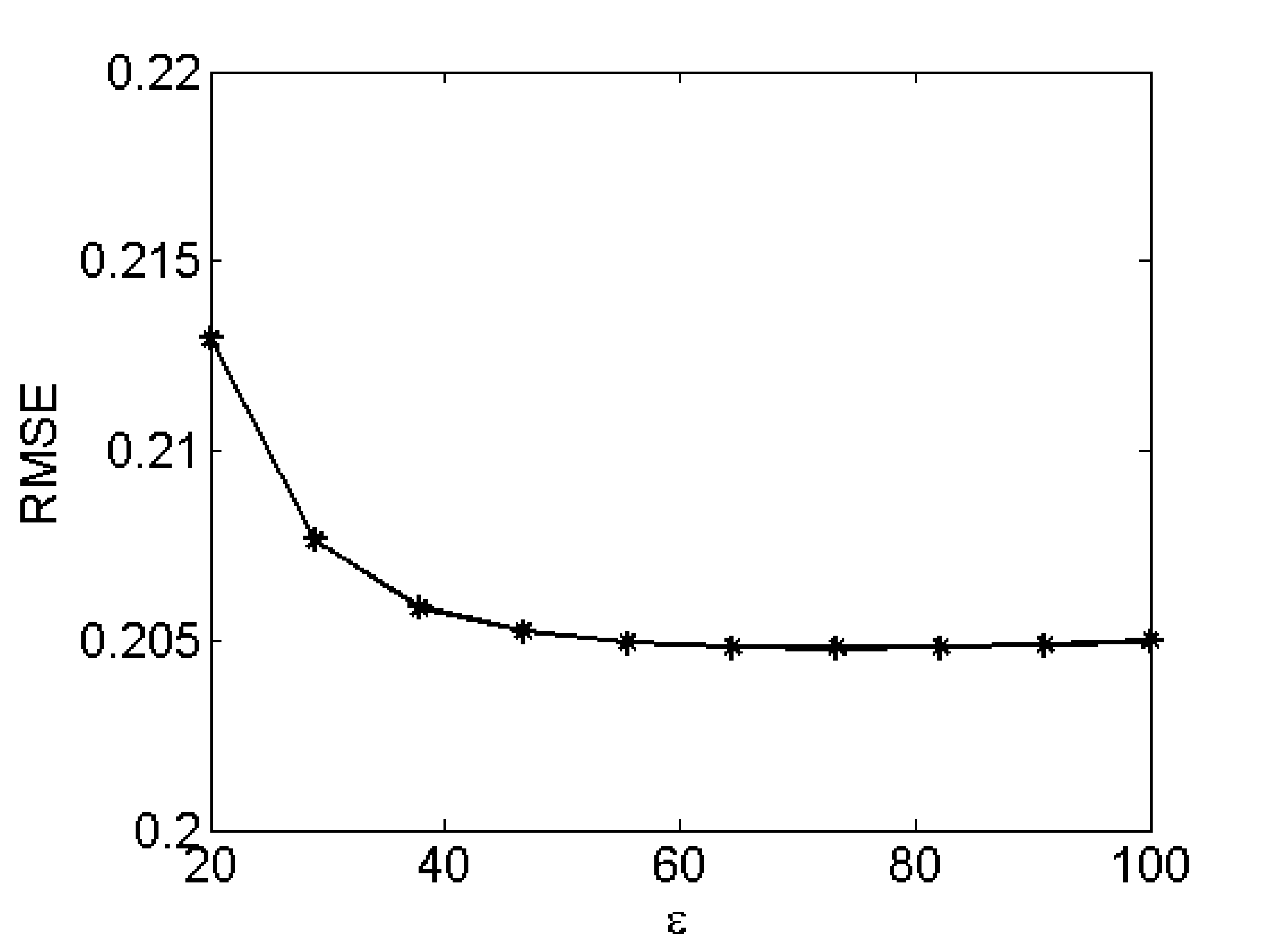}}%,height=0.3\textwidth
\subfigure[ \label{fig: mse_ep_imq2}]%
{\includegraphics[width=0.49\textwidth]{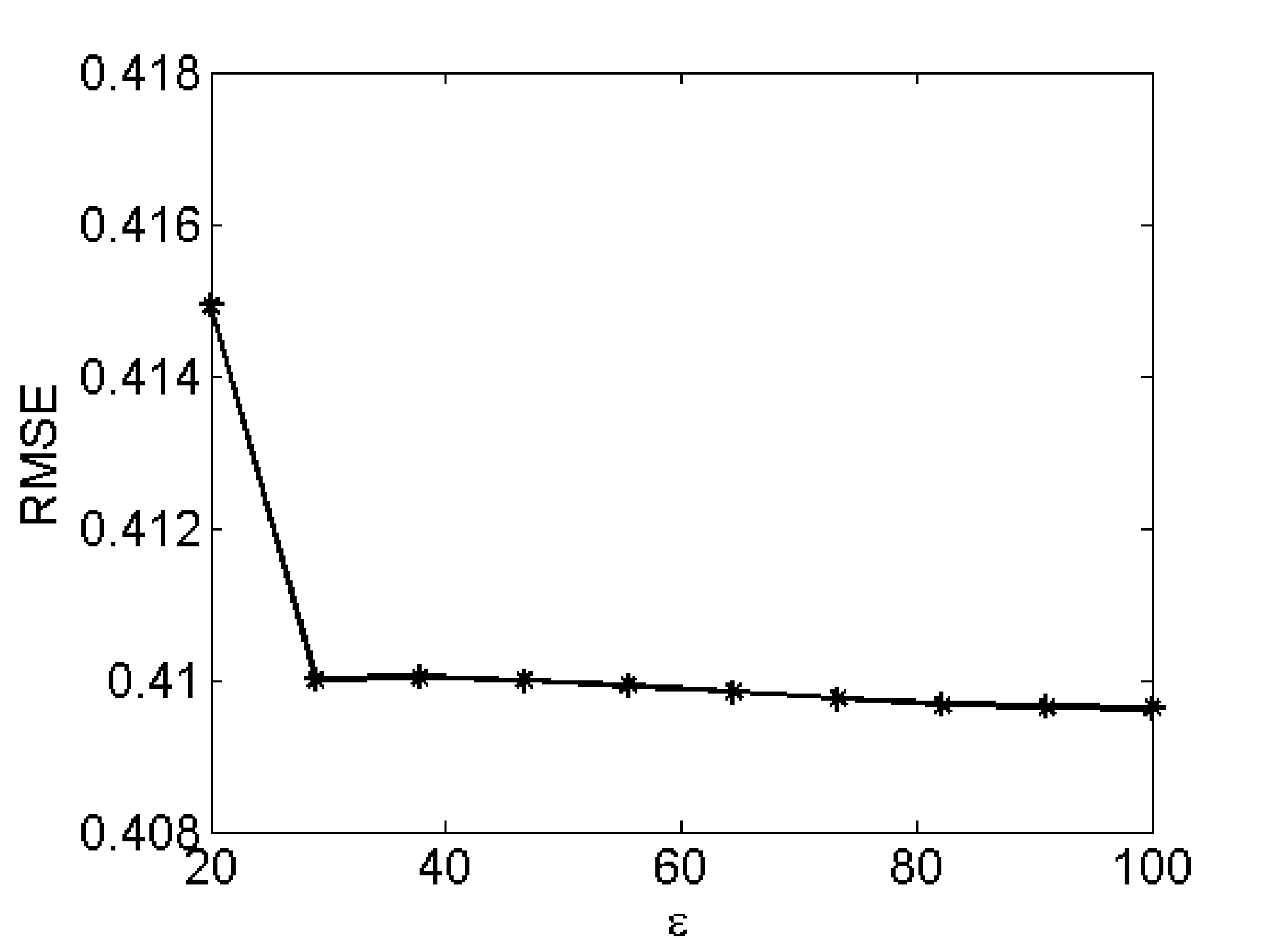}}
\caption{$RMSE$ of Gaussian and inverse multiquadric reconstruction as a function of the parameter $\varepsilon$. (a) Gaussian kernel, crescent-shaped phantom, $N=30,\ M=20,\ K=64,\ \nu=0.5$; (b) Gaussian kernel, bull's eye phantom, $N=50,\ M=40,\ K=64,\ \nu=0.7$; (c) Inverse multiquadric kernel, crescent-shaped phantom, $N=30,\ M=20,\ K=64,\ L_{1}=10,\ L_{2}=10$; (d) Inverse multiquadric kernel, Shepp-Logan phantom, $N=50,\ M=40,\ K=64,\  L_{1}=10,\ L_{2}=10$; }
\label{fig: mse_ep}
\end{figure}

We now consider the second main shape parameter involved in kernel methods, i.e. the kernel shape parameter $\varepsilon$. We will assume to work with an optimal window function parameter (as discussed in the previous paragraph). We observe that in the following cases the behavior of the $RMSE$ as a function of $\varepsilon$ is similar to that of the $RMSE$ as a function of $\nu$.
\begin{itemize}
\item \emph{Inverse multiquadrics}. There is an optimal value $\varepsilon_{opt}$ such that $RMSE$ is minimum, but for $\varepsilon>\varepsilon_{opt}$ the $RMSE$ increase slowly (see Figures \ref{fig: mse_ep_imq1} and \ref{fig: mse_ep_imq2});
\item \emph{Gaussian}. There is an optimal value $\varepsilon_{opt}$ such that $RMSE$ is minimum and for $\varepsilon>\varepsilon_{opt}$ the $RMSE$ increases fast.
Moreover $\varepsilon_{opt}$ increases if the the total number of data increases (see Figure \ref{fig: mse_ep_gauss1} and \ref{fig: mse_ep_gauss2}). %Empirical results show that, remaining in the range of data $70<n<140$, we can approximate $\varepsilon_{opt}\approx\frac{2}{5}n$
\end{itemize}

\begin{figure}[htbp]
\centering%
\subfigure[ \label{fig: rcond_ep_gauss1}]%
{\includegraphics[width=0.49\textwidth]{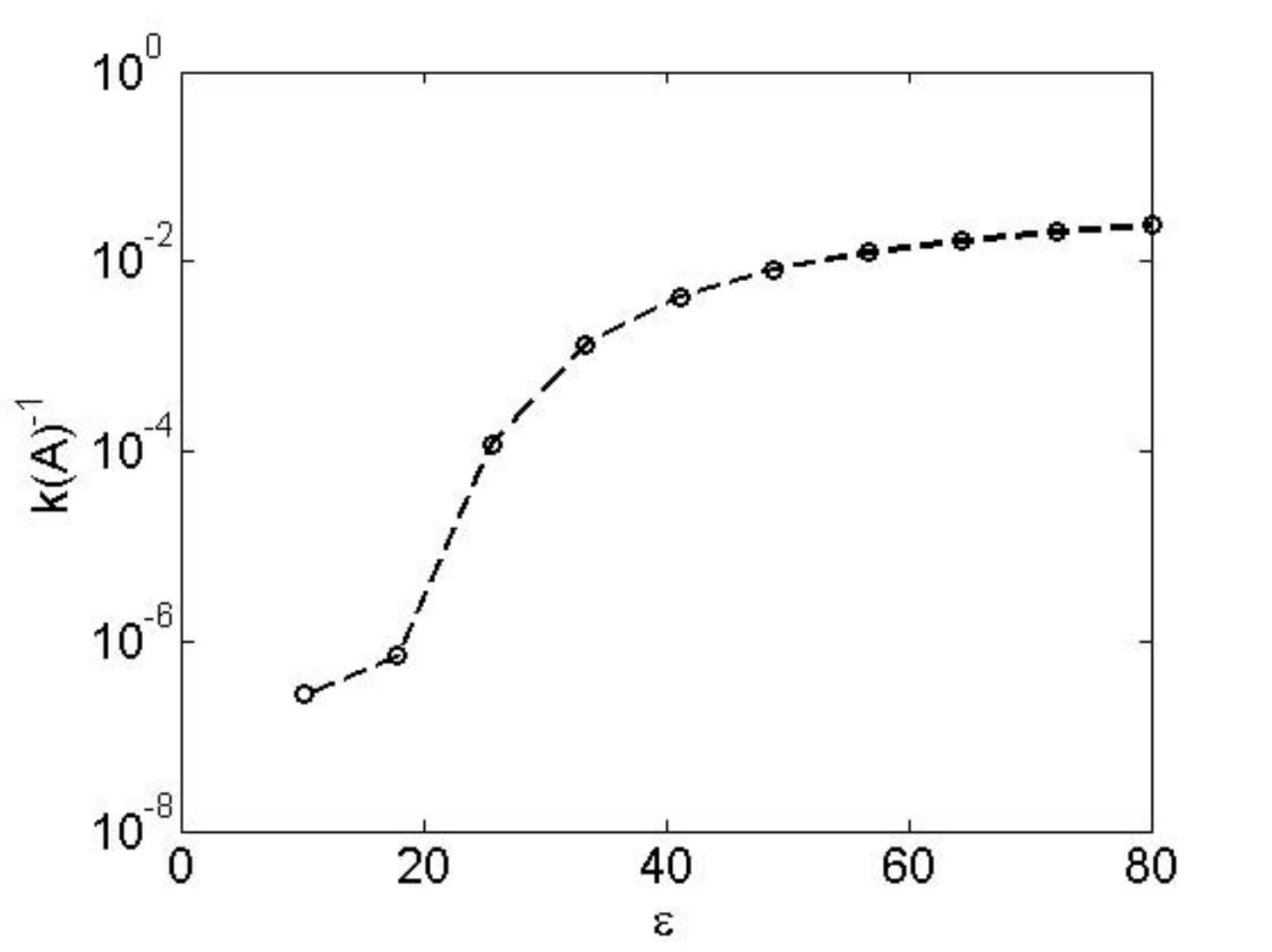}}%,height=0.3\textwidth
\subfigure[ \label{fig: rcond_ep_imq2}]%
{\includegraphics[width=0.49\textwidth]{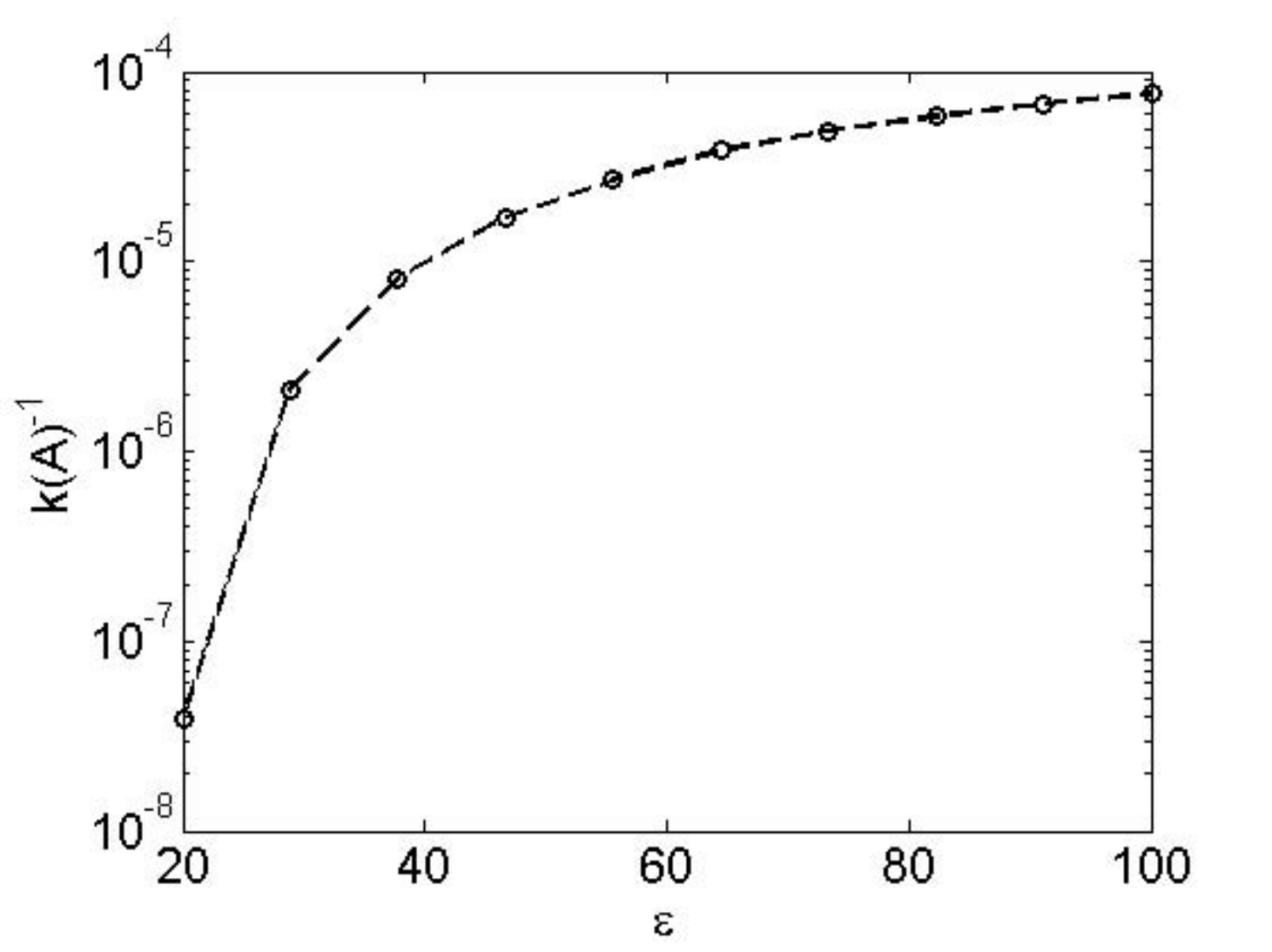}}
\caption{Reciprocal of the condition number of $A$ for Gaussian and inverse multiquadric kernel reconstruction as a function of the parameter $\varepsilon.$ (a) Parameters and data as in Figure \ref{fig: mse_ep_gauss1}; (b)Parameters and data as in Figure \ref{fig: mse_ep_imq2}}
\label{fig: rcond_ep}
\end{figure}

The difference with the window parameter is that now the condition number $k(A)$ is smaller for big value of $\varepsilon$ (equivalently the reciprocal $k(A)^{-1}$ increases with $\varepsilon$ - Figure \ref{fig: rcond_ep}).

\begin{figure}[htbp]
\centering%
\subfigure[$RMSE$ \label{fig: mse_L1_imq1}]%
{\includegraphics[width=0.49\textwidth]{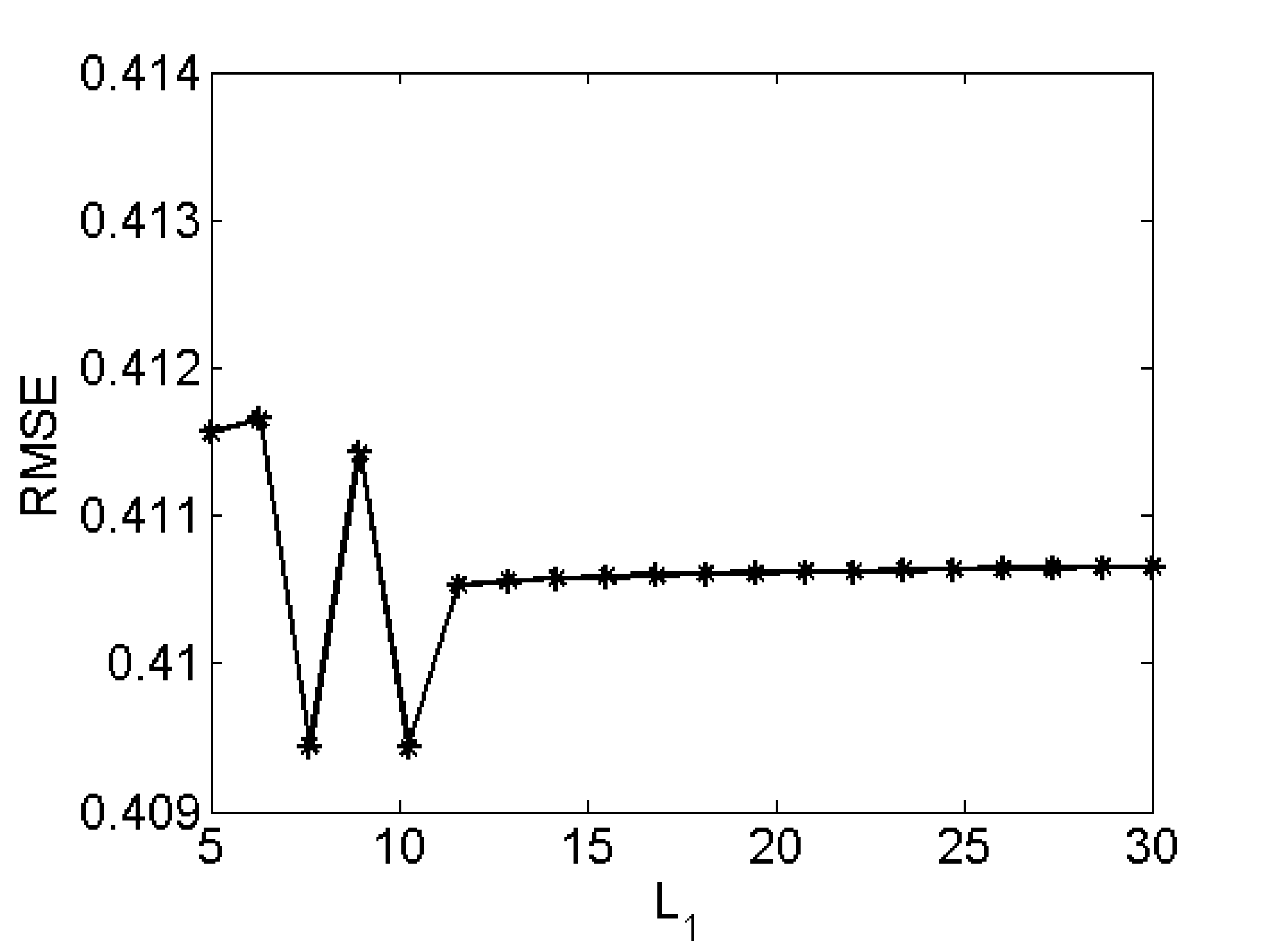}}%,height=0.3\textwidth
\subfigure[$k(A)^{-1}$ \label{fig: rcond_L1_imq1}]%
{\includegraphics[width=0.49\textwidth]{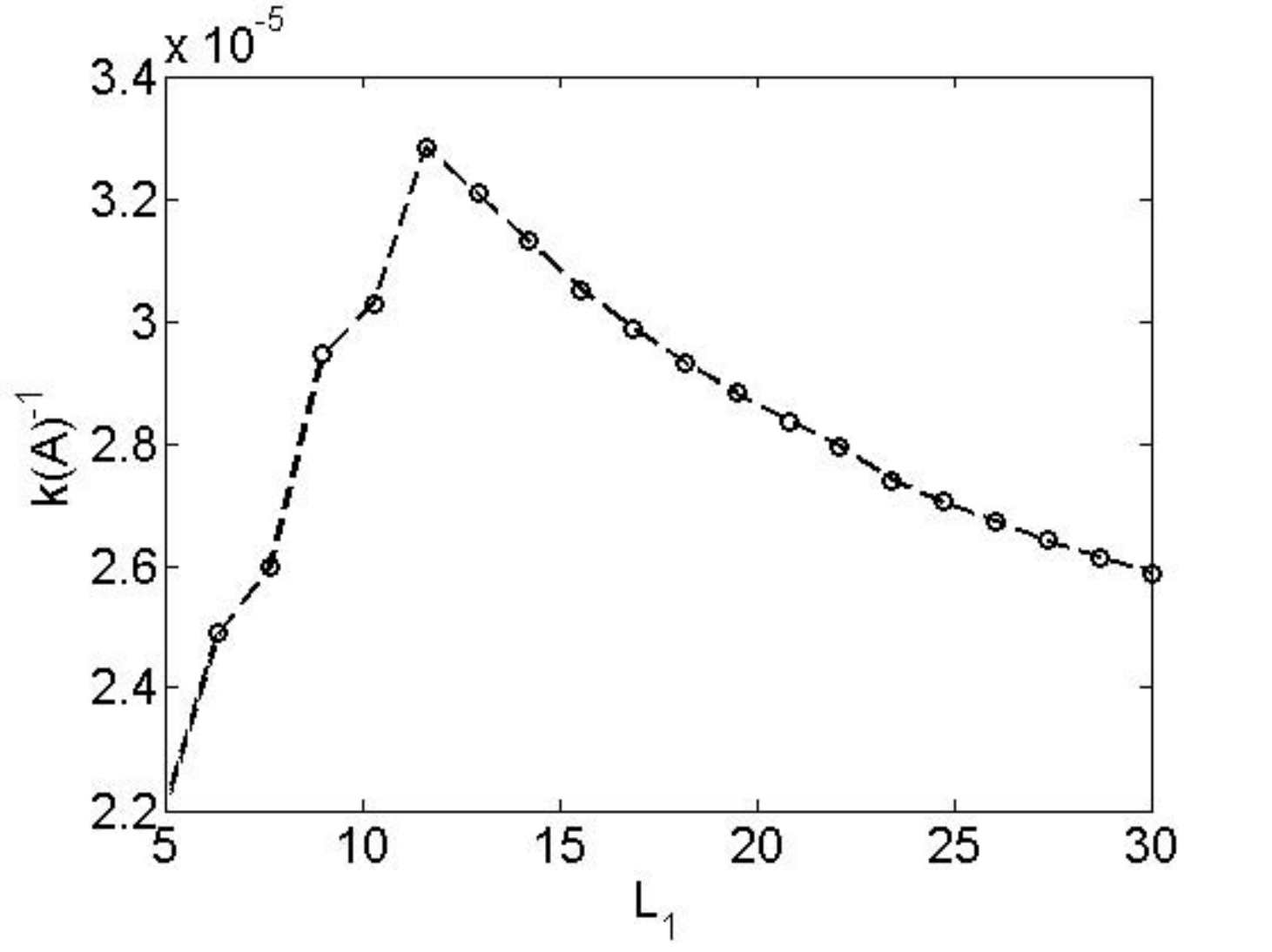}}
\caption{ $RMSE$ and $k^{-1}(A)$ for the inverse multiquadric reconstruction of the Shepp-Logan phantom as a function of $L_{1}$ with $N=36,\ M=25,\ K=60,\ \epsilon=27,\ L_{2}=10.$}
\label{fig: L1_imq}
\end{figure}

Another parameter involved in the inverse multiquadric reconstruction is $L_{1}$. We observe that under a certain threshold the behavior of the error as a function of $L_{1}$ is chaotic, but over this threshold the $RMSE$  remains almost constant. The reciprocal of the condition number is instead decreasing (Figure \ref{fig: L1_imq}).

Consider now the Gaussian-multiquadric kernel 
\begin{equation*}
K(x,y)=\sqrt{1+\rho^{2}\norm{x-y}^{2}}e^{-\varepsilon^{2}\norm{x-y}^{2}}
\end{equation*}. 
In this case we have two different shape parameters $\rho$ and $\varepsilon$. As in the Gaussian case, large values of $\varepsilon$ generate matrix with a better condition number ($\rho$ fixed). For fixed $\varepsilon\gg\rho$ we have $K(x,y)\approx e^{-\varepsilon^{2}\norm{x-y}^{2}}$ a kernel similar to the Gaussian case, but values $\rho\approx 0$ give a bigger conditioning.
\begin{figure}[htbp]
\centering%
\subfigure[$RMSE$ \label{fig: mse_rho_mq2}]%
{\includegraphics[width=0.49\textwidth]{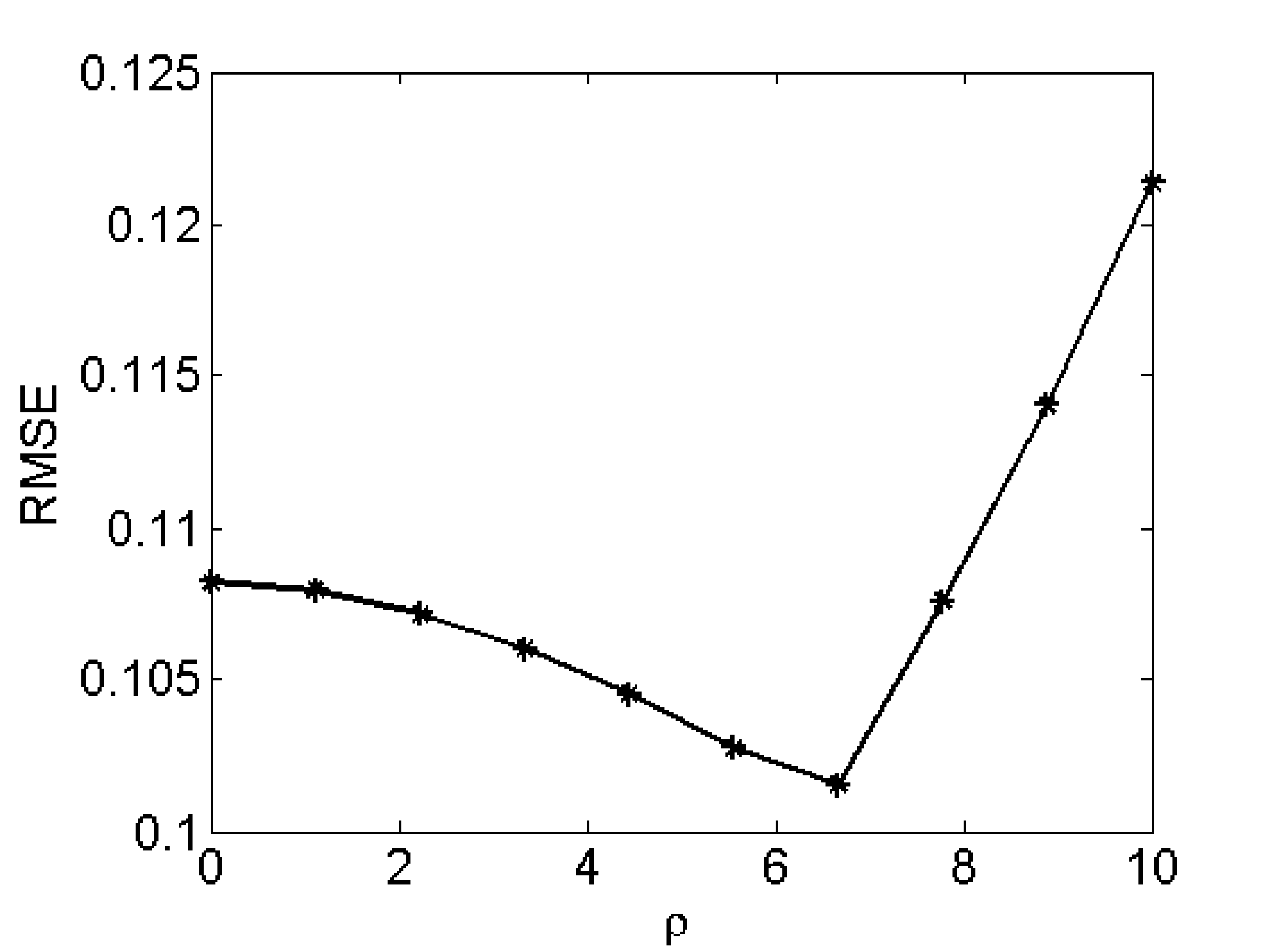}}%,height=0.3\textwidth
\subfigure[$k(A)^{-1}$ \label{fig: rcond_rho_mq2}]%
{\includegraphics[width=0.49\textwidth]{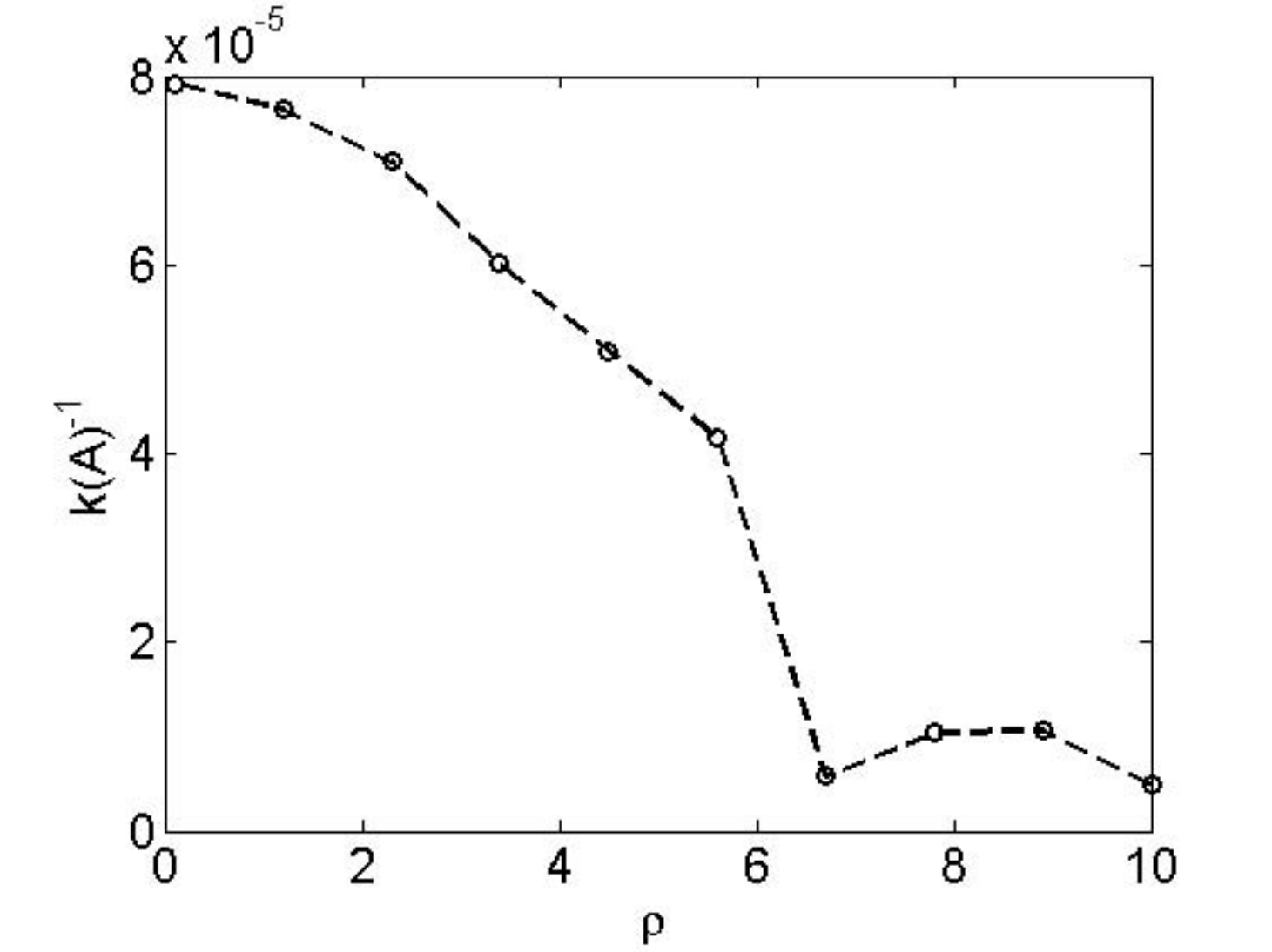}}
\caption{$RMSE$ and $k^{-1}(A)$ for multiquadric reconstruction of the crescent-shaped phantom as a function of $\rho$ with $N=50,\ M=40,\ K=64,\ \epsilon=52,\ \nu=0.5.$}
\label{fig: rho_mq}
\end{figure}

\begin{figure}[htbp]
\centering%
\subfigure[$N=20,\ M=15,\ \varepsilon=20.4$ \label{fig: mse_rho_sl1}]%
{\includegraphics[width=0.49\textwidth]{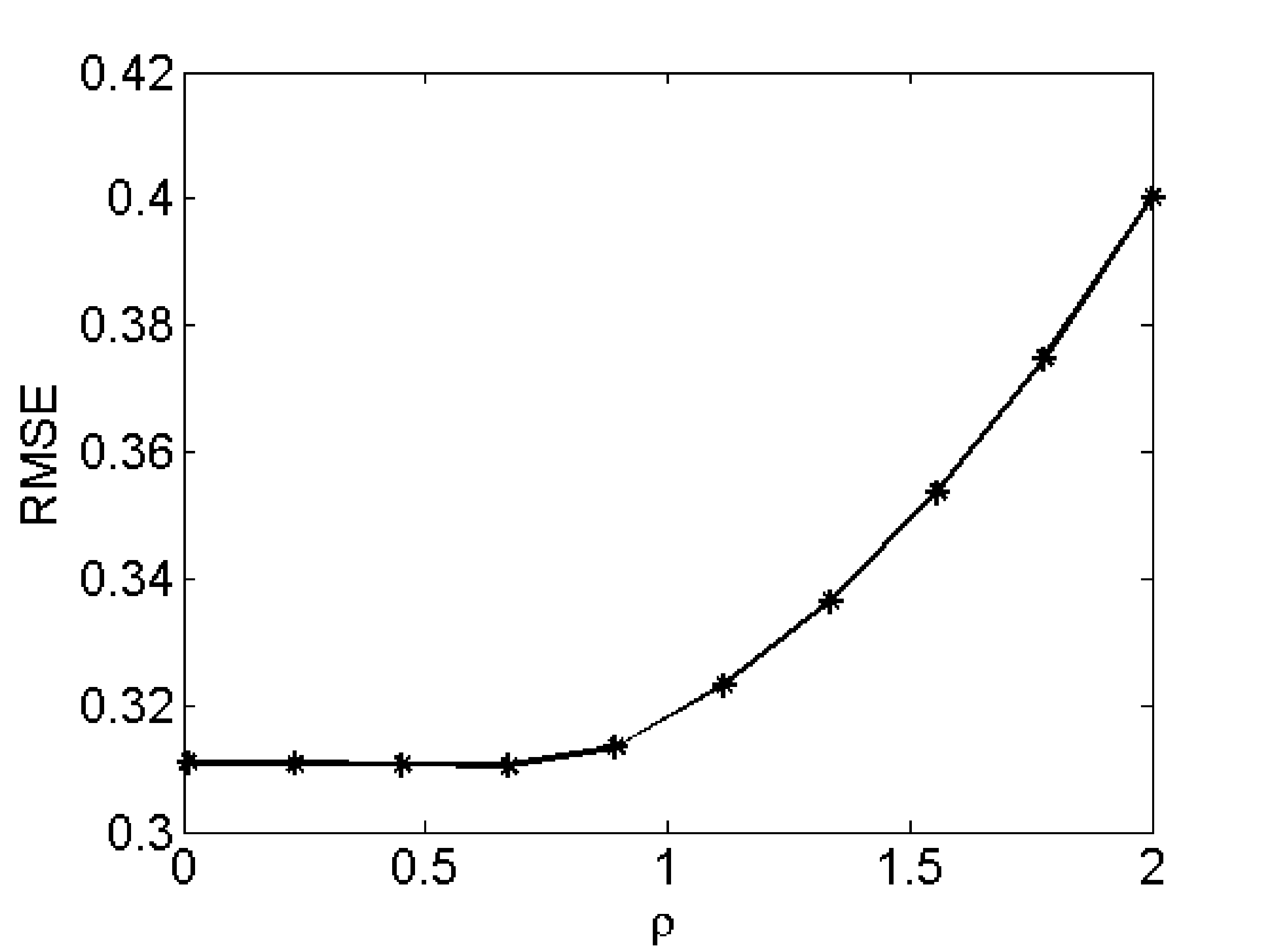}}%,height=0.3\textwidth
\subfigure[$N=50,\ M=40,\ \varepsilon=52.4$ \label{fig: mse_rho_sl2}]%
{\includegraphics[width=0.49\textwidth]{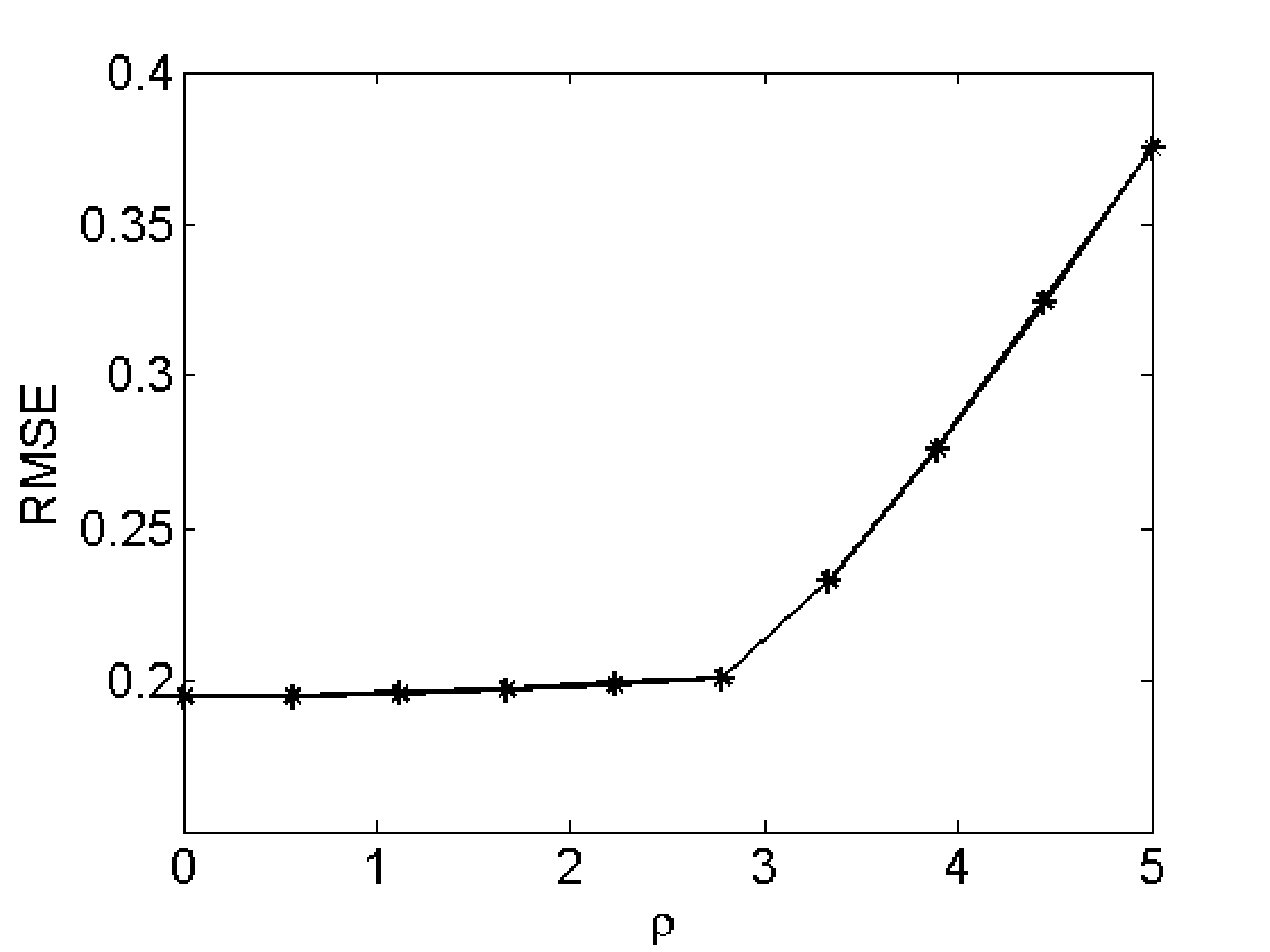}}
\caption{$RMSE$ for multiquadric reconstruction of the Shepp-Logan in function of $\rho$ with $K=256,\ \nu=1.3.$}
\label{fig: rho_mq_shepp}
\end{figure}

If $\rho$ is of the same order of $\varepsilon$ then $k(A)$ is small but $MSE$ becomes larger.
%%%mettere ordini di grandezza dell'errore 
 So we have again a situation with an optimal value for $\rho$ and $\varepsilon$ that arises from the trade-off to have a well conditioned matrix and a good approximation of the non-regularized reconstruction problem (trade-off principle \cite{SCHAB}). 
Moreover, as well as $\varepsilon$ also 	the optimal value of $\rho$ depends on the number of data and on the phantom. For example, using the crescent-shaped phantom, the value of $\rho_{opt}$ varies from $\approx3$ (for $N=20,\ M=15$) to  $\approx7$ (for $N=50,\ M=40$); while considering the Shepp-Logan filter we obtain  $\rho_{opt}<0.5$  (for $N=20,\ M=15$) and $\rho_{opt}<2$  (for $N=50,\ M=40$) (Figure \ref{fig: rho_mq_shepp}).

\begin{figure}[htbp]
\centering%
\subfigure[Crescent-shaped phantom \label{fig: mse_ep_compact_c}]%
{\includegraphics[width=0.33\textwidth]{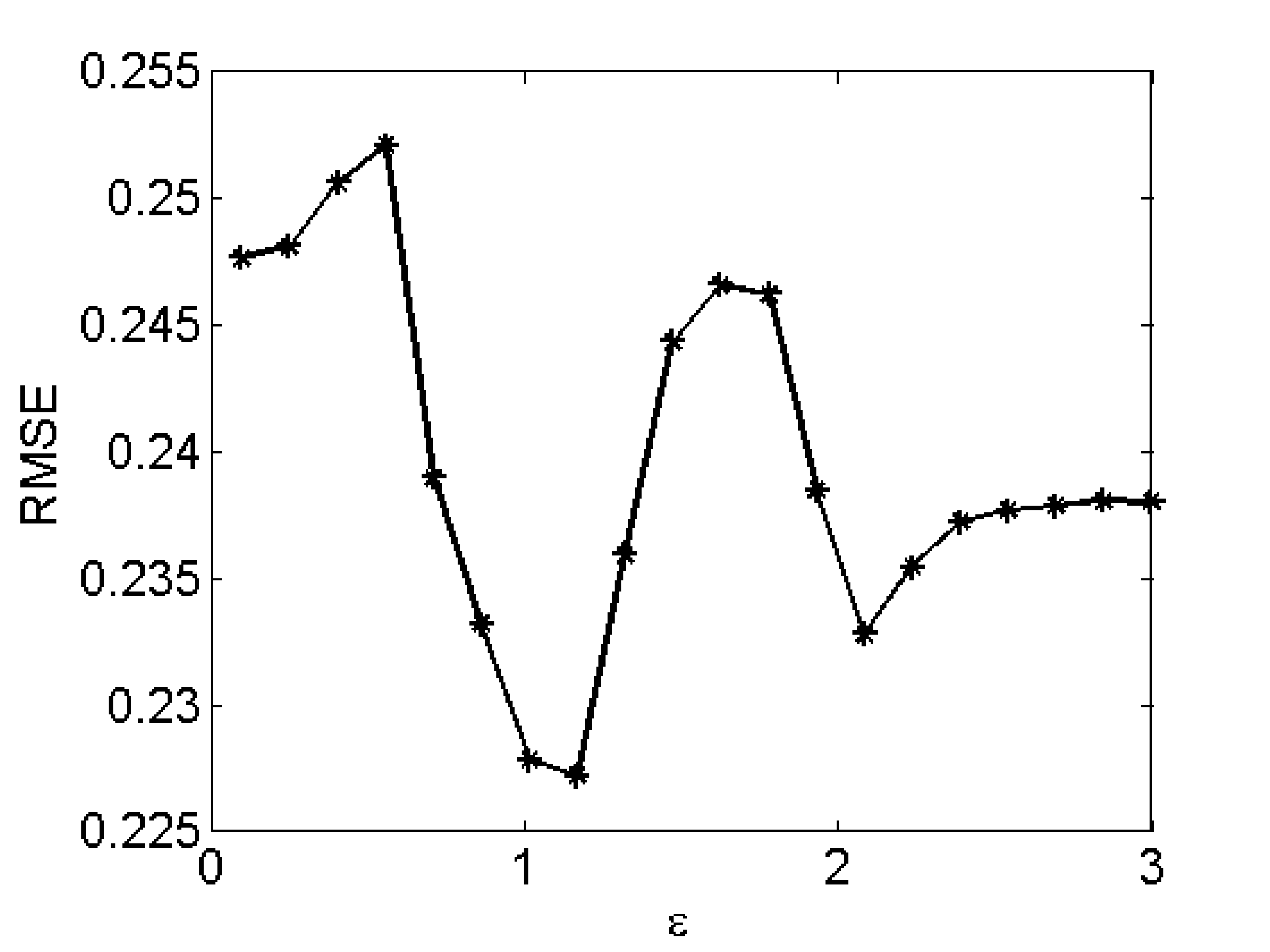}}%,height=0.3\textwidth
\subfigure[Bull's eye phantom \label{fig: mse_ep_compact_b}]%
{\includegraphics[width=0.33\textwidth]{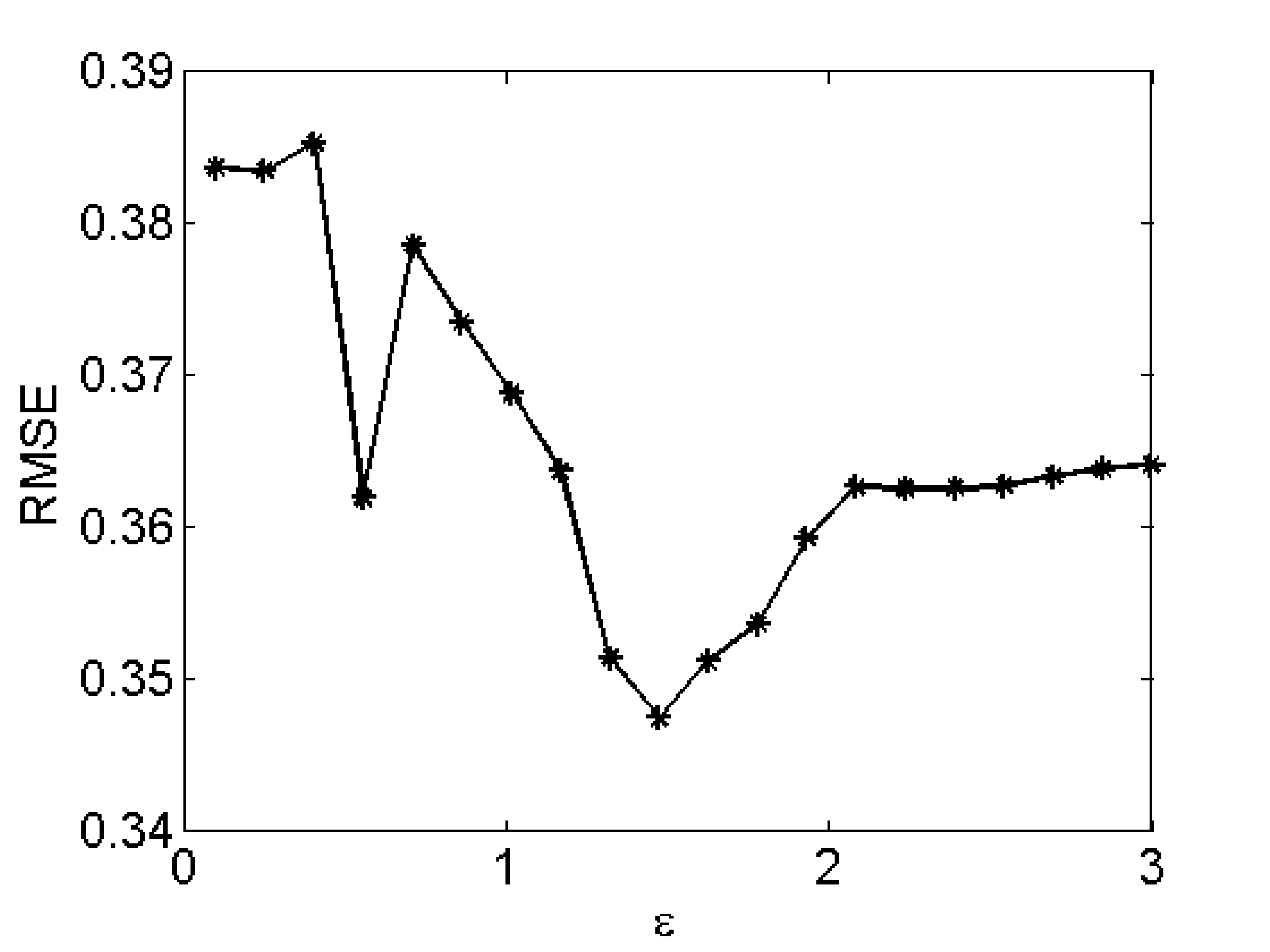}}
\subfigure[Shepp-Logan phantom \label{fig: mse_ep_compact_sl}]%
{\includegraphics[width=0.33\textwidth]{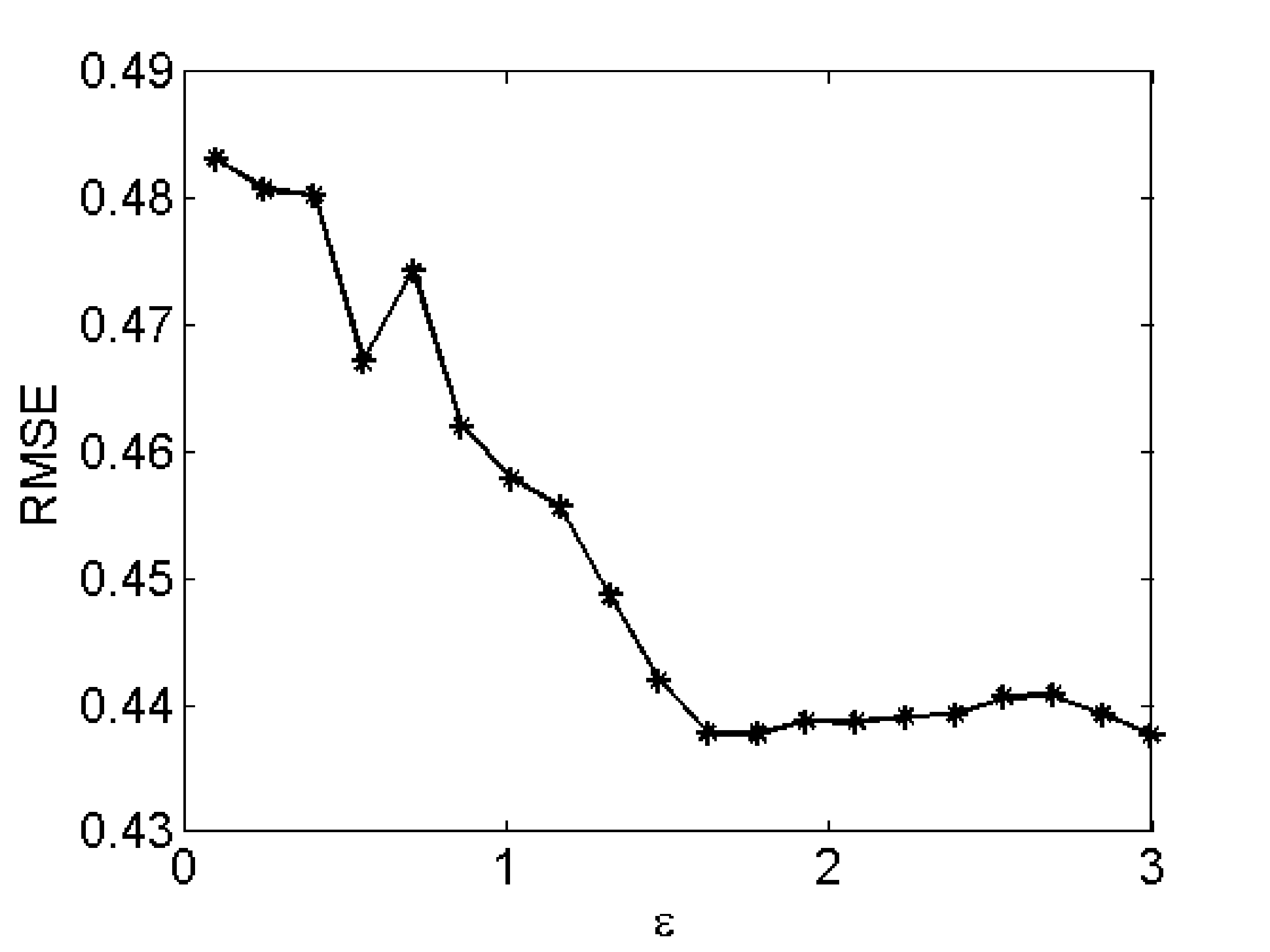}}
\caption{$RMSE$ of compactly supported reconstruction as a function of the parameter $\varepsilon$. $N=30$, $M=20,\ K=64,\ \nu=10^{-6 }.$}
\label{fig: mse_ep_compact}
\end{figure}

\begin{figure}[htbp]
\centering
\includegraphics[width=0.5\textwidth]{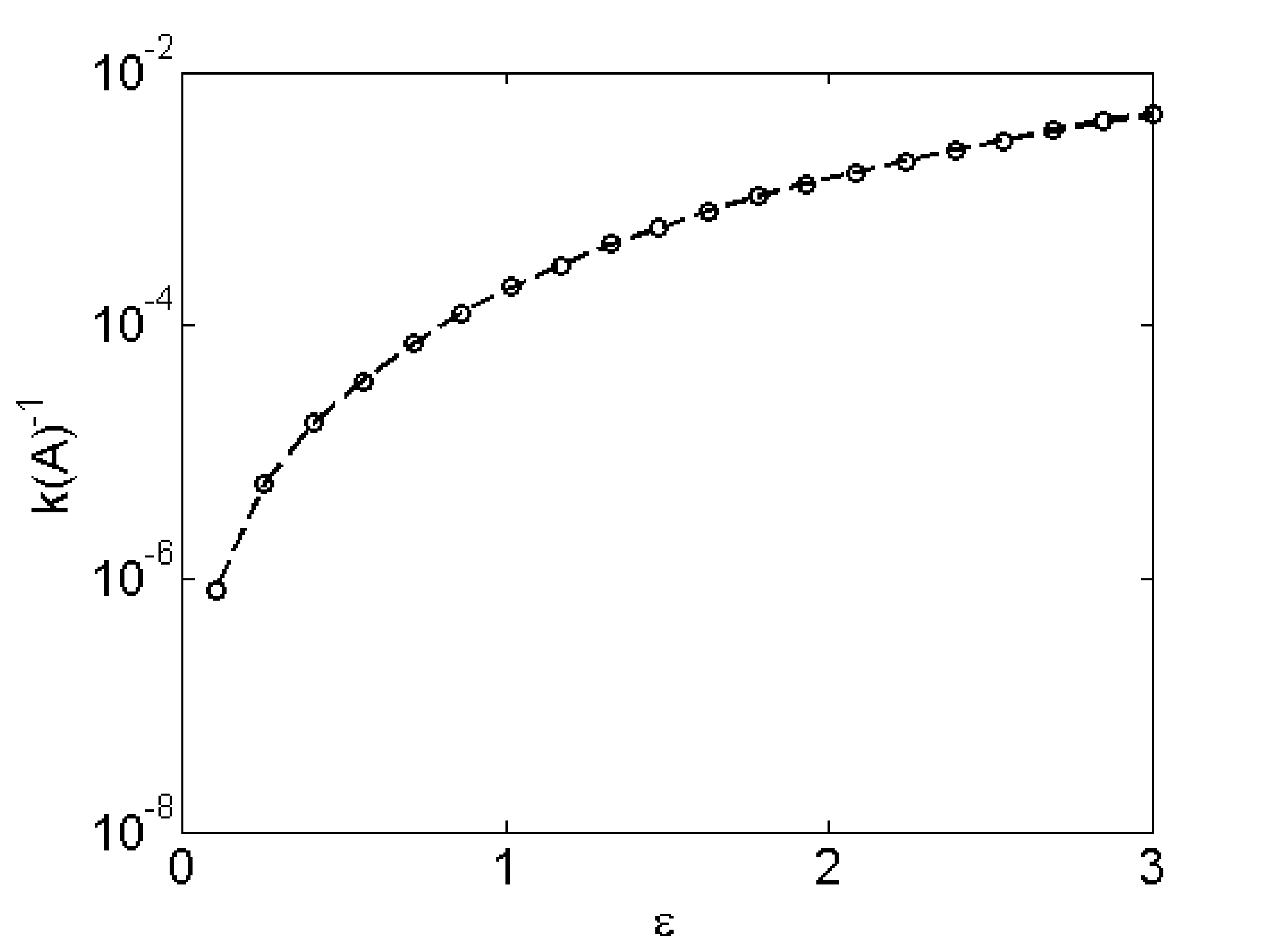}
\caption{Reciprocal of the condition number of the matrix $A$ of compactly supported reconstruction as a function of the parameter $\varepsilon$ with $N=30$, $M=20,\ K=64,\ \varepsilon=1.1.$}
\label{fig: rcond_ep_compact}
\end{figure}

At last we consider again the compactly supported kernel. Varying $\varepsilon$ we see that the $RMSE$ presents a minimum for an optimal value $\varepsilon_{opt}$ that is between 1 and 2 for both the crescent-shaped and the bull's eye phantoms (Figures \ref{fig: mse_ep_compact_c} and \ref{fig: mse_ep_compact_b}). Considering instead the Shepp-Logan phantom the $RMSE$ decreases if $\varepsilon$ increases (Figure \ref{fig: mse_ep_compact_sl}). However, in all the cases, the $RMSE$ remains almost constant for $\varepsilon$ large enough, while the reciprocal of the condition number increases with $\varepsilon$ (Figure \ref{fig: rcond_ep_compact}).

\subsection{Scale parameter}
\begin{figure}[htbp]
\centering%
\subfigure[$RMSE$ \label{fig: mse_h_gauss}]%
{\includegraphics[width=0.49\textwidth]{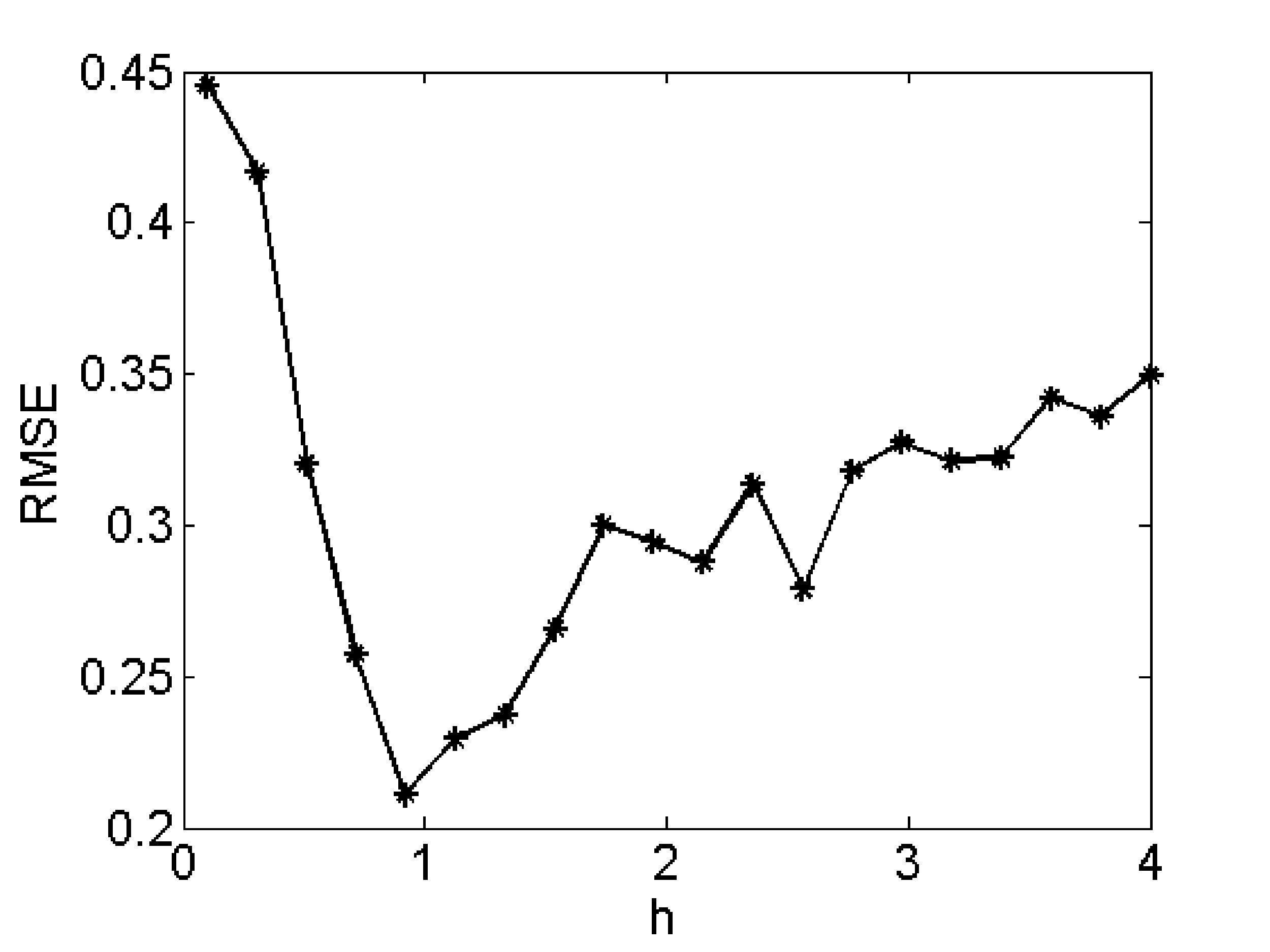}}%,height=0.3\textwidth
\subfigure[$k^{-1}(A)$ \label{fig: rcond_h_gauss}]%
{\includegraphics[width=0.49\textwidth]{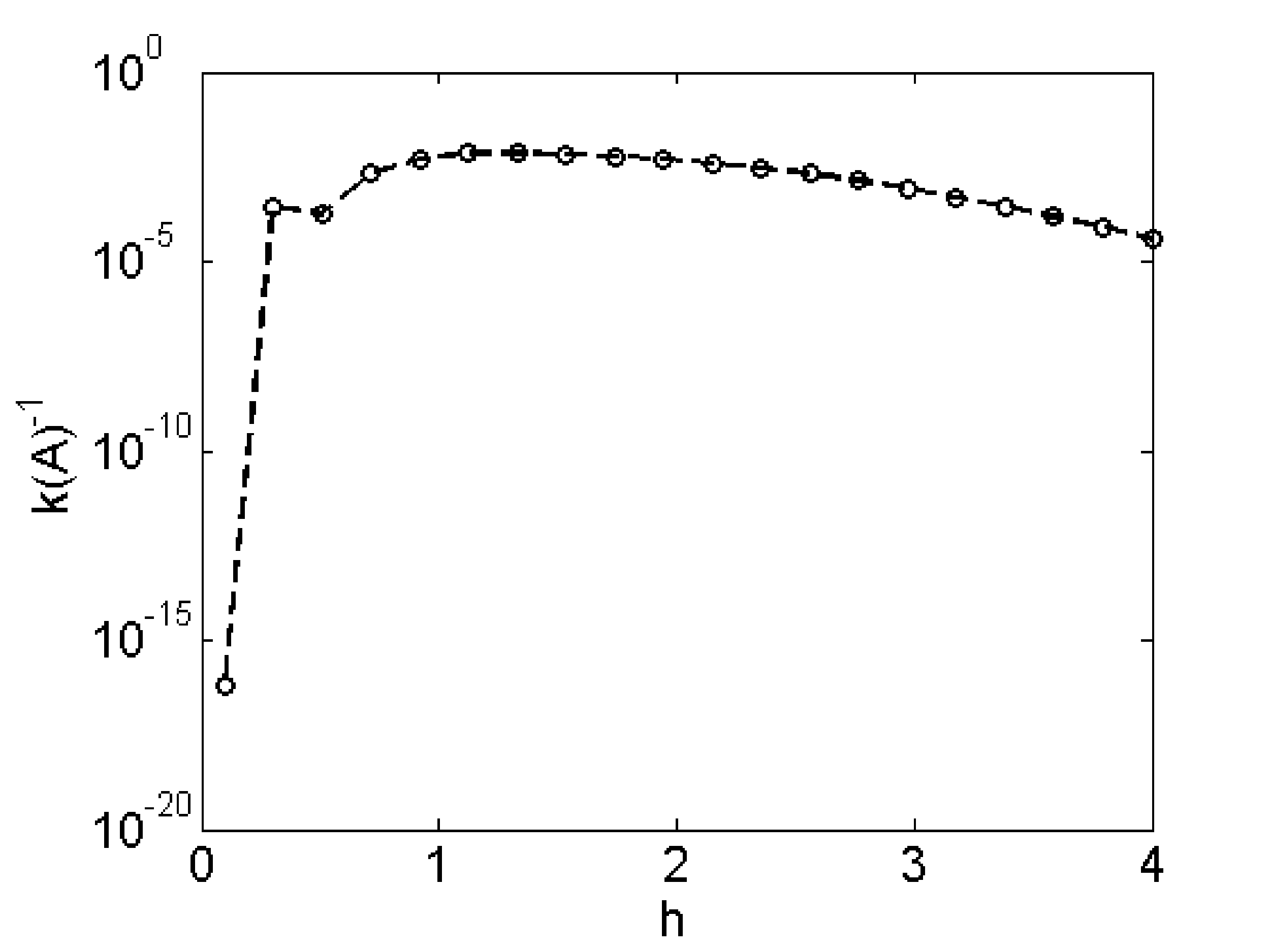}}
\caption{Gaussian kernel reconstruction of the crescent-shaped phantom in function of $h$. Here $\varepsilon=50,\ \nu=0.5, \ N=30,\ M=20,\ K=64.$}
\label{fig: h_gauss}
\end{figure}

Considering the scaled reconstruction problem of section \ref{sec: scaled_problem}, one has also to consider the behavior of the solution depending on the scale parameter $h$. Numerical experiments show that there is an optimal value of $h$ such that the $RMSE$ is minimum and $k^{-1}(A)$ is maximum (Figure \ref{fig: mse_h_gauss} and \ref{fig: rcond_h_gauss}). It turns out that, in most of the cases we examined, the optimal value of $h$ is $h\approx1$. Thus, in our discussion we will consider $h=1$.

\section{Comparison of the methods}
In this section we compare the classical Fourier methods, introduced in chapter \ref{chap: fourier_methods}, with the kernel-based methods of chapter \ref{chap: kernelRec}. In this second case we assume the use of optimal shape parameters.

We compare the solutions of different algorithms varying the phantom and also testing their behavior when introducing some noise in the data. Again we use the $RMSE$ to measure how much the solutions differ from the original phantom.

\begin{figure}[htbp]
\centering%
\subfigure[ \label{fig: fou_rec1}]%
{\includegraphics[width=0.49\textwidth]{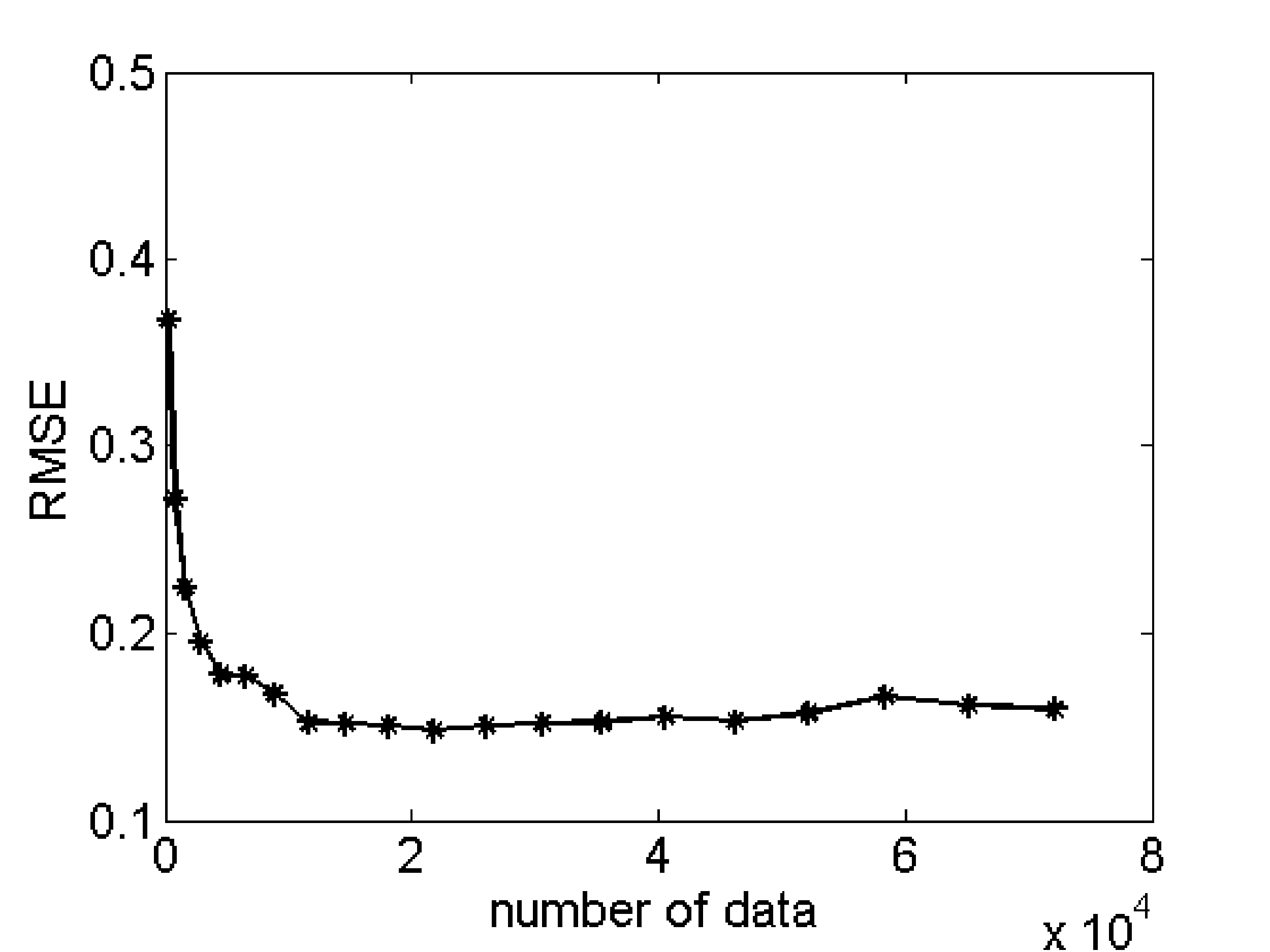}}%,height=0.3\textwidth
\subfigure[ \label{fig: fou_noise1}]%
{\includegraphics[width=0.49\textwidth]{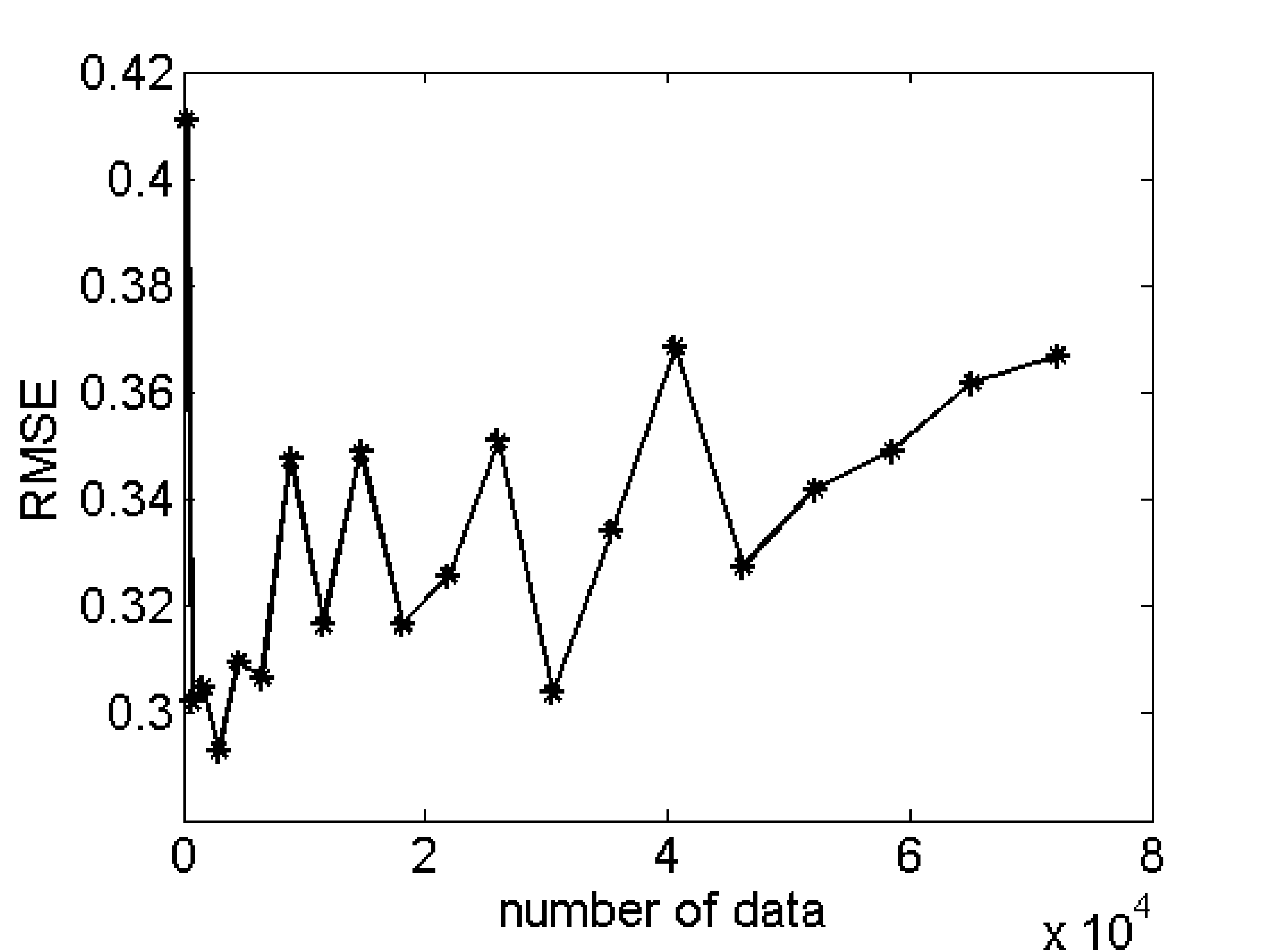}}
\subfigure[ \label{fig: fou_shepp}]%
{\includegraphics[width=0.49\textwidth]{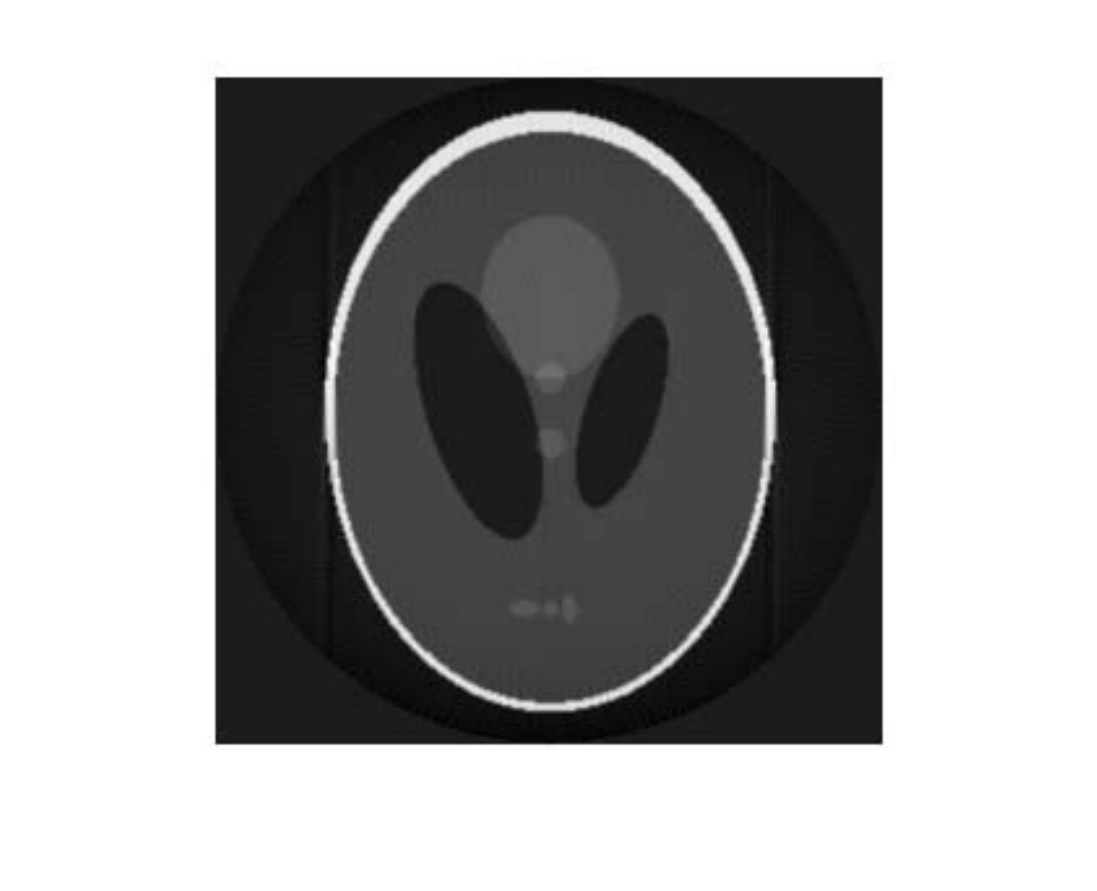}}%,height=0.3\textwidth
\subfigure[ \label{fig: fou_shepp_noise}]%
{\includegraphics[width=0.49\textwidth]{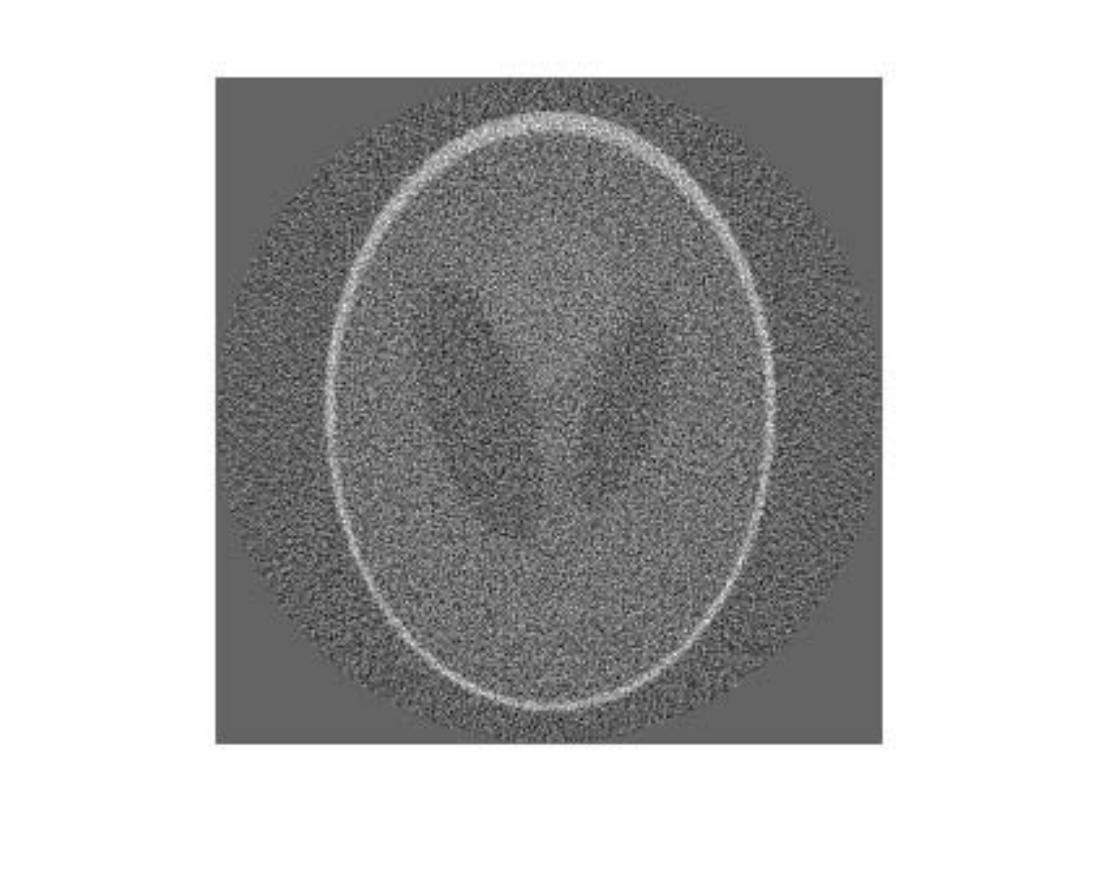}}
\caption{Reconstruction of the Shepp-Logan filter using the back projection formula with Shepp-Logan filter and linear interpolation. (a) $RMSE$ as a function of the number of data (without noise); (b) $RMSE$ as a function of the number of data (with Gaussian noise with mean $\mu=0.01$ and variance $\sigma=0.01$); (c) Reconstruction without noise with $N=180,\ M=200$; (d) Reconstruction with Gaussian noise with $N=180,\ M=200$.}
\label{fig: fou_rec_time}
\end{figure}

We start considering the behavior of the $RMSE$ in function of the number of the data $n$ (where again data are supposed to be taken using a parallel beam geometry). As one would aspect, with both kernel and Fourier-based methods, the $RMSE$ decreases when $n$ increases. In particular is interesting to notice that in the Fourier reconstruction, the $RMSE$ decreases with an exponentially rate and so, for large $n$, there is no big improvement of the solution. For example in the case of the Shepp-Logan phantom, for $n>18090$ there is a variation of the $RMSE$ lower than $5.67\cdot 10^{-3}$ (Figure \ref{fig: fou_rec1}). 

\begin{figure}[htbp]
\centering%
\subfigure[ \label{fig: fou_mq_rec1}]%
{\includegraphics[width=0.49\textwidth]{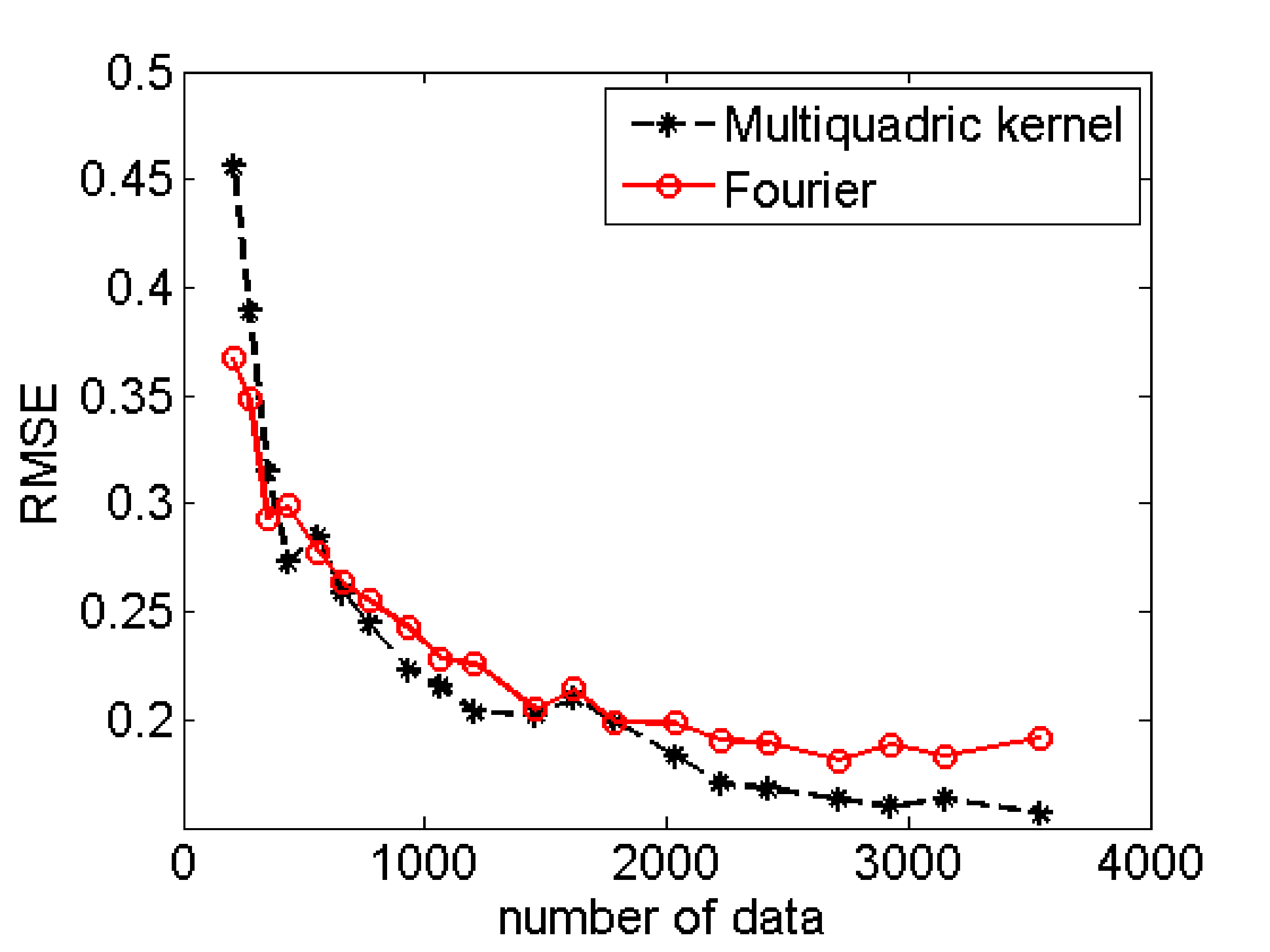}}%,height=0.3\textwidth
\subfigure[ \label{fig: mq_rcond1}]%
{\includegraphics[width=0.49\textwidth]{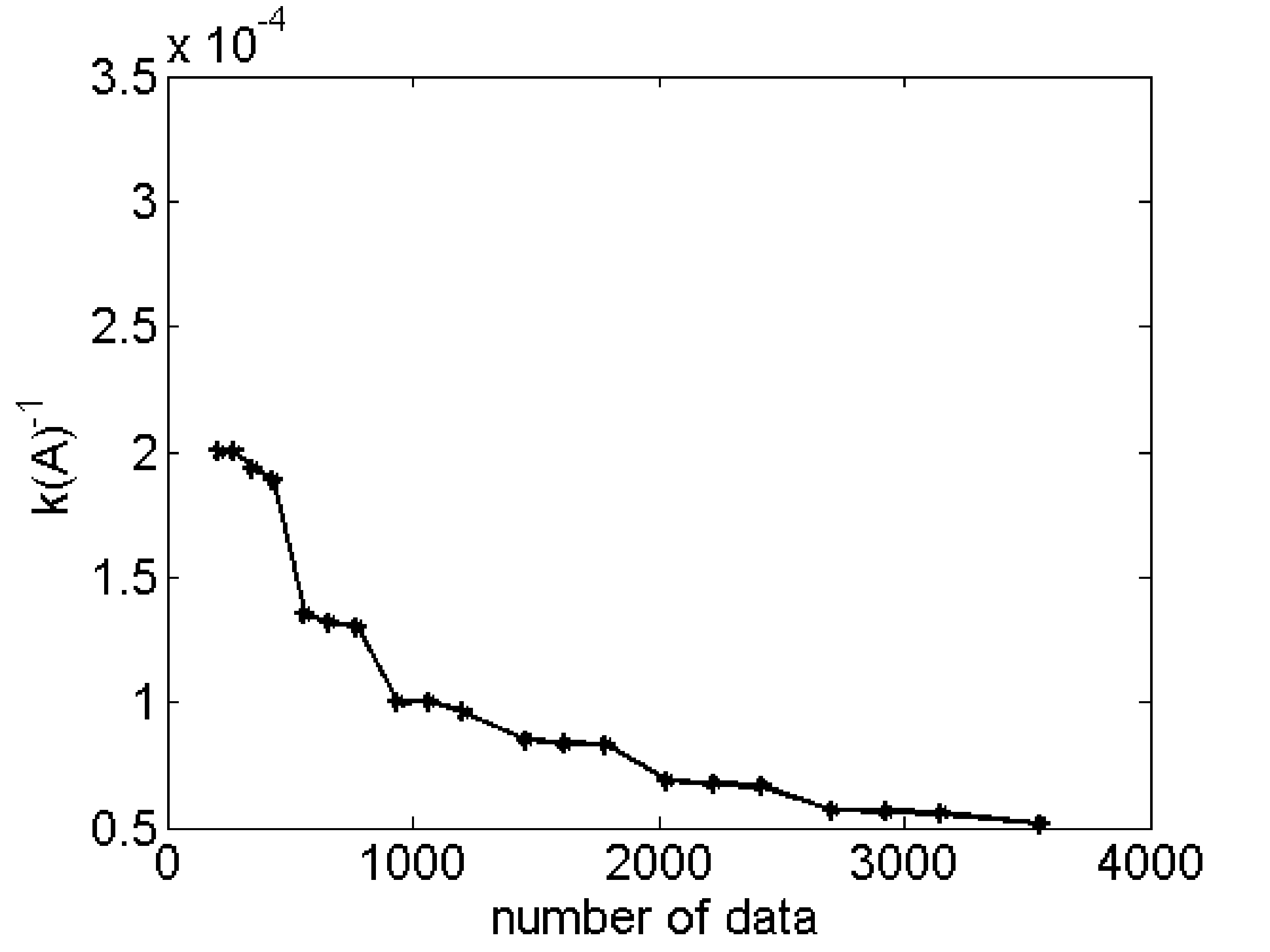}}
\subfigure[ \label{fig: fou_time2}]%
{\includegraphics[width=0.49\textwidth]{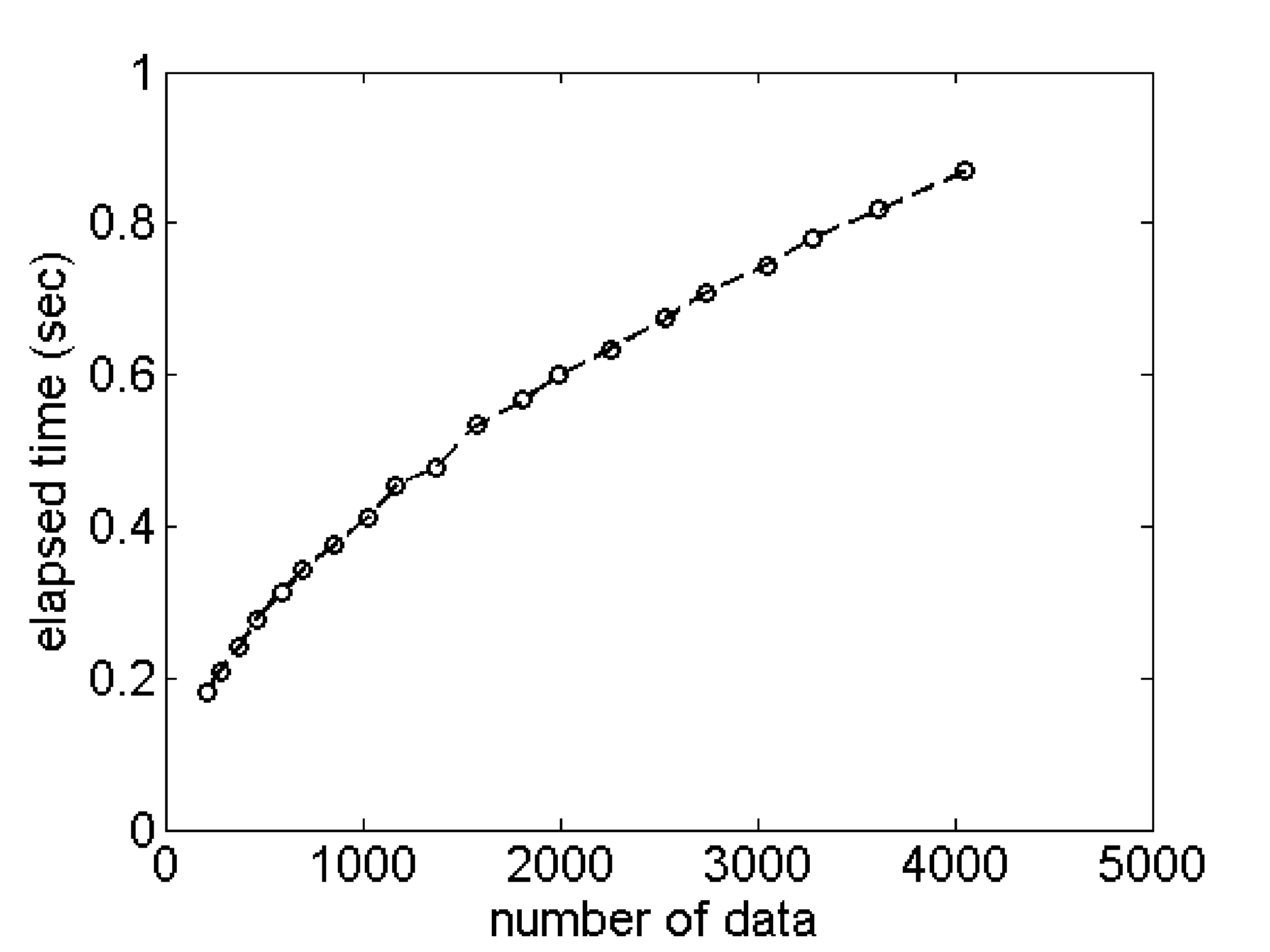}}%,height=0.3\textwidth
\subfigure[ \label{fig: mq_time1}]%
{\includegraphics[width=0.49\textwidth]{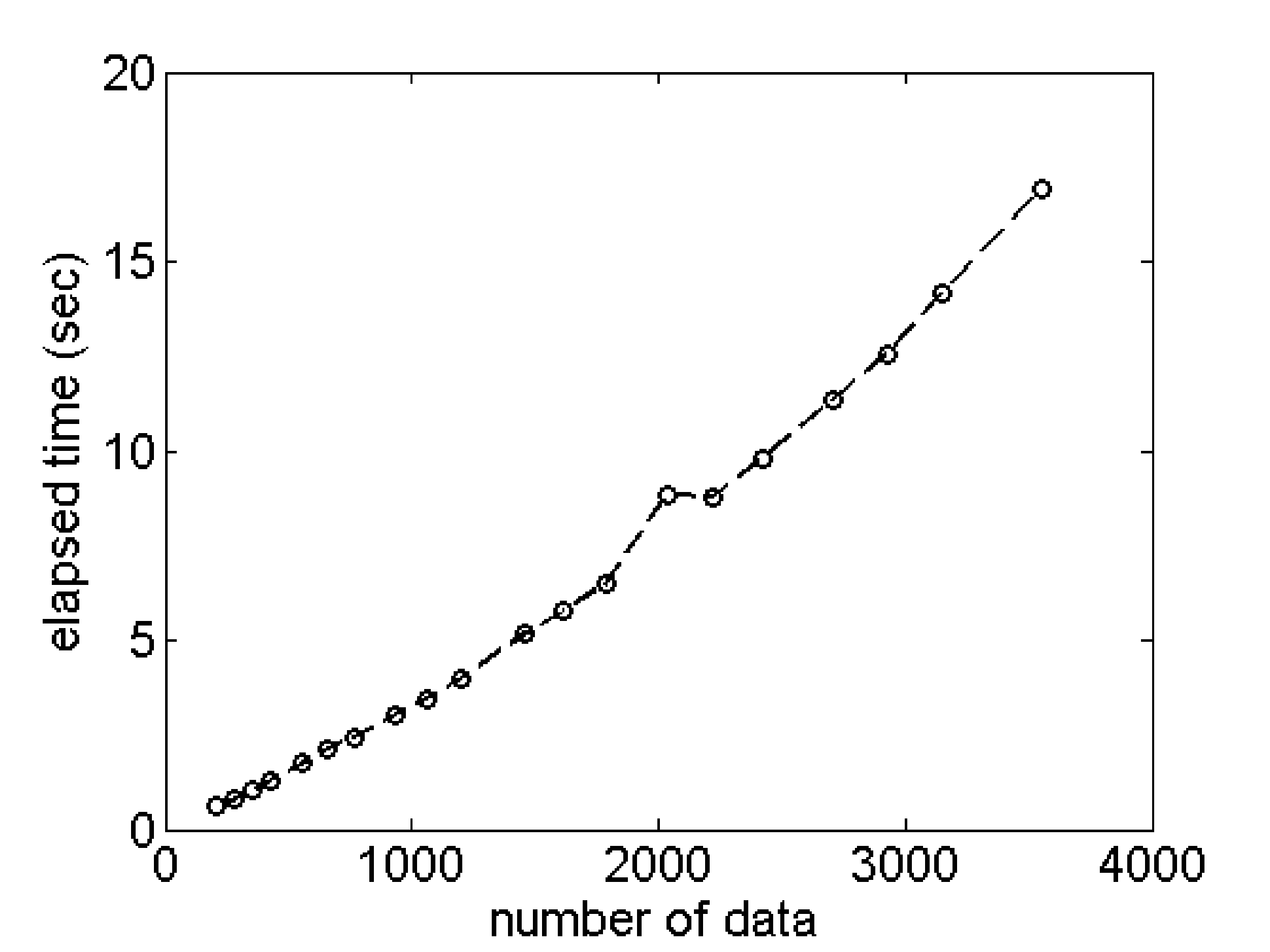}}
\caption{Comparison of the Fourier and multiquadric kernel reconstruction of the Shepp-Logan phantom: (a) Root mean square error; (b) Reciprocal of the condition number of matrix $A$ of the kernel method; (c) Elapsed time for Fourier based method; (d) Elapsed time for multiquadric reconstruction method.}
\label{fig: fou_mq}
\end{figure}

In the case of kernel reconstruction the use of big amount of data should consider CPU limits\footnote{All the computation shown in this chapter are made using a computer with a CPU core i5, 2,53 GHz and a RAM of 4 GB.}
. Indeed the matrix $A$ used to compute the coefficient $c$ of the solution belongs to the space $Mat(n,n)$, moreover, to evaluate the solution on a grid of $K\times K$ pixels, we must multiply $c$ for a matrix $B\in Mat(K^{2},n)$ representing the the basis functions of the space our solution belongs to. Thus, for example, with $K=256,\ N=50,\ M=40$, one obtains two (non-sparse) matrix one with $(N*(2M+1))^{2}=16.402.500$ elements and the other with $K^{2}\cdot(N*(2M+1))=265.420.800$ elements. Furthermore we notice that if $n$ increases, also $k(A)$ increases (see Figure \ref{fig: mq_rcond1}).

\begin{figure}[htbp]
\centering%
\subfigure[ \label{fig: fou_gau_rec1}]%
{\includegraphics[width=0.6\textwidth]{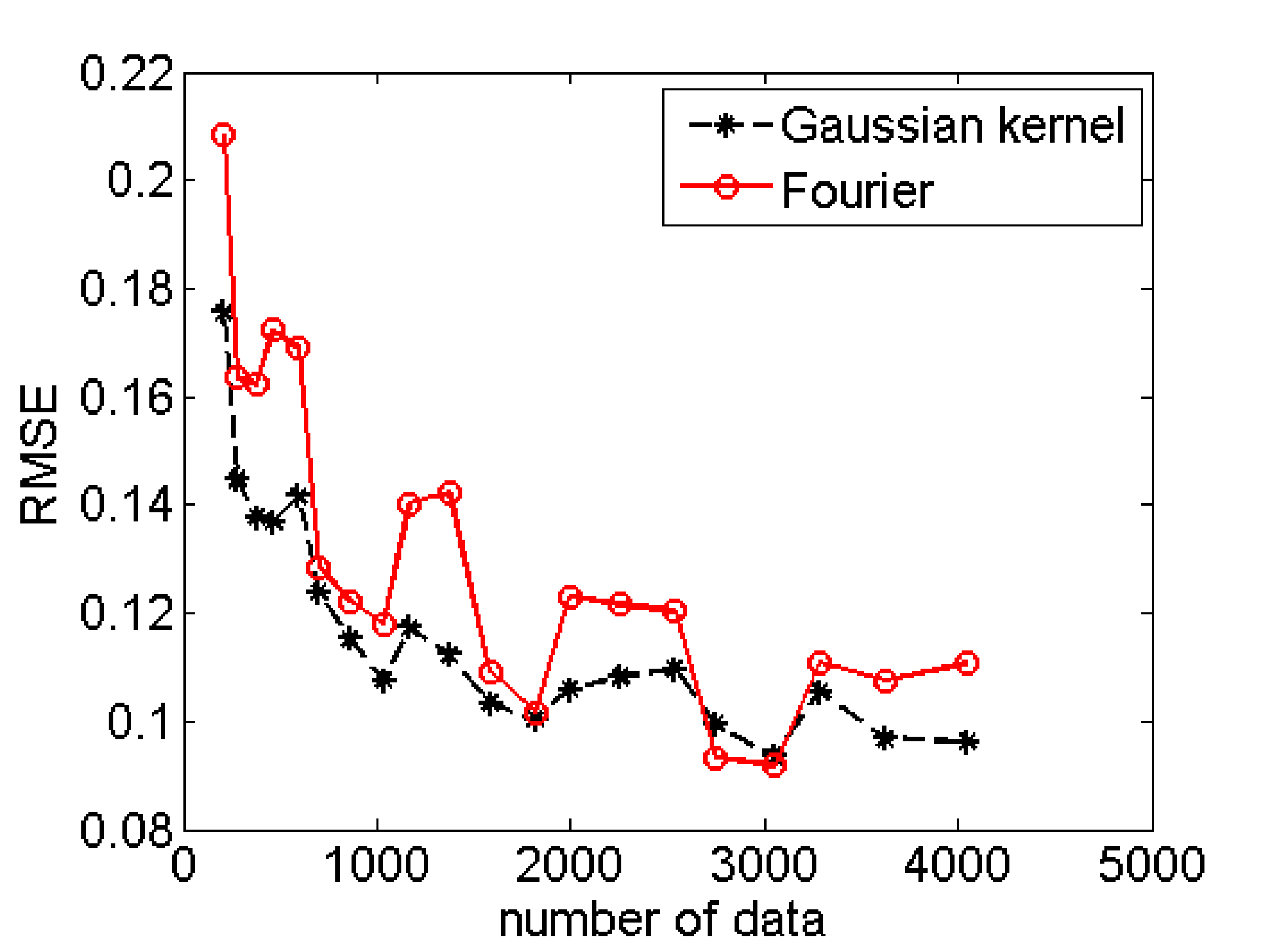}}%,height=0.3\textwidth
\\
\subfigure[ \label{fig: fou_cshape1}]%
{\includegraphics[width=0.49\textwidth]{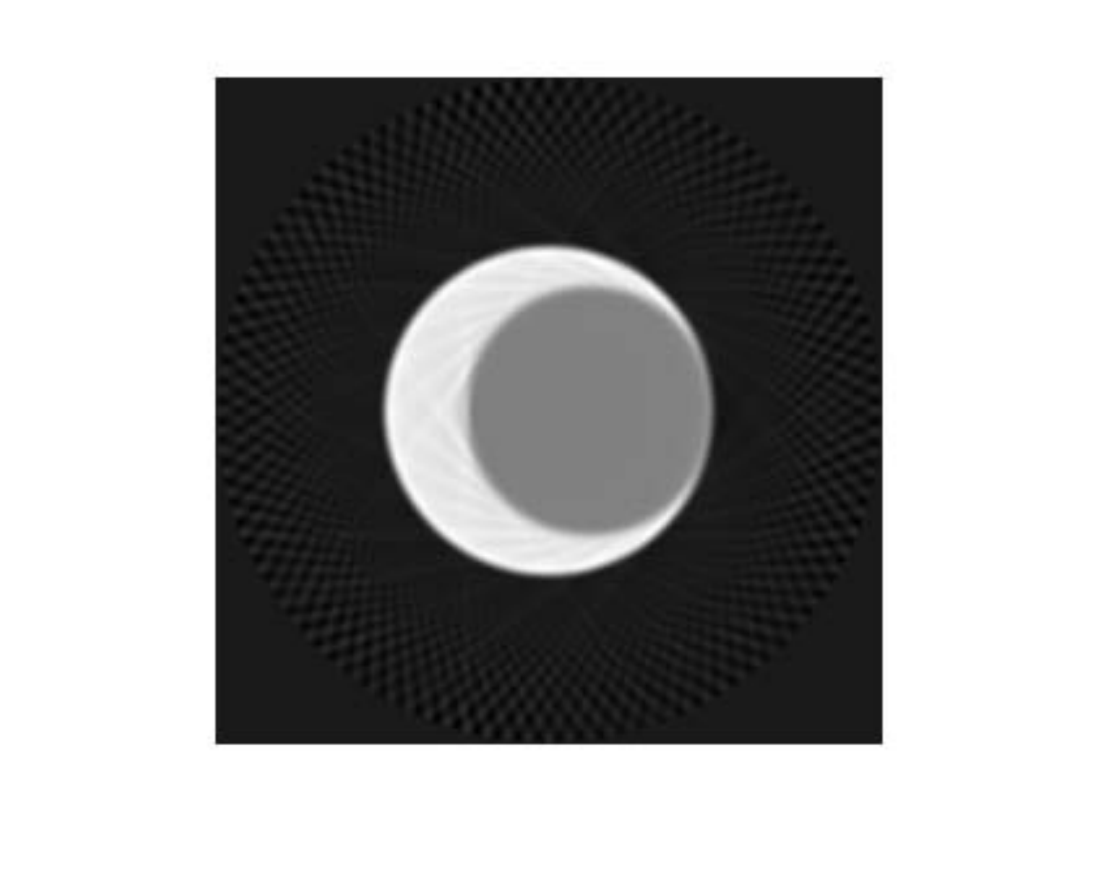}}
\subfigure[ \label{fig: gau_cshape}]%
{\includegraphics[width=0.49\textwidth]{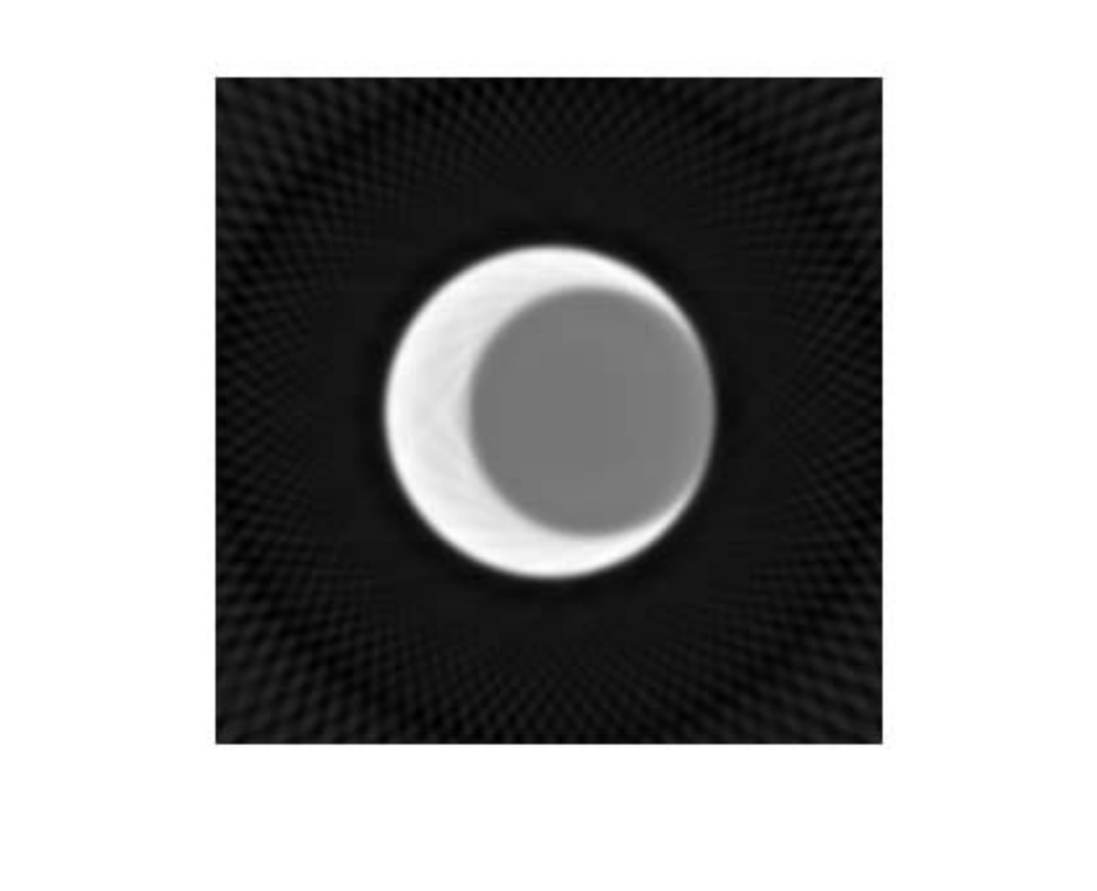}}%,height=0.3\textwidth
\caption{Comparison of Fourier and Gaussian kernel reconstruction methods of the crescent-shaped phantom: (a) Root mean square error; (b) Fourier method $N=50,\ M=40;$ (c) Gaussian kernel $N=50,\ M=40$.}
\label{fig: cshape_fou_gau_rec}
\end{figure}

Considering a problem with reasonable values of $n$ and $K$ (e.g. $n<130,\ K<256$), we observe that the $RMSE$ of the kernel methods behaves in the same way as the Fourier based methods and has a comparable magnitude (Figure \ref{fig: fou_mq_rec1}).
In particular, Figure \ref{fig: cshape_fou_gau_rec} shows how the Gaussian kernel reconstruction applied to the crescent-shaped phantom gives better result w.r.t. to the back projection formula.

Thus we conclude that kernel based methods can be useful in the context of limited number of available data.

\begin{figure}[htbp]
\centering%
\subfigure[ \label{fig: g_mq__imq_c_rec1}]%
{\includegraphics[width=1\textwidth]{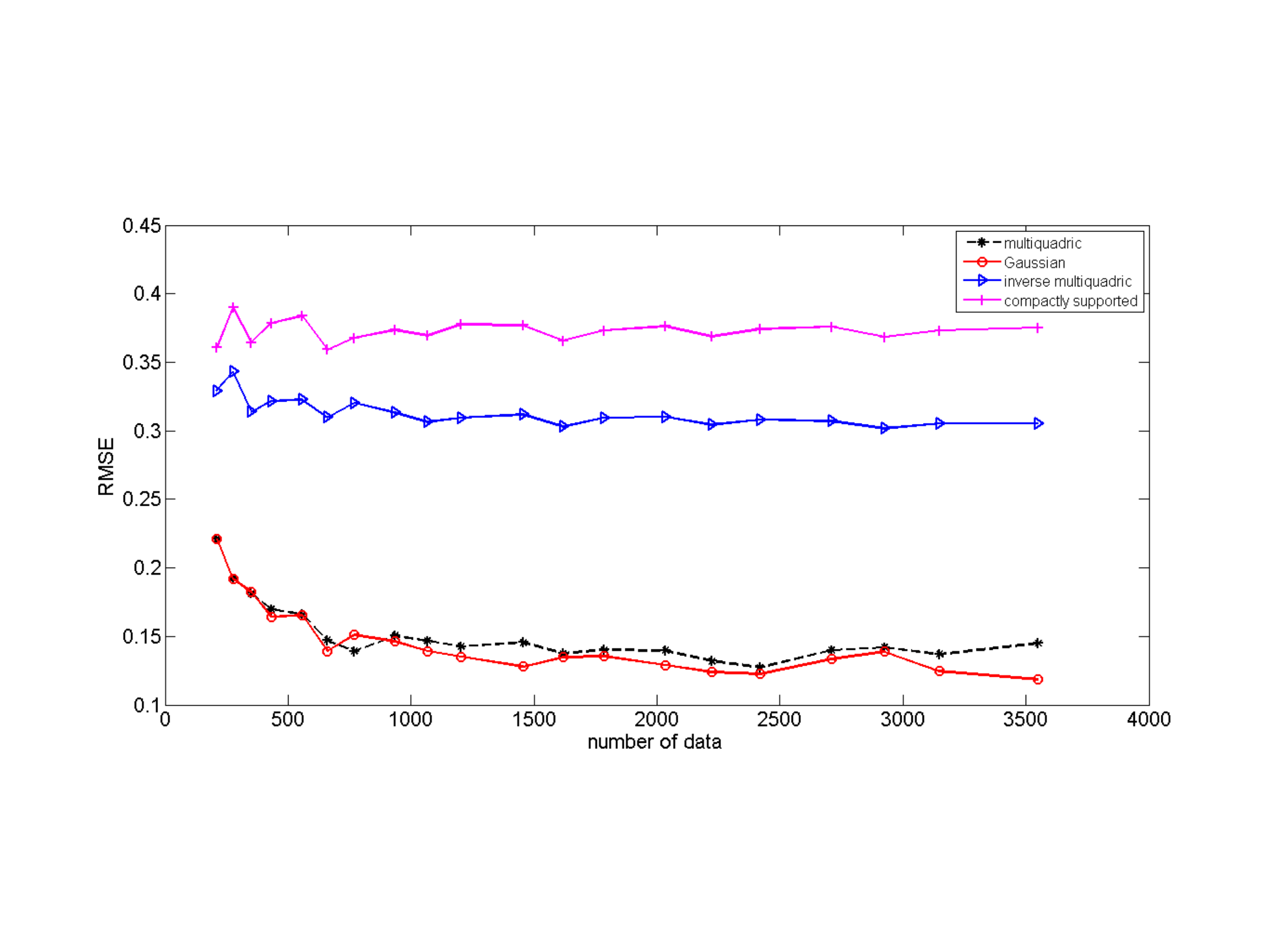}}%,height=0.3\textwidth
\\
\subfigure[ \label{fig: imq_bull1}]%
{\includegraphics[width=0.49\textwidth]{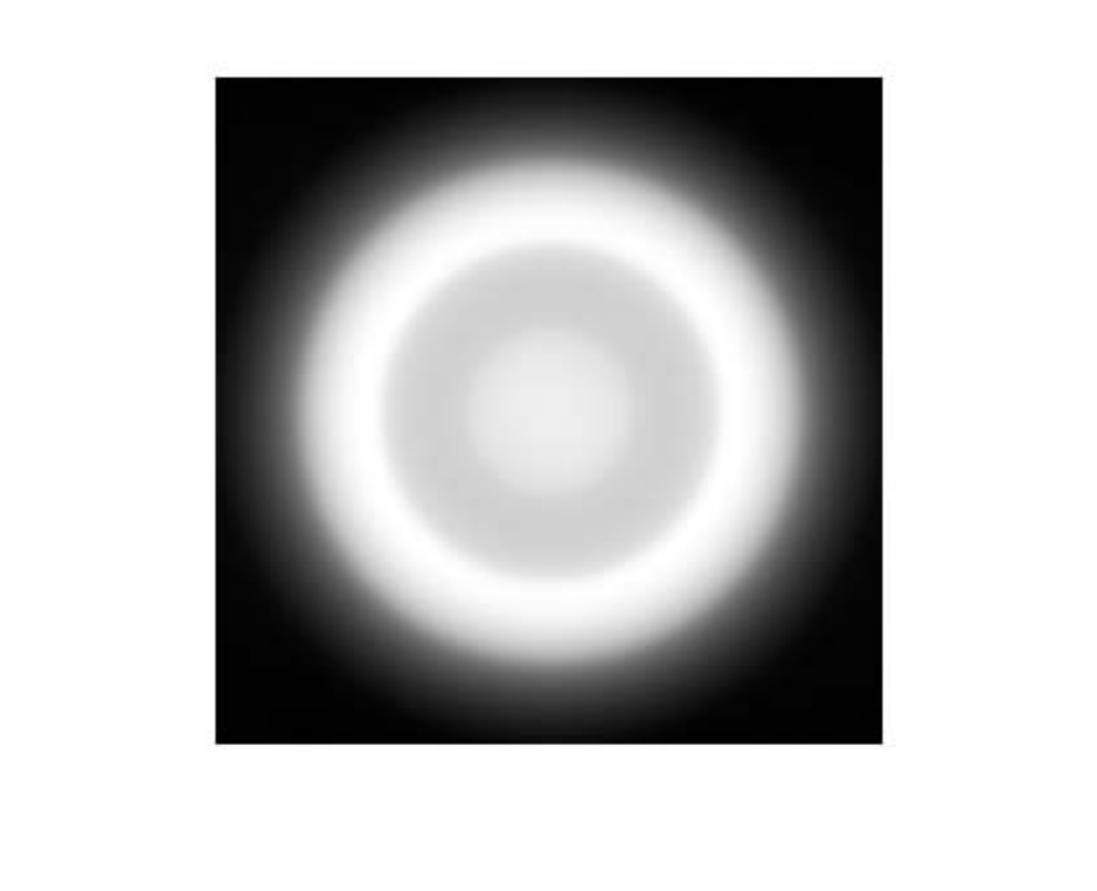}}
\subfigure[ \label{fig: gau_bull1}]%
{\includegraphics[width=0.49\textwidth]{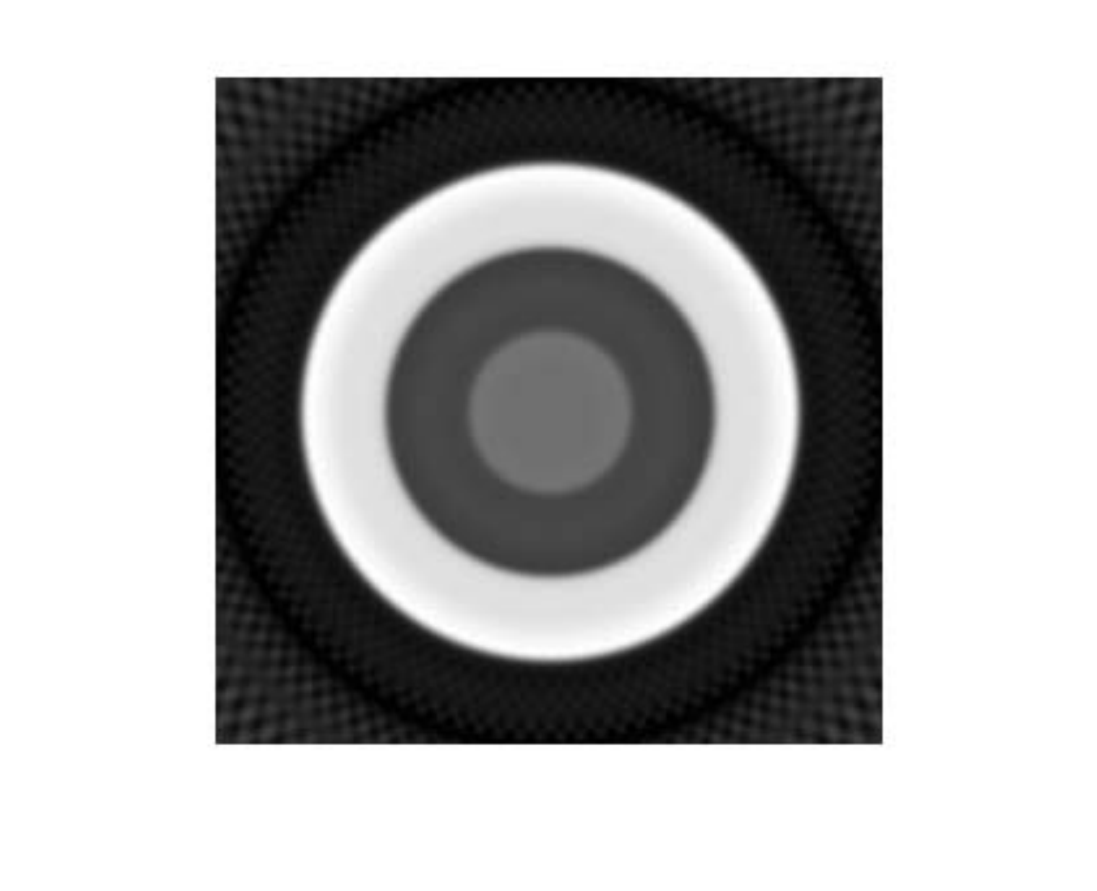}}%,height=0.3\textwidth
\caption{Comparison of kernel reconstruction methods of the bull's eye phantom using different kernel functions: (a) Root mean square error; (b) Reconstruction with inverse multiquadrics kernel $N=50,\ M=35;$ (c) Reconstruction with Gaussian kernel $N=50,\ M=40$.}
\label{fig: bull_ker_rec}
\end{figure}

Comparing the reconstruction of a phantom using different kernels, we see that the Gaussian-multiquadric kernel and the Gaussian kernel give similar results while the truncated-multiquadric kernel and the compactly supported kernel are less accurate (Figure \ref{fig: bull_ker_rec}).

The biggest problem in the kernel methods, a part from the memory limitation, is the computational time. Referring to Figure \ref{fig: fou_time2} and \ref{fig: mq_time1}, we can see how the elapsed time (in sec) during the execution of the algorithm grows exponentially with the number of data, while using Fourier techniques the time depend linearly on the dimension of the problem. We belive that this is due to the implementation in MATLAB of Fourier transform by the FFTW algorithm \cite{FFTW}.

\begin{figure}[htbp]
\centering%
\subfigure[ \label{fig: fou_cfr_noise1}]%
{\includegraphics[width=0.8\textwidth]{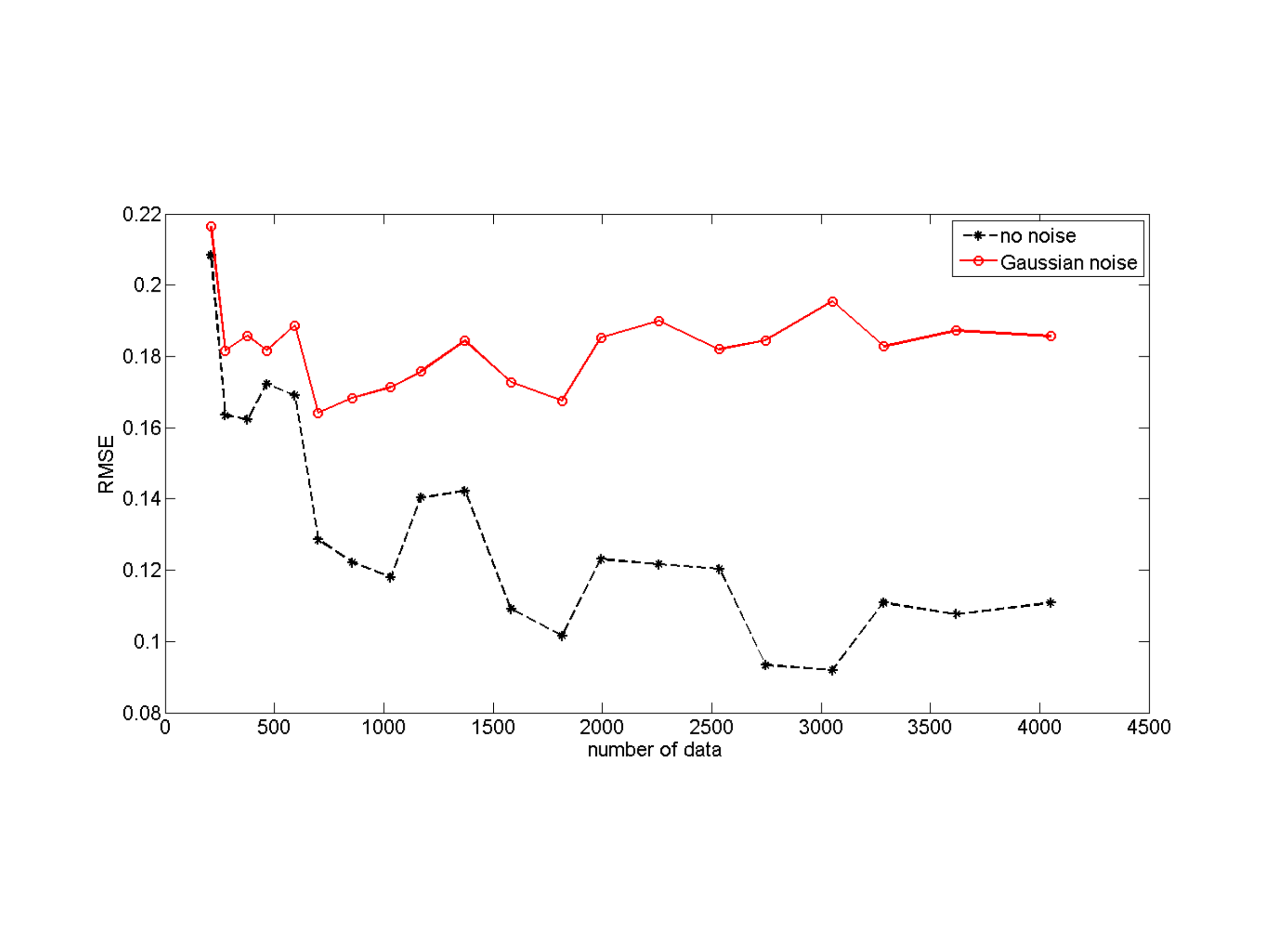}}%,height=0.3\textwidth

\subfigure[ \label{fig: gau_cfr_noise1}]%
{\includegraphics[width=0.8\textwidth]{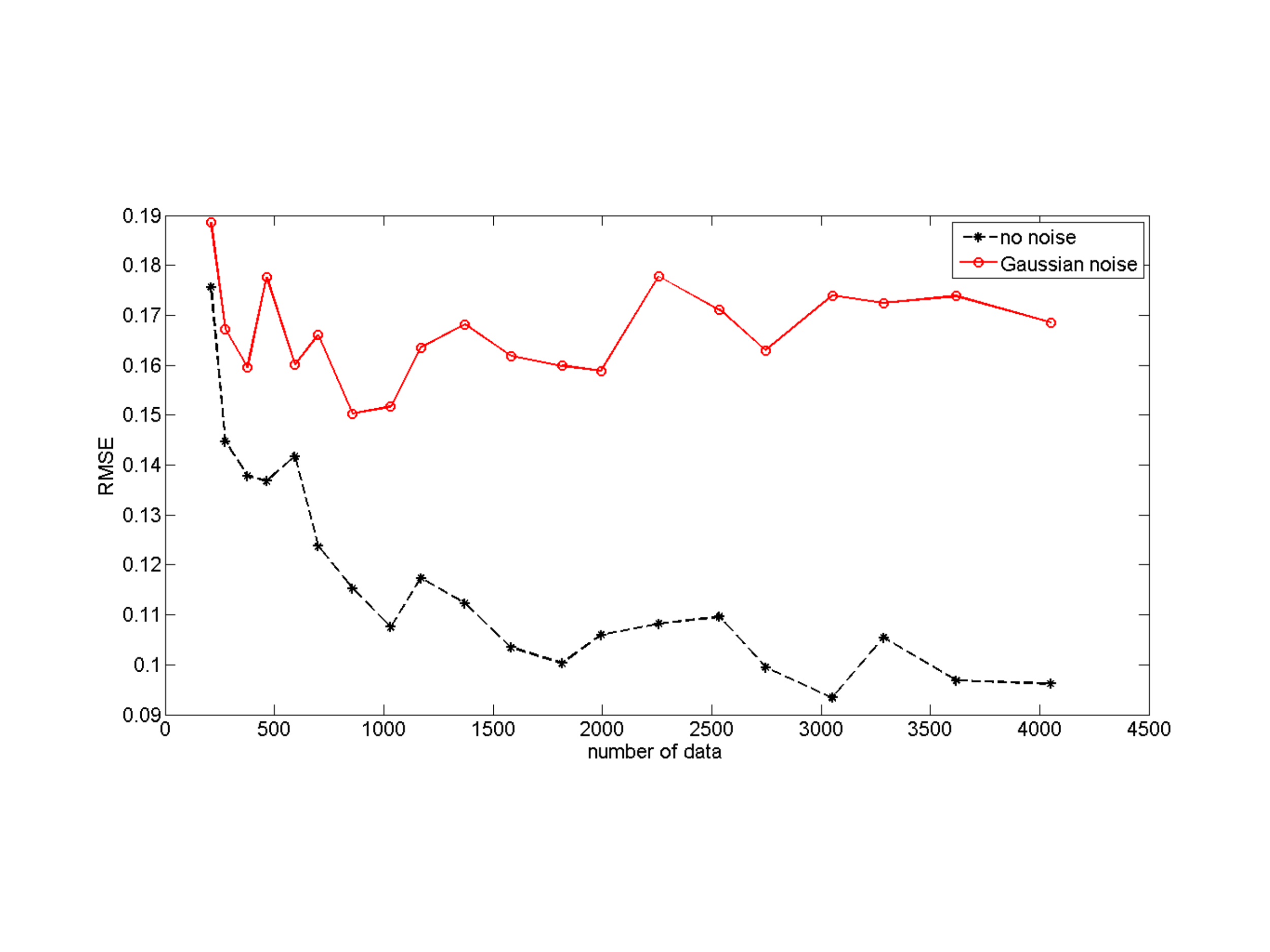}}
\caption{Comparison of the Fourier and Gaussian kernel reconstruction of the crescent-shaped phantom with no noise data and Gaussian noise data with 0.001 mean and 0.001 variance: (a) Fourier method; (b) Gaussian kernel method.}
\label{fig: cfr_noise}
\end{figure}

\begin{figure}[htbp]
\centering%
\subfigure[ \label{fig: fou_cshape_noise}]%
{\includegraphics[width=0.49\textwidth]{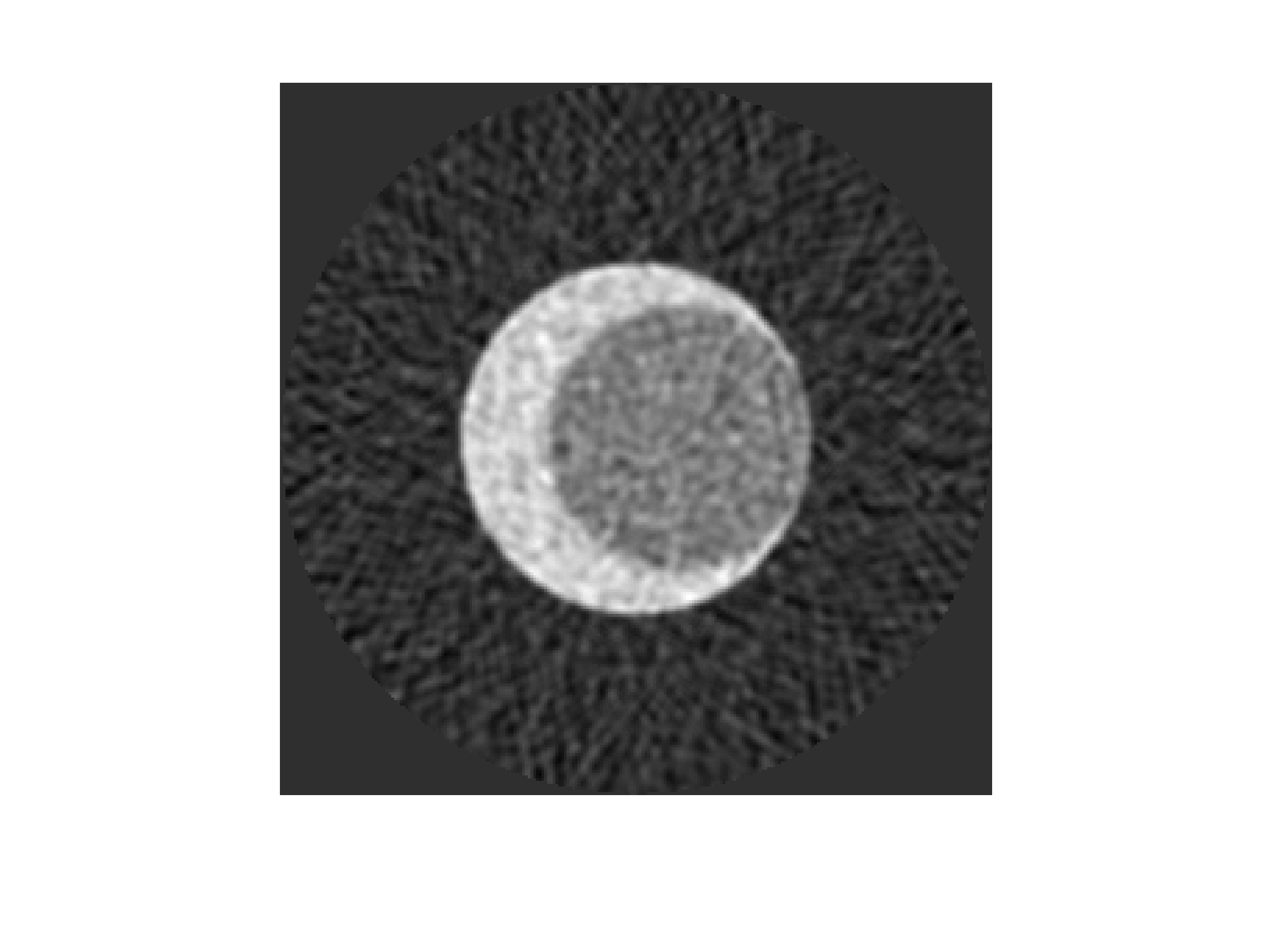}}%,height=0.3\textwidth
\subfigure[ \label{fig: gau_cshape_noise}]%
{\includegraphics[width=0.49\textwidth]{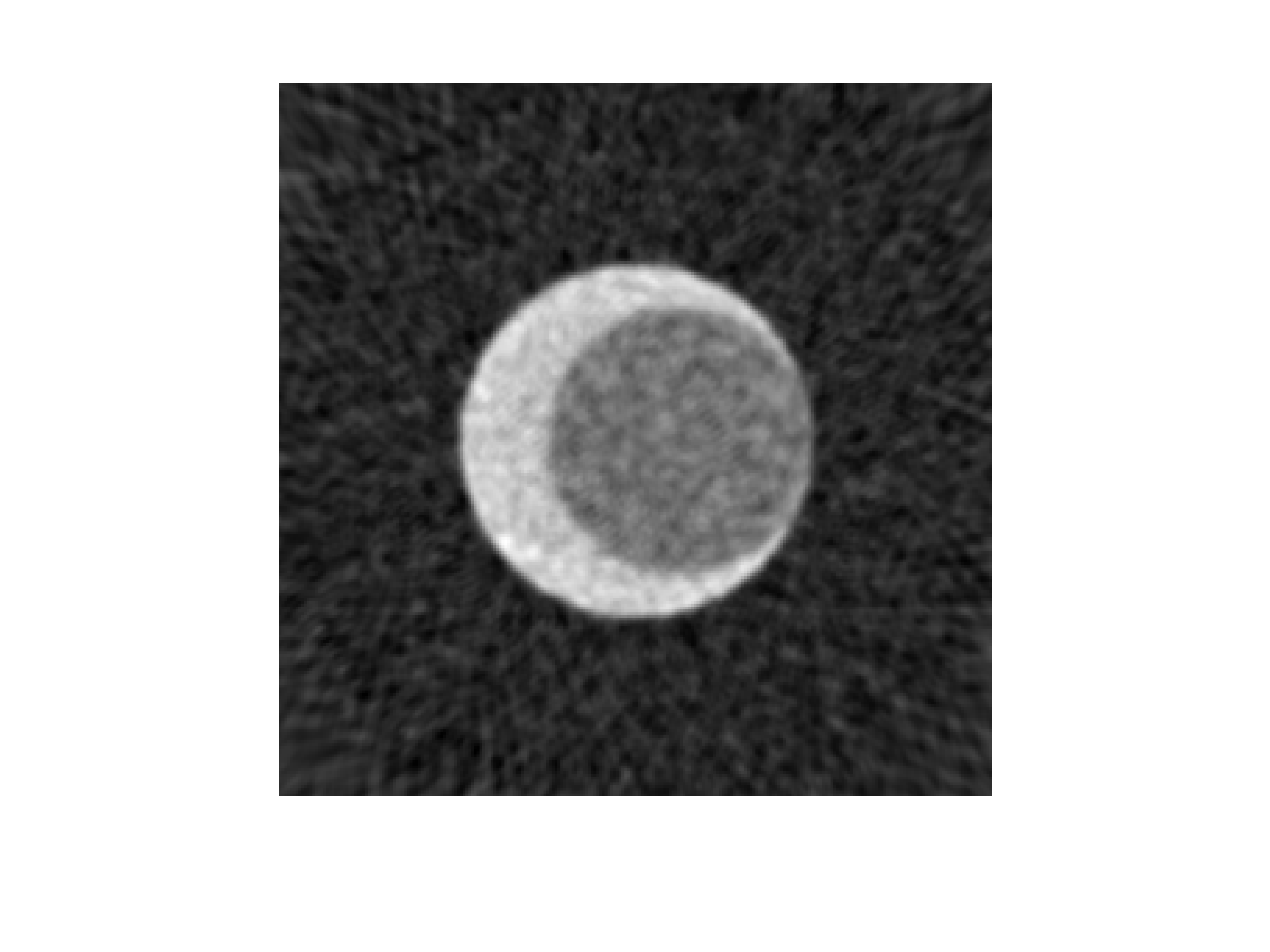}}
\caption{Comparison of the Fourier and Gaussian kernel reconstruction of the crescent-shaped phantom with Gaussian noise data (mean $\mu= 0.001$ and variance $\sigma=0.001$: (a) Fourier method; (b) Gaussian kernel method.}
\label{fig: cfr_cshape_noise}
\end{figure}

Finally we test our methods after the introduction of noise in the data. We first notice that for large $n$ the $RMSE$ of the Fourier methods increase (Figure \ref{fig: fou_noise1}).
Considering a kernel based method we observe (in the range of acceptable $n$) a behavior similar to the Fourier case, where the $RMSEs$ computed with different methods have the same order of magnitude (Figure \ref{fig: cfr_noise} and \ref{fig: cfr_cshape_noise}).

\section{Graphical user interface}
In order to test the various methods on different phantoms we developed a graphical user interface (GUI) allowing the user to choose the options of the reconstruction and the parameters using the mouse and the keyboard and to access to all output information of the method.

\begin{figure}[htbp]
\centering
\includegraphics[width=\textwidth]{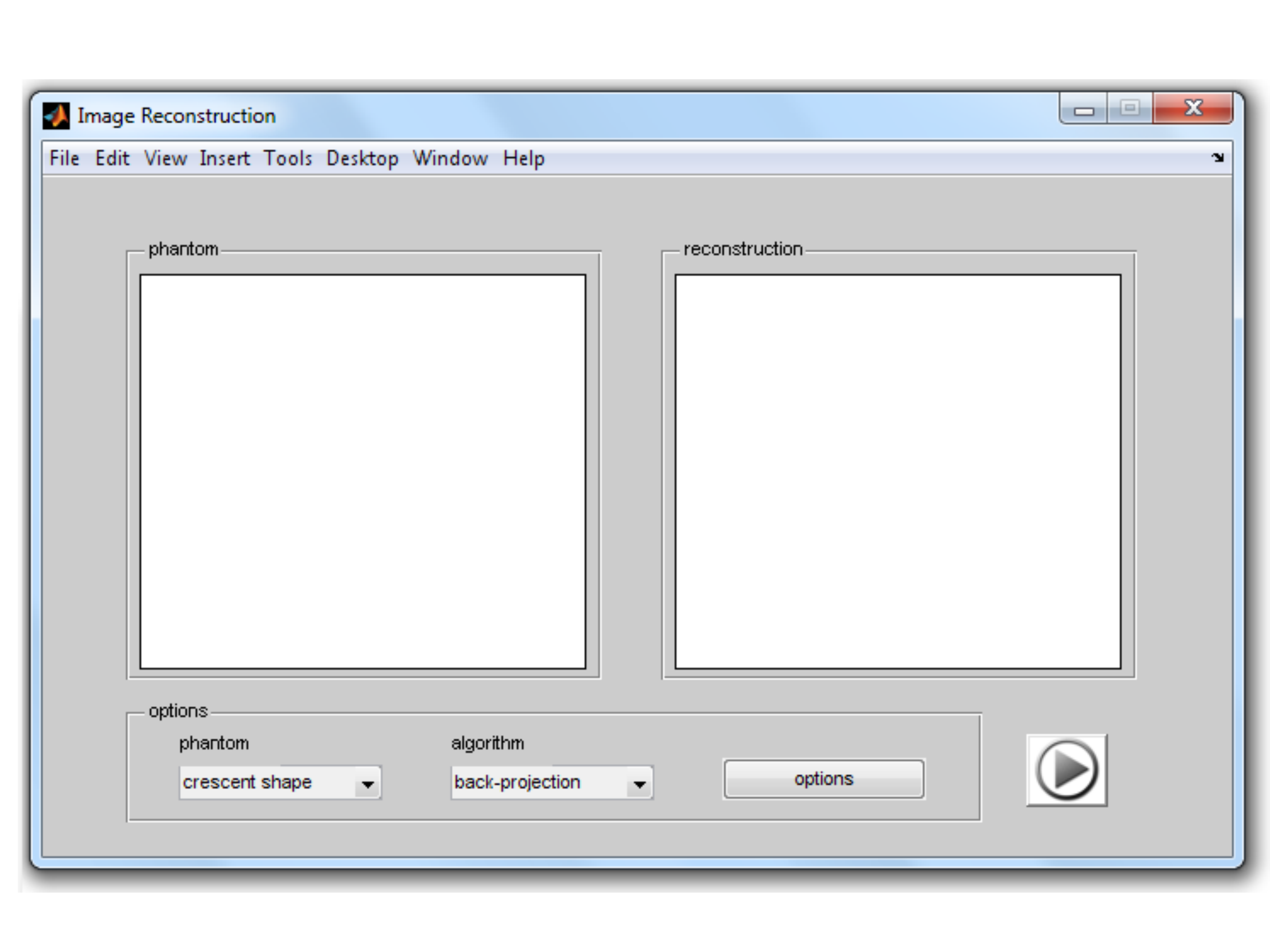} 
\caption{Graphical user interface}
\label{fig: gui}
\end{figure}

\begin{figure}[htbp]
\centering
\includegraphics[width=0.5\textwidth]{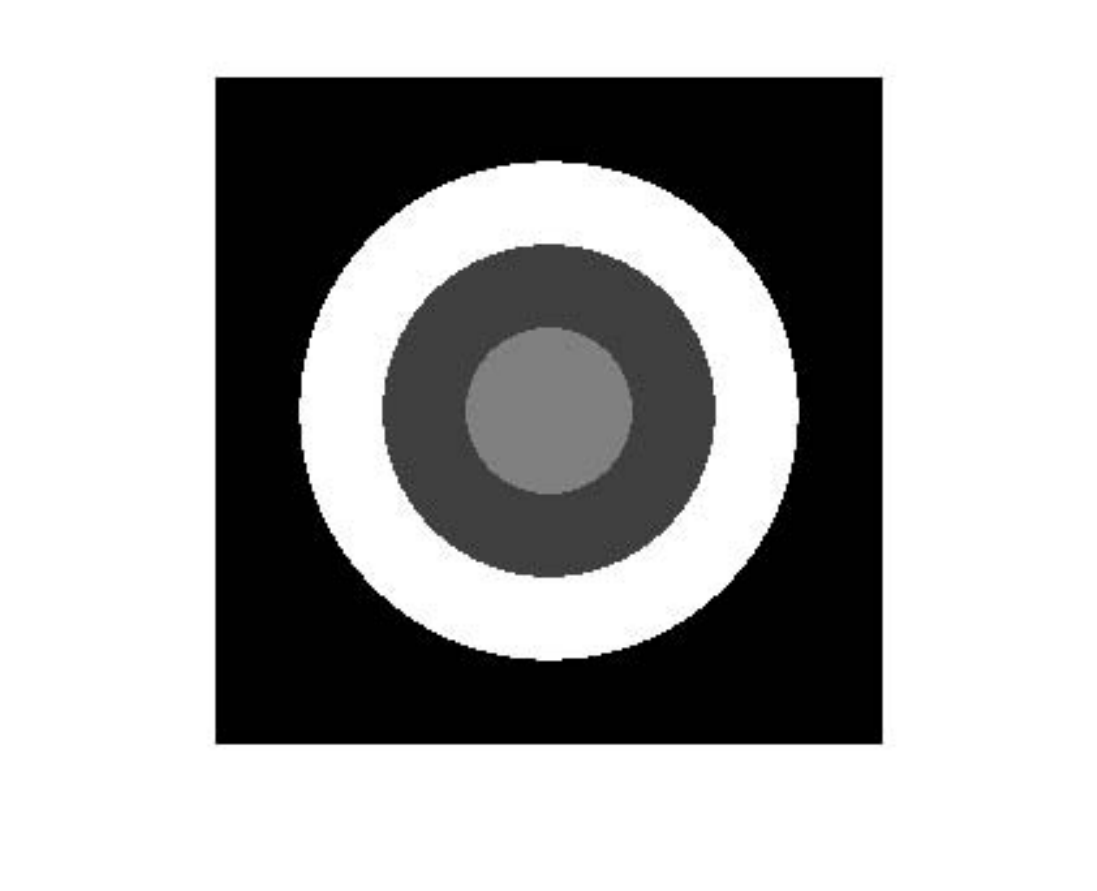} 
\caption{Bull's eye phantom}
\label{fig: bull_eye}
\end{figure}

\begin{figure}[htbp]
\centering
\includegraphics[width=0.6\textwidth]{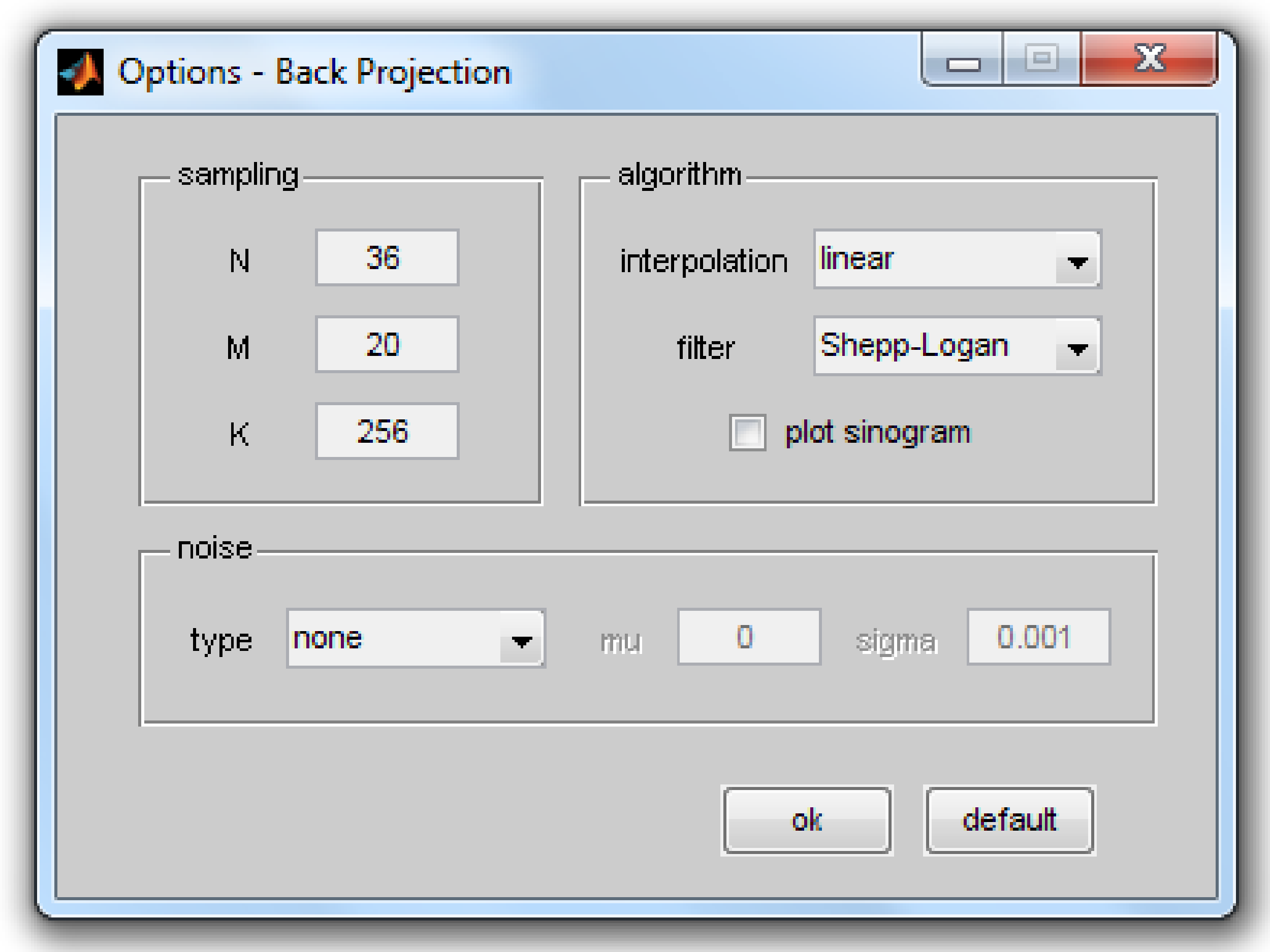} 
\caption{Options window}
\label{fig: option_figure}
\end{figure}

\begin{figure}[htbp]
\centering
\includegraphics[width=0.8\textwidth]{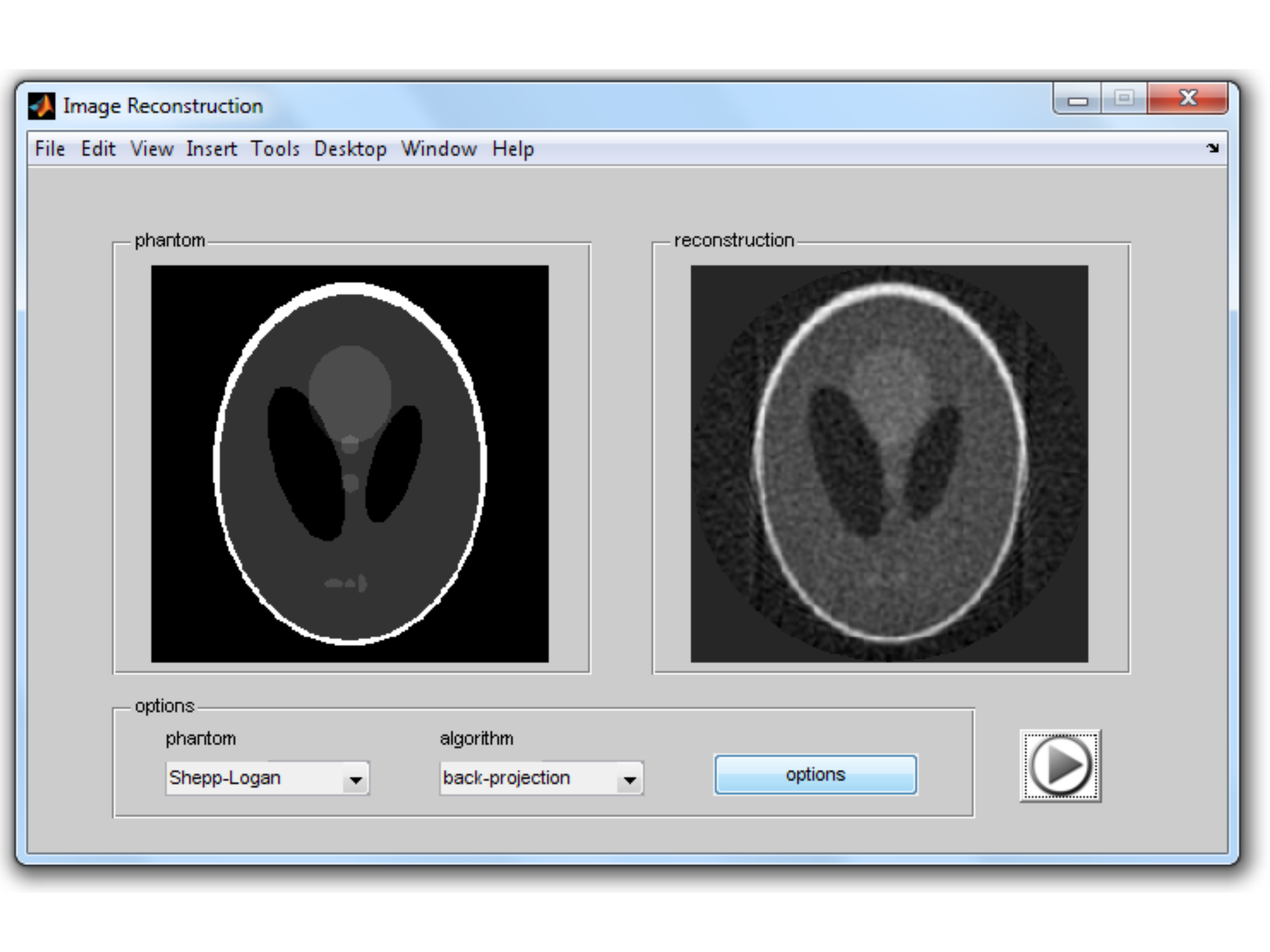} 
\caption{Shepp-Logan phantom reconstruction with back-projection formula and Gaussian noise-data (mean $\mu=0$, variance $\sigma=0.001$).}%N=120, M=50 K=256, linear shepp log, mu=0 sigma=0.001
\label{fig: gui_shepp_noise}
\end{figure}

The GUI has been realized in MATLAB (version 7.7.0 R2008b). Figure \ref{fig: gui} shows the main window of the GUI when it is started. Referring always to Figure \ref{fig: gui} we can see it presents two windows where the original phantom and the reconstructed image of the phantom will be displayed. 
On the bottom of the window there is a panel that, thanks to pop-up menus, allows to choose the phantom and the reconstruction algorithm. There are three available phantoms: the crescent-shaped phantom (introduced in section \ref{subsec: crescent_shape}), the bull's eye phantom (Figure \ref{fig: bull_eye}) and the Shepp-Logan phantom (Figure \ref{fig: shepp_phantom} in section \ref{sec: phantoms}).
Clicking on the options button a second window will be opened (figura \ref{fig: option_figure}).
Thanks to this options window it is possible to modify the number of sampled angle $N$, the number of samples on the $t$ axis ($2M+1$) and the dimension of the output image ($K^{2}$). Moreover, depending on the selected algorithm, it is possible to change the predefined parameters used in the methods. For example, when using the back projection formula, one can choose both interpolation technique (nearest neighbor, linear, cubic) and the low pass filter (Ram-Lak, Shepp-Logan, cosine filter). Finally one can add a certain amount of noise to the Radon data to test the robustness of a method under the action of noise. Possible choices of noise are the Gaussian noise (with mean and variance decided by the user), Poisson noise and shot (or "salt and pepper") noise. Figure \ref{fig: gui_shepp_noise} shows the result of applying the back projection formula to the Shepp-Logan phantom adding Gaussian noise to data.

\begin{figure}[htbp]
\centering
\includegraphics[width=0.8\textwidth]{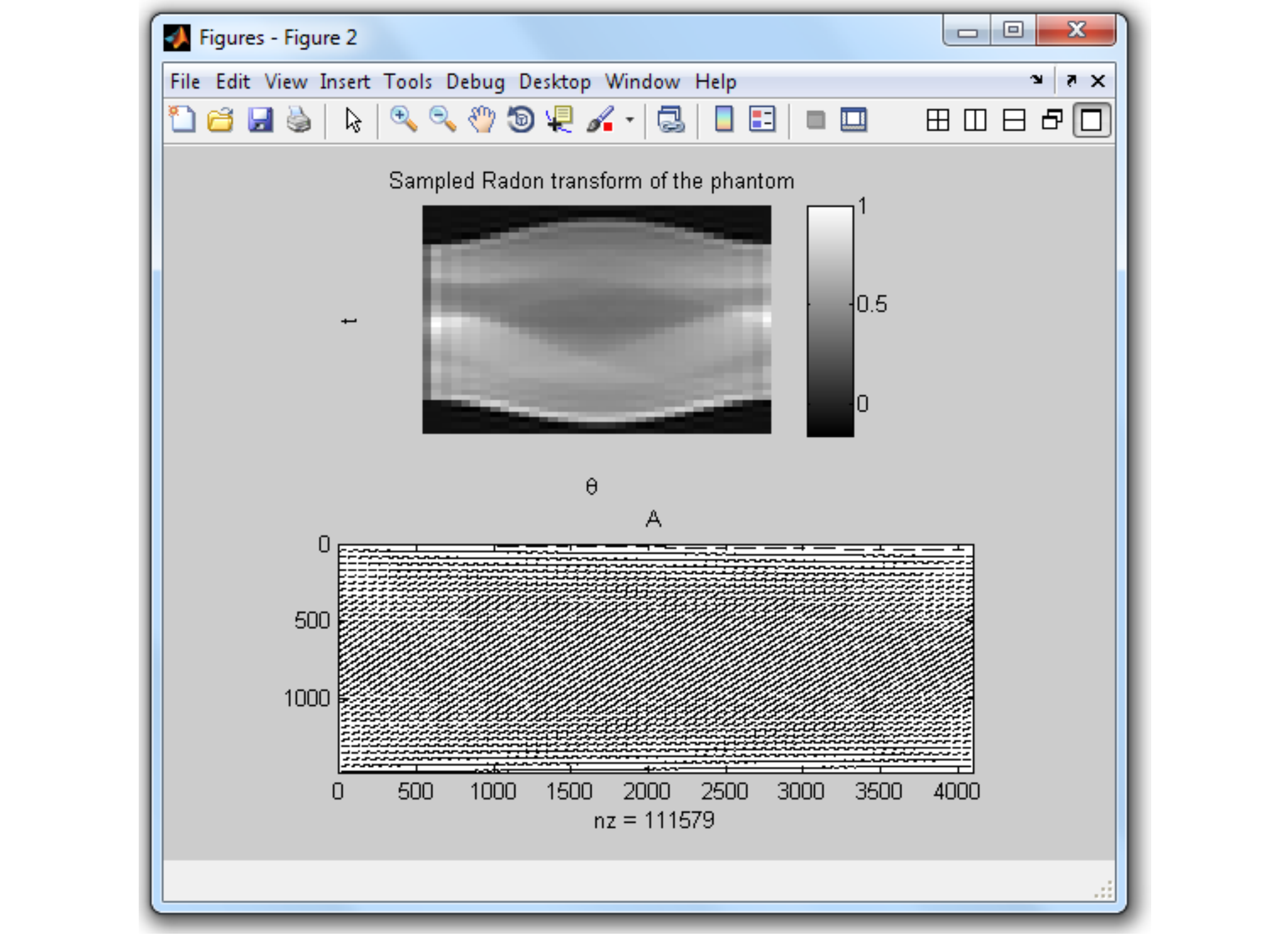} 
\caption{Sinogram and sparsisty of the matrix system of the Kaczmarz method.}%shepp logan phantom default options
\label{fig: sinogram_kaczmarz}
\end{figure}

It is also always possible to decide if to plot in a new figure the sinogram of the phantom, i.e. the sampled Radon transform used for the reconstruction. Figure \ref{fig: sinogram_kaczmarz} shows the case of Kaczmarz's method.
In addition, depending on the algorithm a further plot will be displayed:
\begin{itemize}
\item Using the back projection formula, the interpolation of the convolution between the low pass filter and the Radon transform is shown;
\item With the Kaczmarz's method, the sparsity of the system matrix (figura \ref{fig: sinogram_kaczmarz});
\item In the case of kernel methods, the image of the sysem matrix.
\end{itemize}

\begin{figure}[htbp]
\centering%
\subfigure[\label{fig: gui_cshape_mq}]%
{\includegraphics[width=0.49\textwidth]{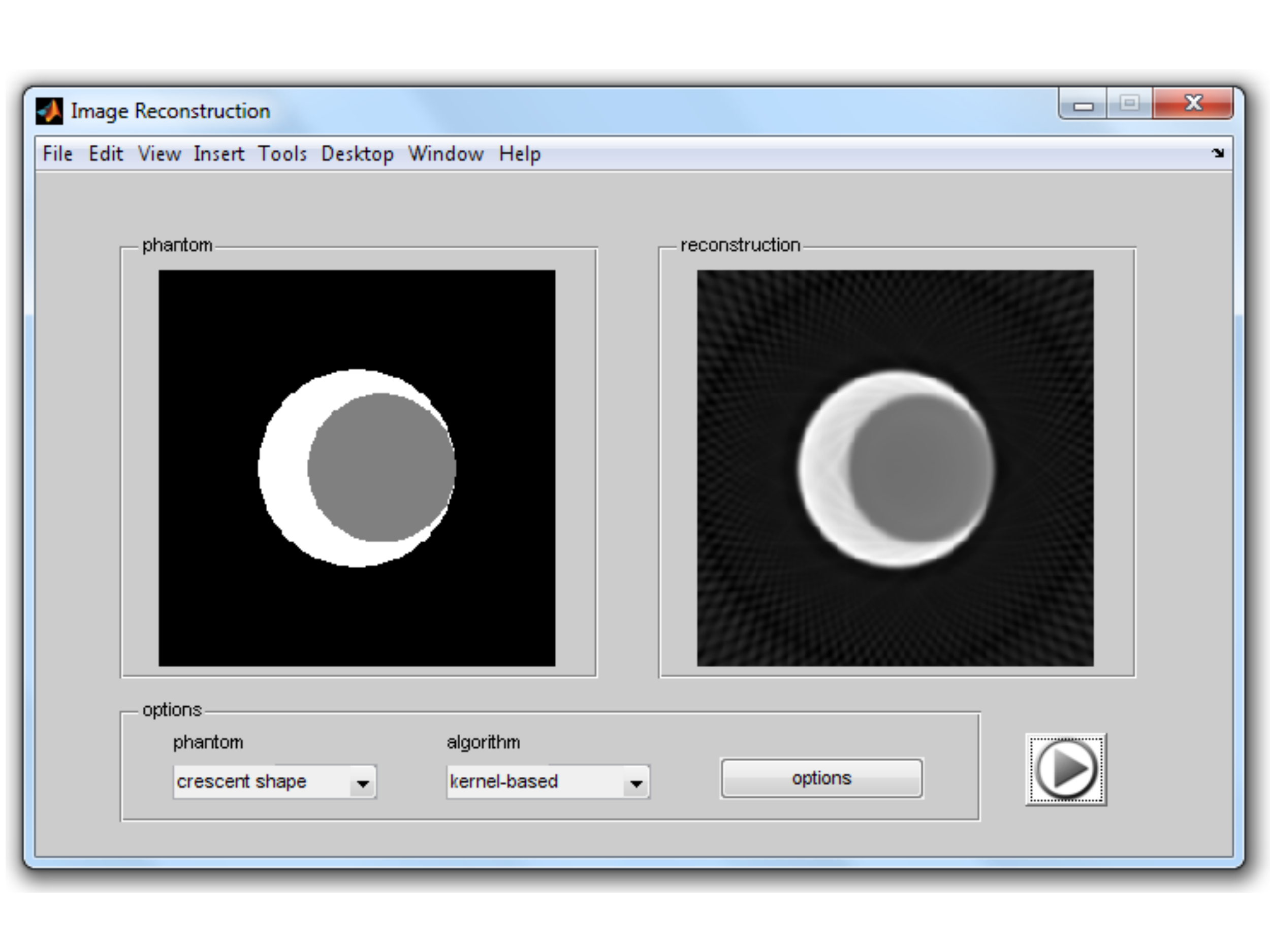}}%,height=0.3\textwidth
\subfigure[\label{fig: gui_bulleye_art}]%
{\includegraphics[width=0.49\textwidth]{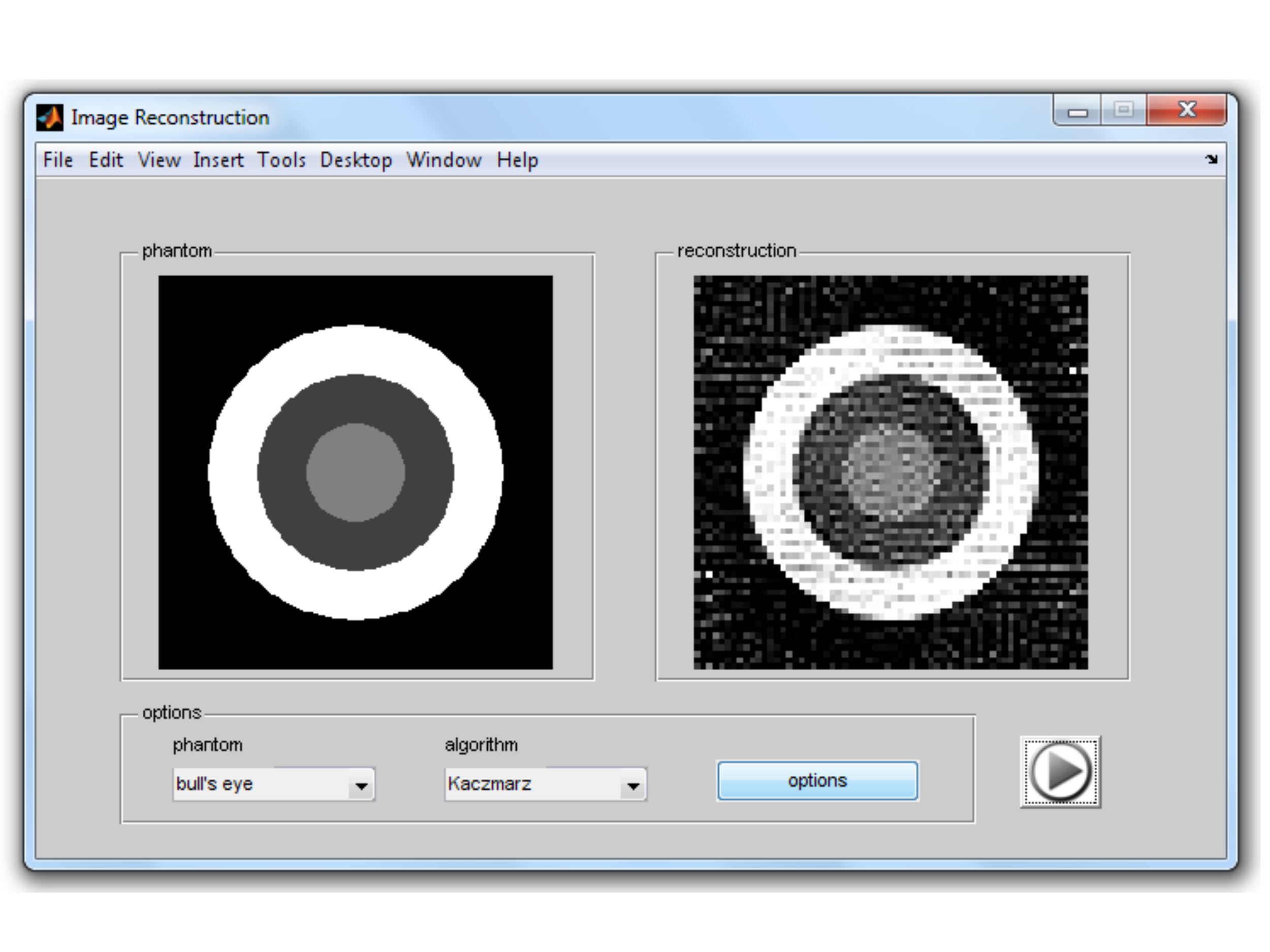}}
\caption{Example of reconstruction using the graphical user interface. (a) Multiquadric reconstruction of the crescent-shaped phantom; (b) Kaczmarz's reconstruction of the bull's eye phantom.}
\label{fig: gui_examples}
\end{figure}

The reconstruction begins pushing the start button. Figures \ref{fig: gui_cshape_mq} and \ref{fig: gui_bulleye_art} show two different applications of the GUI. 
During the computation of the solution, in the MATLAB command window, it is displayed the status of the process, e.g. for the Kaczmarz's method
\begin{verbatim}
Start reconstruction..
Radon transform computed..
computing ART system..
applying Kaczmarz's method..
residual: 
    0.5895

done.
\end{verbatim}
Where residual indicates the norm of the difference $Ax_{k}-b$ of the solution at the final iteration $k$.

Closing the GUI window or saving the workspace from the menu, one can find all information about the solution and the algorithm used in the output structure \texttt{out}. This structure contains the following fields:
\begin{itemize}
\item \texttt{radon}: sampled Radon transform of the phantom ($N\times (2M+1)$ matrix);
\item \texttt{reconstruction}: Reconstructed image of the phantom ($K\times K$ matrix );
\item \texttt{phantom}: name of the used phantom (string);
\item \texttt{algorithm}: name of the used algorithm (string);
\item \texttt{options}: options used in the reconstruction (structure depending on the algorithm).
\end{itemize}
The \texttt{options} structure contains information about the sampling ($N,M,K$), on the noise (type of noise and mean and variance in the Gaussian case) and a logical value that determines if the sinogram is plotted or not. Additional parameters depend on the algorithm:
\begin{itemize}
\item Back projection formula: interpolation technique and the low pass filter;
\item Kaczmarz's method: relaxed parameter $\lambda$ and a maximum number of iterations and a tolerance that decide when to stop the reconstruction process; 
\item Kernel methods: the name of the kernel used (Gaussian, multiquadric, inverse multiquadric), the shape parameter of the kernel and which the window function has been used.%AGGIUNGERE COMPACT SUPPORT
\end{itemize}

The graphical user interface is designed in a way that is easy to extend its functionalities for example adding new phantoms or new reconstruction methods. Another improvement that can be added is the possibility to directly compare the results of a reconstruction using two different methods or different parameters.

\backmatter
%Conclusions

\chapter{Conclusions}
In this thesis we studied the problem of clinical image reconstruction. The problem was faced both from an analytical point of view, considering the mathematical aspect of the problem and the classical methods used to solve it, and also with a numerical approach, implementing new algorithms to solve it and comparing the behavior of the different methods.

In the first part we focused on classical Fourier based methods. These methods are founded on the back projection formula and its discretization . We saw that in a discrete context it is possible to obtain only an approximated solution because of the presence of noise and the constrain to have only a finite amount of data.

In the second part we introduced a different approach for solving the image reconstruction problem, called ART. With this kind of methods the solution is obtained solving a linear system. In particular positive definite kernel can be used to this aim, provided a regularization of the Radon transform functional.

The regularization we used is to multiply the kernel function by a window function so that the Radon transform of the product function is finite. Then, we realized these algorithms using particular kernel and window functions and studied their behavior in function of shape parameters.

In the last part of the thesis we compared kernel based with Fourier based methods. We saw that the quality of the reconstruction of the two methods is similar, also in the case of noisy data. The main limits of kernel based methods are the computational time, that grows exponentially with the size of the problem, and a bound for the number of data usable, indeed the linear system involved in the problem can become huge. 

Possible improvement and further works consist in using other kind of kernels and window functions. In this case the main difficulties can be finding the analytical expression of the Radon transform, or try other regularization techniques for the Radon transform integral. Implementing faster algorithms for solving the linear system, for example generating structured or sparse matrix, can be also another improvement. Moreover, an accurate study of the approximation error can give useful informations, e.g. in the determination of optimal shape parameters. Finally, one can introduce a polynomial term in the expression of the solution and then use conditionally positive definite kernels. 

The big number of applications and the vastness of possibilities that can be followed in using the kernel based approach show why this research field has become so important in the last years. 

\appendix
%Appendix A
\chapter{Appendix A: Inverse multiquadrics kernel matrix}\label{app: invMul}
We compute the elements of the matrix $A$ of the inverse multiquadrics reconstruction problem (section \ref{sec: imqRec}). We recall that in that case 
\begin{align*}
a_{k,j}=\frac{2}{\varepsilon^{2}a}\int_{c_{1}}^{c_{2}}{\text{asinh}\left(\sqrt{\frac{\varepsilon^{2}L^{2}-u^{2}}{1+u^{2}}}\right)\,du},
\end{align*}
with
\begin{align*}
&c_{1}=\varepsilon\max{(-L,-|a|\sqrt{H^{2}-r^{2}}+b)} & &c_{2}=\varepsilon\min{(L,|a|\sqrt{H^{2}-r^{2}}+b)}.
\end{align*}
Hence, all we have to do is the compute 
\begin{equation*}
I=\int{\text{asinh}\left( \sqrt{\frac{M^{2}-u^{2}}{1+u^{2}} }\right)\,du}, \qquad M>0.
\end{equation*}
Using the logarithmic representation of asinh we can write 
\begin{align*}
I=&\int{\log{\left( \sqrt{M^{2}-u^{2}}+\sqrt{1+M^{2}}\right)}\,du}-\frac{1}{2}\int{\log{(1+u^{2})}\,du}=\\
&=I_{1}-\left(\text{atan}{u}+\frac{u}{2}\log{(1+u^{2})}-u\right).
\end{align*}
Integrating $I_{1}$ by parts
\begin{align*}
I_{1}&=u\log{\left( \sqrt{M^{2}-u^{2}}+\sqrt{1+M^{2}}\right)}+\\
&\qquad+\int{\frac{u^{2}}{(M^{2}-u^{2})-\sqrt{M^{2}-u^{2}}\sqrt{M^{2}+1}}\,du}=\\
&=u\log{\left( \sqrt{M^{2}-u^{2}}+\sqrt{1+M^{2}}\right)}+I_{2}.
\end{align*}
Adding and subtracting $M^{2}$ in the numerator, we get
\begin{align*}
I_{2}&=-\int{\frac{M^{2}-u^{2}}{(M^{2}-u^{2})-\sqrt{M^{2}-u^{2}}\sqrt{M^{2}+1}}\,du}+\\
&+M^{2}\int{\frac{1}{(M^{2}-u^{2})-\sqrt{M^{2}-u^{2}}\sqrt{M^{2}+1}}\,du}=\\
&=I_{3}+M^{2}I_{4}.
\end{align*}
Setting $\alpha=\text{acos}{\frac{u}{M}}$ and $c=\sqrt{1+\frac{1}{M^{2}}}$ in $I_{3}$, we have
\begin{align*}
I_{3}&=-\int{\frac{\sqrt{M^{2}-u^{2}}}{\sqrt{M^{2}-u^{2}}+\sqrt{M^{2}+1}}\,du}=M\int{\frac{\sin^{2}{\alpha}}{\sin{\alpha}+c}\,d\alpha}=\\
&=M\int{\frac{\sin^{2}{\alpha}-c^{2}}{\sin{\alpha}+c}\,d\alpha}+Mc^{2}\int{\frac{1}{\sin{\alpha}+c}\,d\alpha}=\\
&=M(-\cos{\alpha}-c\alpha)+Mc^{2}\frac{2\text{atan}\left( \frac{\cos{\alpha}}{\sin{\alpha}+(\sqrt{c-1}+\sqrt{c+1})^{2}} \right) +\alpha}{\sqrt{c-1}\sqrt{c+1}}=\\
&=-u-\left(\sqrt{M^{2}+1}\right)\text{acos}\left( \frac{u}{M}\right)+\\
&\qquad+(M^{2}+1)\left[2\text{atan}\left( \frac{u}{\sqrt{M^{2}-u^{2}}+\sqrt{M^{2}+1}+1}\right)+ \text{acos}\left( \frac{u}{M}\right)\right].
\end{align*}
Finally $I_{4}$:
\begin{align*}
I_{4}=\int{\frac{1}{(M^{2}-u^{2})-\sqrt{M^{2}-u^{2}}\sqrt{M^{2}+1}}\,du}=\text{atan}\left(u\sqrt{\frac{M^{2}+1}{M^{2}-u^{2}}} \right)-\text{atan}u
\end{align*}
Putting together the results, since $I=I_{1}-I_{3}-M^{2}I_{4}$, we obtain
\begin{align*}
&\int{\text{asinh}\left( \sqrt{\frac{M^{2}-u^{2}}{1+u^{2}} }\right)\,du}=\frac{u}{2}\text{asinh}\left( \sqrt{\frac{M^{2}-u^{2}}{1+u^{2}} }\right)-(1+M^{2})\text{atan}u+\\
&\qquad+\sqrt{M^{2}+1}\left(\sqrt{M^{2}+1}-1\right)\text{acos}\left( \frac{u}{M}\right)+M^{2}\text{atan}\left(u\sqrt{\frac{M^{2}+1}{M^{2}-u^{2}}} \right)+\\
&\qquad+2(M^{2}+1)\text{atan}\left( \frac{u}{\sqrt{M^{2}-u^{2}}+\sqrt{M^{2}+1}+1}\right).
\end{align*}
Where, of course, this formula is valid for $|u|<M$.

%\appendix B
\chapter{Appendix B: Compactly supported kernel matrix}\label{app: compSupp}
We compute the elements of the matrix $A$ of the compactly supported function reconstruction problem (section \ref{subsec: compSupp}). We recall that in that case   
\begin{equation*}
a_{k,j}=\int_{\R}{b(x_{s})(1-\nu^{2}\norm{x_{s}})_{+}\,ds},
\end{equation*}
where
\begin{align*}
b(x)=\left\{
\begin{aligned}
&g(t_j-x\cdot v_j) & &\text{if}\ |t_j-x\cdot v_j|\leq\frac{1}{\varepsilon}\\
&0 & &\text{if}\ |t_j-x\cdot v_j|>\frac{1}{\varepsilon}
\end{aligned}
\right.
\end{align*}
and
\begin{align*}
g(t)&=\left\{
\begin{aligned}
&\frac{2}{\varepsilon}\left[\frac{\sqrt{1-\varepsilon^2t^2}}{3}(2\varepsilon^2t^2+1)-\varepsilon^2t^2\text{acosh}\left(\frac{1}{\varepsilon|t|}\right)\right] & &\text{if}\ t\neq0\\
&\frac{2}{3\varepsilon} & &\text{if}\ t=0.
\end{aligned}
\right.
\end{align*}

Since $\norm{x_{s}}^{2}=r^2+s^2$, one obtains
\begin{align*}
a_{kj}&=(1-\nu^{2}r^2)\int_{\nu\sqrt{r^2+s^2}\leq1}{b(x_{s})\,ds}-\nu^2\int_{\nu\sqrt{r^2+s^2}\leq1}{b(x_{s})s^2\,ds}=\\
&=\left\{
\begin{aligned}
&(1-\nu^{2}r^2)I_{1}-\nu^2I_{2} & &\text{if}\ |r|\leq\frac{1}{\nu}\\
&0 & &\text{if}\ |r|>\frac{1}{\nu}
\end{aligned}
\right.
\end{align*}
Setting $t-x_{s}\cdot v=as+b$ and $D_{s}=\left\{s:\ |s|\leq\sqrt{\frac{1}{\nu^2}-r^2},\ \varepsilon|as+b|\leq1\right\}$,
\begin{align*}
I_{1}&=\frac{2}{\varepsilon}\int_{D_{s}}{\left[\frac{\sqrt{1-\varepsilon^2(as+b)^2}}{3}(2\varepsilon^2(as+b)^2+1)\right]\,ds}+\\
&-\frac{2}{\varepsilon}\int_{D_{s}}{\varepsilon^2(as+b)^2\text{acosh}\left(\frac{1}{\varepsilon|as+b|}\right)\,ds}\\
I_{2}&=\frac{2}{\varepsilon}\int_{D_{s}}{\left[\frac{\sqrt{1-\varepsilon^2(as+b)^2}}{3}(2\varepsilon^2(as+b)^2+1)\right]s^2\,ds}+\\
&-\frac{2}{\varepsilon}\int_{D_{s}}{\varepsilon^2(as+b)^2\text{acosh}\left(\frac{1}{\varepsilon|as+b|}\right)s^2\,ds}.
\end{align*}
We distinguish the cases $a=0$ and $a\neq0$:
\begin{itemize}
\item If $a=0$
\begin{equation*}
I_{1}=\left\{
\begin{aligned}
&\frac{4}{\varepsilon}\sqrt{\frac{1}{\nu^2}-r^2}\left[\frac{\sqrt{1-\varepsilon^2b^2}(2\varepsilon^2b^2+1)}{3}-2\varepsilon^2b^2\text{acosh}\left(\frac{1}{\varepsilon|b|}\right)\right] & &\text{if} \ |b|\leq\frac{1}{\varepsilon}\\
&0 & &\text{if} \ |b|>\frac{1}{\varepsilon}
\end{aligned}
\right.
\end{equation*}
\begin{equation*}
I_{2}=\left\{
\begin{aligned}
&\frac{4}{3\varepsilon}\left(\frac{1}{\nu^2}-r^2\right)^{\frac{3}{2}}\left[\frac{\sqrt{1-\varepsilon^2b^2}(2\varepsilon^2b^2+1)}{3}-2\varepsilon^2b^2\text{acosh}\left(\frac{1}{\varepsilon|b|}\right)\right] & &\text{if} \ |b|\leq\frac{1}{\varepsilon}\\
&0 & &\text{if} \ |b|>\frac{1}{\varepsilon}
\end{aligned}
\right.
\end{equation*}
thus we conclude that
\begin{itemize}
\item if $\nu|r|\leq1$, $\varepsilon|b|\leq1$ and $b\neq0$
\begin{equation*}
a_{kj}=\frac{8}{3}\frac{(1-\nu^2r^2)^{3/2}}{\varepsilon\nu}\left[\frac{\sqrt{1-\varepsilon^2b^2}(2\varepsilon^2b^2+1)}{3}-2\varepsilon^2b^2\text{acosh}\left(\frac{1}{\varepsilon|b|}\right)\right];
\end{equation*}
\item if $\nu|r|\leq1$ and $b=0$
\begin{equation*}
a_{kj}=\frac{8}{9\varepsilon}\left(\frac{1}{\nu^2}-r^2\right)^{\frac{2}{3}}
\end{equation*}
\item if $\nu|r|>1$, $\varepsilon|b|>1$, then $a_{k,j}=0$
\end{itemize}
\item If $a\neq0$
\begin{align*}
I_1=\int_{c_{1}}^{c_2}{\frac{2}{\varepsilon}\left[\frac{\sqrt{1-u^2}}{3}(2u^2+1)\right]\,\frac{du}{\varepsilon a}}-\int_{c_1}^{c_2}{\frac{2}{\varepsilon}u^2\text{acosh}\left(\frac{1}{|u|}\right)\,\frac{du}{\varepsilon a}},
\end{align*}
where $u=\varepsilon(as+b)$ and 
\begin{align*}
&c_1=\max{\left(-1,b\varepsilon-\varepsilon |a|\sqrt{\frac{1}{\nu^2}-r^2}\right)}, &c_2=\min{\left(1,b\varepsilon+\varepsilon |a|\sqrt{\frac{1}{\nu^2}-r^2}\right)}
\end{align*}
thus
\begin{align*}
I_1&=\frac{1}{6\varepsilon^2 a}\left[3\arcsin{u}+u\sqrt{1-u^2}(2u^2+1)\right]_{c_1}^{c_2}+\\
&-\frac{2}{\varepsilon^2 a}\left[\frac{1}{6}\arcsin{u}+\frac{u^3}{3}\text{acosh}\left(\frac{1}{|u|}\right)-\frac{u}{6}\sqrt{1-u^2}\right]_{c_1}^{c_2}=\\
&=\frac{1}{3\varepsilon^2 a}\left[\frac{1}{2}\arcsin{u}+u\sqrt{1-u^2}(u^2+\frac{3}{2})-2u^3\text{acosh}\left(\frac{1}{|u|}\right)\right]_{c_1}^{c_2}
\end{align*}
and the second integral becomes
\begin{align*}
I_2&=\int_{c_{1}}^{c_2}{\frac{2}{\varepsilon}\left[\frac{\sqrt{1-u^2}}{3}(2u^2+1)\left(\frac{u-\varepsilon b}{\varepsilon a}\right)^2\right]\,\frac{du}{\varepsilon a}}+\\
&-\int_{c_1}^{c_2}{\frac{2}{\varepsilon}u^2\text{acosh}\left(\frac{1}{|u|}\right)\left(\frac{u-\varepsilon b}{\varepsilon a}\right)^2\,\frac{du}{\varepsilon a}}=\\
&=\frac{2}{3\varepsilon^4 a^3}\left[\frac{1}{12}\left(b^2\varepsilon^2+\frac{1}{10}\right)\arcsin{u}+\frac{2}{9}b\varepsilon (1-u^2)^{\frac{3}{2}}+\right.\\
&-\frac{u^3}{30}(6u^2-15b\varepsilon u+10b^2\varepsilon^2)\text{acosh}\left(\frac{1}{|u|}\right)+\\
&+\frac{\sqrt{1-u^2}}{60}\left(\frac{20}{3}u^5-16b\varepsilon u^4+(10b^2\varepsilon^2+\frac{19}{3})u^3+\right.\\
&\left.\left.-\frac{14}{3}b\varepsilon u^2+(15b^2\varepsilon^2-\frac{1}{2})u-\frac{28}{3}b\varepsilon\right)
\right]_{c_1}^{c_2}.
\end{align*}
Finally, if $\nu|r|\leq1$,  we have: 
\begin{align}
a_{k,j}&=\frac{1}{3\varepsilon^2 a }\left[ 
\frac{1}{2}\arcsin{u}\left( 1-\nu^2r^2-\frac{\nu^2}{\varepsilon^2 a^2}(b^2\varepsilon^2+\frac{1}{10})\right)+\right.\\
&+\sqrt{1-u^2}\left( u(u^2+\frac{3}{2})(\-\nu^2r^2)-\frac{\nu^2q_2(u)}{10\varepsilon^2 a^2}-\frac{4\nu^2b\varepsilon}{3\varepsilon^2 a^2}(1-u^2)\right)+\\
&\left.+u^3\text{acosh}\left(\frac{1}{|u|}\right)\left( \frac{\nu^2q_1(u)}{5\varepsilon^2 a^2}-2(1-\nu^2r^2)\right)
\right]_{c_1}^{c_2}
\label{eq: matrix_cs}
\end{align}
where
\begin{align*}
&q_{1}(u)=6u^2-15b\varepsilon u+10b^2\varepsilon^2\\
&q_{2}(u)=\frac{20}{3}u^5-16b\varepsilon u^4+(10b^2\varepsilon^2+\frac{19}{3})u^3+\\
&-\frac{14}{3}b\varepsilon u^2+(15b^2\varepsilon^2-\frac{1}{2})u-\frac{28}{3}b\varepsilon.
\end{align*}
\end{itemize}
We observe that if $\varepsilon|b|\leq1$, then $c_{1}<c_{2}$ always holds.

While for $\nu|r|>1$, $a_{k,j}=0$.

At last we observe that because of the term $\text{acosh}(|u|^{-1})$ in \eqref{eq: matrix_cs}, we have to consider apart the cases $c_{1}=0,\ c_{2}=0$. In these cases, it is easy to see that, because of continuity, it is sufficient to consider the limit of \eqref{eq: matrix_cs} for $u\rightarrow0$, so that $u^3\text{acosh}(|u|^{-1})\rightarrow0$.

%Bibliography
\cleardoublepage 
\addcontentsline{toc}{chapter}{Bibliography}
\bibliographystyle{plain}
%\bibliography{thesisbib}

\end{document}